\documentclass[10pt,english,reqno]{amsart} 

\calclayout

\usepackage[usenames,dvipsnames]{color}
\newcommand{\yifei}[1]{{\color{orange}   [#1]}}

\usepackage[hidelinks]{hyperref}
\usepackage{graphicx}
\usepackage{amssymb}
\usepackage{epstopdf}
\usepackage{enumerate}
\usepackage{mathabx}
\usepackage{tikz-cd}
\usepackage{MnSymbol}
\usepackage[mathscr]{eucal}
\usepackage[multiple]{footmisc}

\DeclareMathSymbol{A}{\mathalpha}{operators}{`A}
\DeclareMathSymbol{B}{\mathalpha}{operators}{`B}
\DeclareMathSymbol{C}{\mathalpha}{operators}{`C}
\DeclareMathSymbol{D}{\mathalpha}{operators}{`D}
\DeclareMathSymbol{E}{\mathalpha}{operators}{`E}
\DeclareMathSymbol{F}{\mathalpha}{operators}{`F}
\DeclareMathSymbol{G}{\mathalpha}{operators}{`G}
\DeclareMathSymbol{H}{\mathalpha}{operators}{`H}
\DeclareMathSymbol{I}{\mathalpha}{operators}{`I}
\DeclareMathSymbol{J}{\mathalpha}{operators}{`J}
\DeclareMathSymbol{K}{\mathalpha}{operators}{`K}
\DeclareMathSymbol{L}{\mathalpha}{operators}{`L}
\DeclareMathSymbol{M}{\mathalpha}{operators}{`M}
\DeclareMathSymbol{N}{\mathalpha}{operators}{`N}
\DeclareMathSymbol{O}{\mathalpha}{operators}{`O}
\DeclareMathSymbol{P}{\mathalpha}{operators}{`P}
\DeclareMathSymbol{Q}{\mathalpha}{operators}{`Q}
\DeclareMathSymbol{R}{\mathalpha}{operators}{`R}
\DeclareMathSymbol{S}{\mathalpha}{operators}{`S}
\DeclareMathSymbol{T}{\mathalpha}{operators}{`T}
\DeclareMathSymbol{U}{\mathalpha}{operators}{`U}
\DeclareMathSymbol{V}{\mathalpha}{operators}{`V}
\DeclareMathSymbol{W}{\mathalpha}{operators}{`W}
\DeclareMathSymbol{X}{\mathalpha}{operators}{`X}
\DeclareMathSymbol{Y}{\mathalpha}{operators}{`Y}
\DeclareMathSymbol{Z}{\mathalpha}{operators}{`Z}

\newcommand{\integers}{\mathbf Z}
\newcommand{\rationals}{\mathbf Q}

\newcommand{\complexes}{\mathbf C}
\newcommand{\base}{\textnormal{\textsf{k}}}
\newcommand{\coeff}{\textnormal{\textsf{e}}}
\newcommand{\characteristic}{\operatorname{char}}
\newcommand{\Gal}{\operatorname{Gal}}
\newcommand{\Art}{\mathrm{Art}}

\newcommand{\PrL}{\textnormal{\textsf{Pr}}^L}
\newcommand{\PrR}{\textnormal{\textsf{Pr}}^R}

\newcommand{\Cat}{\textnormal{\textsf{Cat}}}
\newcommand{\Spc}{\textnormal{\textsf{Spc}}}
\newcommand{\Sptr}{\textnormal{\textsf{Sptr}}}
\newcommand{\Mod}{\textnormal{-\textsf{Mod}}}
\newcommand{\Comod}{\textnormal{-\textsf{Comod}}}
\newcommand{\Shv}{\textnormal{\textsf{Shv}}}
\newcommand{\CAlg}{\textnormal{\textsf{CAlg}}}
\newcommand{\Alg}{\textnormal{\textsf{Alg}}}
\newcommand{\Ring}{\textnormal{\textsf{Ring}}}

\newcommand{\kernel}{\textnormal{\textsf{K}}}
\newcommand{\Span}{\operatorname{Span}}
\newcommand{\Ind}{\operatorname{Ind}}
\newcommand{\LKE}{\mathrm{LKE}}
\newcommand{\Gp}{\textnormal{\textsf{Gp}}}
\newcommand{\com}{\mathrm{com}}

\newcommand{\opposite}{\mathrm{op}}
\newcommand{\Ho}{\operatorname{Ho}}
\newcommand{\Hom}{\mathrm{Hom}}
\newcommand{\Ext}{\mathrm{Ext}}
\newcommand{\SExt}{\mathscr Ext}
\newcommand{\Maps}{\mathrm{Maps}}
\newcommand{\End}{\operatorname{End}}
\newcommand{\Fun}{\mathrm{Fun}}
\newcommand{\Fib}{\mathrm{Fib}}
\newcommand{\Cofib}{\mathrm{Cofib}}
\newcommand{\profin}{\mathrm{profin}}

\newcommand{\Aut}{\operatorname{Aut}}
\newcommand{\id}{\mathrm{id}}
\newcommand{\ev}{\mathrm{ev}}
\newcommand{\unit}{\mathrm{unit}}
\newcommand{\counit}{\mathrm{counit}}

\newcommand{\deloop}{\textnormal{\textsf{B}}}

\DeclareMathOperator*\colim{colim}

\newcommand{\Stk}{\textnormal{\textsf{Stk}}}
\newcommand{\Gpd}{\textnormal{\textsf{Gpd}}}
\newcommand{\Sch}{\textnormal{\textsf{Sch}}}
\newcommand{\Perf}{\textnormal{\textsf{Perf}}}
\newcommand{\derived}{\textnormal{\textsf{D}}}
\newcommand{\Lis}{\textnormal{\textsf{Lis}}}
\newcommand{\SLis}{\mathscr{L}is}
\newcommand{\etale}{\textnormal{\'et}}
\newcommand{\Zar}{\textnormal{Zar}}

\newcommand{\Spec}{\operatorname{Spec}}
\newcommand{\tr}{\operatorname{tr}}

\newcommand{\Div}{\textnormal{\textsf{Div}}}
\newcommand{\Fr}{\operatorname{Fr}}
\newcommand{\Tr}{\operatorname{Tr}}
\newcommand{\Nm}{\operatorname{Nm}}
\newcommand{\res}{\operatorname{res}}

\newcommand{\SHom}{\mathscr{H}om}
\newcommand{\SMaps}{\mathscr{M}aps}
\newcommand{\SQuad}{\mathscr{Q}uad}
\newcommand{\strict}{\mathrm{st}}
\newcommand{\tors}{\textnormal{-}\mathrm{tors}}
\newcommand{\AJ}{\mathrm{AJ}}
\newcommand{\Pic}{\mathrm{Pic}}
\newcommand{\red}{\mathrm{red}}
\newcommand{\fact}{\mathrm{fact}}
\newcommand{\specialization}{\operatorname{sp}}

\newcommand{\Ktheory}{\textnormal{\textsf{K}}}
\newcommand{\Tate}{\mathrm{Tate}}

\newcommand{\Rep}{\textnormal{\textsf{Rep}}}
\newcommand{\SRep}{\mathscr{R}ep}
\newcommand{\Ran}{\mathrm{Ran}}
\newcommand{\Bun}{\mathrm{Bun}}
\newcommand{\Hec}{\mathrm{Hec}}
\newcommand{\Gr}{\mathrm{Gr}}
\newcommand{\Fl}{\mathrm{Fl}}
\newcommand{\Sht}{\mathrm{Sht}}
\newcommand{\disj}{\mathrm{disj}}
\newcommand{\IC}{\mathrm{IC}}
\newcommand{\aff}{\mathrm{aff}}
\newcommand{\translate}{\textnormal{\textsf{t}}}
\newcommand{\CT}{\operatorname{CT}}
\newcommand{\Av}{\operatorname{Av}}
\newcommand{\swap}{\mathrm{sw}}
\newcommand{\cusp}{\mathrm{cusp}}

\newcommand{\sconn}{\mathrm{sc}}
\newcommand{\der}{\mathrm{der}}
\newcommand{\adjoint}{\mathrm{ad}}
\newcommand{\abelian}{\mathrm{ab}}
\newcommand{\ord}{\operatorname{ord}}
\newcommand{\mon}{\textnormal{-mon}}
\newcommand{\geom}{\mathrm{geom}}
\newcommand{\CExt}{\operatorname{CExt}}
\newcommand{\Hf}{\mathrm{Hf}}
\newcommand{\Torel}{\mathrm{Torel}}

\newcommand{\SL}{\mathrm{SL}}
\newcommand{\GL}{\mathrm{GL}}
\newcommand{\PGL}{\mathrm{PGL}}
\newcommand{\Sp}{\mathrm{Sp}}

\newcommand{\Sat}{\textnormal{\textsf{Sat}}}
\newcommand{\SSat}{\mathscr{S}at}

\newcommand{\loo}[1]{(\!(#1)\!)}
\newcommand{\arc}[1]{[\![#1]\!]}

\makeatletter
\renewcommand{\@secnumfont}{\bfseries}
\makeatother

\newtheorem{thm}[subsubsection]{Theorem}
\newtheorem*{thm*}{Theorem}
\newtheorem{thmx}{Theorem}

\newtheorem{prop}[subsubsection]{Proposition}
\newtheorem{lem}[subsubsection]{Lemma}
\newtheorem{conj}[subsubsection]{Conjecture}
\newtheorem{cor}[subsubsection]{Corollary}

\theoremstyle{definition}
\newtheorem{defn}[subsubsection]{Definition}

\newtheorem{eg}[subsubsection]{Example}
\newtheorem{rem}[subsubsection]{Remark}

\numberwithin{equation}{section}
\allowdisplaybreaks

\newtheoremstyle{void}
    {}{}{}{}
    {\bfseries}{.}{ }
    {\thmname{#1}\thmnumber{#2}\thmnote{ {\mdseries \textit{#3}}}}
\theoremstyle{void}

\newtheorem{void}[subsubsection]{}

\title[Spectral decomposition of genuine cusp forms]{Spectral decomposition of genuine cusp forms over global function fields}

\author{Yifei Zhao}
\date{\today}
\begin{document}

\begin{abstract}
We prove a version of the twisted geometric Satake equivalence and extend the Langlands parametrization of V.~Lafforgue to certain covers of reductive groups.
\end{abstract}

\maketitle

This is an updated version of the same-named article in \emph{Compos.~Math.}~\textbf{160}, 1194--1260, containing several corrections and improvements.

\setcounter{tocdepth}{2}
\tableofcontents

%%%%%%%%%%%%%%

\newpage

\section*{Introduction}

In the seminal work \cite{MR3787407}, V.~Lafforgue constructed the Langlands parametrization of cusp forms for reductive groups over global function fields.

The goal of this article is to extend this parametrization to a large class of covers of reductive groups, following the strategy indicated in \cite[\S14]{MR3787407}. The main new ingredient is a version of the twisted geometric Satake equivalence, extending earlier works of Finkelberg, Lysenko, Reich, and Gaitsgory (\emph{cf.}~\cite{MR2684259, MR2956088, MR3769731}).

The class of covers we treat includes the ones considered by Brylinski and Deligne in \cite{MR1896177}, so our result contributes to the Langlands--Weissman program for Brylinski--Deligne covers, as formulated by Weissman, Gan, and Gao (\emph{cf.}~\cite{MR3802418, MR3802419}). We refer the reader to \cite{MR3802417} for a survey of this program.

\subsection{Spectral decomposition}

\begin{void}
\label{void-global-covering-group-intro}
Let $F$ be a global field of positive characteristic $p$. Denote by $\mathbb A_F$ the topological ring of ad\`eles of $F$. Let $\ell \neq p$ be a prime and fix an algebraic closure $\overline{\rationals}_{\ell}$ of $\rationals_{\ell}$, with a chosen half-integer Tate twist $\overline{\rationals}_{\ell}(\frac{1}{2})$.

Our group-theoretic input is a pair $(G, \mu)$, where $G$ is a split reductive group over $F$ and $\mu$ is a ``parameter" for covers of (the ad\`elic points of) $G$.

The precise meaning of $\mu$ will be explained in \S\ref{void-nature-of-etale-level}. For now, let us take for granted that $\mu$ gives rise to a topological central extension
\begin{equation}
\label{eq-global-covering-group-intro}
1 \rightarrow A \rightarrow \widetilde G \rightarrow G(\mathbb A_F) \rightarrow 1,
\end{equation}
for some finite abelian group $A$, equipped with a canonical splitting over $G(F)$.

Fix an injective character $\zeta : A \hookrightarrow \overline{\rationals}{}_{\ell}^{\times}$. The Langlands--Weissman program, in its global function field incarnation, studies the $\overline{\rationals}_{\ell}$-vector space of \emph{$\zeta$-genuine automorphic forms} on $\widetilde G$, \emph{i.e.}~locally constant functions $f : G(F) \backslash \widetilde G \rightarrow \overline{\rationals}_{\ell}$ which are $A$-equivariant against $\zeta$.
\end{void}

\begin{void}
As in the case for reductive groups, $\zeta$-genuine automorphic forms admit a notion of cuspidality, defined by the vanishing of constant terms.

Furthermore, one may constrain the action of the connected part $Z^{\circ}$ of the center of $G$, by demanding our $\zeta$-genuine automorphic forms to be invariant under the action of a cocompact lattice $\Xi$ in $Z^{\sharp}(F) \backslash \widetilde Z^{\sharp}$, for a torus $Z^{\sharp}$ isogenous to $Z^{\circ}$ with induced cover $\widetilde Z^{\sharp}$. The torus $Z^{\sharp}$ appears because the induced cover of $Z^{\circ}$ may \emph{not} be abelian, while $\widetilde Z^{\sharp}$ is abelian and its image in $\widetilde G$ is central by construction.

Imposing these two conditions leads us to the $\overline{\rationals}_{\ell}$-vector space
\begin{equation}
\label{eq-cusp-form-introduction}
	\Fun_{\cusp, \zeta}(G(F) \backslash \widetilde G/\Xi, \overline{\rationals}_{\ell}).
\end{equation}
The main goal of this article is to decompose \eqref{eq-cusp-form-introduction} according to spectral data.
\end{void}

\begin{void}
\label{void-ell-group-intro}
Let us now describe the spectral data involved. For this, we fix an algebraic closure $\bar F$ of $F$ and denote by $\Gal_F$ the Galois group of $\bar F / F$.

To the pair $(G, \mu)$, one may attach certain combinatorial data. The ``classical" part of these data consists of the root data of $G$ together with a certain quadratic form. Using them, one defines the ``metaplectic dual group" $H$: This is the dual group found in \cite{MR1227098, MR2684259, MR2963537, MR3769731}, for increasingly general kinds of $(G, \mu)$.

Contrary to the case of reductive groups, the combinatorial package associated to $(G, \mu)$ also includes certain $2$-categorical data. Following Weissman's approach in \cite{MR3802418}, we obtain from them an extension of topological groups
\begin{equation}
\label{eq-ell-group-intro}
	1 \rightarrow H(\overline{\rationals}_{\ell}) \rightarrow {}^LH_F \rightarrow \Gal_F \rightarrow 1,
\end{equation}
which we shall refer to as the \emph{$L$-group} of $(G, \mu)$. Then we proceed as in the case of reductive groups, defining an \emph{$L$-parameter} to be an $H(\overline{\rationals}_{\ell})$-conjugacy class of continuous sections $\sigma : \Gal_F \rightarrow {}^LH_F$ of the surjection in \eqref{eq-ell-group-intro}.

Let us mention that the $L$-group \eqref{eq-ell-group-intro} is \emph{not} the only way to package the dual data of $(G, \mu)$. In fact, in the main body of the text, we will consider a refinement of \eqref{eq-ell-group-intro}, which we call ``metaplectic dual data" following Gaitsgory and Lysenko (\emph{cf.}~\cite{MR3769731}). For the moment, however, let us use \eqref{eq-ell-group-intro} to formulate our main result.
\end{void}

\begin{thmx}
\label{thmx-spectral-decomposition}
There is a canonical decomposition
\begin{equation}
\label{eq-spectral-decomposition-intro}
	\Fun_{\cusp, \zeta}(G(F) \backslash \widetilde G/\Xi, \overline{\rationals}_{\ell}) \simeq \bigoplus_{[\sigma]} \mathbf H_{[\sigma]},
\end{equation}
indexed by $L$-parameters $[\sigma]$.
\end{thmx}

\begin{void}
Theorem \ref{thmx-spectral-decomposition}, as stated above, is nearly meaningless. A precise version of it appears as Theorem \ref{thm-spectral-decomposition} in the main body of the text. It includes the compatibility with Hecke action expected of the Langlands parametrization, as well as other information on the $L$-parameters $[\sigma]$ appearing in the index set of \eqref{eq-spectral-decomposition-intro}.

The main body of the text also treats the case where $G$ is not necessarily split. In that case, the left-hand-side of \eqref{eq-spectral-decomposition-intro} must be enlarged to account for different forms of $G$.

Theorem \ref{thmx-spectral-decomposition} is an addendum to a long series of works on the ``automorphic-to-Galois" direction of the Langlands program for global function fields, due to Drinfeld, L.~Lafforgue, and V.~Lafforgue (\emph{cf.}~\cite{MR0498489, MR918745, MR902291, MR1875184, MR3787407}). Contrary to these highly original works, the main ideas of the present text have already been laid out in \cite{MR3787407} and \cite{MR3769731}. The only contribution from the author is a new treatment of the twisted geometric Satake equivalence, discussed in \S\ref{sec-satake-equivalence-intro} below.
\end{void}

\begin{void}[The meaning of $\mu$]
\label{void-nature-of-etale-level}
Let us now say what $\mu$ stands for. This is the only place where our formulation of the Langlands parametrization differs from the usual one found in the Langlands--Weissman program.

Namely, in the latter context, one parametrizes covers of $G(\mathbb A_F)$ by pairs $(E, N)$, where $E$ is a Brylinski--Deligne cover of $G$ and $N\ge 1$ is an integer such that $\mu_N(F)$ has cardinality $N$ (\emph{cf.}~\cite{MR3802418}). However, this formalism leads to technical difficulties in the proof of the twisted geometric Satake equivalence.\footnote{See \S\ref{sec-integral-vs-etale-levels} for a more detailed discussion.}

Fortunately, an alternative parametrization of covers found by Deligne (\emph{cf.}~\cite{MR1441006}) allows us to circumvent these difficulties. Formulated in modern language, a parameter of \cite{MR1441006} is a morphism of pointed (higher) \'etale stacks
\begin{equation}
\label{eq-etale-level-intro}
	\deloop_F G \rightarrow \deloop_F^4 A(1),
\end{equation}
where $\deloop_F$ denotes the deloop functor for \'etale stacks over $F$, and $A(1)$ is the Tate twist of the finite abelian group $A$. In this text, we let $\mu$ be a morphism \eqref{eq-etale-level-intro}, which we shall refer to as an ($A$-valued) \emph{\'etale level}.\footnote{The terminology is borrowed from the quantum geometric Langlands program (which is in turn motivated by 2D conformal field theory), where an analogous notion in de Rham cohomology is refered to as a ``level". In \cite{zhao2022metaplectic}, an  \'etale level is called an ($A$-valued) \'etale metaplectic cover.}

The idea of parametrizing covering groups by \'etale levels entered the Langlands program through Gaitsgory and Lysenko's work \cite{MR3769731}.
\end{void}

\begin{rem}
In \cite[\S2]{MR1441006}, \'etale levels are described by cocycles. One may wonder why we choose not to use this explicit description.

The reason is that in order to define metaplectic dual data, it is necessary to work with additional structures on \'etale levels, which are natural from the point of view of higher algebra but become quite cumbersome when written in terms of cocycles. For this reason, we systematically adopt the tools of higher algebra (\emph{cf.}~\cite{MR2522659, lurie2017higher}) in this text.
\end{rem}

\subsection{Cohomology of Shtukas}

\begin{void}
\label{void-cohomology-of-shtukas-context-intro}
Let us sketch the proof of Theorem \ref{thmx-spectral-decomposition}, essentially following \cite[\S14]{MR3787407}. We begin by formulating an integral version of the problem.

Let $X$ be a smooth, proper, geometrically connected curve over a finite field $\base$, with field of fractions $F$. Let $D \subset X$ be a $\base$-finite closed subscheme and $\mathring X$ be its complement. We shall assume that our \'etale level comes from a morphism of pointed \'etale stacks over $\mathring X$:
$$
\deloop_{\mathring X} G \rightarrow \deloop^4_{\mathring X} A(1).
$$
where $\deloop_{\mathring X}$ stands for the deloop functor for \'etale stacks over $\mathring X$.

Denote by $\Bun_{G, D}$ the moduli stack of $G$-bundles over $X$ trivialized along $D$. The \'etale level $\mu$ determines an (\'etale) $A$-gerbe $\mathscr G_{\Bun_{G, D}}$ over $\Bun_{G, D}$, which ``categorifies" \eqref{eq-global-covering-group-intro} in the following sense: A procedure akin to taking the trace of Frobenius yields an $A$-torsor
\begin{equation}
\label{eq-global-torsor-intro}
\widetilde{\Bun}_{G, D} \rightarrow \Bun_{G, D}(\base),
\end{equation}
whose pullback along the uniformization map $G(\mathbb A_F) \rightarrow \Bun_{G, D}(\base)$ recovers $\widetilde G$.

Instead of $\zeta$-genuine functions on $G(F) \backslash\widetilde G$, we shall study those defined on $\widetilde{\Bun}_{G, D}$.
\end{void}

\begin{void}
\label{void-local-hecke-stack-intro}
Let us introduce, somewhat informally, two more geometric objects: the local Hecke stack and the moduli stack of Shtukas of Drinfeld and Varshavsky (\emph{cf.}~\cite{MR902291, MR2061225}). The formal definitions are provided in \S\ref{void-local-hecke-stack-definition} and \S\ref{void-moduli-of-shtukas}.

For a finite set $I$, the local Hecke stack $\Hec_{G, I}$ parametrizes a point $x^I$ of $\mathring X^I$, together with a modification at $x^I$ of $G$-bundles over the formal disk $D_{x^I}$, represented as
$$
P^0 \overset{x^I}{\sim} P^1.
$$

The moduli stack of Shtukas $\Sht_{G, D}^I$ parametrizes a point $x^I$ of $\mathring X^I$, together with a modification at $x^I$ of $G$-bundles over $X$ equipped with trivializations along $D$:
\begin{equation}
\label{eq-shtuka-point-intro}
(P^0, \phi^0) \overset{x^I}{\sim} (P^1, \phi^1),
\end{equation}
as well as an isomorphism $(P^1, \phi^1) \simeq ({}^{\tau}P^0, {}^{\tau}\phi^0)$, where ${}^{\tau}$ denotes the Frobenius twist.

Restricting \eqref{eq-shtuka-point-intro} to $D_{x^I}$ defines a morphism of stacks over $\mathring X^I$:
\begin{equation}
\label{eq-shtuka-to-local-hecke-intro}
r : \Sht_{G, D}^I \rightarrow \Hec_{G, I}.
\end{equation}
\end{void}

\begin{void}
In addition to $\mathscr G_{\Bun_{G, D}}$, the \'etale level $\mu$ defines an $A$-gerbe $\mathscr G_{\Hec_{G, I}}$ over $\Hec_{G, I}$. The key observation is that the pullback of $\mathscr G_{\Hec_{G, I}}$ along \eqref{eq-shtuka-to-local-hecke-intro} is \emph{canonically} trivial.

This observation has the following practical consequence: Given a ``$(\mathscr G_{\Hec_{G, I}}, \zeta)$-twisted constructible complex" $\mathscr A$ of $\overline{\rationals}_{\ell}$-sheaves over $\Hec_{G, I}$, a notion we define in Appendix \ref{sec-twisted-ell-adic-sheaves}, its pullback $r^*\mathscr A$ is \emph{untwisted}, \emph{i.e.}~a usual constructible complex of $\overline{\rationals}_{\ell}$-sheaves. In particular, one may take the (compactly supported) cohomology of $r^*\mathscr A$ along the structural morphism $\Sht_{G, D}^I \rightarrow \mathring X^I$ to obtain ``cohomology of Shtukas with coefficients in $\mathscr A$."\footnote{More precisely, we shall take compactly supported cohomology along $\Sht_{G, D}^I/\Xi \rightarrow \mathring X^I$, which is a relative ind-Deligne--Mumford stack. The result is then an ind-constructible complex over $\mathring X^I$.} For the kind of objects $\mathscr A$ we consider, $r^*\mathscr A$ will furthermore be perverse relative to $\mathring X^I$.

In the special case $I = \emptyset$, the moduli stack $\Sht_{G, D}^{\emptyset}$ coincides with the discrete stack $\Bun_{G, D}(\base)$ and $r^* \mathscr G_{\Hec_{G, \emptyset}}$ is the trivial $A$-gerbe. However, the canonical trivialization on $r^*\mathscr G_{\Hec_{G, \emptyset}}$ is \emph{not} the identity map. Rather, it corresponds to the $A$-torsor \eqref{eq-global-torsor-intro}. The cohomology of $\Sht_{G, D}^{\emptyset}$ with coefficients in $\overline{\rationals}_{\ell}$, taken in the above sense, thus encodes $\zeta$-genuine functions on $\widetilde{\Bun}_{G, D}$.
\end{void}

\begin{void}
The source of coefficients on $\Sht_{G, D}^I$ is provided by the twisted geometric Satake equivalence, explained in more detail in \S\ref{sec-satake-equivalence-intro} below.

For now, it suffices to say that this equivalence supplies a family of functors
\begin{equation}
\label{eq-satake-functors-intro}
	\Rep(({}^LH_{\mathring X})^I) \rightarrow \derived_{\mathscr G, \zeta}(\Hec_{G, I})
\end{equation}
natural in the finite set $I$. Here, the source is an appropriately defined category of representation of $I$-copies of an integral form ${}^LH_{\mathring X}$ of the $L$-group \eqref{eq-ell-group-intro}, and the target is the derived category of $(\mathscr G_{\Hec_{G, I}}, \zeta)$-twisted constructible complexes of $\overline{\rationals}_{\ell}$-sheaves over $\Hec_{G, I}$.

In fact, we shall also need a version of \eqref{eq-satake-functors-intro} for \emph{iterated} Hecke stacks, which allows us to endow the cohomology of Shtukas with equivariance structure with respect to the partial Frobenius endomorphisms of $\mathring X^I$.
\end{void}

\begin{void}
With the above preparation, we are in a position to essentially repeat V.~Lafforgue's construction of the spectral decomposition. This is done in \S\ref{sec-spectral-decomposition} of the article. (It involves some confusing computations with gerbes, but no serious difficulties.)

We emphasize that all of the above is already explained in \cite[\S14]{MR3787407}, with one minor difference: The gerbes considered there are the gerbes parametrizing $N$th roots of unity of line bundles (where $N \mid q - 1$, for $q$ the cardinality of $\base$). The gerbes we consider are more general, so we explain in \S\ref{sec-global-function-fields-preparation} how the relevant pieces of structure on them arise.

One may summarize the situation as follows: In generalizing \cite{MR3787407} to covers, the only new feature is that we consider more general coefficients in the cohomology of Shtukas.
\end{void}

\subsection{Geometric Satake equivalence}
\label{sec-satake-equivalence-intro}

\begin{void}
Let us now discuss the twisted geometric Satake equivalence, whose proof occupies the bulk of this article.

For this result, we may take $X$ to be a smooth curve over an arbitrary field $\base$.\footnote{For the application to global function fields, $X$ will be the possibly open curve $\mathring X$ of \S\ref{void-cohomology-of-shtukas-context-intro}.} Let $G$ be a split reductive group equipped with an $A$-valued \'etale level $\mu$ defined over $X$. To $\mu$, one associates an $A$-gerbe $\mathscr G_{\Hec_G}$ over the local Hecke stack $\Hec_G$. (As opposed in \S\ref{void-local-hecke-stack-intro}, we assume for simplicity that $I = \{1\}$ here and omit it from the notation.) We fix an injective character $\zeta : A \hookrightarrow \overline{\rationals}_{\ell}$ as in \S\ref{void-global-covering-group-intro}.

In this context, we may consider the full subcategory
$$
\Sat_{\mathscr G, \zeta}(\Hec_G) \subset \derived_{\mathscr G, \zeta}(\Hec_G)
$$
consisting of objects whose pullback to the affine Grassmannian $\Gr_G$ are perverse and universally locally acyclic (ULA) relative to $X$.

Using several pieces of canonical structure on $\mathscr G_{\Hec_G}$, one may endow $\Sat_{\mathscr G, \zeta}(\Hec_G)$ with a symmetric monoidal structure. Furthermore, one ``tweaks" the commutativity constraint on $\Sat_{\mathscr G, \zeta}(\Hec_G)$ using the sum of positive roots as for the usual Satake category (\emph{cf.}~\cite[\S6]{MR2342692}). This results in a symmetric monoidal category
\begin{equation}
\label{eq-twisted-satake-category-intro}
{}^+\Sat_{\mathscr G, \zeta}(\Hec_G).
\end{equation}
\end{void}

\begin{rem}
The above construction of the twisted Satake category \eqref{eq-twisted-satake-category-intro} is due to Gaitsgory and Lysenko (\emph{cf.}~\cite{MR3769731}), based on Reich's work \cite{MR2956088}. It differs from Finkelberg and Lysenko's approach to the twisted geometric Satake equivalence (\emph{cf.}~\cite{MR2265675, MR2684259, lysenko2014twisted}) by removing factorization line bundles over $\Gr_G$ from the picture.

Since \cite{MR3769731} is more focused on statements than proofs, we shall supply a self-contained construction of the symmetric monoidal category \eqref{eq-twisted-satake-category-intro} in \S\ref{sec-the-satake-category}.
\end{rem}

\begin{void}
Let us now turn to the spectral side.

To the pair $(G, \mu)$, we shall attach another pair $(H, \nu)$, where $H$ is the metaplectic dual group (\emph{cf.}~\S\ref{void-ell-group-intro}) and $\nu$ is a morphism of \'etale sheaves of $\mathbb E_{\infty}$-monoids over $X$:
\begin{equation}
\label{eq-metaplectic-dual-morphism-intro}
\nu : \hat Z_H \rightarrow \deloop^2_X A.
\end{equation}
Here, $\hat Z_H$ is the character group of the center $Z_H$ of $H$, viewed as a constant \'etale sheaf.

The pair $(H, \nu)$ is our version of the \emph{metaplectic dual data}. Contrary to the same-named notion introduced in \cite{MR3769731}, our construction of $(H, \nu)$ is of group-theoretic nature and remains valid over an arbitrary base scheme $S$ over which $A$ has invertible order.\footnote{Our original motivation for giving this construction is that it is uniform for number fields and function fields (\emph{cf.}~\cite{zhao2022metaplectic}). In the present context, it also has the practical consequence of removing an assumption on the characteristic of $\base$ from \cite{MR3769731} (\emph{cf.}~the discussion in \S\ref{sec-reich}).}
\end{void}

\begin{void}
To formulate our version of the twisted geometric Satake equivalence, we need to add to \eqref{eq-metaplectic-dual-morphism-intro} a term having to do with $\vartheta$-characteristics over $X$.

More precisely, we shall define another $\mathbb E_{\infty}$-monoidal morphism
\begin{equation}
\label{eq-metaplectic-dual-morphism-theta-twist-intro}
\nu + \vartheta : \hat Z_H \rightarrow \deloop^2_X A.
\end{equation}
which coincides with $\nu$ whenever a $\vartheta$-characteristic over $X$ is chosen (\emph{cf.}~\S\ref{void-theta-shift}).

This modification is closely related to a phenomenon known in the Langlands--Weissman program: The $L$-group arises from the Baer sum of two extensions of the Galois group by $Z_H(\overline{\rationals}_{\ell})$ (\emph{cf.}~\cite{MR3802418}). The addition of the $\vartheta$-term in \eqref{eq-metaplectic-dual-morphism-theta-twist-intro} corresponds to the twist by Weissman's meta-Galois group (also known as ``the first twist").

In fact, we shall prove a precise result to this effect: Our Theorem \ref{thm-meta-galois-twist} asserts that in the function field context, Weissman's meta-Galois group is the fundamental group of the $\{\pm 1\}$-gerbe of $\vartheta$-characteristics.
\end{void}

\begin{void}
Let us write $\Rep_H$ for the symmetric monoidal category of $H$-representations on $\overline{\rationals}_{\ell}$-local systems over $X$. It is of \'etale local nature over $X$ and admits a $\hat Z_H$-grading defined by the action of $Z_H$.

Using a procedure explained in \S\ref{sec-monoidal-twist}, we may use the $\mathbb E_{\infty}$-monoidal morphism \eqref{eq-metaplectic-dual-morphism-theta-twist-intro} to ``twist" $\Rep_H$, obtaining a new symmetric monoidal category
\begin{equation}
\label{eq-dual-category-intro}
\Rep_{H, \nu + \vartheta}.
\end{equation}

This is a version of the representation category of the $L$-group. More precisely, the integral form of the $L$-group ${}^LH_X$ is defined so that $\Rep({}^LH_X)$ is \emph{monoidally} equivalent to \eqref{eq-dual-category-intro}. However, we caution the reader that this equivalence is generally \emph{incompatible} with the commutativity constraints.

From the perspective of the twisted geometric Satake equivalence, it is the symmetric monoidal category \eqref{eq-dual-category-intro} which appears naturally on the dual side.
\end{void}

\begin{thmx}
\label{thmx-satake-equivalence}
There is a canonical equivalence of symmetric monoidal categories
\begin{equation}
\label{eq-satake-equivalence-intro}
	{}^+\Sat_{\mathscr G, \zeta}(\Hec_G) \simeq \Rep_{H, \nu + \vartheta}.
\end{equation}
\end{thmx}

\begin{void}
A more general version of Theorem \ref{thmx-satake-equivalence} appears as Theorem \ref{thm-satake-equivalence}. In the same subsection, we state various properties of our Satake equivalence which are needed for applications to global function fields.

Let us point out that the statement of Theorem \ref{thm-satake-equivalence} has already appeared in \cite[\S9]{MR3769731}, at least if we assume that our metaplectic dual data coincide with those of \cite{MR3769731}. However, \emph{op.cit.}~does not provide a proof of this statement. By inspecting the existing literature \cite{MR2684259, lysenko2014twisted, MR2956088}, we have identified at least one place where new ideas are needed for the proof: the construction of the fiber functor.\footnote{We provide a detailed discussion on the status of the existing literature in \S\ref{sec-finkelberg-lysenko}-\ref{sec-reich}.}

Let us briefly explain the difficulty involved.
\end{void}

\begin{void}[The trouble with the fiber functor]
\label{void-fiber-functor-trouble}
First, we note that the $A$-gerbe $\mathscr G_{\Hec_G}$ is nontrivial along the fibers over $X$ in general, so there is no analogue of the global cohomology functor in the twisted setting. As a substitute, we consider the constant term functor
\begin{equation}
\label{eq-constant-term-intro}
{}^+\Sat_{\mathscr G, \zeta}(\Hec_G) \rightarrow \Sat_{\mathscr G, \zeta}(\Hec_T),
\end{equation}
where $T$ is the universal Cartan of $G$.\footnote{One needs to choose a Borel subgroup $B \subset G$ to obtain an \'etale level $\mu_T$ for $T$. However, the resulting category $\Sat_{\mathscr G, \zeta}(\Hec_T)$, as well as the functor \eqref{eq-constant-term-intro}, is independent of the choice of $B$.} Let us also take for granted the twisted geometric Satake equivalence for $T$, which is a symmetric monoidal equivalence
\begin{equation}
\label{eq-satake-torus-intro}
\Sat_{\mathscr G, \zeta}(\Hec_T) \simeq \Rep_{T_H, \nu + \vartheta},
\end{equation}
for $T_H$ the canonical maximal torus of $H$.

As Reich points out, $\Rep_{T_H, \nu + \vartheta}$ may \emph{not} admit a fiber functor due to the $(\nu + \vartheta)$-twist, and therefore one cannot directly apply the Tannakian formalism to ${}^+\Sat_{\mathscr G, \zeta}(\Hec_G)$ (\emph{cf.}~\cite[\S V.1]{MR2956088}). However, one may try to resolve the issue as follows: By ``untwisting" both \eqref{eq-constant-term-intro} and the equivalence \eqref{eq-satake-torus-intro}, one obtains the fiber functor as the composition of symmetric monoidal functors
\begin{align*}
	{}^+\Sat_{\mathscr G, \zeta}(\Hec_G)_{-(\nu + \vartheta)} &\rightarrow \Sat_{\mathscr G, \zeta}(\Hec_T)_{-(\nu + \vartheta)} \\
	& \simeq (\Rep_{T_H, \nu + \vartheta})_{-(\nu + \vartheta)} \simeq \Rep_{T_H} \rightarrow \Lis(X),
\end{align*}
where $\Lis(X)$ denotes the category of $\overline{\rationals}_{\ell}$-local systems on $X$. It is then possible to apply the Tannakian formalism to recognize ${}^+\Sat_{\mathscr G, \zeta}(\Hec_G)_{-(\nu + \vartheta)}$ as $\Rep_H$, and ultimately deduce \eqref{eq-satake-equivalence-intro} from this equivalence by restoring the twist.
\end{void}

\begin{void}
To realize this strategy, however, one must construct a $\hat Z_H$-grading on the twisted Satake category ${}^+\Sat_{\mathscr G, \zeta}(\Hec_G)$, compatible with its symmetric monoidal structure.

With respect to the natural maps
$$
\hat T_H \twoheadrightarrow \hat Z_H \rightarrow \pi_1 G,
$$
this $\hat Z_H$-grading is supposed to be a \emph{refinement} of the $\pi_1 G$-grading defined by the connected components of $\Hec_G$, as well as a \emph{coarsening} of the $\hat T_H$-grading that one may obtain from the semisimplicity of the (pointwise) Satake category. Note, however, that the $\hat T_H$-grading is \emph{incompatible} with the symmetric monoidal structure, so it is insufficient for the construcion of the desired $\hat Z_H$-grading.

Let us make this discussion more concrete by considering a classical example.
\end{void}

\begin{eg}
\label{eg-metaplectic-double-cover}
Take $G = \SL_2$ and $A = \{\pm 1\}$, under the assumption $\characteristic \base \neq 2$. As for the \'etale level $\mu$, we shall take the canonical lift of the mod $2$ \'etale Chern class
$$
c_2 \text{ mod }2 \in H^4_{\etale}(\deloop \SL_2, \{\pm 1\}^{\otimes 2}).
$$
The existence of this canonical lift is justified by the calculation of the reduced cohomology of $\deloop\SL_2$ in degrees $\le 4$ (\emph{cf.}~\cite{MR1441006}). It induces the metaplectic double cover of $\SL_2(F)$, for $F$ a local field of characteristic $\neq 2$.

For this choice of $\mu$, the dual group $H$ is identified with $\SL_2$, endowed with its canonical pinning. In particular, we have $\hat T_H \simeq \integers$ and $\hat Z_H \simeq \integers/2$, while $\pi_1 G$ vanishes.

Furthermore, if we fix a $\vartheta$-characteristic over $X$, then $\Rep_{T_H, \nu + \vartheta}$ is canonically identified with the symmetric monoidal category of $\integers$-graded $\overline{\rationals}_{\ell}$-local systems over $X$, whose commutativity constraint is given by \emph{the sign rule}. In particular, the functor of summing up its $\integers$-graded components is \emph{not} symmetric monoidal. The discussion of \S\ref{void-fiber-functor-trouble} is thus necessary already in this classical example.
\end{eg}

\begin{void}[The $\hat Z_H$-grading]
\label{void-center-grading-intro}
Let us now explain the geometric origin of the $\hat Z_H$-grading on the twisted Satake category, established in \S\ref{sec-grading-by-virtual-connected-components} of the present article.\footnote{An earlier version of this article contains a less direct construction of the $\hat Z_H$-grading. The present construction makes use of the ``canonical quadratic structure" established more recently in \cite{shi2025extendedpureinnerforms}.}

To do so, we shall view objects of ${}^+\Sat_{\mathscr G, \zeta}(\Hec_G)$ as $(\mathscr G_{\Gr_G}, \zeta)$-twisted perverse sheaves over $\Gr_G$ (relative to $X$) equivariant with respect to the arc group $L^+G$. Here, $\mathscr G_{\Gr_G}$ is the pullback of $\mathscr G_{\Hec_G}$ to $\Gr_G$, viewed as an $L^+G$-equivariant $A$-gerbe.

The key observation is that on each connected component $\Gr_G^{\theta}$ of $\Gr_G$, corresponding to $\theta\in\pi_1 G$, the $A$-gerbe $\mathscr G_{\Gr_G}$ admits \emph{a family} of $L^+G_{\adjoint}$-equivariance structures parametrized by lifts of $\theta$ to $\hat Z_H$, for $G_{\adjoint}$ the adjoint form of $G$. Distinct lifts of $\theta$ \emph{always} give rise to distinct $L^+G_{\adjoint}$-equivariance structures.

The $\hat Z_H$-grading on ${}^+\Sat_{\mathscr G, \zeta}(\Hec_G)$ is now defined as follows: For any $\xi \in \hat Z_H$ with image $\theta$ in $\pi_1 G$, the $\xi$-graded component of ${}^+\Sat_{\mathscr G, \zeta}(\Hec_G)$ consists of $L^+G_{\adjoint}$-equivariant $(\mathscr G_{\Gr_G}, \zeta)$-twisted perverse sheaves over $\Gr_G^{\theta}$, with respect to the $L^+G_{\adjoint}$-equivariance structure of $\mathscr G_{\Gr_G}$ \emph{corresponding to $\xi$}. It is straightforward to verify that the $\hat Z_H$-grading on ${}^+\Sat_{\mathscr G, \zeta}(\Hec_G)$ obtained this way is compatible with the symmetric monoidal structure.

The fact that, under the Satake equivalence \eqref{eq-satake-equivalence-intro}, the natural $\hat Z_H$-grading on $\Rep_{H, \nu + \vartheta}$ corresponds to this $\hat Z_H$-grading on ${}^+\Sat_{\mathscr G, \zeta}(\Hec_G)$ seems interesting in itself.
\end{void}

\begin{void}
A few other components of our proof of the twisted geometric Satake equivalence have not appeared in the existing literature. They are in:
\begin{enumerate}
	\item the proof of semisimplicity of the pointwise twisted Satake category;
	\item the construction of the equivalence for tori;
	\item the canonical identification of the Tannaka dual group with $H$.
\end{enumerate}

Let us discuss them in turn.
\end{void}

\begin{void}[Semisimplicity]
Point (1) is the assertion that the twisted Satake category associated to a geometric point $x$ of $X$ is semisimple (\emph{cf.}~Corollary \ref{cor-semisimplicity}).

As in the untwisted setting, this assertion follows from the parity vanishing of intersection cohomology sheaves on Schubert varieties (\emph{cf.}~\cite[proof of Proposition 1 \& \S A.7]{MR1826370}). However, the proof of parity vanishing there does \emph{not} apply in the twisted setting.\footnote{In an earlier version of this paper, it is \emph{incorrectly} asserted that the semisimplicity of the twisted Satake category can be proved as in the untwisted case. This error is pointed out by Gaitsgory.}

That being said, at least for $G$ simple and simply connected, the required parity vanishing is known by Lusztig's work \cite[\S5]{MR1276910}, predating all works on the twisted geometric Satake equivalence. In \S\ref{sec-pointwise-studies} of this article, we establish the parity vanishing following the argument of \cite[Proposition 3.12]{MR4108915}.

I thank Gurbir Dhillon and Michael Finkelberg for the references \cite{MR1276910, MR4108915}. It must also be mentioned that the recent preprint \cite{dhillon2025endoscopymetaplecticaffinehecke} establishes \emph{significantly stronger} results concerning $\derived_{\mathscr G, \zeta}(\Hec_{G, x})$ than those in this article.

Another proof of semisimplicity, avoiding parity vanishing altogether, is suggested by Gaitsgory and reproduced in \S\ref{sec-gaitsgory-semisimplicity} of this article.
\end{void}

\begin{void}[Equivalence for tori]
Unlike \cite{MR3769731}, our metaplectic dual data are \emph{not} defined using the factorization structure of the affine Grassmannian. As a result, the twisted geometric Satake equivalence for tori is less tautological for us.

We perform the necessary calculation in \S\ref{sec-tori}. The results there are in turn used to relate Weissman's meta-Galois group to $\vartheta$-characteristics (\emph{cf.}~Theorem \ref{thm-meta-galois-twist}).
\end{void}

\begin{void}[Canonicity of $H$]
We not only compute the root data of the Tannaka dual group in this text, but also construct its isomorphism with $H$. This is needed for treating covers of not necessarily split reductive groups.

In the untwisted setting, such an isomorphism is constructed in \cite[VI.11]{fargues2021geometrization}. There are some intricacies involved in the twisted setting: Essentially, if we view \eqref{eq-metaplectic-dual-morphism-intro} as an $\mathbb E_{\infty}$-monoidal morphism $\hat T_H \rightarrow \deloop^2_X A$ trivialized over the root lattice of $H$, then the canonicity of \eqref{eq-satake-equivalence-intro} dictates which trivialization we must take. In this text, we specify the trivialization which makes \eqref{eq-satake-equivalence-intro} canonical.
\end{void}

\subsection{Convention}

\begin{void}
Let us note some convention and notation used throughout this article.
\end{void}

\begin{void}
\label{void-linear-spaces}
Denote by $\Spc$ the $\infty$-category of $\infty$-groupoids, whose objects are referred to as \emph{spaces}. Denote by $\Spc_*$ the $\infty$-category of pointed spaces.

Denote by $\Sptr$ the $\infty$-category of spectra and by $\Omega^{\infty} : \Sptr \rightarrow \Spc_*$ the forgetful functor. For any (classical) ring $R$, we write $R\Mod$ for the $\infty$-category of $HR$-module spectra.

Since the restriction of $\Omega^{\infty}$ to \emph{connective} $HR$-module spectra is conservative, we may view a connective $HR$-module spectrum as a space with additional structure, and refer to it simply as an \emph{$R$-linear space}. Similarly, we often call a morphim of connective $HR$-module spectra \emph{$R$-linear} to distinguish it from a morphism of the underlying spaces.
\end{void}

\begin{void}
\label{void-bar-construction}
Given a site $\mathscr C$ and an $\infty$-category $\mathscr O$ admitting finite limits, we write $\Shv(\mathscr C, \mathscr O)$ for the $\infty$-category of $\mathscr O$-valued sheaves over $\mathscr C$.

We shall extensively use the iterated Bar construction for the $\infty$-category $\Shv(\mathscr C, \Spc)$ endowed with the Cartesian symmetric monoidal structure (\emph{cf.}~\cite[\S 5.2.3]{lurie2017higher}). For any integer $n\ge 0$, we view the $n$th iterated Bar construction as a functor
$$
\deloop^n : \Shv(\mathscr C, \Alg_{\mathbb E_n}(\Spc)) \rightarrow \Shv(\mathscr C, \Spc_*).
$$
We shall also refer to $\deloop$ as the \emph{deloop} functor.

For a $\Sptr$-valued sheaf $\mathscr A$ over $\mathscr C$, there is a canonical isomorphism in $\Shv(\mathscr C, \Spc_*)$:
$$
\deloop^n \mathscr A \simeq \Omega^{\infty}(\mathscr A[n]).
$$
In other words, the value of $\deloop^n \mathscr A$ at any $c\in\mathscr C$ coincides with $\Omega^{\infty}\Gamma(c, \mathscr A[n])$.
\end{void}

\begin{void}
Given a site $\mathscr C$ and a sheaf of abelian groups $\mathscr A$ over $\mathscr C$, we refer to a section of $\deloop^2\mathscr A$ at $c\in\mathscr C$ as an \emph{$\mathscr A$-gerbe} over $c$.

Contrary to most literature on banded gerbes, we denote the symmetric monoidal product on $A$-gerbes \emph{additively}. We find this convention more natural in view of its relation to the $\integers$-linear structure on $\deloop^2\mathscr A$.
\end{void}

\begin{void}
Denote by $\Ring$ the category of (unital, commutative) rings. A \emph{prestack} is a functor of $\infty$-categories $\Ring \rightarrow \Spc$.

By an \emph{indscheme}, we shall mean a prestack $Z$ which admits a presentation $Z \simeq \colim_{\alpha} Z_{\alpha}$ as a filtered colimit, where each $Z_{\alpha}$ is a quasi-compact quasi-separated scheme and each transition map $Z_{\alpha} \rightarrow Z_{\alpha'}$ is a closed immersion.

Given a ring $R$ (respectively a scheme $S$), an \emph{$R$-prestack} is a prestack over $\Spec R$ (respectively over $S$). We employ the same convention for \emph{$R$-schemes}, \emph{$R$-indschemes}, \emph{etc}.
\end{void}

\subsection{Acknowledgements}

\begin{void}
This article, as well as its prequel \cite{zhao2022metaplectic}, grew out of my attempt to fill in the proof of the twisted geometric Satake equivalence formulated in \cite{MR3769731}. I benefited tremedously from conversations with Dennis Gaitsgory on this topic.

In addition, I owe a large part of my understanding of covering groups to conversations with Wee Teck Gan, Sergey Lysenko, and Sam Raskin. I also thank Claudius Heyer for illuminating conversations during the preparation of this article.

I thank K\c{e}stutis \v{C}esnavi\v{c}ius, Aron Heleodoro, and Cong Xue for offering me the opportunities to present parts of this work.
\end{void}

\begin{void}
Since the first version of this article appeared, several errors were found and improvements were discovered, both by myself and by the community. This article is essentially rewritten to take all of them into account.

In this process, I received help from many mathematicians: Lin Chen, Gurbir Dhillon, Tianyi Feng, Tony Feng, Dennis Gaitsgory, Naoki Imai, and Luozi Shi. I am deeply grateful to all of them.
\end{void}

\newpage

\part{Geometric Satake equivalence}

\section{The Satake category}
\label{sec-the-satake-category}

In this section, we construct the twisted Satake category. We begin by constructing a gerbe $\mathscr G_{\Hec_G}$ over the local Hecke stack $\Hec_G$ in \S\ref{sec-local-hecke-stack-gerbe}. This gerbe has a natural factorization structure compatible with the convolution structure of $\Hec_G$. This is precisely formulated and constructed in \S\ref{sec-factorization}. Using $\mathscr G_{\Hec_G}$, we construct the twisted Satake category in \S\ref{sec-satake-category}, subject to some technical results whose proofs are supplied in \S\ref{sec-outer-convolution-diagram}. In \S\ref{sec-constant-term-functors}, we define constant term functors and study their basic properties.

Our construction of the twisted Satake category is largely parallel to \cite{MR3769731}, except that we offer a more formal treatment of factorization.

\subsection{The local $A$-gerbe}
\label{sec-local-hecke-stack-gerbe}

\begin{void}[Context]
\label{void-local-hecke-stack-context}
Let $\base$ be a field and $X$ be a smooth curve over $\base$, \emph{i.e.}~$X$ is a quasi-compact and separated $\base$-scheme which is smooth of relative dimension $1$.

Let $G$ be a smooth affine group $X$-scheme. Let $A$ be a finite abelian group whose order is invertible in $\base$. Let $\mu$ be an $A$-valued \emph{\'etale level} of $G$, \emph{i.e.}~a morphism of pointed \'etale $X$-stacks:
\begin{equation}
\label{eq-etale-level}
\mu : \deloop_X G \rightarrow \deloop^4_X A(1),
\end{equation}
where $\deloop_X$ denotes the Bar construction for sheaves on the big \'etale site of $X$ (\emph{cf.}~\S\ref{void-bar-construction}).
\end{void}

\begin{void}
\label{void-ran-space-definition}
Denote by $\Ran$ the $\base$-prestack sending a $\base$-algebra $R$ to the set of finite subsets of $X(R)$. It is called the \emph{Ran space} of $X$. For an $R$-point $\underline x$ of $\Ran$, we write $D_{\underline x}$ for the formal completion of $X_R := \Spec R \times X$ along the union of the graphs $\Gamma_{\underline x} := \bigcup \Gamma_{x^i}$ over $x^i \in \underline x$.

We shall view $D_{\underline x}$ as an affine scheme and refer to it as the \emph{formal disk} around $\underline x$. The \emph{formal punctured disk} around $\underline x$ is defined to be the affine scheme $\mathring D_{\underline x} := D_{\underline x} \setminus \Gamma_{\underline x}$.

%For a finite set $I$ and an $R$-point $x^I = (x^i)_{i\in I}$ of $X^I$, we also write $\Gamma_{x^I}$, $D_{x^I}$, and $\mathring D_{x^I}$ for the affine schemes defined for $\underline x := \bigcup_{i\in I} x^i$ as a finite subset of $X(R)$.
\end{void}

\begin{void}
\label{void-local-hecke-stack-definition}
Denote by $L^+G$ (respectively $L^+G$) the group $\Ran$-prestack whose $R$-points are pairs $(\underline x, g)$, where $\underline x$ is an $R$-point of $\Ran$ and $g : D_{\underline x} \rightarrow G$ (respectively $g : \mathring D_{\underline x} \rightarrow G$) is a morphism of affine $X$-schemes.

Denote by $\Hec_G$ the $\Ran$-prestack whose $R$-points are quadruples $(\underline x, P^0, P^1, \alpha)$, where $\underline x$ is an $R$-point of $\Ran$, $P^0$, $P^1$ are $G$-bundles over $D_{\underline x}$, and $\alpha$ is an isomorphism of their restrictions to $\mathring D_{\underline x}$. We shall represent an $R$-point of $\Hec_G$ by the expression
\begin{equation}
\label{eq-hecke-modification-expression}
P^0 \overset{\underline x}{\sim} P^1,
\end{equation}
which we refer to as a \emph{modification} of $G$-bundles over $D_{\underline x}$ along $\underline x$.

Denote by $\Gr_G$ the $\Ran$-prestack whose $R$-points are modifications \eqref{eq-hecke-modification-expression}, where $P^0$ is the trivial $G$-bundle.

The $\Ran$-prestacks $L^+ G$, $L G$, $\Hec_G$, $\Gr_G$ are called the \emph{arc group}, the \emph{loop group}, the \emph{local Hecke stack}, and the \emph{affine Grassmannian}, respectively. For any morphism $\underline x : S \rightarrow \Ran$ of $\base$-prestacks, we write $L_{\underline x}^+G$, $L_{\underline x}G$, $\Hec_{G, \underline x}$, $\Gr_{G, \underline x}$ for their base changes along $\underline x$.
\end{void}

\begin{void}
We recall the well-known representability statements for the above $\Ran$-prestacks, referring to \cite[\S3.1]{MR3752460} for details.

The arc group $L^+G$ is an affine group $\Ran$-scheme (\emph{i.e.}~the structural morphism $L^+G \rightarrow \Ran$ is representable in affine group schemes). The loop group $LG$ is an ind-affine group $\Ran$-indscheme.

The affine Grassmannian $\Gr_G$ is a $\Ran$-indscheme of ind-finite type, and the map $LG \rightarrow \Gr_G$, sending $(\underline x, g)$ to the modification of the trivial $G$-bundle given by multiplication by $g$, realizes $\Gr_G$ as a quotient of \'etale $\Ran$-stacks\footnote{As above, this means that for a fixed $R$-point $\underline x$ of $\Ran$, the base change $\Gr_{G, \underline x}$ is the quotient $L_{\underline x}G/L^+_{\underline x}G$ of \'etale $R$-stacks. Note that $\Ran$ itself is \emph{not} an \'etale stack, \emph{cf.}~\cite[Warning 2.4.4]{lurie2014weil}.}:
$$
L G / L^+G \simeq \Gr_G.
$$

Likewise, $\Hec_G$ is identified as a quotient of \'etale $\Ran$-stacks:
\begin{equation}
\label{eq-local-hecke-stack-as-quotient}
\Hec_G \simeq L^+G \backslash LG / L^+G.
\end{equation}
\end{void}

\begin{rem}
\label{rem-hecke-stack-ind-presentation}
Although $\Hec_G$ is not of ind-finite type over $\Ran$, it can be approximated by finite type algebraic stacks as follows: For an $R$-point $\underline x$ of $\Ran$, there is an indscheme presentation $\Gr_{G, \underline x} \simeq \colim_{\alpha} Z_{\alpha}$, where each $Z_{\alpha}$ is an $L^+_{\underline x}G$-stable $R$-scheme of finite type, on which the $L_{\underline x}^+G$-action factors through a finite type quotient $H_{\beta}$, so we have
\begin{equation}
\label{eq-hecke-stack-ind-presentation}
\Hec_{G, \underline x} \simeq \colim_{\alpha} \lim_{\beta \gg \alpha} H_{\beta} \backslash Z_{\alpha},
\end{equation}
where $\beta\gg \alpha$ means that $\beta$ ranges over a cofinal subset of indices depending on $\alpha$. In other words, $\Hec_{G, \underline x}$ is an ind-pro-algebraic stack.
\end{rem}

\begin{void}[Convolution]
\label{void-local-hecke-stack-convolution}
The $\Ran$-prestack $\Hec_G$ occurs as the $1$-simplices of a groupoid $\Ran$-prestack $\Hec_G^{[\cdot]}$. Indeed, for each integer $n\ge 0$, we define an $R$-point of $\Hec_G^{[n]}$ to be a chain of modifications
\begin{equation}
\label{eq-hecke-modification-chain-expression}
P^0 \overset{\underline x}{\sim} P^1 \overset{\underline x}{\sim} \cdots \overset{\underline x}{\sim} P^n.
\end{equation}
Given a morphism $f : [m] \rightarrow [n]$ in the simplicial category, the induced morphism $\Hec_G^{[n]} \rightarrow \Hec_G^{[m]}$ sends \eqref{eq-hecke-modification-chain-expression} to the chain of modifications
$$
P^{f^{-1}(0)} \overset{\underline x}{\sim} P^{f^{-1}(1)} \overset{\underline x}{\sim} \cdots \overset{\underline x}{\sim} P^{f^{-1}(n)},
$$
where each modification is the composition of the corresponding segment in \eqref{eq-hecke-modification-chain-expression}. The resulting simplicial $\Ran$-prestack $\Hec_G^{[\cdot]}$ is a groupoid (\emph{i.e.}~it obeys Segal's axioms, \emph{cf.}~\cite[Definition 6.1.2.7]{MR2522659}).

We refer to this groupoid structure on $\Hec_G$ as the \emph{convolution structure}. Note that $\Hec_G^{[\cdot]}$ is the \v{C}ech nerve of $L^+(\deloop G) \rightarrow L(\deloop G)$, where $L^+(\deloop G)$ (respectively $L(\deloop G)$) parametrizes morphisms $D_{\underline x} \rightarrow \deloop G$ (respectively $\mathring D_{\underline x} \rightarrow \deloop G$) with $\underline x\in\Ran(R)$.
\end{void}

\begin{void}[The $A$-gerbe $\mathscr G_{\Hec_G}$]
\label{void-local-hecke-stack-gerbe}
The goal of this subsection is to construct an $A$-gerbe $\mathscr G_{\Hec_G}$ over $\Hec_G$ using the \'etale level $\mu$, \emph{cf.}~\eqref{eq-etale-level}.

Indeed, in \S\ref{void-local-integration-map-construction} below, we shall construct a morphism of spaces
\begin{equation}
\label{eq-local-integration-map}
\int_{\mathring D} : \Maps_*(\deloop_X G, \deloop^4_XA(1)) \rightarrow \Maps(\Hec_G, \deloop^2 A),
\end{equation}
where $\Maps_*(\cdot, \cdot)$ denotes the mapping space of pointed $X$-stacks. Then we will set $\mathscr G_{\Hec_G}$ to be the image of $\mu$ under \eqref{eq-local-integration-map}.

This requires the formalism of trace maps, which we now turn to.
\end{void}

\begin{void}[The local trace map]
\label{void-trace-map}
Given an $R$-point $\underline x$ of $\Ran$, we shall construct a morphism in the pro-category of $H\integers$-module spectra
\begin{equation}
\label{eq-trace-map}
\tr_{\underline x} : \Gamma(D_{\underline x} \text{ mod }\mathring D_{\underline x}, \hat{\integers}(1)[2]) \rightarrow \Gamma(\Spec R, \hat{\integers}),
\end{equation}
where $\hat{\integers}$ is the pro-(abelian group) $``\lim_{\characteristic\base \nmid n}"\integers/n$ and the notation $\Gamma(\mathscr Y\text{ mod }\mathscr Z, \cdot)$, for any morphism $\mathscr Z \rightarrow \mathscr Y$ of $\base$-prestacks, means the fiber of the pullback map $\Gamma(\mathscr Y, \cdot) \rightarrow \Gamma(\mathscr Z, \cdot)$. We shall refer to \eqref{eq-trace-map} as the \emph{trace map}.

Indeed, write $\iota : \Gamma_{\underline x} \rightarrow X_R$ for the closed immersion. By the formal base change theorem (\emph{cf.}~\cite[Corollary 6.6.4]{MR1360610}, \cite[Theorem 6.11]{MR4278670}), the source of \eqref{eq-trace-map} is identified as
\begin{align}
\notag
\Gamma(D_{\underline x} \text{ mod }\mathring D_{\underline x}, \hat{\integers}(1)[2]) & \simeq \Gamma(X_R \text{ mod }X_R\setminus\Gamma_{\underline x}, \hat{\integers}(1)[2]) \\
\label{eq-trace-map-formal-base-change}
& \simeq \Gamma(\Gamma_{\underline x}, \iota^! \hat{\integers}(1)[2])
\end{align}
via pullback along $D_{\underline x} \rightarrow X_R$. Since the structural morphism $X_R \rightarrow \Spec R$ is smooth of relative dimension $1$, we may identify $\iota^! \hat{\integers}(1)[2]$ with $\pi^!\hat{\integers}$, for $\pi : \Gamma_{\underline x} \rightarrow \Spec R$ the projection. We then have a morphism
\begin{align}
\notag
\Gamma(\Gamma_{\underline x}, \iota^! \hat{\integers}(1)[2]) &\simeq \Gamma(\Gamma_{\underline x}, \pi^! \hat{\integers}) \\
\label{eq-trace-map-adjunction}
& \simeq \Gamma(\Spec R, \pi_* \pi^! \hat{\integers}) \rightarrow \Gamma(\Spec R, \hat{\integers}),
\end{align}
where the last morphism is the co-unit of the adjunction from properness of $\pi$. The desired morphism \eqref{eq-trace-map} is the composition of \eqref{eq-trace-map-formal-base-change} and \eqref{eq-trace-map-adjunction}.
\end{void}

\begin{rem}
\label{rem-trace-map-projection-formula}
In the construction of \eqref{eq-trace-map}, we worked with coefficients in $\hat{\integers}$. We could instead perform the construction for any complex $\mathscr A$ of torsion \'etale sheaves of invertible order over $\base$ (the case of interest being $\mathscr A = A[2]$) and obtain
\begin{equation}
\label{eq-trace-map-with-coefficients}
\tr_{\underline x} : \Gamma(D_{\underline x}\text{ mod }\mathring D_{\underline x}, \mathscr A(1)[2]) \rightarrow \Gamma(\Spec R, \mathscr A).
\end{equation}

Note that \eqref{eq-trace-map} and \eqref{eq-trace-map-with-coefficients} are related by a version of the projection formula, asserting that the diagram below commutes
$$
\begin{tikzcd}[column sep = 1em]
	\Gamma(\Spec R, \mathscr A) \otimes \Gamma(D_{\underline x}\text{ mod }\mathring D_{\underline x}, \hat{\integers}(1)[2]) \ar[r, "\otimes"]\ar[d, "\id\otimes \tr_{\underline x}"] & \Gamma(D_{\underline x}\text{ mod }\mathring D_{\underline x}, \mathscr A(1)[2]) \ar[d, "\tr_{\underline x}"] \\
	\Gamma(\Spec R, \mathscr A) \otimes \Gamma(\Spec R, \hat{\integers}) \ar[r, "\otimes"] & \Gamma(\Spec R, \mathscr A)
\end{tikzcd}
$$
where the horizontal arrows are defined by multiplication of coefficients. This holds because the top horizontal arrow is induced from the canonical morphism
$$
\pi^*\mathscr A \otimes \pi^!\hat{\integers} \rightarrow \pi^!\mathscr A,
$$
which corresponds, under adjunction, to the composition
$$
\pi_*(\pi^*\mathscr A \otimes \pi^!\hat{\integers}) \simeq \mathscr A \otimes \pi_*\pi^!\hat{\integers} \rightarrow \mathscr A \otimes \hat{\integers} \xrightarrow{\otimes} \mathscr A.
$$
\end{rem}

\begin{rem}
\label{rem-kummer-class-of-tautological-line-bundle}
Consider the case where $\underline x$ is defined by an $R$-point $x$ of $X$ and write $D_x$ instead of $D_{\underline x}$, \emph{etc}. Denote by
$$
\Psi : \deloop \mathbb G_m \rightarrow \deloop^2 \hat{\integers}(1)
$$
the Kummer morphism, \emph{i.e.}~the deloop of the compatible family of boundary maps $\mathbb G_m \rightarrow \deloop(\integers/n(1))$ given by the degree-$n$ covers with $\characteristic\base \nmid n$.

Consider the line bundle $\mathscr O(x)$ over $D_x$, equipped with its canonical trivialization over $\mathring D_x$. Then $\Psi(\mathscr O(x))$ defines an object of (the underlying space of) $\Gamma(D_{\underline x}\text{ mod }\mathring D_{\underline x}, \hat{\integers}(1)[2])$. Its image under $\tr_x$ is the constant section $1$.

Indeed, this can be checked over $\base$-points of $R$, where it follows from cohomological purity of the closed immersion $\iota : \Spec \base \simeq \Gamma_x \rightarrow X$.
\end{rem}

\begin{void}[Construction of \eqref{eq-local-integration-map}]
\label{void-local-integration-map-construction}
Using the cohomological interpretation of $\deloop^n A$ (\emph{cf.}~\S\ref{void-bar-construction}), it suffices to construct a map of $H\integers$-module spectra
\begin{equation}
\label{eq-local-integration-map-linear}
\Gamma(\deloop_X G \text{ mod } X, A(1)[4]) \rightarrow \Gamma(\Hec_G, A[2]).
\end{equation}

To construct \eqref{eq-local-integration-map-linear}, we proceed as follows: Let $P^0\overset{\underline x}{\sim} P^1$ be an $R$-point of $\Hec_G$, which we view as two morphisms $D_{\underline x} \rightarrow \deloop_X G$ with an identification over $\mathring D_{\underline x}$. Pulling back along them and subtracting, we obtain a morphism of $H\integers$-module spectra
\begin{equation}
\label{eq-local-integration-map-linear-pullback-subtraction}
 (P^1)^* - (P^0)^* : \Gamma(\deloop_X G \text{ mod } X, A(1)[4]) \rightarrow \Gamma(D_{\underline x} \text{ mod }\mathring D_{\underline x}, A(1)[4]).
\end{equation}

The desired morphism \eqref{eq-local-integration-map-linear} is obtained as the composition of \eqref{eq-local-integration-map-linear-pullback-subtraction} with the trace map \eqref{eq-trace-map-with-coefficients} (for $\mathscr A := A[2]$), using naturality of the construction in the $R$-point $P^0\overset{\underline x}{\sim} P^1$.

As explained in \S\ref{void-local-hecke-stack-gerbe}, this concludes the construction of $\mathscr G_{\Hec_G}$.
\end{void}

\begin{rem}
The construction of $\mathscr G_{\Hec_G}$ also appears in \cite[\S7.3]{MR3769731}. As explained in \emph{loc.cit.}, this construction may be used to prove that every factorization $A$-gerbe on $\Gr_G$ canonically descends to $\Hec_G$, a result claimed in \cite[Theorem III.2.10]{MR2956088} but not justified adequately. On the other hand, our approach to the twisted geometric Satake equivalence does \emph{not} need this result.
\end{rem}

\begin{void}
\label{void-local-hecke-stack-gerbe-multiplicative}
Next, we shall explain in what sense $\mathscr G_{\Hec_G}$ is ``multiplicative" with respect to the convolution structure on $\Hec_G$ (\emph{cf.}~\S\ref{void-local-hecke-stack-convolution}).

Indeed, the group structure on $\deloop^2 A$ realizes it as the $1$-simplices of a groupoid prestack $\deloop^2 A^{[\cdot]}$. By pullback, we obtain a groupoid $\Ran$-prestack $\deloop^2_{\Ran} A^{[\cdot]}$. We define a \emph{multiplicative $A$-gerbe} over $\Hec_G$ to be a morphism of groupoid $\Ran$-prestacks
$$
\Hec_G^{[\cdot]} \rightarrow \deloop^2_{\Ran} A^{[\cdot]}.
$$
\end{void}

\begin{lem}
\label{lem-local-hecke-stack-gerbe-multiplicative}
The $A$-gerbe $\mathscr G_{\Hec_G}$ canonically lifts to a multiplicative $A$-gerbe.
\end{lem}

\begin{proof}
Given an $R$-point $\underline x$ of $\Ran$, the groupoid $\Hec_G(R)\times_{\Ran(R)} \{\underline x\}$ of lifts of $\underline x$ to $\Hec_G$ is the \v{C}ech nerve of the morphism in $\Spc$ given by pullback:
\begin{equation}
\label{eq-hecke-groupoid-point-as-cech-nerve}
\Maps(D_{\underline x}, \deloop G) \rightarrow \Maps(\mathring D_{\underline x}, \deloop G).
\end{equation}
The construction of \S\ref{void-local-integration-map-construction} supplies a morphism in $\Spc$:
\begin{equation}
\label{eq-hecke-groupoid-multiplicative-gerbe-plain}
\Hec_G(R) \times_{\Ran(R)} \{\underline x\} \rightarrow \Gamma(\Spec R, \deloop^2 A)
\end{equation}
It remains to lift \eqref{eq-hecke-groupoid-multiplicative-gerbe-plain} to a morphism of groupoids naturally in $\underline x$.

The \'etale level $\mu$ induces a morphism from \eqref{eq-hecke-groupoid-point-as-cech-nerve} to
\begin{equation}
\label{eq-hecke-groupoid-point-etale-level-as-cech-nerve}
\Maps(D_{\underline x}, \deloop^4A(1)) \rightarrow \Maps(\mathring D_{\underline x}, \deloop^4A(1)).
\end{equation}
The morphism \eqref{eq-hecke-groupoid-point-etale-level-as-cech-nerve} naturally lifts to a morphism of $H\integers$-module spectra, thanks to the $\integers$-linear structure on $\deloop^4 A(1)$ (\emph{cf.}~\S\ref{void-bar-construction}).

Note that given any morphism $f : a \rightarrow b$ in a stable $\infty$-category $\mathscr C$, there is a canonical morphism from the \v{C}ech nerve of $f$ to $\Fib(f)$ as groupoid objects in $\mathscr C$: It is induced from the canonical map from $f$ to $0 \rightarrow \Cofib(f)$. Note, furthermore, that on $1$-simplices, this morphism specializes to the subtraction map
$$
a\times_b a \rightarrow \Fib(f).
$$

The \v{C}ech nerve of \eqref{eq-hecke-groupoid-point-etale-level-as-cech-nerve} thus maps to $\Gamma(D_{\underline x}\text{ mod }\mathring D_{\underline x}, \deloop^4A(1))$, as groupoid objects in $H\integers$-module spectra. Composing with $\mu$, we see that the morphism
\begin{align}
\label{eq-hecke-groupoid-maps-to-linear-groupoid}
\Hec_G(R) \times_{\Ran(R)} \{\underline x\} &\rightarrow \Gamma(D_{\underline x}\text{ mod }\mathring D_{\underline x}, \deloop^4A(1)) \\
\notag
(P^0\overset{\underline x}{\sim} P^1) & \mapsto \mu(P^1) - \mu(P^0)
\end{align}
lifts to a morphism of groupoid objects in $\Spc$. The morphism \eqref{eq-hecke-groupoid-multiplicative-gerbe-plain} is the composition of \eqref{eq-hecke-groupoid-maps-to-linear-groupoid} with the trace map \eqref{eq-trace-map-with-coefficients}, the latter being a morphism of groupoid objects as it comes from a morphism of $H\integers$-module spectra.
\end{proof}

\subsection{Factorization}
\label{sec-factorization}

\begin{void}
In \S\ref{sec-local-hecke-stack-gerbe}, we have constructed a $\base$-prestack $\Hec_G$ together with an $A$-gerbe $\mathscr G_{\Hec_G}$. The pair $(\Hec_G, \mathscr G_{\Hec_G})$ furthermore admits a multiplicative structure (\emph{cf.}~Lemma \ref{lem-local-hecke-stack-gerbe-multiplicative}).

In this subsection, we shall endow $(\Hec_G, \mathscr G_{\Hec_G})$ with an additional structure, which is a precise formulation of its ``unital factorization structure" compatible with the given multiplicative structure (\emph{cf.}~Proposition \ref{prop-local-hecke-stack-gerbe-factorization}).
\end{void}

\begin{void}
\label{void-ran-operad}
We begin by lifting $\Ran$ to a $\base$-presheaf of colored operads: For a $\base$-algebra $R$, the set $\Ran(R)$ admits the structure of a colored operad, with operations defined by
$$
\Maps_{\Ran(R)}(\{\underline x{}_i\}_{i\in I}, \underline y) :=
\begin{cases}
	* & \text{$\underline x{}_i$ and $\underline x{}_j$ are disjoint if $i\neq j$ and $\bigsqcup_{i\in I} \underline x_i \subset \underline y$} \\
	\emptyset & \text{otherwise}
\end{cases}
$$
Here, two $R$-points $\underline x$, $\underline y$ of $\Ran$ are called \emph{disjoint} if $\Gamma_{\underline x} \cap \Gamma_{\underline y} = \emptyset$ as closed subsets of $X_R$.

For any $\base$-presheaf of symmetric monoidal $\infty$-categories $\mathscr O$, we define the $\infty$-category of \emph{$\Ran$-algebras} in $\mathscr O$ as the mapping space
$$
\Alg_{\Ran}(\mathscr O) := \Maps(\Ran, \mathscr O)
$$
taken in the $\infty$-category of $\base$-presheaves valued in the $\infty$-category of $\infty$-operads. (We view symmetric monoidal $\infty$-categories as $\infty$-operads via the tautological forgetful functor.)
\end{void}

\begin{rem}
\label{rem-ran-algebra-factorization-morphisms}
Informally, a $\Ran$-algebra in $\mathscr O$ is a compatible family of $\Ran(R)$-algebras $\mathscr A_R$ in $\mathscr O(R)$, for each $\base$-algebra $R$.

The datum of a $\Ran(R)$-algebra in $\mathscr O(R)$ consists of an object $\mathscr A_{\underline x} \in \mathscr O(R)$ for each $\underline x \in \Ran(R)$, along with structural morphisms
\begin{equation}
\label{eq-factorization-algebra-morphisms}
f_{\{\underline x{}_i\} \rightarrow y} : \bigotimes_{i\in I} \mathscr A_{\underline x{}_i} \rightarrow \mathscr A_{\underline y}
\end{equation}
for any finite pairwise disjoint set $\underline x{}_i$ ($i\in I$) of elements of $\Ran(R)$ contained in $\underline y$. Furthermore, the maps \eqref{eq-factorization-algebra-morphisms} admit homotopy coherent data with respect to compositions of such containments.
\end{rem}

\begin{rem}
\label{rem-ran-algebra-truncated}
For the purpose of this article, the $\Ran$-algebras we shall consider take values in the full subcategory of $1$-truncated objects of $\mathscr O(R)$ for every $\base$-algebra $R$.

In this case, homotopy coherence for the structural morphisms \eqref{eq-factorization-algebra-morphisms} can be stated explicitly: For a composition of operations $\{\underline y{}_j\}_{j\in J} \rightarrow \underline z$ and $\{\underline x{}_i\}_{i\in I_j} \rightarrow \underline y{}_j$ ($j\in J$), there is $2$-isomorphism rendering the diagram
$$
\begin{tikzcd}[column sep = 1.5em]
\bigotimes_{j\in J} (\bigotimes_{i\in I_j} \mathscr A_{\underline x{}_i}) \ar[d, "\bigotimes_{j\in J} f_{\{\underline x{}_i\}\rightarrow\underline y{}_j}"]\ar[r, "\simeq"] & \bigotimes_{\substack{j\in J \\ i\in I_j}} \mathscr A_{\underline x{}_i}  \ar[d, "f_{\{\underline x{}_i\} \rightarrow \underline z}"] \\
\bigotimes_{j\in J} \mathscr A_{\underline y{}_j} \ar[r, "f_{\{\underline y{}_j\}\rightarrow \underline z}"] & \mathscr A_{\underline z}
\end{tikzcd}
$$
commute and obeying the coherence condition for triple compositions.
\end{rem}

\begin{void}
\label{void-factorization-algebra}
Let $\mathscr O$ be a $\base$-presheaf of symmetric monoidal $\infty$-categories.

A $\Ran$-algebra $\mathscr A$ in $\mathscr O$ is called a \emph{(unital) factorization algebra} if for any disjoint $R$-points $\underline x$, $\underline y$ of $\Ran$, the structural morphism
$$
f_{\{\underline x, \underline y\} \rightarrow \underline x \sqcup \underline y} :
\mathscr A_{\underline x} \otimes \mathscr A_{\underline y} \rightarrow \mathscr A_{\underline x \sqcup \underline y}
$$
is an isomorphism (\emph{cf.}~Remark \ref{rem-ran-algebra-factorization-morphisms}).
\end{void}

\begin{rem}
There is another homotopy-coherent formulation of (unital) factorization structures due to Raskin (\emph{cf.}~\cite{raskin2014chiral}).

The above formulation has the advantage of fitting into the formalism of \emph{generalized $\infty$-operads} (\emph{cf.}~\cite[\S2.3.2]{lurie2017higher}). It is independently found by Lin Chen. In forthcoming joint work, we aim to develop the theory of factorization categories from this point of view.
\end{rem}

\begin{void}
\label{void-local-hecke-stack-factorization}
We shall apply the above paradigm to the $\base$-presheaf $\Span(\Stk)$ assigning to each $\base$-algebra $R$ the symmetric monoidal $\infty$-category of spans of \'etale $R$-stacks, and endow $\Hec_G$ with the structure of a factorization algebra in $\Span(\Stk)$.

Indeed, given a $\base$-algebra $R$ and an operation $\{\underline x{}_i\}_{i\in I} \rightarrow \underline y$ in $\Ran(R)$, we have natural morphisms of formal disks
\begin{equation}
\label{eq-formal-disk-factorization-morphism}
\bigsqcup_{i\in I} D_{\underline x{}_i} \simeq D_{\underline x} \rightarrow D_{\underline y} \quad (\underline x := \bigsqcup_{i\in I}\underline x{}_i),
\end{equation}
and compatibly, a cospan of formal punctured disks
\begin{equation}
\label{eq-formal-punctured-disk-factorization-morphism}
\bigsqcup_{i\in I} \mathring D_{\underline x{}_i} \simeq \mathring D_{\underline x} \rightarrow D_{\underline y} \setminus \Gamma_{\underline x} \leftarrow \mathring D_{\underline y}.
\end{equation}
The assignment $\underline x \mapsto D_{\underline x}$ (respectively, $\underline x\mapsto \mathring D_{\underline x}$) defines a $\Ran(R)$-algebra in the category (respectively, the cospan category) of affine schemes over $X$ equipped with the disjoint union symmetric monoidal structure.

By mapping into $\deloop G$, we endow $\Hec_G$ with the structure of a $\Ran$-algebra in $\Span(\Stk)$. Note that given an operation $\{\underline x{}_i\} \rightarrow \underline y$ in $\Ran(R)$, the structural morphism on $\Hec_G$ is the span of \'etale $R$-stacks
\begin{equation}
\label{eq-local-hecke-stack-factorization-morphism}
\begin{tikzcd}[column sep = 1.5em]
	\Hec_{G, \{\underline x{}_i\}\rightarrow \underline y} \ar[r, "p_{\underline y}"]\ar[d, "\prod_{i\in I} p_{\underline x{}_i}"] & \Hec_{G, \underline y} \\
	\prod_{i\in I} \Hec_{G, \underline x{}_i}
\end{tikzcd}
\end{equation}
where $\Hec_{G, \{\underline x{}_i\}\rightarrow \underline y}$ parametrizes modifications $P^0 \overset{\underline x}{\sim} P^1$ of $G$-bundles \emph{over $D_{\underline y}$ along $\underline x$}.

Since both morphisms in \eqref{eq-formal-punctured-disk-factorization-morphism} are isomorphisms for $\underline x = \underline y$, the $\Ran$-algebra $\Hec_G$ in $\Span(\Stk)$ is in fact a factorization algebra.
\end{void}

\begin{rem}
\label{rem-affine-grassmannian-factorization}
We may repeat the construction of \S\ref{void-local-hecke-stack-factorization} for the affine Grassmannian $\Gr_G$ (\emph{cf.}~\S\ref{void-local-hecke-stack-definition}), which lifts $\Gr_G$ to a factorization algebra in $\Stk$, the $\base$-presheaf assigning to $R$ the symmetric monoidal $\infty$-category of \'etale $R$-stacks (as opposed to $\Span(\Stk)$).

Indeed, this is because a $G$-bundle over $D_{\underline y}$ equipped with a trivialization off $\Gamma_{\underline x}$ is equivalent to a $G$-bundle over $D_{\underline x}$ equipped with a trivialization off $\Gamma_{\underline x}$, by the Beauville--Laszlo lemma (\emph{cf.}~\cite{MR1320381}).

Under the forgetful functor $\Stk \rightarrow \Span(\Stk)$, we may realize the structural morphism $\Gr_G \rightarrow \Hec_G$ as a morphism of factorization algebras in $\Span(\Stk)$.
\end{rem}

\begin{void}[Variant: convolution]
\label{void-local-hecke-stack-factorization-convolution}
The construction of \S\ref{void-local-hecke-stack-factorization} has an obvious variant when we take the convolution structure on $\Hec_G$ into account (\emph{cf.}~\S\ref{void-local-hecke-stack-convolution}).

Namely, denote by $\Gpd(\Stk)$ the $\base$-presheaf assigning to $R$ the symmetric monoidal $\infty$-category of groupoid \'etale $R$-stacks. Then $\Hec_G^{[\cdot]}$ has the natural structure of a factorization algebra in $\Span(\Gpd(\Stk))$.

Note that the structural morphism \eqref{eq-local-hecke-stack-factorization-morphism} for $\Hec_G^{[\cdot]}$ involves the \'etale $R$-stack $\Hec_{G, \{\underline x{}_i\} \rightarrow \underline y}^{[\cdot]}$ parametrizing chains of modifications of $G$-bundles over $D_{\underline y}$ along $\underline x$.
\end{void}

\begin{void}
Finally, we shall involve the multiplicative $A$-gerbe $\mathscr G_{\Hec_G}$.

Denote by $\Span(\Gpd(\Stk))_{/\deloop^2 A}$ the $\base$-presheaf assigning to a $\base$-algebra $R$ the $\infty$-category of spans of groupoid \'etale $R$-stacks over $\deloop^2 A^{[\cdot]}$.

Concretely, an object of $\Span(\Gpd(\Stk))_{/\deloop^2 A}(R)$ is a groupoid \'etale $R$-stack $\mathscr Y^{[\cdot]}$ equipped with a multiplicative $A$-gerbe $\mathscr G : \mathscr Y^{[\cdot]} \rightarrow \deloop^2 A^{[\cdot]}$, and a morphism $(\mathscr Y_1^{[\cdot]}, \mathscr G_2) \rightarrow (\mathscr Y_2^{[\cdot]}, \mathscr G_2)$ is a correspondence
$$
\mathscr Y_1^{[\cdot]} \xleftarrow{p_1} \mathscr Y_{12}^{[\cdot]} \xrightarrow{p_2} \mathscr Y_2^{[\cdot]}
$$
of groupoid \'etale $R$-stacks along with an isomorphism of multiplicative $A$-gerbes
$$
p_1^*\mathscr G_1 \simeq p_2^*\mathscr G_2.
$$

The symmetric monoidal structure on $\Span(\Gpd(\Stk))_{/\deloop^2 A}$ is given by Cartesian product of groupoid stacks along with external sum of multiplicative $A$-gerbes.
\end{void}

\begin{prop}
\label{prop-local-hecke-stack-gerbe-factorization}
The pair $(\Hec_G, \mathscr G_{\Hec_G})$ canonically lifts to a factorization algebra in $\Span(\Gpd(\Stk))_{/\deloop^2 A}$.
\end{prop}

\begin{proof}
For each $\base$-algebra $R$, the proof of Lemma \ref{lem-local-hecke-stack-gerbe-multiplicative} supplies a morphism
\begin{equation}
\label{eq-gerbe-ran-algebra-morphism}
\mathscr G_{\Hec_G} : \Hec_G(R) \rightarrow \Gamma(\Spec R, \deloop^2 A),
\end{equation}
where the source $\Hec_G(R)$ lifts to a $\Ran(R)$-algebra in $\Gpd(\Spc)$ (\emph{cf.}~\S\ref{void-local-hecke-stack-factorization-convolution}) while the target $\Gamma(\Spec R, \deloop^2 A)$ is a commutative algebra in $\Spc$, hence a $\Ran(R)$-algebra. It remains to lift \eqref{eq-gerbe-ran-algebra-morphism} to a morphism of $\Ran(R)$-algebras in $\Gpd(\Spc)$, naturally in $R$.

Recall that \eqref{eq-gerbe-ran-algebra-morphism} is the composition of two maps: the map induced from \eqref{eq-hecke-groupoid-maps-to-linear-groupoid} as $\underline x$ varies in $\Ran(R)$
\begin{equation}
\label{eq-hecke-groupoid-maps-to-linear-groupoid-ran-algebra}
\Hec_G(R) \rightarrow \Gamma(D_{(\cdot)} \text{ mod }\mathring D_{(\cdot)}, \deloop^4A(1))
\end{equation}
and the trace map \eqref{eq-trace-map-with-coefficients}. The morphism \eqref{eq-hecke-groupoid-maps-to-linear-groupoid-ran-algebra} naturally lifts to a morphism of $\Ran(R)$-algebras in $\Span(\Gpd(\Spc))$, because the $\Ran(R)$-algebra structures are induced from those of $\mathring D_{(\cdot)}$ and $D_{(\cdot)}$ (\emph{cf.}~\S\ref{void-local-hecke-stack-factorization}). We shall argue that the target of \eqref{eq-hecke-groupoid-maps-to-linear-groupoid-ran-algebra} lifts to a $\Ran(R)$-algebra in $H\integers\Mod$ (with Cartesian symmetric monoidal structure) and that the trace map
\begin{equation}
\label{eq-trace-map-ran-algebra}
\tr_{(\cdot)} : \Gamma(D_{(\cdot)}\text{ mod }\mathring D_{(\cdot)}, \mathscr A(1)[2]) \rightarrow \Gamma(\Spec R, \mathscr A)
\end{equation}
lifts to a map of $\Ran(R)$-algebras in $H\integers\Mod$.

Consider the category $\Div^+(X_R/R)$ of effective Cartier divisors of $X_R$ relative to $R$, with morphisms given by inclusions and symmetric monoidal structure given by sums. The association $\underline x\mapsto \Gamma_{\underline x}$ defines a $\Ran(R)$-algebra $\Gamma_{(\cdot)}$ in $\Div^+(X_R/R)$. The functor
\begin{align*}
\Div^+(X_R / R) &\rightarrow H\integers\Mod \\
Z &\mapsto \Gamma(X_R \text{ mod }X_R \setminus Z, \mathscr A(1)[2])
\end{align*}
is lax symmetric monoidal, so it carries the $\Ran(R)$-algebra $\Gamma_{(\cdot)}$ to the $\Ran(R)$-algebra $\Gamma(X_R \text{ mod }X_R\setminus\Gamma_{(\cdot)}, \mathscr A(1)[2])$ in $H\integers\Mod$. The latter provides the desired lift of the target of \eqref{eq-hecke-groupoid-maps-to-linear-groupoid-ran-algebra} by formal base change (\emph{cf.}~\S\ref{void-trace-map}).

It remains to lift the morphism
$$
\Gamma(X_R \text{ mod }X_R\setminus \Gamma_{(\cdot)}, \mathscr A(1)[2]) \simeq \Gamma(\Gamma_{(\cdot)}, \pi^!\mathscr A) \rightarrow \Gamma(\Spec R, \mathscr A),
$$
provided by adjunction for $\pi : \Gamma_{\underline x} \rightarrow \Spec R$ (over $\underline x \in \Ran(R)$), to a morphism of $\Ran(R)$-algebras in $H\integers\Mod$. The compatibility with structural morphisms are supplied by the naturality of the adjunction map
$$
\Gamma(Z, \pi^!\mathscr A) \rightarrow \Gamma(\Spec R, \mathscr A)
$$
for any $Z \in \Div^+(X_R/R)$ with projection $\pi : Z \rightarrow \Spec R$, with respect to $Z$.
\end{proof}

\subsection{The Satake category}
\label{sec-satake-category}

\begin{void}[Coefficients]
\label{void-satake-category-coefficients}
Let $\ell$ be a prime invertible in $\base$ and $\coeff$ be a finite extension of $\rationals_{\ell}$. We assume that $A$ is a subgroup of $\coeff^{\times}$ and write $\zeta : A \hookrightarrow \coeff^{\times}$ for the inclusion.

We shall use the functor \eqref{eq-twisted-sheaves-functor} assigning to an object $(X, \mathscr G) \in \Sch_{/\deloop^2 A}$ the $\infty$-category $\derived_{\mathscr G, \zeta}(X)$ of $(\mathscr G, \zeta)$-twisted constructible complexes over $X$.
\end{void}

\begin{void}
\label{void-constructible-complexes-indschemes}
For an indscheme $Z$ equipped with an $A$-gerbe $\mathscr G$, we write $\derived_{\mathscr G, \zeta}(Z)$ for the colimit of the $\infty$-categories $\derived_{\mathscr G, \zeta}(X)$ indexed by the poset of closed subschemes $X \hookrightarrow Z$ along pushforward functors.

In particular, a given indscheme presentation $Z \simeq \colim Z_{\alpha}$ realizes $\derived_{\mathscr G, \zeta}(Z)$ as the colimit of $\infty$-categories
$$
\derived_{\mathscr G, \zeta}(Z) \simeq \colim_{\alpha} \derived_{\mathscr G, \zeta}(Z_{\alpha}).
$$

Given any $S$-point $\underline x$ of $\Ran$ (for $S \in \Sch$), we may apply this construction to the indscheme $\Gr_{G, \underline x}$ equipped with the pullback $\mathscr G_{\Gr_G}$ of the $A$-gerbe $\mathscr G_{\Hec_G}$ (\emph{cf.}~\S\ref{void-local-hecke-stack-gerbe}) along the structural map $\Gr_{G, \underline x} \rightarrow \Hec_{G, \underline x}$. This yields the $\infty$-category $\derived_{\mathscr G, \zeta}(\Gr_{G, \underline x})$ of $(\mathscr G_{\Gr_G}, \zeta)$-twisted constructible complexes over $\Gr_{G, \underline x}$.
\end{void}

\begin{void}
\label{void-twisted-sheaves-hecke-stack}
We give a slightly \emph{ad hoc} definition of $\derived_{\mathscr G, \zeta}(\Hec_{G, \underline x})$, sufficient for our purposes.

Namely, given an $S$-point $\underline x$ of $\Ran$ (for $S \in \Sch$), we write $\derived_{\mathscr G, \zeta}(\Hec_{G, \underline x})$ for the $\infty$-category of $L_{\underline x}^+G$-equivariant objects of $\derived_{\mathscr G, \zeta}(\Gr_{G, \underline x})$, formed with respect to the canonical $L_{\underline x}^+G$-equivariance structure of $\mathscr G_{\Gr_G}$.

In terms of the presentation \eqref{eq-hecke-stack-ind-presentation}, we have
$$
\derived_{\mathscr G, \zeta}(\Hec_{G, \underline x}) \simeq \colim_{\alpha} \colim_{\beta \gg \alpha} \derived_{\mathscr G, \zeta}(Z_{\alpha})^{H_{\beta}},
$$
where the colimit over $\beta$ is taken along forgetful functors. (It stabilizes over the cofinal set of $\beta$ where $\ker(L_{\underline x}^+G \rightarrow H_{\beta})$ is pro-unipotent and acts trivially on $Z_{\alpha}$.)
\end{void}

\begin{void}
In the remainder of this subsection, we assume that \emph{$G$ is reductive}.

This hypothesis guarantees that $\Gr_G \rightarrow \Ran$ is ind-proper (\emph{cf.}~\cite[Theorem 3.1.3]{MR3752460}). This in turn implies that each boundary map
$$
\Hec_G^{[n]} \rightarrow \Hec_G^{[m]}
$$
of the groupoid $\Ran$-prestack $\Hec_G$ (\emph{cf.}~\S\ref{void-local-hecke-stack-convolution}) is ind-proper.
\end{void}

\begin{void}[Convolution]
\label{void-convolution-product}
Fix $S \in \Sch$ and an $S$-point $\underline x$ of $\Ran$. We shall construct a monoidal operation on $\derived_{\mathscr G, \zeta}(\Hec_{G, \underline x})$, called the \emph{convolution product}.

Indeed, given $\mathscr A_1, \cdots, \mathscr A_n \in \derived_{\mathscr G, \zeta}(\Hec_{G, \underline x})$, their convolution product is given by
\begin{equation}
\label{eq-convolution-product}
\mathscr A_1 \circ \cdots \circ \mathscr A_n := m_! (p_1^*\mathscr A_1 \otimes \cdots \otimes p_n^*\mathscr A_n),
\end{equation}
where $m$, $p_1, \cdots, p_n$ are the face maps
\begin{equation}
\label{eq-convolution-diagram}
\begin{tikzcd}[column sep = 1.5em]
	\Hec_{G, \underline x}^{[n]} \ar[r, "\prod p_j"]\ar[d, "m"] & \prod_{j = 1}^n \Hec_{G, \underline x} \\
	\Hec_{G, \underline x}
\end{tikzcd}
\end{equation}
sending $P^0 \overset{\underline x}{\sim} \cdots \overset{\underline x}{\sim} P^n$ to $P^0\overset{\underline x}{\sim} P^n$, respectively the segments $P^0\overset{\underline x}{\sim} P^1, \cdots, P^{n-1}\overset{\underline x}{\sim} P^n$.

The formation of \eqref{eq-convolution-product} appeals to the identification of $A$-gerbes
\begin{equation}
\label{eq-gerbe-multiplicativity-isomorphism}
m^*\mathscr G_{\Hec_{G, \underline x}} \simeq p_1^*\mathscr G_{\Hec_{G, \underline x}} + \cdots + p_n^*\mathscr G_{\Hec_{G, \underline x}}
\end{equation}
supplied by Lemma \ref{lem-local-hecke-stack-gerbe-multiplicative}.

The monoidal unit is $e_!(\coeff)$, where $e : \deloop L_{\underline x}^+G \simeq \Hec_{G, \underline x}^{[0]} \rightarrow \Hec_{G, \underline x}$ is the unit section, formed with respect to the trivialization of $e^*(\mathscr G_{\Hec_{G, \underline x}})$ supplied by Lemma \ref{lem-local-hecke-stack-gerbe-multiplicative}.
\end{void}

\begin{rem}
Certainly, we expect the above structure to lift to a monoidal structure on the $\infty$-category $\derived_{\mathscr G, \zeta}(\Hec_{G, \underline x})$. To give a formal construction, however, we need to view $(\Hec_{G, \underline x}, \mathscr G_{\Hec_{G, \underline x}})$ as a monoid in the $\infty$-category of spans of $\Stk_{/\deloop^2 A}$ and apply a ``twisted" $6$-functor formalism, the latter being unavailable at the time of writing.
\end{rem}

\begin{void}
\label{void-satake-subcategory}
Given an $S$-point $\underline x$ of $\Ran$, we define
$$
\Sat_{\mathscr G, \zeta}(\Hec_{G, \underline x}) \subset \derived_{\mathscr G, \zeta}(\Hec_{G, \underline x})
$$
for the full subcategory consisting of objects $\mathscr A$ whose pullback to $\Gr_{G, \underline x}$ are universally locally acyclic (ULA) and perverse relative to $S$ (\emph{cf.}~Definition \ref{defn-universal-local-acyclicity}, \S\ref{void-relative-perversity}).
\end{void}

\begin{rem}
\label{rem-satake-category-forgetful-to-grassmannian}
Pullback along $\Gr_{G, \underline x} \rightarrow \Hec_{G, \underline x}$ realizes $\Sat_{\mathscr G, \zeta}(\Hec_{G, \underline x})$ as the full subcategory of $\derived_{\mathscr G, \zeta}(\Gr_{G, \underline x})$ consisting of objects which are ULA and perverse relative to $S$ and constant along $L^+_{\underline x}G$-orbits.

Indeed, this holds because $L^+_{\underline x}G\rightarrow S$ is pro-smooth with connected geometric fibers: It is an inverse limit of Weil restrictions of $G$ along finite locally free morphisms.
\end{rem}

\begin{prop}
\label{prop-convolution-structure-satake}
Given $S \in \Sch$ and an $S$-point $\underline x$ of $\Ran$,
\begin{enumerate}
	\item the monoidal unit $e_!(\coeff)$ of $\derived_{\mathscr G, \zeta}(\Hec_{G, \underline x})$ belongs to $\Sat_{\mathscr G, \zeta}(\Hec_{G, \underline x})$;
	\item the convolution product of any $\mathscr A_1, \mathscr A_2 \in \Sat_{\mathscr G, \zeta}(\Hec_{G, \underline x})$ belongs to $\Sat_{\mathscr G, \zeta}(\Hec_{G, \underline x})$.
\end{enumerate}
\end{prop}

\begin{void}
\label{void-satake-category-general-base-monoidal}
The proof of Proposition \ref{prop-convolution-structure-satake} will be supplied in \S\ref{void-convolution-fusion-well-defined}.

It implies that $\Sat_{\mathscr G, \zeta}(\Hec_{G, \underline x})$ inherits a monoidal structure from $\derived_{\mathscr G, \zeta}(\Hec_{G, \underline x})$, which we shall refer to as the \emph{convolution monoidal structure}. We shall refer to the monoidal category $(\Sat_{\mathscr G, \zeta}(\Hec_{G, \underline x}), \circ)$ as the \emph{Satake category} at $\underline x$.

Next, we shall show that when $\underline x$ is defined by a finite product of $X$, the Satake category can be upgraded to a symmetric monoidal category.
\end{void}

\begin{rem}
\label{rem-satake-category-finite-set-functoriality}
Let $I$ be a finite set. It determines a morphism $X^I \rightarrow \Ran$, sending $(x^i)_{i\in I}$ to $\underline x := \bigcup_{i\in I} x^i$, hence an \'etale $X^I$-stack $\Hec_{G, I}$ by pullback (and likewise $\Gr_{G, I}$, $L^+_I G$, \emph{etc.}).

The factorization structure on $\Hec_G$ includes, in particular, the following span of \'etale stacks associated to each map of finite sets $\varphi : I \rightarrow J$ (\emph{cf.}~\eqref{eq-local-hecke-stack-factorization-morphism})
\begin{equation}
\label{eq-local-hecke-stack-unital-structure}
\begin{tikzcd}[column sep = 1em]
	\Hec_{G, I \rightarrow J} \ar[r, "\iota"]\ar[d, "\pi"] & \Hec_{G, J} \\
	\Hec_{G, I}
\end{tikzcd}
\end{equation}
where $\Hec_{G, I \rightarrow J}$ parametrizes a point $y^J$ of $X^J$ (with image $x^I$ in $X^I$) and a modification of $G$-bundles over $D_{\underline y}$ along $\underline x$. (Here, $\underline x$, $\underline y$ denote the images of $x^I$, $y^J$ in $\Ran$.) Note that $\Hec_{G, I \rightarrow J}$ may alternatively be presented as the quotient $L^+_{J} G \backslash (\Gr_{G, I} \times_{X^I} X^J)$ of \'etale stacks, showing that $\iota$ is a closed immersion.

Moreover, the factorization structure on $(\Hec_G, \mathscr G_{\Hec_G})$ (\emph{cf.}~Proposition \ref{prop-local-hecke-stack-gerbe-factorization}) allows us to lift \eqref{eq-local-hecke-stack-unital-structure} to a span in $\Stk_{/\deloop^2 A}$. Thus, $\iota_! \pi^*$ restricts to a monoidal functor
\begin{equation}
\label{eq-satake-category-finite-set-functoriality}
	\varphi_! : \Sat_{\mathscr G, \zeta}(\Hec_{G, I}) \rightarrow \Sat_{\mathscr G, \zeta}(\Hec_{G, J})
\end{equation}

This shows that the assignment of $\Sat_{\mathscr G, \zeta}(\Hec_{G, I})$ to $I$ is functorial.
\end{rem}

\begin{void}[Fusion product]
\label{void-fusion-product}
For any integer $n\ge 1$, denote by $X^{\bigsqcup_n I, \disj} \subset X^{\bigsqcup_n I}$ the open subscheme consisting of $I$-tuples $(x^i_1)_{i\in I}, \cdots, (x_n^i)_{i\in I}$ whose images $\underline x{}_1, \cdots, \underline x{}_n$ in $\Ran$ are pairwise disjoint (\emph{cf.}~\S\ref{void-ran-operad}) and by $\Hec_{G, \bigsqcup_n I}^{\disj}$ the base change of $\Hec_{G, \bigsqcup I_n}$ to $X^{\bigsqcup_n I, \disj}$.

We have natural open immersions of \'etale $X^{\bigsqcup_n I}$-stacks
\begin{equation}
\label{eq-hecke-stack-factorization-diagram}
\begin{tikzcd}[column sep = 1em]
	\Hec_{G, \bigsqcup_n I}^{\disj} \ar[r, "f"]\ar[d, "j"] & (\Hec_{G, I})^n \\
	\Hec_{G, \bigsqcup_n I}
\end{tikzcd}
\end{equation}
where the map $f$ comes from the factorization structure on $\Hec_G$ (\emph{cf.}~\S\ref{void-local-hecke-stack-factorization}). Pulling back the external tensor product of $\mathscr A_1, \cdots, \mathscr A_n \in \Sat_{\mathscr G, \zeta}(\Hec_{G, I})$ along $f$ yields
\begin{equation}
\label{eq-satake-sheaves-external-tensor-product}
f^*(\mathscr A_1 \boxtimes \cdots \boxtimes \mathscr A_n) \in \Sat_{\mathscr G, \zeta}(\Hec_{G, \bigsqcup_n I}^{\disj})
\end{equation}
whose formation uses the isomorphism of $A$-gerbes
\begin{equation}
\label{eq-gerbe-factorization-isomorphism}
	f^*(\mathscr G_{\Hec_{G, I}} \boxplus \cdots \boxplus \mathscr G_{\Hec_{G, I}}) \simeq j^*\mathscr G_{\Hec_{G, \bigsqcup_n I}}
\end{equation}
obtained from the lift of $f$ to a morphism in $\Stk_{/\deloop^2 A}$ (\emph{cf.}~Proposition \ref{prop-local-hecke-stack-gerbe-factorization}). Proposition \ref{prop-satake-category-disjoint-pullback} below shows that \eqref{eq-satake-sheaves-external-tensor-product} canonically extends along $j$ to an object
\begin{equation}
\label{eq-external-fusion-product}
\mathscr A_1 \star_{X^I} \cdots \star_{X^I} \mathscr A_n \in \Sat_{\mathscr G, \zeta}(\Hec_{G, \bigsqcup_n I}),
\end{equation}
called the \emph{external fusion product} of $\mathscr A_1, \cdots, \mathscr A_n$.

The pullback of \eqref{eq-external-fusion-product} to the diagonal $X^I \rightarrow X^{\bigsqcup_n I}$ is called the \emph{fusion product}
\begin{equation}
\label{eq-fusion-product}
\mathscr A_1 \star\cdots \star\mathscr A_n \in \Sat_{\mathscr G, \zeta}(\Hec_{G, I}).
\end{equation}
\end{void}

\begin{prop}
\label{prop-satake-category-disjoint-pullback}
For any finite set $I$ and integer $n\ge 1$,
\begin{enumerate}
	\item the pullback functor
	\begin{equation}
	\label{eq-satake-category-disjoint-pullback}
		j^* : \Sat_{\mathscr G, \zeta}(\Hec_{G, \bigsqcup_n I}) \rightarrow \Sat_{\mathscr G, \zeta}(\Hec_{G, \bigsqcup_n I}^{\disj})
	\end{equation}
	is fully faithful;
	
	\item the object \eqref{eq-satake-sheaves-external-tensor-product} belongs to the essential image of \eqref{eq-satake-category-disjoint-pullback}.
\end{enumerate}
\end{prop}

\begin{void}
The proof of Proposition \ref{prop-satake-category-disjoint-pullback} will be supplied in \S\ref{void-convolution-fusion-well-defined}. For now, let us explain how the fusion product \eqref{eq-fusion-product} defines a symmetric monoidal category $(\Sat_{\mathscr G, \zeta}(\Hec_{G, I}), \star)$.

Indeed, we use \eqref{eq-fusion-product} as the $n$-ary operation on $\Sat_{\mathscr G, \zeta}(\Hec_{G, I})$. To see its $\Sigma_n$-equivariance, it suffices to construct an isomorphism in $\Sat_{\mathscr G, \zeta}(\Hec_{G, \bigsqcup_n I}^{\disj})$ for each $\sigma \in \Sigma_n$:
\begin{equation}
\label{eq-fusion-product-symmetric-group-equivariance}
\sigma^*f^*(\mathscr A_1\boxtimes \cdots \boxtimes \mathscr A_n) \simeq f^*(\mathscr A_{\sigma(1)} \boxtimes \cdots \mathscr A_{\sigma(n)}),
\end{equation}
compatible with compositions. The isomorphism \eqref{eq-fusion-product-symmetric-group-equivariance} comes from the $\Sigma_n$-equivariance of $f$ as a morphism in $\Stk_{/\deloop^2 A}$.

Next, we argue that $e_!(\coeff)$ (\emph{cf.}~\S\ref{void-convolution-product}) is a unit for $\Sat_{\mathscr G, \zeta}(\Hec_{G, I})$ also with respect to the fusion product. This follows from the canonical isomorphism in $\Sat_{\mathscr G, \zeta}(\Hec_{G, I \sqcup I}^{\disj})$
$$
f^*(\mathscr A \boxtimes e_!(\coeff)) \simeq j^*(\varphi_! \mathscr A),
$$
where $\varphi : I \rightarrow I \sqcup I$ is the inclusion of the first summand (\emph{cf.}~Remark \ref{rem-satake-category-finite-set-functoriality}), as $\varphi_! \mathscr A$ restricts to $\mathscr A$ along the diagonal.
\end{void}

\begin{prop}[``Convolution $=$ fusion"]
\label{prop-convolution-equals-fusion}
For any finite set $I$, there is a canonical equivalence of monoidal categories
\begin{equation}
\label{eq-convolution-equals-fusion}
(\Sat_{\mathscr G, \zeta}(\Hec_{G, I}), \circ) \simeq (\Sat_{\mathscr G, \zeta}(\Hec_{G, I}), \star).
\end{equation}
\end{prop}

\begin{proof}
We shall argue that the convolution monoidal structure lifts $(\Sat_{\mathscr G, \zeta}(\Hec_{G, I}), \star)$ to a monoid in the (2-)category $\CAlg(\Cat)$ of symmetric monoidal categories.

Concretely, this means that given $\mathscr A_1, \mathscr A_2, \mathscr B_1, \mathscr B_2 \in \Sat_{\mathscr G, \zeta}(\Hec_{G, I})$, we have commutation of fusion and convolution products
$$
(\mathscr A_1 \circ \mathscr A_2) \star (\mathscr B_1 \circ \mathscr B_2) \simeq (\mathscr A_1 \star \mathscr B_1) \circ (\mathscr A_2 \star \mathscr B_2).
$$
This follows from the construction of the fusion product and the fact that $(\Hec_G, \mathscr G_{\Hec_G})$ is a factorization algebra in $\Span(\Gpd(\Stk))_{/\deloop^2 A}$ (\emph{cf.}~Proposition \ref{prop-local-hecke-stack-gerbe-factorization}).

The equivalence \eqref{eq-convolution-equals-fusion} now follows from a categorical analogue of the Eckmann--Hilton argument (\emph{cf.}~\cite[Proposition 2.4.3.9, Proposition 3.2.4.7]{lurie2017higher}).
\end{proof}

\begin{void}[Normalization]
\label{void-normalized-satake-category}
By Proposition \ref{prop-convolution-equals-fusion}, we may regard $\Sat_{\mathscr G, \zeta}(\Hec_{G, I})$ as endowed with either the convolution or the fusion symmetric monoidal structure.

Next, we shall ``normalize" its commutativity constraint to ensure that the constant term functors are symmetric monoidal (\emph{cf.}~Lemma \ref{lem-constant-term-symmetric-monoidal}). There is a locally constant function
\begin{equation}
\label{eq-hecke-stack-parity-function}
d_G : \Hec_{G, I} \rightarrow \integers/2
\end{equation}
defined as follows: \'Etale locally over $X^I$, we may assume that $G$ is split, with a Borel subgroup $B$ and maximal quotient torus $B\twoheadrightarrow T$. Connected components of $\Hec_{G, I}$ are in bijection with $\pi_1 G$, which is a quotient of the cocharacter lattice of $T$. We define \eqref{eq-hecke-stack-parity-function} by
$$
d_G(\theta) := \langle 2\check{\rho}, \theta\rangle \text{ mod }2,
$$
where $2\check{\rho}$ is the sum of positive roots. This expression is independent of the choice of $B$ (as it varies in a connected family), so we obtain \eqref{eq-hecke-stack-parity-function} globally over $X^I$.

Given $\mathscr A_1, \mathscr A_2 \in \Sat_{\mathscr G, \zeta}(\Hec_{G, I})$ supported on connected components $\theta_1, \theta_2$ of $\Hec_{G, I}$,  we modify their commutativity constraint $c_{\mathscr A_1, \mathscr A_2}$ in $\Sat_{\mathscr G, \zeta}(\Hec_{G, I})$ by
\begin{equation}
\label{eq-modified-commutativity-constraint}
\tilde c_{\mathscr A_1, \mathscr A_2} := (-1)^{d_G(\theta_1) d_G(\theta_2)} c_{\mathscr A_1, \mathscr A_2}.
\end{equation}
Denote by ${}^+\Sat_{\mathscr G, \zeta}(\Hec_{G, I})$ the symmetric monoidal category $\Sat_{\mathscr G, \zeta}(\Hec_{G, I})$ equipped with the \emph{modified} commutativity constraint \eqref{eq-modified-commutativity-constraint}.

We shall refer to ${}^+\Sat_{\mathscr G, \zeta}(\Hec_{G, I})$ as the \emph{normalized Satake category}.
\end{void}

\begin{rem}
\label{rem-satake-category-normalization-monoidal-equivalence}
By definition, we have an equivalence of \emph{monoidal} categories
$$
{}^+\Sat_{\mathscr G, \zeta}(\Hec_{G, I}) \simeq \Sat_{\mathscr G, \zeta}(\Hec_{G, I})
$$
which is \emph{incompatible} with the symmetric monoidal structures.
\end{rem}

\subsection{Outer convolution diagram}
\label{sec-outer-convolution-diagram}

\begin{void}
We remain in the context of \S\ref{void-local-hecke-stack-context} and \S\ref{void-satake-category-coefficients}, and assume that $G$ is reductive.

In the previous subsection, we have defined the Satake category subject to two unproven assertions: Proposition \ref{prop-convolution-structure-satake} and Proposition \ref{prop-satake-category-disjoint-pullback}. The goal of this subsection is to supply their proofs.

In fact, we shall deduce both assertions from the ``outer convolution diagram". This is an adaptation of an argument due to Gaitsgory \cite{MR1826370} to the twisted setting.
\end{void}

\begin{void}
\label{void-outer-convolution-diagram-construction}
Let $n\ge 0$ be an integer. Denote by $\widetilde{\Hec}{}^{[n]}_{G}$ the prestack parametrizing $n$ points $\underline x{}_1, \cdots, \underline x{}_n$ of $\Ran$, together with a chain of modifications
\begin{equation}
\label{eq-outer-convolution-point}
P^0 \overset{\underline x{}_1}{\sim} P^1 \overset{\underline x{}_2}{\sim} \cdots \overset{\underline x{}_n}{\sim} P^n
\end{equation}
of $G$-bundles over $D_{\underline x{}_1\cup\cdots\cup\underline x{}_n}$.

The structural morphism $\widetilde{\Hec}{}_G^{[n]} \rightarrow \Ran^n$ realizes $\widetilde{\Hec}{}_G^{[n]}$ as a $\Ran^n$-stack. Given $n$ finite sets $I_1, \cdots, I_n$, we may base change $\widetilde{\Hec}{}_G^{[n]}$ along the product $X^{I_1\sqcup \cdots \sqcup I_n} \rightarrow \Ran^n$ to obtain a $\base$-stack $\widetilde{\Hec}{}_G^{I_1, \cdots, I_n}$. For our applications, only the case $I_1 = \cdots = I_n := I$ for some finite set $I$ is relevant. Let us display the structural morphisms
\begin{equation}
\label{eq-outer-convolution-diagram}
\begin{tikzcd}[column sep = 1.5em]
	\widetilde{\Hec}{}_G^{I,\cdots, I} \ar[r, "\prod\widetilde p_j"]\ar[d, "\widetilde m"] & \prod_{j = 1}^n \Hec_{G, I} \\
	\Hec_{G, \bigsqcup_n I}
\end{tikzcd}
\end{equation}
where $\widetilde m$, $\widetilde p_1, \cdots, \widetilde p_n$ send \eqref{eq-outer-convolution-point} to $P^0 \overset{\underline x{}_1\cup\cdots \underline x{}_n}{\sim} P^n$, respectively the segments $P^0 \overset{\underline x{}_1}{\sim} P^1, \cdots, P^{n-1} \overset{\underline x{}_n}{\sim} P^n$. It is known that $\widetilde m$ is ind-proper (\emph{cf.}~\cite[Corollary 2.10]{MR3178249}).

Note that the pullback of \eqref{eq-outer-convolution-diagram} along the diagonal $\Delta : X^I \rightarrow X^{\bigsqcup_n I}$ recovers \eqref{eq-convolution-diagram}, whereas pulling back $\widetilde{\Hec}{}_G^{I, \cdots, I}$ to the disjoint locus $X^{\bigsqcup_n I, \disj} \subset X^{\bigsqcup_n I}$ in \eqref{eq-outer-convolution-diagram} yields \eqref{eq-hecke-stack-factorization-diagram}.
\end{void}

\begin{lem}
\label{lem-outer-convolution-gerbe-isomorphism}
There is a canonical isomorphism of $A$-gerbes over $\widetilde{\Hec}{}_G^{I, \cdots, I}$:
\begin{equation}
\label{eq-outer-convolution-gerbe-isomorphism}
\widetilde m^* \mathscr G_{\Hec_{G, I}} \simeq \widetilde p_1^* \mathscr G_{\Hec_{G, I}} + \cdots + \widetilde p_n^*\mathscr G_{\Hec_{G, I}}
\end{equation}
with the following properties:
\begin{enumerate}
	\item its pullback to the diagonal $X^I$ is identified with \eqref{eq-gerbe-multiplicativity-isomorphism};
	\item its pullback to the disjoint locus $X^{\bigsqcup_n I, \disj}$ is identified with \eqref{eq-gerbe-factorization-isomorphism}.
\end{enumerate}
\end{lem}

\begin{proof}
For ease of notation, we will treat the case $n = 2$, the general case being similar.

Given an $R$-point of $\widetilde{\Hec}{}_G^{I, I}$, corresponding to $R$-points $x_1^I$, $x_2^I$ of $X^I$ and modifications $P^0 \overset{\underline x{}_1}{\sim} P^1 \overset{\underline x{}_2}{\sim} P^2$, we interpret $P^0$, $P^1$, $P^2$ as morphisms $D_{\underline x{}_1\cup\underline x{}_2} \rightarrow \deloop G$ and form
\begin{equation}
\label{eq-outer-convolution-multiplicative-gerbe-value}
\mu(P^2) - \mu(P^0) \in \Gamma(D_{\underline x{}_1\cup\underline x{}_2} \text{ mod }\mathring D_{\underline x{}_1\cup\underline x{}_2}, \deloop^4A(1)).
\end{equation}
By definition, its image under $\tr_{\underline x{}_1 \cup \underline x{}_2}$ yields the value of $\widetilde m{}^*\mathscr G_{\Hec_{G, I}}$ (\emph{cf.}~\S\ref{void-local-integration-map-construction}).

Let us express \eqref{eq-outer-convolution-multiplicative-gerbe-value} as the sum
\begin{align*}
\mu(P^2) - \mu(P^0) \simeq \mu(P^2) &- \mu(P^1) \\
&+ \mu(P^1) - \mu(P^0),
\end{align*}
where the two summands come from $\Gamma(D_{\underline x{}_i} \text{ mod }\mathring D_{\underline x{}_i}, \deloop^4A(1))$ for $i = 2, 1$, respectively. The desired isomorphism \eqref{eq-outer-convolution-gerbe-isomorphism} now follows from the additivity of the trace map.

Propertes (1) and (2) are immediate consequences of the construction.
\end{proof}

\begin{void}
\label{void-outer-convolution-product}
Given $\mathscr A_1, \cdots, \mathscr A_n \in \Sat_{\mathscr G, \zeta}(\Hec_{G, I})$, we may use the isomorphism \eqref{eq-outer-convolution-gerbe-isomorphism} to form the \emph{outer convolution product}
\begin{equation}
\label{eq-outer-convolution-product}
\mathscr A_1 \circ_{X^I} \cdots \circ_{X^I} \mathscr A_n := \widetilde m_!(\widetilde p_1^*\mathscr A_1 \otimes \cdots \otimes \widetilde p_n^* \mathscr A_2)
\end{equation}
as an object of the $\infty$-category $\derived_{\mathscr G, \zeta}(\Hec_{G, \bigsqcup_n I})$.

Using Lemma \ref{lem-outer-convolution-gerbe-isomorphism}, we obtain isomorphisms
\begin{align}
	\label{eq-outer-convolution-recovers-convolution}
	\Delta^*(\mathscr A_1 \circ_{X^I} \cdots \circ_{X^I} \mathscr A_n) &\simeq \mathscr A_1 \circ \cdots \circ \mathscr A_n, \\
	\label{eq-outer-convolution-recovers-fusion}
	j^*(\mathscr A_1 \circ_{X^I} \cdots \circ_{X^I} \mathscr A_n) &\simeq f^*(\mathscr A_1 \boxtimes\cdots\boxtimes \mathscr A_n),
\end{align}
where the morphisms $j$ and $f$ are as in \S\ref{void-fusion-product}.
\end{void}

\begin{lem}
\label{lem-outer-convolution-product-satake}
The complex $\mathscr A_1\circ_{X^I} \cdots \circ_{X^I} \mathscr A_n$ belongs to $\Sat_{\mathscr G, \zeta}(\Hec_{G, \bigsqcup_n I})$.
\end{lem}

\begin{proof}
By definition, we need to show that $\mathscr A_1 \circ_{X^I} \cdots \circ_{X^I} \mathscr A_n$ is ULA and perverse relative to $X^{\bigsqcup_n I}$ (\emph{cf.}~\S\ref{void-satake-subcategory}). The ULA property holds because $\widetilde m$ is ind-proper.

To show that $\mathscr A_1 \circ_{X^I} \cdots \circ_{X^I} \mathscr A_n$ is perverse relative to $X^{\bigsqcup_n I}$, we first use the isomorphism \eqref{eq-outer-convolution-recovers-fusion} to see that its pullback along $j$ is perverse relative to $X^{\bigsqcup_n I, \disj}$. Then we conclude using its ULA property relative to $X^{\bigsqcup_n I}$ and the fact that the nearby cycles functor is perverse $t$-exact (\emph{cf.}~\cite[Corollaire 4.5]{MR1293970}).
\end{proof}

\begin{void}
\label{void-convolution-fusion-well-defined}
We now prove Propositions \ref{prop-convolution-structure-satake} and \ref{prop-satake-category-disjoint-pullback} simultaneously.

\begin{proof}[Proof of Proposition \ref{prop-satake-category-disjoint-pullback}]
Assertion (1) follows from \cite[Theorem 6.8]{MR4630128}. Assertion (2) follows from Lemma \ref{lem-outer-convolution-product-satake} and the isomorphism \eqref{eq-outer-convolution-recovers-fusion}.
\end{proof}

\begin{proof}[Proof of Proposition \ref{prop-convolution-structure-satake}]
Assertion (1) is immediate.

To prove assertion (2), we note that $\mathscr A_1 \circ \cdots \circ \mathscr A_n$ is ULA relative to $S$ because $m$ is ind-proper. The statement that it is perverse relative to $S$ can be proved over geometric points, so we reduce to the case where $S$ is the spectrum of an algebraically closed field.

In this case, the $S$-point $\underline x$ of $\Ran$ is of the form $\underline x = \bigcup_{i\in I} x_i$, where $\{x_i\}_{i\in I}$ are pairwise disjoint $S$-points of $X$. The factorization structure of $(\Hec_G, \mathscr G_{\Hec_G})$ (\emph{cf.}~Proposition \ref{prop-local-hecke-stack-gerbe-factorization}) allows us to further reduce to the case of a single $S$-point $x$ of $X$.

Now, any object $\mathscr A \in \Sat_{\mathscr G, \zeta}(\Hec_{G, x})$ extends to some \'etale neighborhood of $x$. Thus the result follows from Lemma \ref{lem-outer-convolution-product-satake} and the isomorphism \eqref{eq-outer-convolution-recovers-convolution} (for $I = \{1\}$).
\end{proof}
\end{void}

\subsection{Constant term functors}
\label{sec-constant-term-functors}

\begin{void}
\label{void-constant-term-context}
We remain in the context of \S\ref{void-local-hecke-stack-context} and \S\ref{void-satake-category-coefficients} and assume that $G$ is reductive. Furthermore, we assume the existence of, and choose a square root $\coeff(\frac{1}{2})$ of the Tate twist $\coeff(1)$ over $\Spec \base$. (This is done for convenience, \emph{cf.}~Remark \ref{rem-tate-twist-renormalization-removal}.)

Let $P$ be a parabolic subgroup of $G$ with Levi quotient $P\twoheadrightarrow M$. For any $S$-point $\underline x$ of $\Ran$ ($S \in \Sch$), we have a diagram of indschemes:
$$
\begin{tikzcd}[column sep = 0em]
	& \Gr_{P, \underline x} \ar[dl, swap, "q"]\ar[dr, "p"] \\
	\Gr_{G, \underline x} & & \Gr_{M, \underline x}
\end{tikzcd}
$$
\end{void}

\begin{void}[Definition of $\CT_P$]
\label{void-constant-term-functor}
Since the pullback of $\mu$ along $P \subset G$ canonically descends to $M$, we have an isomorphism of $A$-gerbes over $\Gr_{P, \underline x}$ (already over $\Hec_{P, \underline x}$, in fact):
\begin{equation}
\label{eq-constant-term-gerbe-compatibility}
q^*\mathscr G_{\Gr_G} \simeq p^*\mathscr G_{\Gr_M},
\end{equation}
which we may use to form the functor
\begin{equation}
\label{eq-constant-term-naive}
	q_! p^* : \derived_{\mathscr G, \zeta}(\Gr_{G, \underline x}) \rightarrow \derived_{\mathscr G, \zeta}(\Gr_{M, \underline x}).
\end{equation}

We shall ``normalize" \eqref{eq-constant-term-naive} as follows. Over an \'etale cover of $S$, we may split $M$ and choose a Borel subgroup $B_M \subset M$, with induced Borel $B := P\times_M B_M$ of $G$. Denote by $T$ the maximal quotient torus of $B_M$, so $T$ is the abstract Cartan of both $G$ and $M$. Denote by $2\check{\rho}_G$ (respectively $2\check{\rho}_M$) the sum of positive roots of $G$ (respectively $M$). We obtain a locally constant function
\begin{equation}
\label{eq-levi-grassmannian-dimension-function}
d_{G, M} : \Gr_{M, \underline x} \rightarrow \integers
\end{equation}
taking value $\langle 2\check{\rho}_G - 2\check{\rho}_M, \theta\rangle$ over the connected component corresponding to $\theta \in \pi_1 M$. This expression is independent of the choice of $B_M$, so \eqref{eq-levi-grassmannian-dimension-function} descends to $S$.

We define the \emph{constant term functor} to be
\begin{equation}
\label{eq-constant-term-functor-definition}
\CT_P := q_! p^* ( \frac{d_{G, M}}{2} ) [d_{G, M}],
\end{equation}
using the fixed half-integral Tate twist (\emph{cf.}~\S\ref{void-constant-term-context}).
\end{void}

\begin{rem}
\label{rem-constant-term-composition}
The formation of \eqref{eq-constant-term-functor-definition} respects composition in the following sense. Given a parabolic subgroup $P_1 \subset M$ with Levi quotient $P_1 \twoheadrightarrow M_1$, we have an isomorphism
$$
\CT_{P\times_M P_1} \simeq \CT_{P_1} \circ \CT_P,
$$
where $P\times_M P_1$ is viewed as a parabolic subgroup of $G$. This follows from base change.
\end{rem}

\begin{rem}
\label{rem-constant-term-arc-equivariance}
The functor \eqref{eq-constant-term-functor-definition} carries $L_{\underline x}^+G$-equivariant objects to $L_{\underline x}^+M$-equivariant objects. Indeed, it suffices to show that $p_!$ carries $L_{\underline x}^+P$-equivariant objects to $L_{\underline x}^+M$-equivariant objects, and this holds because the kernel of $L_{\underline x}^+P \rightarrow L_{\underline x}^+M$ is pro-unipotent.

Therefore, \eqref{eq-constant-term-functor-definition} induces a functor
$$
\CT_P : \derived_{\mathscr G, \zeta}(\Hec_{G, \underline x}) \rightarrow \derived_{\mathscr G, \zeta}(\Hec_{M, \underline x}).
$$
\end{rem}

\begin{void}
Next, we prove some basic properties of the constant term functor \eqref{eq-constant-term-functor-definition}. Our main ingredient is a twisted version of Braden's hyperbolic localization theorem, which we establish over a general base scheme in \S\ref{sec-hyperbolic-localization}, following Richarz's strategy in \cite{MR3912059}.
\end{void}

\begin{prop}
\label{prop-constant-term-reflective-properties}
Let $\mathscr A$ be an $L^+_{\underline x}G$-equivariant object of $\derived_{\mathscr G, \zeta}(\Gr_{G, \underline x})$.
\begin{enumerate}
	\item $\mathscr A$ vanishes if and only if $\CT_P(\mathscr A)$ does;
	\item $\mathscr A$ is ULA relative to $S$ if and only if $\CT_P(\mathscr A)$ is;
	\item $\mathscr A$ is connective (respectively, coconnective) in the perverse $t$-structure relative to $S$ if and only if $\CT_P(\mathscr A)$ is.
\end{enumerate}
\end{prop}

\begin{proof}
All statements are of \'etale local nature over $S$, so we may fix a split maximal torus and a Borel subgroup $T \subset B \subset G$, realizing $P$ as a standard parabolic of $G$. Since $\CT_P$ respects composition (\emph{cf.}~Remark \ref{rem-constant-term-composition}), we further reduce to the case $P = B$.

Consider the $\mathbb G_m$-action on $G$ by conjugation by a regular dominant cocharacter. It has attractor $B$, repeller $B^-$ (the opposite Borel), and fixed point locus $T$ (\emph{cf.}~\S\ref{void-hyperbolic-localization-notions}). The induced $\mathbb G_m$-action on $\Gr_{G, \underline x}$ has attractor $\Gr_{B, \underline x}$, repeller $\Gr_{B^-, \underline x}$, and fixed point locus $\Gr_{T, \underline x}$. In particular, \eqref{eq-hyperbolic-localization-diagram} leads to the following diagram of $S$-indschemes
\begin{equation}
\begin{tikzcd}[column sep = 0em]
	& \Gr_{B, \underline x} \ar[dl, swap, "q^+"]\ar[dr, shift left = 0.5ex, "p^+"] \\
	\Gr_{G, \underline x} & & \Gr_{T, \underline x} \ar[ul, shift left = 0.5ex, "i^+"] \ar[dl, shift right = 0.5ex, swap, "i^-"] \\
	& \Gr_{B^-, \underline x} \ar[ul, "q^-"] \ar[ur, swap, shift right = 0.5ex, "p^-"]
\end{tikzcd}
\end{equation}
and we have a canonical isomorphism in $\derived_{\mathscr G, \zeta}(\Gr_{T, \underline x})$ (\emph{cf.}~Theorem \ref{thm-hyperbolic-localization})
\begin{equation}
\label{eq-constant-term-hyperbolic-localization}
\CT_B(\mathscr A) \simeq (p^-)_*(q^-)^! \mathscr A(\frac{d_{G, T}}{2}) [d_{G, T}].
\end{equation}

To prove statement (1), we may use Lemma \ref{lem-hyperbolic-localization-base-change} to reduce to the case $S = \Spec \bar{\base}$. Using the factorization structure of $(\Gr_G, \mathscr G_{\Gr_G})$ (\emph{cf.}~Proposition \ref{prop-local-hecke-stack-gerbe-factorization}), we may further assume that $\underline x$ factors through a $\bar{\base}$-point $x$ of $X$. The implication $\CT_B(\mathscr A) \simeq 0 \Rightarrow \mathscr A \simeq 0$ holds, because every $L_x^+G$-orbit in $\Gr_{G, x}$ meets the image of \emph{some} connected component of $\Gr_{B, x}$ at a single $\bar{\base}$-point (\emph{cf.}~\cite[Proof of Theorem 3.2]{MR2342692}).

To prove statement (2), we may use Proposition \ref{prop-hyperbolic-localization-ula-preservation} to handle the ``$\Rightarrow$" direction. For the ``$\Leftarrow$" direction, we assume that $\CT_B(\mathscr A)$ is ULA and prove that the canonical morphism below is an isomorphism:
\begin{equation}
\label{eq-constant-term-ula-dualization-map}
\mathbf D(\mathscr A) \boxtimes \mathscr A \rightarrow \SHom(\overset{\leftarrow}{\pi}{}^* \mathscr A, \overset{\rightarrow}{\pi}{}^!\mathscr A).
\end{equation}
Here, $\mathbf D$ stands for Verdier duality over $\Gr_{G, \underline x}$ relative to $S$ and $\overset{\leftarrow}{\pi}$, $\overset{\rightarrow}{\pi}$ are the two projections from $\Gr_{G, \underline x}\times_S \Gr_{G, \underline x}$ to $\Gr_{G, \underline x}$ (\emph{cf.}~\S\ref{void-verdier-duality}).

By statement (1), it suffices to prove that \eqref{eq-constant-term-ula-dualization-map} becomes an isomorphism after applying $\CT_{B^-\times B}$, the constant term functor associated to $G\times G$, the \'etale level $(-\mu, \mu)$, and the parabolic subgroup $B^-\times B$. By Lemma \ref{lem-hyperbolic-localization-base-change} and Proposition \ref{prop-hyperbolic-localization-verdier-duality}, the image of \eqref{eq-constant-term-ula-dualization-map} under $\CT_{B^-\times B}$ may be identified with
$$
\mathbf D(\CT_B\mathscr A) \boxtimes \CT_B(\mathscr A) \rightarrow \SHom(\overset{\leftarrow}{\pi}{}^* \CT_B\mathscr A, \overset{\rightarrow}{\pi}{}^!\CT_B\mathscr A),
$$
which is an isomorphism because $\CT_B(\mathscr A)$ is ULA.

To prove statement (3), we observe that the ``$\Rightarrow$" direction follows from \eqref{eq-constant-term-hyperbolic-localization}, as the left-hand-side is right $t$-exact while the right-hand-side is left $t$-exact, and the ``$\Leftarrow$" follows from the ``$\Rightarrow$" direction and statement (1).
\end{proof}

\begin{void}
\label{void-constant-term-satake-category}
By Proposition \ref{prop-constant-term-reflective-properties}, \eqref{eq-constant-term-functor-definition} restricts to a functor on the Satake categories:
\begin{equation}
\label{eq-constant-term-satake-category}
\CT_P : \Sat_{\mathscr G, \zeta}(\Hec_{G, \underline x}) \rightarrow \Sat_{\mathscr G, \zeta}(\Hec_{M, \underline x}).
\end{equation}

When $\underline x$ is given by $X^I \rightarrow \Ran$ for a finite set $I$, let us incorporate the symmetric monoidal structure (\emph{cf.}~\S\ref{void-normalized-satake-category}) and write \eqref{eq-constant-term-satake-category} as a functor
\begin{equation}
\label{eq-constant-term-satake-category-finite-set}
	\CT_P : {}^+\Sat_{\mathscr G, \zeta}(\Hec_{G, I}) \rightarrow {}^+\Sat_{\mathscr G, \zeta}(\Hec_{M, I}).
\end{equation}
\end{void}

\begin{lem}
\label{lem-constant-term-symmetric-monoidal}
The functor \eqref{eq-constant-term-satake-category-finite-set} is naturally symmetric monoidal.
\end{lem}

\begin{proof}
This follows from the classical argument (\emph{cf.}~\cite[Proposition 6.4]{MR2342692}), as the cohomology shift $[d_{G, M}]$ in the definition of \eqref{eq-constant-term-functor-definition} is precisely accounted for by the normalization of the commutativity constraints (\emph{cf.}~\S\ref{void-normalized-satake-category}).
\end{proof}

\subsection{Verdier duality}

\begin{void}
We remain in the context of \S\ref{void-local-hecke-stack-context} and \S\ref{void-satake-category-coefficients} and assume that $G$ is reductive.

Fix an $S$-point $\underline x$ of $\Ran$ ($S \in \Sch$). The goal of this subsection is to use the Verdier duality functor over $\Gr_{G, \underline x}$ relative to $S$ to construct monoidal duals in $\Sat_{\mathscr G, \zeta}(\Hec_{G, \underline x})$ with respect to the convolution monoidal structure (\emph{cf.}~\S\ref{void-satake-category-general-base-monoidal}).

This construction is known in the untwisted setting (\emph{cf.}~\cite[\S11]{MR2342692}, \cite[Proposition VI.8.2]{fargues2021geometrization}) and its adaptation to the twisted setting is straightforward.
\end{void}

\begin{void}
First, we note that Verdier duality (\emph{cf.}~\S\ref{void-verdier-duality}) yields an equivalence
\begin{equation}
\label{eq-verdier-duality-grassmannian}
	\mathbf D : \derived_{\mathscr G, \zeta}(\Gr_{G, \underline x}) \simeq \derived_{-\mathscr G, \zeta}(\Gr_{G, \underline x}).
\end{equation}
preserving perverse ULA objects. After restricting to the full subcategories of perverse ULA objects, it also preserves the property of being constant along $L^+_{\underline x}G$-orbits. These properties also hold for the inverse of \eqref{eq-verdier-duality-grassmannian}, being itself a Verdier duality functor.

Therefore, \eqref{eq-verdier-duality-grassmannian} restricts to an equivalence (\emph{cf.}~Remark \ref{rem-satake-category-forgetful-to-grassmannian})
\begin{equation}
\label{eq-verdier-duality-satake-category}
	\mathbf D : \Sat_{\mathscr G, \zeta}(\Hec_{G, \underline x}) \simeq \Sat_{-\mathscr G, \zeta}(\Hec_{G, \underline x}).
\end{equation}
\end{void}

\begin{void}
On the other hand, consider the \emph{swap} automorphism
\begin{align*}
\swap : \Hec_{G, \underline x} &\rightarrow \Hec_{G, \underline x}, \\
(P^0 \overset{\underline x}{\sim} P^1) & \mapsto (P^1 \overset{\underline x}{\sim} P^0).
\end{align*}

By construction of $\mathscr G_{\Hec_G}$, the pullback $\swap^*(-\mathscr G_{\Hec_G})$ is canonically identified with $\mathscr G_{\Hec_G}$, so $\swap^*$ induces an equivalence
\begin{equation}
\label{eq-swap-satake-autoequivalence}
	\swap^* : \Sat_{-\mathscr G, \zeta}(\Hec_{G, \underline x}) \simeq \Sat_{\mathscr G, \zeta}(\Hec_{G, \underline x}).
\end{equation}
\end{void}

\begin{prop}
\label{prop-convolution-monoidal-dual}
For any $\mathscr A \in \Sat_{\mathscr G, \zeta}(\Hec_{G, \underline x})$, the object $\swap^*\mathbf D(\mathscr A)$ is a left dual of $\mathscr A$ with respect to the convolution monoidal structure.
\end{prop}

\begin{proof}
We need to exhibit the unit and co-unit morphisms
\begin{align}
\label{eq-convolution-dual-unit}
	\mathbf 1 & \rightarrow \swap^*\mathbf D(\mathscr A) \circ \mathscr A \\
\label{eq-convolution-dual-counit}
	\swap^*\mathbf D(\mathscr A) \circ \mathscr A & \rightarrow \mathbf 1,
\end{align}
where $\mathbf 1 := e_!(\coeff)$ is the monoidal unit (\emph{cf.}~Proposition \ref{prop-convolution-structure-satake}), and show that they obey the adjunction property.

In order to do so, we shall use the Cartesian square:
\begin{equation}
\label{eq-hecke-stack-multiplication-diagonal-restriction}
\begin{tikzcd}[column sep = 1em]
	\Hec_{G, \underline x} \ar[r, "\delta"]\ar[d, "t"] & \Hec_{G, \underline x}^{[2]} \ar[d, "m"] \\
	L^+_{\underline x}\deloop G \ar[r, "e"] & \Hec_{G, \underline x}
\end{tikzcd}
\end{equation}
where $m$ is the multiplication map (\emph{cf.}~\S\ref{void-convolution-product}) and $t$, $\delta$ sends $P^0\overset{\underline x}{\sim} P^1$ to $P^1$, respectively the concatenation with its own inverse $P^1\overset{\underline x}{\sim} P^0\overset{\underline x}{\sim} P^1$. Note that if we express $\Hec_{G, \underline x}^{[2]}$ as the fiber product $\Hec_{G, \underline x} \times_{L_{\underline x}^+\deloop G} \Hec_{G, \underline x}$ along the projections $p_1$, $p_2$, then $\delta$ corresponds to the morphism $(\swap, \id)$.

To construct \eqref{eq-convolution-dual-unit}, we express the right-hand-side as $m_*(p_1^*\swap^*\mathbf D(\mathscr A) \otimes p_2^*\mathscr A)$. By adjunction and base change along \eqref{eq-hecke-stack-multiplication-diagonal-restriction}, it suffices to construct a global section of the $!$-restriction of $\mathbf D(\mathscr A) \boxtimes \mathscr A$ along the diagonal of the morphism $i : \Hec_{G, \underline x} \rightarrow L_{\underline x}^+\deloop G$ sending $P^0\overset{\underline x}{\sim} P^1$ to $P^0$. However, our assumption on $\mathscr A$ implies that it is ULA with respect to $i$. Thus this $!$-restriction is isomorphic to $\SHom(\mathscr A, \mathscr A)$ (\emph{cf.}~Remark \ref{rem-internal-hom-via-duality}). The required global section is the identity.

To construct \eqref{eq-convolution-dual-counit}, we perform a similar analysis and reduce to constructing a morphism $\mathbf D(\mathscr A) \otimes \mathscr A \rightarrow i^!(\coeff)$, which is given by the co-unit of Verdier duality. We omit verifying the adjunction property satisfied by \eqref{eq-convolution-dual-unit} and \eqref{eq-convolution-dual-counit}.
\end{proof}

\medskip

\section{The Satake equivalence}

The goal of this section is to formulate the twisted geometric Satake equivalence (\emph{cf.}~Theorem \ref{thm-satake-equivalence}). We define metaplectic dual data in \S\ref{sec-metaplectic-dual-data} and use it to construct the dual category of ``twisted representations" in \S\ref{sec-twisted-representations}. We state the equivalence in \S\ref{sec-satake-equivalence-statement} and list a few of its properties, which will be established in the course of its proof.

Prior works on the twisted geometric Satake equivalence include \cite{MR2265675, MR2684259, lysenko2014twisted, MR2956088, MR3769731}. We compare our statement to the literature in Appendix \ref{sec-prior-works-satake}.

\subsection{The dual data}
\label{sec-metaplectic-dual-data}

\begin{void}
We remain in the context of \S\ref{void-local-hecke-stack-context} and assume in addition that $G$ is reductive. We shall construct a pair $(H, \nu)$, which we call \emph{metaplectic dual data}, where
\begin{enumerate}
	\item $H$ is a locally constant \'etale sheaf over $X$ of pinned split reductive group $\integers$-schemes;
	\item $\nu : \hat Z_H \rightarrow \deloop_X^2 A$ is morphism of \'etale sheaves over $X$ of $\mathbb E_{\infty}$-monoids.
\end{enumerate}

Here, $Z_H$ denotes the center of $H$ and $\hat Z_H$ its character group, viewed as a (locally constant) \'etale sheaf of abelian groups over $X$.

The construction of $(H, \nu)$ is insensitive to the geometry of $X$, so up to \S\ref{void-metaplectic-dual-morphism}, we shall replace $X$ by an arbitrary scheme $S$ over which $A$ has invertible order.
\end{void}

\begin{void}
\label{void-reductive-group-notation}
We first fix some notation asssociated to the reductive group $G$.
\begin{enumerate}
	\item $G_{\sconn}$ (respectively $G_{\adjoint}$) is the simply connected (respectively adjoint) form $G$. The kernel of $G \rightarrow G_{\adjoint}$ is the center $Z$, and the stack quotient $G/G_{\sconn}$ is the cocenter $G_{\abelian}$ of $G$, \emph{cf.}~\cite[\S1.1.5]{shi2025extendedpureinnerforms}.
	\item $T$ (respectively $T_{\sconn}$, $T_{\adjoint}$) is the universal Cartan of $G$ (respectively $G_{\sconn}$, $G_{\adjoint}$). The locally constant \'etale sheaf of \emph{cocharacters} of $T$ (respectively $T_{\sconn}$, $T_{\adjoint}$) is denoted by $\Lambda$ (respectively $\Lambda_{\sconn}$, $\Lambda_{\adjoint}$), with dual $\check{\Lambda}$ (respectively $\check{\Lambda}_{\sconn}$, $\check{\Lambda}_{\adjoint}$).
	
	In particular, we have
	\begin{align*}
		Z &\simeq \Fib(\Lambda \rightarrow \Lambda_{\adjoint}) \otimes_{\integers} \mathbb G_m \\
		G_{\abelian} &\simeq \pi_1 G \otimes_{\integers} \mathbb G_m \text{ for }\pi_1 G := \Lambda/\Lambda_{\sconn}.
	\end{align*}
	\item $\Delta \subset \Lambda$ (respectively $\check{\Delta} \subset \check{\Lambda}$) is the subsheaf of simple coroots (respectively simple roots). Then $\Delta$ spans $\Lambda_{\sconn}$ and $\check{\Delta}$ spans $\check{\Lambda}_{\adjoint}$.
\end{enumerate}
\end{void}

\begin{void}
\label{void-etale-level-fiber-sequence}
Recall that we have a canonical fiber sequence of \'etale sheaves of $\integers$-linear spaces\footnote{Recall that a ``$\integers$-linear space" means a connective $H\integers$-module spectrum (\emph{cf.}~\S\ref{void-linear-spaces}).} over $S$ (\emph{cf.}~\cite[Proposition 5.1.11]{zhao2022metaplectic})
\begin{equation}
\label{eq-etale-level-fiber-sequence}
	\SMaps_{\integers}(\pi_1 G, \deloop^2 A) \rightarrow \SMaps_*(\deloop G, \deloop^4 A(1)) \rightarrow \SQuad(\Lambda, A(-1))_{\strict},
\end{equation}
where $\SMaps_{\integers}(\cdot, \cdot)$ (respectively $\SMaps_*(\cdot, \cdot)$) is the \'etale sheaf of mapping spaces of sheaves of $\integers$-linear spaces (respectively pointed spaces), and $\SQuad(\Lambda, A(-1))_{\strict}$ denotes the \'etale sheaf of \emph{strictly Weyl-invariant} quadratic forms $\Lambda \rightarrow A(-1)$ (to be recalled in \S\ref{void-etale-level-invariants}).

The first map in \eqref{eq-etale-level-fiber-sequence} is given by tensoring with $\mathbb G_m$ and pulling back along $G \rightarrow G_{\abelian}$. The second map in \eqref{eq-etale-level-fiber-sequence} is induced from the canonical isomorphism
$$
\pi_0 \SMaps_*(\deloop G, \deloop^4A(1)) \simeq \SQuad(\Lambda, A(-1))_{\strict}.
$$
\end{void}

\begin{rem}
\label{rem-etale-level-fiber-sequence-multiplicative-group}
For $G = \mathbb G_m$, the fiber sequence \eqref{eq-etale-level-fiber-sequence} reads
\begin{equation}
\label{eq-etale-level-fiber-sequence-multiplicative-group}
\deloop^2 A \rightarrow \SMaps_*(\deloop\mathbb G_m, \deloop^4A(1)) \rightarrow A(-1).
\end{equation}

This fiber sequence admits a canonical splitting, sending $a \in A(-1)$ to $a\otimes (\Psi\otimes\Psi)$, where $\Psi : \deloop\mathbb G_m \rightarrow \deloop^2\hat{\integers}(1)$ is the Kummer morphism (\emph{cf.}~\cite[Remark 4.2.8]{zhao2022metaplectic}).
\end{rem}

\begin{void}
\label{void-etale-level-invariants}
By definition, the \'etale level $\mu$ is a global section of $\SMaps_*(\deloop G, \deloop^4 A(1))$. Let us use \eqref{eq-etale-level-fiber-sequence} to fix notation associated to $\mu$.

\begin{enumerate}
	\item $Q$ is the strictly Weyl-invariant quadratic form $\Lambda \rightarrow A(-1)$ associated to $\mu$, \emph{i.e.}~its image under the second map of \eqref{eq-etale-level-fiber-sequence};
	\item $b$ is the symmetric form on $\Lambda$ associated to $Q$, defined by the identity
	$$
	b(\lambda_1, \lambda_2) := Q(\lambda_1 + \lambda_2) - Q(\lambda_1) - Q(\lambda_2).
	$$
\end{enumerate}

The condition of strict Weyl-invariance on $Q$ means that the equality below holds for any $\alpha \in \Delta$ and $\lambda \in \Lambda$:
\begin{equation}
\label{eq-strict-weyl-invariance}
	b(\alpha, \lambda) = Q(\alpha) \langle\check{\alpha}, \lambda \rangle.
\end{equation}
\end{void}

\begin{void}[Construction of $H$]
\label{void-metaplectic-dual-group}
Denote by $\Lambda^{\sharp} \subset \Lambda$ the kernel of $b$ and by $\check{\Lambda}^{\sharp}$ its dual, which we may view as a lattice inside $\check{\Lambda}\otimes_{\integers}\rationals$. For each $\alpha \in \Delta$, we write
\begin{align*}
\alpha^{\sharp} &:= \ord(Q(\alpha)) \cdot \alpha, \\
\check{\alpha}^{\sharp} &:= \ord(Q(\alpha))^{-1} \cdot \check{\alpha},
\end{align*}
where $\ord(Q(\alpha))$ is the order of $Q(\alpha) \in A(-1)$.

It follows from strict Weyl-invariance \eqref{eq-strict-weyl-invariance} that $\alpha^{\sharp}$ belongs to $\Lambda^{\sharp}$ and $\check{\alpha}^{\sharp}$ belongs to $\check{\Lambda}^{\sharp}$. Furthermore, the quadruple $(\{\alpha^{\sharp}\} \subset \Lambda^{\sharp}, \{\check{\alpha}^{\sharp}\} \subset \check{\Lambda}^{\sharp})$, together with the correspondence $\alpha^{\sharp} \leftrightarrow \check{\alpha}^{\sharp}$ for each $\alpha\in\Delta$, defines a locally constant \'etale sheaf of root data.

We write $H$ for the locally constant \'etale sheaf of pinned split reductive group $\integers$-schemes with root data $(\{\alpha^{\sharp}\} \subset \Lambda^{\sharp}, \{\check{\alpha}^{\sharp}\} \subset \check{\Lambda}^{\sharp})$. (Here, $\{\alpha^{\sharp}\}$ are the simple \emph{roots} of $H$, \emph{etc.})
\end{void}

\begin{rem}
\label{rem-metaplectic-dual-center-as-quotient}
By definition of $H$, the \'etale sheaf $\hat Z_H$ may be realized as the quotient
\begin{equation}
\label{eq-metaplectic-dual-center-characters}
	\hat Z_H \simeq \Lambda^{\sharp} / \Lambda^{\sharp}_{\sconn},
\end{equation}
where $\Lambda^{\sharp}_{\sconn}$ is the $\integers$-linear span of $\alpha^{\sharp}$ ($\alpha \in \Delta$).
\end{rem}

\begin{void}[The form $b_1$]
\label{void-auxiliary-form-primary}
For later purposes, we shall also need the bilinear form
\begin{equation}
\label{eq-auxiliary-form-primary}
b_1 : \Lambda_{\sconn} \otimes \Lambda_{\adjoint} \rightarrow A(-1)
\end{equation}
defined as the $\integers$-linear extension of the expression \eqref{eq-strict-weyl-invariance} for $\alpha \in \Delta$ and $\lambda \in \Lambda_{\adjoint}$, where we interpret $\langle \cdot, \cdot\rangle$ as the canonical pairing between $\Lambda_{\adjoint}$ and the root lattice.
\end{void}

\begin{rem}
\label{rem-metaplectic-root-lattice-as-kernel}
The form $b_1$ provides us another interpretation of the sublattice $\Lambda^{\sharp}_{\sconn} \subset \Lambda_{\sconn}$ (\emph{cf.}~Remark \ref{rem-metaplectic-dual-center-as-quotient}). Namely, it is the (one-sided) kernel of $b_1$.

To see this, we write an arbitrary element of $\Lambda_{\sconn}$ as a sum $\sum_{\alpha\in\Delta} d_{\alpha} \alpha$ for $d_{\alpha} \in \integers$. This element belongs to the kernel of $b_1$ if and only if its pairing with every fundamental coweight $\omega_{\alpha} \in \Lambda_{\adjoint}$ ($\alpha \in \Delta$) vanishes. However, this pairing is computed by $Q(\alpha) d_{\alpha}$, so it vanishes if and only if $d_{\alpha}$ is a multiple of $\ord Q(\alpha)$.
\end{rem}

\begin{void}
\label{void-sharp-torus-symmetric-monoidal-level}
For the moment, let us assume that $G$ is split and choose a Borel subgroup $B \subset G$ as well as a splitting $T\subset B$ of the projection $B \twoheadrightarrow T$.

Pulling back $\mu$ along $T^{\sharp} \rightarrow T \rightarrow G$, for $T^{\sharp} := \Lambda^{\sharp}\otimes_{\integers} \mathbb G_m$, yields an $\mathbb E_{\infty}$-monoidal morphism $\mu_{T^{\sharp}} : \deloop T^{\sharp} \rightarrow \deloop^4A(1)$ (\emph{cf.}~\cite[Proposition 4.6.2]{zhao2022metaplectic}). Taking sections over $\deloop\mathbb G_m$ then yields an $\mathbb E_{\infty}$-monoidal morphism
\begin{align}
\notag
\Lambda^{\sharp} &\simeq \SMaps_*(\deloop\mathbb G_m, \deloop T^{\sharp}) \\
\label{eq-commutative-cover-section}
& \rightarrow \SMaps_*(\deloop \mathbb G_m, \deloop^4A(1)) \simeq \deloop^2 A \oplus A(-1),
\end{align}
where the last isomorphism is given by the canonical splitting of \eqref{eq-etale-level-fiber-sequence-multiplicative-group}. Projecting \eqref{eq-commutative-cover-section} onto the first factor yields an $\mathbb E_{\infty}$-monoidal morphism
\begin{equation}
\label{eq-commutative-cover-section-first-factor}
	\Lambda^{\sharp} \rightarrow \deloop^2 A,
\end{equation}
while projecting it onto the second factor yields the restriction of $Q$ to $\Lambda^{\sharp}$.

In Lemma \ref{lem-etale-level-weyl-equivariance} below, we shall endow \eqref{eq-commutative-cover-section-first-factor} with a canonical equivariance structure with respect to the Weyl-action on $\Lambda^{\sharp}$.
\end{void}

\begin{void}[$G$-equivariance of $\mu$]
\label{void-conjugation-equivariance}
Consider the $G$-action on itself by inner automorphisms. It induces a $G$-action on the pointed stack $\deloop G$. This construction is functorial in sheaves of $\mathbb E_1$-monoids. Applying to the loop space $\Omega \mu : G \rightarrow \deloop^3 A(1)$ of $\mu$, we see that
\begin{equation}
\label{eq-etale-level-conjugation-equivariance}
\mu : \deloop G \rightarrow \deloop^4 A(1)
\end{equation}
is equivariant with respect to the $G$-, and $\deloop^3A(1)$-actions via $\Omega\mu$.

Moreover, using the $\mathbb E_{\infty}$-monoid structure on $\deloop^3 A(1)$, we may trivialize the $\deloop^3A(1)$-action on itself by inner automorphisms, obtaining a trivialization of the induced $\deloop^3 A(1)$-action on $\deloop^4 A(1)$. This shows that \eqref{eq-etale-level-conjugation-equivariance} factors through a morphism
\begin{equation}
\label{eq-etale-level-equivariance-extension}
	\widetilde{\mu} : (\deloop G)/G \rightarrow \deloop^4A(1).
\end{equation}

In other words, $\mu$ admits a canonical $G$-equivariance structure.\footnote{Informally, if one thinks of $\mu$ as a ``central extension" $E$ of $G$ by $\deloop^2 A(1)$, then this $G$-equivariance structure corresponds to an extension of the conjugation action of $G$ along $E \rightarrow G$.}
\end{void}

\begin{rem}
\label{rem-conjugation-equivariance-trivialization}
Let us point out a confusing aspect of \eqref{eq-etale-level-equivariance-extension}. Indeed, the $G$-action on $\deloop G$ as an \'etale stack is canonically trivial: Given a $G$-bundle $P$ and its twist $g(P) := P\times^G G$ by a point $g$ of $G$ (where we use the conjugation map $G \rightarrow G$, $h\mapsto g h g^{-1}$ in the formation of $P \times^G G$), we have the canonical isomorphism of $G$-bundles
\begin{align}
\label{eq-conjugation-action-classifying-stack-trivialization}
	P & \simeq g(P), \\
\notag
	p & \mapsto (p, g).
\end{align}
This is \emph{not} a trivialization of the $G$-action on $\deloop G$ as a \emph{pointed} stack.

Nevertheless, let us use this trivialization to express the isomorphism $\mu(g(P)) \simeq \mu(P)$ provided by \eqref{eq-etale-level-equivariance-extension}. Namely, by identifying its source with $\mu(P)$ using \eqref{eq-conjugation-action-classifying-stack-trivialization}, then this becomes the automorphism of $\mu(P)$ defined by the section $(\Omega \mu)(g)$ of $\deloop^3A(1)$.
\end{rem}

\begin{void}[Weyl-equivariance]
\label{void-weyl-equivariance}
Denote by $N_G(T)$ the normalizer of $T$ in $G$ and by $W := N_G(T)/T$ the Weyl group. The natural $W$-action on $\Lambda$ induces a $W$-action on $\Lambda^{\sharp}$. The restriction of \eqref{eq-etale-level-equivariance-extension} along $T^{\sharp} \rightarrow G$ and $N_G(T) \subset G$ yields an extension of $\mu_{T^{\sharp}}$:
\begin{equation}
\label{eq-etale-level-sharp-torus-equivariance-extension}
\widetilde{\mu}_{T^{\sharp}} : (\deloop T^{\sharp})/N_G(T) \rightarrow \deloop^4A(1).
\end{equation}

The following result shows that \eqref{eq-etale-level-sharp-torus-equivariance-extension} equips $\mu_{T^{\sharp}}$ with a $W$-equivariance structure.
\end{void}

\begin{lem}
\label{lem-etale-level-weyl-equivariance}
The morphism \eqref{eq-etale-level-sharp-torus-equivariance-extension} canonically factors through $(\deloop T^{\sharp})/W$.
\end{lem}

\begin{proof}
Consider the Cartesian diagram of \'etale stacks
\begin{equation}
\label{eq-weyl-group-classifying-stack}
\begin{tikzcd}[column sep = 1.5em]
	\deloop T \ar[r]\ar[d] & \deloop N_G(T) \ar[d, "p"] \\
	S \ar[r] & \deloop W
\end{tikzcd}
\end{equation}
and the \'etale sheaf $\SMaps_*(\deloop T^{\sharp}, \deloop^4A(1))$ over $\deloop W$. The unit map
$$
\SMaps_*(\deloop T^{\sharp}, \deloop^4A(1)) \rightarrow p_*p^*\SMaps_*(\deloop T^{\sharp}, \deloop^4A(1))
$$
induces a fiber sequence by taking global sections:
\begin{equation}
\label{eq-normalizer-equivariance-vs-weyl-equivariance}
\Maps_*((\deloop T^{\sharp})/W, \deloop^4A(1)) \rightarrow \Maps_*((\deloop T^{\sharp})/N_G(T), \deloop^4A(1)) \rightarrow \Hom(\Lambda \otimes \Lambda^{\sharp}, A(-1))^W.
\end{equation}
Here, the third term arises by base change along \eqref{eq-weyl-group-classifying-stack}: The complex corresponding to $\SMaps_*(\deloop T^{\sharp}, \deloop^4 A(1))$ is concentrated in degrees $\ge -2$, with cohomology $\SHom(\Lambda^{\sharp}, A)$ in the lowest degree (\emph{cf.}~\S\ref{void-etale-level-fiber-sequence}), while reduced cohomology of $\deloop T$ is concentrated in degrees $\ge 2$, with $\Hom(\Lambda, \cdot)$ in the lowest degree.

From \eqref{eq-normalizer-equivariance-vs-weyl-equivariance}, we see that $W$-equivariant pointed morphisms $\deloop T^{\sharp} \rightarrow \deloop^4 A(1)$ form a full subcategory of $N_G(T)$-equivariant ones. It remains to show that the image of \eqref{eq-etale-level-sharp-torus-equivariance-extension} in $\Hom(\Lambda \otimes \Lambda^{\sharp}, A(-1))^W$ vanishes. For this, we may restrict \eqref{eq-etale-level-sharp-torus-equivariance-extension} along $\deloop T \rightarrow \deloop N_G(T)$ and identify the resulting morphism
\begin{equation}
\label{eq-etale-level-sharp-torus-equivariance-extension-torus}
\widetilde{\mu}_{T^{\sharp}} : (\deloop T^{\sharp}) / T \rightarrow \deloop^4 A(1)
\end{equation}
as a family of pointed morphisms $\deloop T^{\sharp} \rightarrow \deloop^4A(1)$ parametrized by $\deloop T$. By the commutator computation of \cite[\S5.2]{zhao2022metaplectic}, \eqref{eq-etale-level-sharp-torus-equivariance-extension-torus} is the sum of the constant family $\mu_{T^{\sharp}}$ with the morphism
$$
b\otimes (\Psi \boxtimes \Psi) : \deloop T\times \deloop T^{\sharp} \rightarrow \deloop^4 A(1),
$$
which vanishes because $\Lambda^{\sharp}$ is the kernel of $b$. This implies that \eqref{eq-etale-level-sharp-torus-equivariance-extension-torus} is the constant family $\mu_{T^{\sharp}}$ over $\deloop T$, so its class in $\Hom(\Lambda \otimes \Lambda^{\sharp}, A(-1))$ vanishes.
\end{proof}

\begin{rem}
\label{rem-commutative-cover-independence-of-borel}
The $W$-equivariant $\mathbb E_{\infty}$-monoidal morphism \eqref{eq-commutative-cover-section-first-factor} is canonically assigned to $G$, \emph{i.e.}~independent of the choice of $B$ and the splitting of $B \twoheadrightarrow T$. While this can be verified directly following \cite[\S 5.2.6]{zhao2022metaplectic}, let us offer a cleaner treatment.

Denote by $\Torel$ the $S$-scheme parametrizing subgroup schemes $T \subset B \subset G$, where $B$ is a Borel subgroup and $T$ is a maximal torus of $G$. Any choice of a base point $e : S \rightarrow \Torel$ realizes $\Torel$ as isomorphic to $G / T$. Thus, the geometric fibers of the structural map
\begin{equation}
\label{eq-torel-structural-morphism}
\pi : \Torel \rightarrow S
\end{equation}
have the cohomology of the flag variety of $G$. From this, we deduce that for any complex of torsion \'etale sheaves of invertible order over $S$ concentrated in degrees $\ge -2$, the cofiber of the pullback map
\begin{equation}
\label{eq-pullback-to-torel-cohomology}
\pi^* : \Gamma(S, \mathscr A) \rightarrow \Gamma(\Torel, \mathscr A)
\end{equation}
is concentrated in degrees $\ge 0$. Furthermore, the choice of a base point $e$ realizes $H^0$ of the cofiber of \eqref{eq-pullback-to-torel-cohomology} as $H^2(\Torel \text{ mod }e, H^{-2}\mathscr A)$.

Let us perform the construction of the $W$-equivariant $\mathbb E_{\infty}$-monoidal morphism \eqref{eq-commutative-cover-section-first-factor} with respect to the \emph{universal} Borel and maximal torus over $\Torel$. This yields a $W$-equivariant $\mathbb E_{\infty}$-monoidal morphism
\begin{equation}
\label{eq-metaplectic-dual-morphism-over-torel}
\mu_{T^{\sharp}} : \deloop_{\Torel} T^{\sharp} \simeq \Torel \times \deloop T^{\sharp} \rightarrow \deloop^4A(1).
\end{equation}

It remains to show that \eqref{eq-metaplectic-dual-morphism-over-torel} descends along \eqref{eq-torel-structural-morphism}. Note that this is a \emph{condition} instead of additional data, by the computation of the cofiber of \eqref{eq-pullback-to-torel-cohomology} above, applied to $\mathscr A$ the complex corresponding to $\SMaps_{\mathbb E_{\infty}}(\deloop T^{\sharp}, \deloop^4A(1))^W$. To check that this condition is met, we may work \'etale locally on $S$ and choose a base point $e$ of $\Torel$, realizing the obstruction as an element of
\begin{align*}
H^2(\Torel \text{ mod }e, H^{-2}\mathscr A) & \simeq H^2(\Torel \text{ mod }e, \Hom(\Lambda^{\sharp}, A)^W) \\
& \simeq \Hom(\Lambda_{\sconn}, \Hom(\Lambda^{\sharp}, A)^W(-1)) \subset \Hom(\Lambda_{\sconn} \otimes \Lambda^{\sharp}, A(-1)).
\end{align*}
However, by the commutator computation of \cite[\S5.2]{zhao2022metaplectic}, this element is the restriction of the symmetric form $b$ along $\Lambda_{\sconn} \otimes \Lambda^{\sharp} \subset \Lambda \otimes \Lambda$, which vanishes by definition of $\Lambda^{\sharp}$.
\end{rem}

\begin{void}[$G_{\adjoint}$-equivariance of $\mu_{G_{\sconn}}$]
\label{void-adjoint-equivariance-extension}
We note a variant of the construction of \S\ref{void-conjugation-equivariance}.

Consider the $G_{\adjoint}$-action on $G_{\sconn}$ by conjugation. It induces a $G_{\adjoint}$-action on $\deloop G_{\sconn}$. Since the pullback $\mu_{G_{\sconn}}$ of $\mu$ to $G_{\sconn}$ is uniquely determined by its quadratic form (\emph{cf.}~\S\ref{void-etale-level-fiber-sequence}), it is naturally $G_{\adjoint}$-equivariant, \emph{i.e.}~$\mu_{G_{\sconn}}$ factors through
\begin{equation}
\label{eq-simply-connected-etale-level-equivariant-extension}
\widetilde{\mu}_{G_{\sconn}} : (\deloop G_{\sconn})/G_{\adjoint} \rightarrow \deloop^4A(1).
\end{equation}

Note that $(\deloop G_{\sconn})/G_{\adjoint}$ is canonically identified with the classifying stack of $G_{\sconn} \rtimes G_{\adjoint}$, so \eqref{eq-simply-connected-etale-level-equivariant-extension} may be viewed as an \'etale level for the reductive group scheme $G_{\sconn} \rtimes G_{\adjoint}$ with universal Cartan $T_{\sconn} \times T_{\adjoint}$.

Let us compute the symmetric form
\begin{equation}
\label{eq-simply-connected-etale-level-extension-bilinear-form}
\widetilde b : (\Lambda_{\sconn} \oplus \Lambda_{\adjoint}) \otimes (\Lambda_{\sconn} \oplus \Lambda_{\adjoint}) \rightarrow A(-1)
\end{equation}
associated to \eqref{eq-simply-connected-etale-level-equivariant-extension}.
\end{void}

\begin{lem}
\label{lem-adjoint-extension-symmetric-form-computation}
The symmetric form \eqref{eq-simply-connected-etale-level-extension-bilinear-form} equals the matrix
\begin{equation}
\label{eq-adjoint-extension-symmetric-form-computation}
\widetilde b =
\begin{pmatrix}
	b_{\sconn} & b_1 \\
	(b_1)^{\dagger} & 0
\end{pmatrix}
\end{equation}
where $b_{\sconn}$ is the restriction of $b$ to $\Lambda_{\sconn}$ and $(b_1)^{\dagger}$ is the transpose of $b_1$.
\end{lem}

\begin{proof}
We may choose a Borel subgroup $B\subset G$ and a splitting of $B \twoheadrightarrow T$, realizing $T_{\sconn}$ (respectively $T_{\adjoint}$) as a maximal torus of $G_{\sconn}$ (respectively $G_{\adjoint}$). The restriction of \eqref{eq-simply-connected-etale-level-equivariant-extension} along $T_{\sconn} \subset G_{\sconn}$, $T_{\adjoint} \subset G_{\adjoint}$ yields
\begin{equation}
\label{eq-simply-connected-etale-level-torus-equivariant-extension}
\widetilde{\mu}_{T_{\sconn}} : (\deloop T_{\sconn})/T_{\adjoint} \rightarrow \deloop^4 A(1),
\end{equation}
which we may view as a family of pointed morphism $\deloop T_{\sconn} \rightarrow \deloop^4A(1)$ parametrized by $\deloop T_{\adjoint}$.

By the commutator computation of \cite[Proposition 5.5.4]{zhao2022metaplectic}, \eqref{eq-simply-connected-etale-level-torus-equivariant-extension} is the sum of the constant family $\mu_{T_{\sconn}}$ with the morphism
$$
b_1 \otimes (\Psi \boxtimes \Psi) : \deloop T_{\sconn} \times \deloop T_{\adjoint} \rightarrow \deloop^4A(1)
$$
where $b_1$ is the form \eqref{eq-auxiliary-form-primary}. The equality \eqref{eq-adjoint-extension-symmetric-form-computation} follows.
\end{proof}

\begin{void}[Trivialization along $\Lambda^{\sharp}_{\sconn}$]
\label{void-whittaker-torsor}
We are now ready to construct the canonical trivialization of \eqref{eq-commutative-cover-section-first-factor} over the sublattice $\Lambda_{\sconn}^{\sharp} \subset \Lambda^{\sharp}$. Since the pullback of $\mu_{T^{\sharp}}$ to $T^{\sharp}_{\sconn} := \Lambda^{\sharp}_{\sconn} \otimes_{\integers} \mathbb G_m$ is canonically $\integers$-linear (\emph{cf.}~\S\ref{void-etale-level-fiber-sequence}), it suffices to trivialize the value of $\nu$ on the basis $\{\alpha^{\sharp}\}_{\alpha\in\Delta}$ of $\Lambda^{\sharp}_{\sconn}$. We shall do so using the extension $\widetilde{\mu}_{G_{\sconn}}$ (\emph{cf.}~\S\ref{void-adjoint-equivariance-extension}).

Indeed, Lemma \ref{lem-adjoint-extension-symmetric-form-computation} implies that $\Lambda_{\sconn}^{\sharp} \oplus \Lambda_{\adjoint}$ belongs to the kernel of the symmetric form $\widetilde b$ attached to $\widetilde{\mu}_{G_{\sconn}}$ (\emph{cf.}~Remark \ref{rem-metaplectic-root-lattice-as-kernel}). The constructions of \S\ref{void-sharp-torus-symmetric-monoidal-level}--\S\ref{void-weyl-equivariance}, applied to $\widetilde{\mu}_{G_{\sconn}}$, yields a Weyl-equivariant $\mathbb E_{\infty}$-monoidal (in fact, $\integers$-linear) morphism
$$
\widetilde{\nu} : \Lambda_{\sconn}^{\sharp} \oplus \Lambda_{\adjoint} \rightarrow \deloop^2 A.
$$

On the other hand, the Weyl group $W$ of $G_{\sconn}$ embeds naturally in the Weyl group of $G_{\sconn} \rtimes G_{\adjoint}$, where its action on $\Lambda_{\sconn} \oplus \Lambda_{\adjoint}$ is determined by the formula
\begin{equation}
\label{eq-weyl-group-extended-action}
w (0, \lambda) - (0, \lambda) = (w(\lambda) - \lambda, 0)
\end{equation}
for any $w\in W$ and $\lambda \in \Lambda_{\adjoint}$. This follows, for example, by realizing $T$ as a maximal torus of $G$ and performing the following computation in $G_{\sconn} \rtimes T_{\adjoint}$ for any $\dot w \in G_{\sconn}$ lifting $w$:
\begin{align*}
(\dot w, 1)(1, t)(\dot w^{-1}, 1) &= (\dot w, 1)(t \dot w^{-1} t^{-1}, 1)(1, t) \\
& = (\dot w t \dot w^{-1} t^{-1}, 1)(1, t) = (w(t)t^{-1}, 1)(1, t).
\end{align*}

Altogether, we arrive at a diagram of $W$-equivariant $\mathbb E_{\infty}$-monoidal morphisms
\begin{equation}
\label{eq-simply-connected-commutative-extension}
\begin{tikzcd}[column sep = 1.5em]
	\Lambda_{\sconn}^{\sharp} \ar[r, phantom, "\subset"]\ar[d, phantom, sloped, "\subset"] & \Lambda_{\sconn}^{\sharp} \oplus \Lambda_{\adjoint} \ar[d, "\widetilde{\nu}"] \\
	\Lambda^{\sharp} \ar[r, "\nu"] & \deloop^2 A
\end{tikzcd}
\end{equation}
Writing $\omega \in \Lambda_{\adjoint}$ for the fundamental coweight associated to $\alpha$ and $d := \ord Q(\alpha)$, we obtain canonical isomorphisms using \eqref{eq-weyl-group-extended-action} and the equivariance of $\widetilde{\nu}$ with respect to the simple reflection $s_{\check{\alpha}}\in W$:
\begin{align*}
	\nu(\alpha^{\sharp}) \simeq \widetilde{\nu}(\alpha^{\sharp}, 0) & \simeq \widetilde{\nu}(s_{\check{\alpha}}(0, -d\omega)) - \widetilde{\nu}(0, -d\omega) \\
	& \simeq \widetilde{\nu}(0, -d\omega) - \widetilde{\nu}(0, -d\omega) \simeq 0.
\end{align*}

This provides the desired trivialization of $\nu(\alpha^{\sharp})$.
\end{void}

\begin{rem}
An \emph{a priori} different trivialization of \eqref{eq-commutative-cover-section-first-factor} over $\Lambda^{\sharp}_{\sconn}$ is constructed in \cite[\S 6.1.5]{zhao2022metaplectic}. We have not compared these trivializations, although it seems likely that they coincide. The advantage of the trivialization constructed in \S\ref{void-whittaker-torsor} is that it can be directly related to the affine Grassmannian (\emph{cf.}~Lemma \ref{lem-whittaker-torsor-canonical-trivialization}). This relation will ultimately be responsible for the \emph{canonicity} of our twisted geometric Satake equivalence.
\end{rem}

\begin{void}[Construction of $\nu$]
\label{void-metaplectic-dual-morphism}
Assuming that $G$ is split, we have constructed the $\mathbb E_{\infty}$-monoidal morphism \eqref{eq-commutative-cover-section-first-factor} together with a trivialization over $\Lambda^{\sharp}_{\sconn}$ (\emph{cf.}~\S\ref{void-whittaker-torsor}). In view of \eqref{eq-metaplectic-dual-center-characters}, the morphism \eqref{eq-commutative-cover-section-first-factor} factors through an $\mathbb E_{\infty}$-monoidal morphism
\begin{equation}
\label{eq-metaplectic-dual-morphism}
\nu : \hat Z_H \rightarrow \deloop_X^2 A.
\end{equation}

The $\mathbb E_{\infty}$-monoidal morphism \eqref{eq-metaplectic-dual-morphism} is functorially assigned to $(G, \mu)$, so for any reductive group $S$-scheme $G$, we obtain $\nu$ by \'etale descent.
\end{void}

\begin{void}[The $\vartheta$-shift]
\label{void-theta-shift}
Finally, we involve the geometry of the base curve $X$.

Let us assume that $A$ admits nontrivial $2$-torsion. By our assumption that $|A|$ is invertible in $\base$, this forces $\characteristic \base \neq 2$.

Since the restriction of $b$ to $\Lambda^{\sharp}$ vanishes, the restriction of $Q$ to $\Lambda^{\sharp}$ is a linear form valued in $A(-1)_{2\tors}$, the subsheaf of $A(-1)$ of $2$-torsion elements. It factors through $\hat Z_H$ via \eqref{eq-metaplectic-dual-center-characters}, defining a character
\begin{equation}
\label{eq-2-torsion-character}
\epsilon : \hat Z_H \rightarrow A(-1)_{2\tors}.
\end{equation}

Denote by $\Omega_X$ the canonical line bundle of $X$ relative to $\base$ and consider the image $\Psi(\Omega_X^{-1})$ of $\Omega_X^{-1}$ under $\Psi$ (\emph{cf.}~Remark \ref{rem-kummer-class-of-tautological-line-bundle}), viewed as a global section of $\deloop^2_X\hat{\integers}(1)$. Taking tensor product with \eqref{eq-2-torsion-character}, we obtain a $\integers$-linear morphism
\begin{equation}
\label{eq-theta-morphism}
\vartheta := \epsilon\otimes \Psi(\Omega_X^{-1}) : \hat Z_H \rightarrow \deloop_X^2(A_{2\tors}).
\end{equation}

We may then add \eqref{eq-theta-morphism} to \eqref{eq-metaplectic-dual-morphism} to obtain the $\mathbb E_{\infty}$-monoidal morphism
\begin{equation}
\label{eq-metaplectic-dual-morphism-theta-shifted}
\nu + \vartheta : \hat Z_H \rightarrow \deloop^2_XA,
\end{equation}
which we call the \emph{$\vartheta$-shifted dual datum}.

By convention, the notation $\nu + \vartheta$ stands for $\nu$ itself when $A$ has trivial $2$-torsion.
\end{void}

\begin{rem}
\label{rem-theta-characteristic-trivializes-first-twist}
Since \eqref{eq-2-torsion-character} is $2$-torsion-valued, the formation of \eqref{eq-theta-morphism} only depends on the reduction of $\Psi(\Omega_X^{-1})$ mod $2$. In particular, we may replace $\Omega_X^{-1}$ by $\Omega_X$ without changing \eqref{eq-theta-morphism}. The fact that we prefer $\Omega_X^{-1}$ is because it arises naturally from the geometry of the affine Grassmannian (\emph{cf.}~the proof of Proposition \ref{prop-sharp-torus-commutative-gerbe-identification}).

Note that the reduction of $\Psi(\Omega_X^{-1})$ mod $2$ is the $\mu_2$-gerbe of $\vartheta$-characteristics over $X$. In particular, a choice of a $\vartheta$-characteristic trivializes \eqref{eq-theta-morphism} and renders the $\vartheta$-shift irrelevant. When the ground field $\base$ is algebraically closed or a finite field, then $\vartheta$-characteristics over $X$ always exist (\emph{cf.}~\cite{MR286136}).
\end{rem}

\subsection{Twisted $H$-representations}
\label{sec-twisted-representations}

\begin{void}
We will now invoke the coefficient data of \S\ref{void-satake-category-coefficients}.

For each $S \in \Sch$, we write $\Lis(S)$ for the full subcategory of $\derived(S)$ consisting of locally constant sheaf of finite-dimensional $\coeff$-vector spaces, viewed as a symmetric monoidal category under tensor product.

Given a coalgebra $\mathscr A$ in $\Ind\Lis(S)$, we write $\mathscr A\Comod(\Lis(S))$ for the symmetric monoidal $\coeff$-linear category of $\mathscr A$-comodules in $\Lis(S)$.

The structure sheaf $\mathscr O_H$ admits the structure of a Hopf algebra in $\Ind\Lis(X)$. Thus, we obtain an \'etale sheaf of symmetric monoidal $\coeff$-linear categories $\SRep_H$ over $X$, sending $S$ to $\mathscr O_H\Comod(\Lis(S))$. We write $\Rep_H$ for its global section.
\end{void}

\begin{void}
\label{void-twisted-representations-split}
For a moment, let us assume that $\hat Z_H$ is constant. (This happens for example when $G$ is split, \emph{cf.}~\S\ref{void-metaplectic-dual-group}.) Then for any $\xi \in \hat Z_H$, we have the full subsheaf $\SRep_H^{\xi} \subset \SRep_H$ consisting of objects on which $Z_H$ acts via the character $\xi$. The assignment $\xi \mapsto \SRep_H^{\xi}$ realizes $\SRep_H$ as a $\hat Z_H$-graded sheaf of $\coeff$-linear symmetric monoidal categories (\emph{cf.}~\S\ref{void-graded-sheaf-symmetric-monoidal-categories}, \S\ref{void-graded-symmetric-monoidal-category-sum}).

Let us apply the twisting construction of \S\ref{void-symmetric-monoidal-twist}, with the $\mathbb E_{\infty}$-monoidal morphism $\hat Z_H \rightarrow \deloop^2 \coeff^{\times}$ given by the composition of \eqref{eq-metaplectic-dual-morphism-theta-shifted} with (the deloop of) $\zeta$. This yields a $\hat Z_H$-graded \'etale sheaf of $\coeff$-linear symmetric monoidal categories $\SRep_{H, \nu + \vartheta}$.

Since \'etale sheaves of symmetric monoidal $\coeff$-linear categories satisfy \'etale descent, the construction of \S\ref{void-twisted-representations-split} generalizes to the case where $\hat Z_H$ is locally constant. Thus, we obtain an \'etale sheaf of symmetric monoidal $\coeff$-linear categories $\SRep_{H, \nu + \vartheta}$, whose global section is denoted by $\Rep_{H, \nu + \vartheta}$.

Objects of $\Rep_{H, \nu + \vartheta}$ are called \emph{$(\nu + \vartheta)$-twisted $H$-representations on $\coeff$-local systems}.
\end{void}

\begin{void}
\label{void-twisted-representations-multiple-point}
Next, we extend the construction of $\SRep_{H, \nu + \vartheta}$ to finite products of copies of $X$. By \'etale descent, it again suffices to treat the case where $\hat Z_H$ is constant.

For any finite set $I$, the external tensor product $\mathscr O_H^{\boxtimes I}$ is a Hopf algebra in $\Ind\Lis(X^I)$. We write $\SRep_{H^{\boxtimes I}}$ for the \'etale sheaf of symmetric monoidal $\coeff$-linear categories over $X^I$, sending $S$ to $\mathscr O_H^{\boxtimes I}\Comod(\Lis(X^I))$. Then $\SRep_{H^{\boxtimes I}}$ lifts to a $(\hat Z_H)^{\boxplus I}$-graded sheaf of symmetric monoidal $\coeff$-linear categories.

Applying the twisting construction of \S\ref{void-symmetric-monoidal-twist} with the $\mathbb E_{\infty}$-monoidal morphism
\begin{align}
\notag
	(\nu + \vartheta)^{\boxplus I} : (\hat Z_H)^{\boxplus I} &\rightarrow (\deloop^2_{X^I} A)^{\oplus I} \\
\label{eq-metaplectic-dual-morphism-external-sum}
	&  \xrightarrow{\sum} \deloop^2_{X^I} A \xrightarrow{\zeta} \deloop^2_{X^I}\coeff^{\times}
\end{align}
we obtain a $(\hat Z_H)^{\boxplus I}$-graded \'etale sheaf of $\coeff$-linear symmetric monoidal categories $\SRep_{H^{\boxtimes I}, (\nu + \vartheta)^{\boxplus I}}$ over $X^I$, with global section $\Rep_{H^{\boxtimes I}, (\nu + \vartheta)^{\boxplus I}}$.
\end{void}

\begin{rem}
\label{rem-twisted-representation-finite-set-functoriality}
Let us explain the functoriality of $\Rep_{H^{\boxtimes I}, (\nu + \vartheta)^{\boxplus I}}$ with respect to the finite set $I$, in parallel with Remark \ref{rem-satake-category-finite-set-functoriality}.

Indeed, given a morphism of finite sets $\varphi : I \rightarrow J$, we have an induced morphism $\Delta_{\varphi} : X^J \rightarrow X^I$, hence a pullback functor $(\Delta_{\varphi})^* : \Lis(X^I) \rightarrow \Lis(X^J)$. There is a natural morphism of Hopf algebras $(\Delta_{\varphi})^*(\mathscr O_H^{\boxtimes I}) \rightarrow \mathscr O_H^{\boxtimes J}$ in $\Ind\Lis(X^J)$, compatible with restrictions of \eqref{eq-metaplectic-dual-morphism-external-sum}, so we obtain a functor of symmetric monoidal $\coeff$-linear categories
$$
(\Delta_{\varphi})^* : \Rep_{H^{\boxtimes I}, (\nu + \vartheta)^{\boxplus I}} \rightarrow \Rep_{H^{\boxtimes J}, (\nu + \vartheta)^{\boxplus J}},
$$
which respects compositions of morphisms of finite sets.
\end{rem}

\begin{void}
\label{void-symmetric-monoidal-dual-data-fiber-sequence}
Next, we shall rewrite $\Rep_{H^{\boxtimes I}, (\nu + \vartheta)^{\boxplus I}}$ in terms of an ``$L$-group". To do so, we need to extract the $\integers$-linear part of $\nu$ (\emph{cf.}~\S\ref{void-metaplectic-dual-morphism}).

Let us perform this construction in a more abstract setting: Let $\Xi$ be an \'etale sheaf of abelian groups over a base scheme $S$ over which $A$ has invertible order. There is a fiber sequence of \'etale sheaves of $\integers$-linear spaces
\begin{equation}
\label{eq-symmetric-monoidal-dual-data-fiber-sequence}
\SMaps_{\integers}(\Xi, \deloop^2 A) \rightarrow \SMaps_{\mathbb E_{\infty}}(\Xi, \deloop^2 A) \rightarrow \SHom(\Xi, A_{2\tors}),
\end{equation}
where the first map is induced from the forgetful functor and the second map is defined as follows: Given an $\mathbb E_{\infty}$-monoidal morphism $\Xi \rightarrow \deloop^2 A$, its fiber is a symmetric monoidal extension $\widetilde{\Xi}$ of $\Xi$ by $\deloop A$. Associating to each $\xi \in \widetilde{\Xi}$ the commutativity constraint $c_{\xi, \xi} : \xi\otimes \xi \simeq \xi\otimes \xi$ defines a character $\Xi \rightarrow A_{2\tors}$.

A key observation, due to Gaitsgory and Lysenko (\emph{cf.}~\cite[\S4.8]{MR3769731}), is that \eqref{eq-symmetric-monoidal-dual-data-fiber-sequence} canonically splits. To construct this splitting, we may assume $A_{2\tors} \simeq \integers/2$ and treat the universal case $\Xi := \integers/2$, where we lift the identity character of $\integers/2$ to the $\mathbb E_{\infty}$-monoidal morphism $\integers/2 \rightarrow \deloop^2 \integers/2$ which is trivial as an $\mathbb E_1$-monoidal morphism, but with commutativity constraint specified by the pairing
$$
\integers/2 \otimes \integers/2 \rightarrow \integers/2,\quad a\otimes b\mapsto ab.
$$
\end{void}

\begin{void}
Let us apply the above construction to $S := X$ and $\Xi := \hat Z_H$. Using the splitting of \eqref{eq-symmetric-monoidal-dual-data-fiber-sequence}, we may project $\nu$ onto a $\integers$-linear morphism
\begin{equation}
\label{void-metaplectic-dual-morphism-linear-component}
{}^0\nu : \hat Z_H \rightarrow \deloop^2_X A.
\end{equation}

We shall also involve the $\vartheta$-shift (\emph{cf.}~\S\ref{void-theta-shift}): Using ${}^0\nu$ instead of $\nu$ in the construction of \eqref{eq-metaplectic-dual-morphism-theta-shifted}, we obtain a $\integers$-linear morphism
\begin{equation}
\label{void-metaplectic-dual-morphism-linear-component-theta-shifted}
{}^0\nu + \vartheta : \hat Z_H \rightarrow \deloop^2_X A.
\end{equation}
\end{void}

\begin{rem}
\label{rem-metaplectic-dual-morphism-self-commutativity-constraint}
The image of $\nu$ along the second map of \eqref{eq-symmetric-monoidal-dual-data-fiber-sequence} is a homomorphism $\hat Z_H \rightarrow A_{2\tors}$. By \cite[Proposition 4.6.6]{zhao2022metaplectic}, this coincides with the homomorphism \eqref{eq-2-torsion-character} under the canonical isomorphism
$$
A(-1)_{2\tors} \simeq A_{2\tors}.
$$
\end{rem}

\begin{void}[The $L$-group]
\label{void-ell-group}
We now assume that $X$ is geometrically connected and affine. This hypothesis implies that $X$ is an algebraic $K(\pi, 1)$-space (\emph{cf.}~\cite[\S1.4]{MR3404650}).

We fix a geometric point $\bar x$ of $X$ and write $Z_{H, \bar x}$ (respectively $H_{\bar x}$) for the fiber of $Z_H$ (respectively $H$) at $\bar x$, which is equipped with a natural $\pi_1^{\etale}(X, \bar x)$-action.

We also fix an algebraic closure $\overline{\rationals}_{\ell}$ of $\rationals_{\ell}$ containing $\coeff$. Inducing \eqref{void-metaplectic-dual-morphism-linear-component-theta-shifted} along $\zeta$, we obtain a $\integers$-linear morphism $\hat Z_H \rightarrow \deloop^2_X\overline{\rationals}{}_{\ell}^{\times}$, or equivalently an \'etale $Z_H(\overline{\rationals}_{\ell})$-gerbe. Choose a trivialization of this $Z_H(\overline{\rationals}_{\ell})$-gerbe over $\bar x$ and denote by ${}^LZ_{H, X}$ its fundamental group, which is well-defined by the $K(\pi, 1)$-property of $X$ (\emph{cf.}~\cite[Theorem 19.6]{MR3802418}). By construction, ${}^LZ_{H, X}$ fits into a short exact sequence
\begin{equation}
\label{eq-ell-group-center}
1 \rightarrow Z_{H, \bar x}(\overline{\rationals}_{\ell}) \rightarrow {}^L Z_{H, X} \rightarrow \pi_1^{\etale}(X, \bar x) \rightarrow 1
\end{equation}
satisfying the following properties:
\begin{enumerate}
	\item the conjugation action of ${}^LZ_{H, X}$ on $Z_{H, \bar x}(\overline{\rationals}_{\ell})$ factors through the natural action of $\pi_1^{\etale}(X, \bar x)$ on $Z_{H, \bar x}(\overline{\rationals}_{\ell})$;
	\item there is a finite extension $\coeff \subset \coeff_1$ contained in $\overline{\rationals}_{\ell}$ and a finite quotient $\pi_1^{\etale}(X, \bar x)\twoheadrightarrow \Gamma$ such that \eqref{eq-ell-group-center} is induced from an extension of $\Gamma$ by $Z_{H, \bar x}(\coeff_1)$.
\end{enumerate}

Finally, we define the \emph{$L$-group} to be the short exact sequence induced from \eqref{eq-ell-group-center} along the $\pi_1^{\etale}(X, \bar x)$-equivariant inclusion $Z_{H, \bar x}(\overline{\rationals}_{\ell}) \hookrightarrow H_{\bar x}(\overline{\rationals}_{\ell})$:
\begin{equation}
\label{eq-ell-group}
	1 \rightarrow H_{\bar x}(\overline{\rationals}_{\ell}) \rightarrow {}^LH_X \rightarrow \pi_1^{\etale}(X, \bar x) \rightarrow 1.
\end{equation}
\end{void}

\begin{rem}
\label{rem-ell-group-baer-sum}
The extension \eqref{eq-ell-group-center} is constructed using \eqref{void-metaplectic-dual-morphism-linear-component-theta-shifted} as input, the latter being the sum of $\integers$-linear morphisms ${}^0\nu$ and $\vartheta$. We may apply the same construction to $\vartheta$ and ${}^0\nu$ individually, obtaining extensions ${}^LZ_{H, X}^{(1)}$, respectively ${}^LZ_{H, X}^{(2)}$, of $\pi_1^{\etale}(X, \bar x)$ by $Z_{H, \bar x}(\overline{\rationals}_{\ell})$. Then \eqref{eq-ell-group-center} is the Baer sum
$$
{}^LZ_{H, X} \simeq {}^LZ_{H, X}^{(1)} + {}^LZ_{H, X}^{(2)}.
$$

We shall explain in \S\ref{sec-meta-galois-twist} that ${}^LZ_{H, X}^{(1)}$ (defined by $\vartheta$) is Weissman's ``meta-Galois twist" (\emph{cf.}~\cite[\S4]{MR3802418}, also called the ``first twist" in \cite{MR3276164}) when $\base$ is a finite field.
\end{rem}

\begin{void}
\label{void-ell-group-multiple-point}
We remain in the context of \S\ref{void-ell-group} and fix a finite set $I$. The product $X^I$ is also an algebraic $K(\pi, 1)$-space. Applying the construction of \emph{loc.cit.}, with \eqref{void-metaplectic-dual-morphism-linear-component-theta-shifted} replaced by the $\integers$-linear morphism
$$
({}^0\nu + \vartheta)^{\boxplus I} : (\hat Z_H)^{\boxplus I} \rightarrow \deloop^2_{X^I} A,
$$
we obtain a short exact sequence
$$
1 \rightarrow H_{\bar x}^I(\overline{\rationals}_{\ell}) \rightarrow {}^LH_{X^I} \rightarrow \pi_1^{\etale}(X^I, \bar x^I) \rightarrow 1.
$$

Denote by $\Rep({}^LH_{X^I})$ the symmetric monoidal category of continuous ${}^LH_{X^I}$-representations on finite-dimensional $\overline{\rationals}_{\ell}$-vector spaces which are algebraic over $H^I_{\bar x}(\overline{\rationals}_{\ell})$. On the other hand, recall the symmetric monoidal category $\Rep_{H^{\boxtimes I}, (\nu + \vartheta)^{\boxplus I}}$ (\emph{cf.}~\S\ref{void-twisted-representations-multiple-point}), with $\coeff$ replaced by $\overline{\rationals}_{\ell}$.
\end{void}

\begin{prop}
For any finite set $I$, there is a canonical equivalence of $\overline{\rationals}_{\ell}$-linear \emph{monoidal} categories
\begin{equation}
\label{eq-twisted-representation-as-ell-group-representation}
\Rep_{H^{\boxtimes I}, (\nu + \vartheta)^{\boxplus I}} \simeq \Rep({}^LH_{X^I}).
\end{equation}
\end{prop}

\begin{proof}
It follows from \cite[Proposition 6.4.8]{zhao2022metaplectic} that $\Rep({}^LH_{X^I})$ is canonically equivalent to $\Rep_{H^{\boxtimes I}, ({}^0\nu + \vartheta)^{\boxplus I}}$ as $\overline{\rationals}_{\ell}$-linear \emph{symmetric} monoidal categories.

It remains to observe that $({}^0\nu + \vartheta)^{\boxplus I}$ and $(\nu + \vartheta)^{\boxplus I}$ are isomorphic as $\mathbb E_1$-monoidal morphisms, but this follows from the fact that the splitting of \S\ref{void-symmetric-monoidal-dual-data-fiber-sequence} is defined by an $\mathbb E_1$-monoidally trivial morphism.
\end{proof}

\begin{rem}
The equivalence \eqref{eq-twisted-representation-as-ell-group-representation} is natural in $I$ in the following sense: Given a map of finite sets $\varphi : I \rightarrow J$, we have a commutative square
\begin{equation}
\label{eq-twisted-representation-as-ell-group-representation-naturality}
\begin{tikzcd}[column sep = 1.5em]
	\Rep_{H^{\boxtimes I}, (\nu + \vartheta)^{\boxplus I}} \ar[r, phantom, "\simeq"] \ar[d, "(\Delta_{\varphi})^*"] & \Rep({}^LH_{X^I}) \ar[d, "(\Delta_{\varphi})^*"] \\
	\Rep_{H^{\boxtimes J}, (\nu + \vartheta)^{\boxplus J}} \ar[r, phantom, "\simeq"] & \Rep({}^LH_{X^J}).
\end{tikzcd}
\end{equation}
where the left vertical arrow is defined as in Remark \ref{rem-twisted-representation-finite-set-functoriality} and the right vertical arrow is restriction along the natural map $\Delta_{\varphi} : {}^LH_{X^J} \rightarrow {}^LH_{X^I}$ corresponding to $\varphi$. Furthermore, \eqref{eq-twisted-representation-as-ell-group-representation-naturality} is compatible with compositions.
\end{rem}

\subsection{Statement of the equivalence}
\label{sec-satake-equivalence-statement}

\begin{void}
\label{void-satake-equivalence-statement-context}
Let $X$ be a smooth curve over a field $\base$ and $G$ be a reductive group $X$ scheme. Let $\coeff$ be a finite extension of $\rationals_{\ell}$, for a prime $\ell$ invertible in $\base$. We assume the existence of a square root $\coeff(\frac{1}{2})$ of $\coeff(1)$ over $\Spec\base$ and fix it. Let $\zeta : A \hookrightarrow \coeff^{\times}$ be a finite subgroup and $\mu$ be an $A$-valued \'etale level of $G$.

For any finite set $I$, we have defined the symmetric monoidal $\coeff$-linear abelian categories ${}^+\Sat_{\mathscr G, \zeta}(\Hec_{G, I})$ (\emph{cf.}~\S\ref{void-normalized-satake-category}) and $\Rep_{H^{\boxtimes I}, (\nu + \vartheta)^{\boxplus I}}$ (\emph{cf.}~\S\ref{void-twisted-representations-multiple-point}).

The statement below is our version of the geometric Satake equivalence.
\end{void}

\begin{thm}
\label{thm-satake-equivalence}
For any finite set $I$, there is a canonical equivalence of symmetric monoidal $\coeff$-linear categories
\begin{equation}
\label{eq-satake-equivalence}
{}^+\Sat_{\mathscr G, \zeta}(\Hec_{G, I}) \simeq \Rep_{H^{\boxtimes I}, (\nu + \vartheta)^{\boxplus I}}.
\end{equation}
\end{thm}

\begin{rem}
\label{rem-tate-twist-renormalization-removal}
In our formulation of Theorem \ref{thm-satake-equivalence}, we sacrificed some generality by invoking the square root $\coeff(\frac{1}{2})$. This is done so that the dual group $H$ appears naturally in the equivalence \eqref{eq-satake-equivalence}.

As in the classical setting, one can remove half-integral Tate twists by \emph{not} incorporating any Tate twist in the definition of the constant term functor (\emph{cf.}~\S\ref{void-constant-term-functor}). This would lead to a different fiber functor, hence a different Tannaka dual $H^{\geom}$. The analysis is parallel to \cite[Appendix A]{MR3346175}, so we will omit it.
\end{rem}

\begin{rem}
We formulated Theorem \ref{thm-satake-equivalence} for $\ell$-adic constructible sheaves. This is motivated by our main application for global function fields.

When the ground field $\base$ is $\complexes$, we may also define ${}^+\Sat_{\mathscr G, \zeta}(\Hec_{G, I})$ using constructible sheaves of $\coeff$-vector spaces (for any field $\coeff$ of characteristic zero) in the classical topology, or using algebraic D-modules. Our proof of Theorem \ref{thm-satake-equivalence} applies in these contexts.
\end{rem}

\begin{void}[Functoriality in $I$]
\label{void-satake-equivalence-finite-set-functoriality}
The equivalence \eqref{eq-satake-equivalence} will be natural in $I$ with respect to the functorialities specified in Remark \ref{rem-satake-category-finite-set-functoriality} and Remark \ref{rem-twisted-representation-finite-set-functoriality}.

More precisely, for any map of finite sets $\varphi : I \rightarrow J$, there is a canonical $2$-isomorphism rendering the following diagram of symmetric monoidal $\coeff$-linear categories commute:
\begin{equation}
\label{eq-satake-equivalence-finite-set-functoriality}
\begin{tikzcd}[column sep = 1.5em]
	{}^+\Sat_{\mathscr G, \zeta}(\Hec_{G, I}) \ar[r, phantom, "\simeq"]\ar[d, "\varphi_!"] & \Rep_{H^{\boxtimes I}, (\nu + \vartheta)^{\boxplus I}} \ar[d, "(\Delta_{\varphi})^*"] \\
	{}^+\Sat_{\mathscr G, \zeta}(\Hec_{G, J}) \ar[r, phantom, "\simeq"] & \Rep_{H^{\boxtimes J}, (\nu + \vartheta)^{\boxplus J}}
\end{tikzcd}
\end{equation}

Furthermore, the $2$-isomorphism in \eqref{eq-satake-equivalence-finite-set-functoriality} respects compositions of maps of finite set.
\end{void}

\begin{void}[Compatibility with constant terms]
\label{void-satake-equivalence-constant-term-compatibility}
Let $P$ be a parabolic subgroup of $G$ with Levi quotient $P \twoheadrightarrow M$. The pullback of the \'etale level $\mu$ along $P \subset G$ canonically descends to an \'etale level $\mu_M$ of $M$. Applying the construction of the metaplectic dual data to $\mu_M$, we obtain $(H_M, \nu_M)$ (\emph{cf.}~\S\ref{sec-metaplectic-dual-data}).

By construction, there is a morphism
\begin{equation}
\label{eq-parabolic-subgroup-dual-embedding}
h_P : H_M \rightarrow H
\end{equation}
of locally constant \'etale sheaves over $X$ of pinned split reductive group $\integers$-schemes. It induces a morphism $\hat h_P : \hat Z_H \rightarrow \hat Z_{H_M}$ of characters of the centers with $\nu_M\circ \hat h_P \simeq \nu$.

The equivalence \eqref{eq-satake-equivalence} renders the following diagram commute
$$
\begin{tikzcd}[column sep = 1.5em]
	{}^+\Sat_{\mathscr G, \zeta}(\Hec_{G, I}) \ar[r, phantom, "\simeq"]\ar[d, "\CT_P"] & \Rep_{H^{\boxtimes I}, (\nu + \vartheta)^{\boxplus I}} \ar[d, "h_P^*"] \\
	{}^+\Sat_{\mathscr G, \zeta}(\Hec_{M, I}) \ar[r, phantom, "\simeq"] & \Rep_{H_M^{\boxtimes I}, (\nu_M + \vartheta)^{\boxplus I}}
\end{tikzcd}
$$
where $\CT_P$ is the constant term functor \eqref{eq-constant-term-satake-category-finite-set} and $h_P^*$ is restriction along \eqref{eq-parabolic-subgroup-dual-embedding}.
\end{void}

\begin{rem}
In the context of \S\ref{void-ell-group}, we have defined the category $\Rep({}^LH_X)$ (and more generally $\Rep({}^LH_{X^I})$ for a finite set $I$). By combining \eqref{eq-twisted-representation-as-ell-group-representation} and \eqref{eq-satake-equivalence} and replacing $\coeff$ by $\overline{\rationals}_{\ell}$, we obtain a canonical equivalence of \emph{monoidal} $\overline{\rationals}_{\ell}$-linear categories
\begin{equation}
\label{eq-satake-equivalence-ell-group}
	\Sat_{\mathscr G, \zeta}(\Hec_{G, I}) \simeq \Rep({}^LH_{X^I})
\end{equation}
natural in $I$. Here, we omitted the normalization of the Satake category as it does not affect the underlying monoidal category (\emph{cf.}~Remark \ref{rem-satake-category-normalization-monoidal-equivalence}).
\end{rem}

\medskip

\section{Pointwise studies of the Satake category}
\label{sec-pointwise-studies}

Throughout this section, we fix a smooth curve $X$ over an \emph{algebraically closed} field $\base$ and a \emph{split} reductive group $X$-scheme $G$. Let $\coeff$ be a finite extension of $\rationals_{\ell}$, for a prime $\ell$ invertible in $\base$. We shall also fix a square root $\coeff(\frac{1}{2})$ of $\coeff(1)$. Let $\zeta : A \hookrightarrow \coeff^{\times}$ be a finite subgroup and $\mu$ be an $A$-valued \'etale level of $G$.

The goal of this section is to study $\coeff$-linear abelian category $\Sat_{\mathscr G, \zeta}(\Hec_{G, x})$ for a $\base$-point $x$ of $X$ (\emph{cf.}~\S\ref{void-satake-subcategory}). Namely, we will prove that it is semisimple and identify its simple objects. In fact, we will prove a stronger assertion where we allow $x$ to be any $S$-point of $X$ (\emph{cf.}~Theorem \ref{thm-semisimplicity}). The proof of this theorem is an adaptation of Lusztig and Yun's method in handling monodromic Hecke categories (\emph{cf.}~\cite{MR4108915}). This idea of using \cite{MR4108915} is suggested by Gurbir Dhillon.\footnote{At least when the \'etale level $\mu$ admits an integral lift (\emph{cf.}~\S\ref{sec-integral-vs-etale-levels}), the recent work \cite{dhillon2025endoscopymetaplecticaffinehecke} goes much further than what we prove in this section.}

\subsection{Schubert cells}

\begin{void}
\label{void-split-reductive-killing-pair-context}
We use the notation \S\ref{void-reductive-group-notation} for notions associated to $G$. Furthermore, we fix a maximal torus and a Borel subgroup $T\subset B\subset G$. (In particular, $T$ is identified with the universal Cartan of $G$.) Denote by $\Lambda^+ \subset \Lambda$ the submonoid of dominant coweights.

We shall work over one copy of $X$ and write $\Hec_G$ for $\Hec_{G, \{1\}}$ (\emph{cf.}~Remark \ref{rem-satake-category-finite-set-functoriality}). We use similar notation $\Gr_G$, $L^+G$, $LG$, \emph{etc.}
\end{void}

\begin{void}
\label{void-schubert-cell-notation}
For any $\lambda \in \Lambda$, we write $\Hec_G^{\lambda}$ (respectively $\Gr_G^{\lambda}$) for the corresponding Schubert cell of $\Hec_G$ (respectively $\Gr_G$).

Denote by $\varpi^{\lambda} : X \rightarrow \Gr_T$ the map sending $x$ to the modification $\mathscr O \overset{x}{\sim} \mathscr O(\lambda \Gamma_x)$ of $T$-bundles. We use the same notation for its composition with the inclusion into $\Gr_G$ (and also for its further composition with the projection onto $\Hec_G$). By definition, the $L^+G$-action on $\varpi^{\lambda} : X \rightarrow \Gr_G$ has orbit $\Gr_G^{\lambda}$ and stabilizer $L^+G \cap \varpi^{\lambda} L^+G \varpi^{-\lambda}$.

Denote by $P^{\lambda} \subset G$ the parabolic subgroup generated by root subgroups $N_{\alpha}$ corresponding to those roots $\check{\alpha}$ with $\langle\check{\alpha}, \lambda\rangle \le 0$.\footnote{If we assume $\lambda \in \Lambda^+$, then $P^{\lambda}$ is \emph{opposite} to the standard parabolic subgroup corresponding to simple roots orthogonal to $\lambda$.} Denote by $M^{\lambda}$ the Levi quotient of $P^{\lambda}$. With the above notation, we have morphisms
\begin{equation}
\label{eq-schubert-cell-natural-morphisms}
\begin{tikzcd}[column sep = 1em]
	X \ar[r, "\varpi^{\lambda}"] & \Hec_G^{\lambda} \ar[r, phantom, "\subset"]\ar[d, "\pi^{\lambda}"] & \Hec_G \\
	& \deloop_X (M^{\lambda})
\end{tikzcd}
\end{equation}
where $\pi^{\lambda}$ is the composition (\emph{cf.}~\cite[\S2.1]{MR3752460})
\begin{equation}
\label{eq-bialynicki-birula-map}
\Hec_G^{\lambda} \simeq \deloop_X (L^+G \cap \varpi^{\lambda} L^+G \varpi^{-\lambda}) \rightarrow \deloop_X(P^{\lambda}) \rightarrow \deloop_X(M^{\lambda}).
\end{equation}
\end{void}

\begin{void}
\label{void-gerbe-schubert-cell-notation}
Denote by $\mathscr G_{\Hec_G^{\lambda}}$ the restriction of the $A$-gerbe $\mathscr G_{\Hec_G}$ (\emph{cf.}~\S\ref{void-local-hecke-stack-gerbe}) to $\Hec_G^{\lambda}$. We shall determine $\mathscr G_{\Hec_G^{\lambda}}$ in relation to the morphisms in \eqref{eq-schubert-cell-natural-morphisms}.

Write $\mathscr G_{\varpi^{\lambda}}$ for the pullback of $\mathscr G_{\Hec_G}$ along $\varpi^{\lambda}$. It is an $A$-gerbe over $X$, so it defines an $A$-gerbe $\mathscr G_{\varpi^{\lambda}}|_{\mathscr Z}$ on any $X$-prestack $\mathscr Z$ by pullback.

On the other hand, the fiber sequence \eqref{eq-etale-level-fiber-sequence} yields an isomorphism
$$
\SHom(\pi_1(M^{\lambda}), A(-1))
\simeq
\SMaps_*(\deloop_X(M^{\lambda}), \deloop^2_X A),
$$
so any character $\chi : \pi_1(M^{\lambda}) \rightarrow A(-1)$ induces a pointed morphism $\chi\otimes\Psi :\deloop_X(M^{\lambda}) \rightarrow \deloop^2_XA$, where $\Psi$ is the Kummer morphism (\emph{cf.}~Remark \ref{rem-kummer-class-of-tautological-line-bundle}).

Let $b$ the symmetric form associated to $\mu$ (\emph{cf.}~\S\ref{void-etale-level-invariants}). By construction, all roots of $M^{\lambda}$ are orthogonal to $\lambda$. The equality \eqref{eq-strict-weyl-invariance} thus implies that the character $b(\lambda, \cdot) : \Lambda \rightarrow A(-1)$ factors through $\pi_1(M^{\lambda})$, so we obtain a pointed morphism
$$
b(\lambda, \cdot)\otimes \Psi : \deloop_X(M^{\lambda}) \rightarrow \deloop^2_XA.
$$

The following result is an analogue of \cite[Lemma 2.4]{MR2684259} for \'etale levels.
\end{void}

\begin{prop}
\label{prop-gerbe-schubert-cell}
There is a canonical isomorphism of $A$-gerbes
\begin{equation}
\label{eq-gerbe-schubert-cell}
\mathscr G_{\Hec_G^{\lambda}} - (\mathscr G_{\varpi^{\lambda}}|_{\Hec_G^{\lambda}}) \simeq (\pi^{\lambda})^*(b(\lambda, \cdot)\otimes\Psi).
\end{equation}
\end{prop}

\begin{void}[The forms $b_2$]
\label{void-additional-bilinear-forms}
For the proof of Proposition \ref{prop-gerbe-schubert-cell}, as well as for many later purposes, we shall recall the canonical quadratic structure on $\mu$ (\emph{cf.}~\cite[Proposition 3.1.3]{shi2025extendedpureinnerforms}).

To state it, we need the bilinear pairing
$$
b_2 : \pi_1 G \otimes \Fib(\Lambda \rightarrow \Lambda_{\adjoint}) \rightarrow A(-1),
$$
characterized by the property that the adjoints of $b_2$, $b$, $b_1$ (\emph{cf.}~\S\ref{void-etale-level-invariants}, \S\ref{void-auxiliary-form-primary}) fit into a map of fiber sequences
\begin{equation}
\label{eq-second-bilinear-form-adjoint-definition}
\begin{tikzcd}[column sep = 1.5em]
	\Fib(\Lambda \rightarrow \Lambda_{\adjoint}) \ar[r, "b_2"]\ar[d] & \SHom(\pi_1G, A(-1)) \ar[d] \\
	\Lambda \ar[r, "b"]\ar[d] & \SHom(\Lambda, A(-1)) \ar[d] \\
	\Lambda_{\adjoint} \ar[r, "b_1"] & \SHom(\Lambda_{\sconn}, A(-1))
\end{tikzcd}
\end{equation}

In particular, by tensoring $b_2$ with the self-tensor product of the Kummer morphism $\Psi^{\otimes 2} : \deloop\mathbb G_m \otimes \deloop\mathbb G_m \rightarrow \deloop^4\hat{\integers}(2)$, we obtain a pairing
\begin{equation}
\label{eq-second-auxiliary-form-induced-pairing}
b_2\otimes \Psi^{\otimes 2} : \deloop G_{\abelian} \otimes \deloop Z \rightarrow \deloop^4 A(1).
\end{equation}
\end{void}

\begin{void}[Canonical quadratic structure]
\label{void-canonical-quadratic-structure}
Consider the $\deloop Z$-action on $\deloop G$ induced from the $Z$-action on $G$. Denote by
$$
a : \deloop G \times \deloop Z \rightarrow \deloop G
$$
the action morphism and by $p_1$, $p_2$ the projections of $\deloop G\times\deloop Z$ onto $\deloop G$, respectively $\deloop Z$.

By \cite[Proposition 3.1.3]{shi2025extendedpureinnerforms}, there is a canonical isomorphism
\begin{equation}
\label{eq-canonical-quadratic-structure}
a^*\mu \simeq (p_1)^*\mu + (p_2)^*\mu_Z + b_2\otimes \Psi^{\otimes 2},
\end{equation}
where $\mu_Z$ denotes the pullback of $\mu$ along $Z \subset G$, and we (slightly abusively) use $b_2\otimes \Psi^{\otimes 2}$ to denote the pullback of \eqref{eq-second-auxiliary-form-induced-pairing} to $\deloop G\times\deloop Z$.

Furthermore, \eqref{eq-canonical-quadratic-structure} is compatible with the canonical trivializations of the two sides along $e \times \deloop Z$, respectively $\deloop G \times e$, and admits natural cocycle data over $\deloop G\times\deloop Z\times\deloop Z$.
\end{void}

\begin{rem}
For $G = T$ a split torus, we have $b_2 = b$ and \eqref{eq-canonical-quadratic-structure} reduces to the canonical quadratic structure constructed in \cite[Proposition 4.7.3]{zhao2022metaplectic}. For the proof of Proposition \ref{prop-gerbe-schubert-cell} (but not for later purposes), this special case is sufficient.
\end{rem}

\begin{void}
\label{void-gerbe-schubert-cell-proof}
Let us now construct the isomorphism \eqref{eq-gerbe-schubert-cell}.

\begin{proof}[Proof of Proposition \ref{prop-gerbe-schubert-cell}]
Recall that $\pi^{\lambda}$ is the quotient of the affine-space fibration $\Gr_G^{\lambda} \rightarrow G/P^{\lambda}$ by the action of $L^+G \rightarrow G$, which has a pro-unipotent kernel. Thus, pullback by $\pi^{\lambda}$ induces an equivalence of (discrete) groupoids
\begin{equation}
\label{eq-bb-pullback-equivalence-of-pointed-gerbes}
(\pi^{\lambda})^* : \Maps_*(\deloop_X(M^{\lambda}), \deloop^2_X A) \simeq \Maps_*(\Hec_G^{\lambda}, \deloop_X^2 A),
\end{equation}
where $\Hec_G^{\lambda}$ is viewed as a pointed $X$-stack via $\varpi^{\lambda}$.

The left-hand-side of \eqref{eq-gerbe-schubert-cell} admits the natural structure of a \emph{pointed} morphism $\Hec_G^{\lambda} \rightarrow \deloop_X^2 A$, so it remains to identify it with the image of $b(\lambda, \cdot)\otimes\Psi$ under \eqref{eq-bb-pullback-equivalence-of-pointed-gerbes}.

The inclusion $T \subset G$ gives rise to a commutative square
$$
\begin{tikzcd}[column sep = 1em]
	\Maps_*(\deloop_X(M^{\lambda}), \deloop_X^2 A) \ar[r, phantom, "\simeq"]\ar[d] & \Maps_*(\Hec_G^{\lambda}, \deloop_X^2 A) \ar[d] \\
	\Maps_*(\deloop_X T, \deloop_X^2 A) \ar[r, phantom, "\simeq"] & \Maps_*(\Hec_T^{\lambda}, \deloop_X^2A)
\end{tikzcd}
$$
where the vertical arrows are fully faithful (\emph{cf.}~\S\ref{void-gerbe-schubert-cell-notation}). Thus, it suffices to construct the isomorphism \eqref{eq-gerbe-schubert-cell} after restriction to $\Hec_T^{\lambda}$, or equivalently to its reduced locus $\deloop_X L^+T$.

Let $P$ be an $R$-point of $\deloop_X L^+T$, viewed as a $T$-bundle over $D_x$ with $x \in X(R)$. It defines the $R$-ponit $P \overset{x}{\sim} P(\lambda x)$ of $\Hec_T^{\lambda}$. The value of $\mathscr G_{\Hec_T}$ at this $R$-point is the image of
$$
\mu(P(\lambda x)) - \mu(P) \in \Gamma(D_x \text{ mod }\mathring D_x, \deloop^4A(1))
$$
under the trace map (\emph{cf.}~\S\ref{void-local-integration-map-construction}). We shall use the canonical quadratic structure \eqref{eq-canonical-quadratic-structure} (for $G = T$) to identify this section:
$$
\mu(P(\lambda x)) - \mu(P) \simeq \mu(\mathscr O(\lambda x)) + (b\otimes \Psi^{\otimes 2})(P, \mathscr O(\lambda x)).
$$
Note that under the trace map, $\mu(\mathscr O(\lambda x))$ yields the value of $\mathscr G_{\varpi^{\lambda}}|_{\Hec_T}$.

It thus remains to identify the image of $(b\otimes\Psi^{\otimes 2})(P, \mathscr O(\lambda x))$ under the trace map with the image of $P|_{\Gamma_x}$ under $b(\lambda, \cdot)\otimes\Psi$. To do so, we shall write
$$
(b\otimes\Psi^{\otimes 2})(P, \mathscr O(\lambda x)) \simeq (b(\lambda, \cdot)\otimes\Psi)(P) \otimes \Psi(\mathscr O(x)),
$$
where $(b(\lambda, \cdot)\otimes\Psi)(P)$ is viewed as an object of $\Gamma(D_x, \deloop^2 A)$ and $\Psi(\mathscr O(x))$ is viewed as an object of $\Gamma(D_x \text{ mod }\mathring D_x, \deloop^2 \hat{\integers}(1))$.

Note that $(b(\lambda, \cdot)\otimes\Psi)(P)$ canonically descends to $(b(\lambda, \cdot)\otimes\Psi)(P|_{\Gamma_x})$ along $D_x\rightarrow \Spec R$. We thus obtain the desired identification
\begin{align*}
	\tr_x ((b\otimes\Psi^{\otimes 2})(P, \mathscr O(\lambda x))) & \simeq \tr_x ((b(\lambda, \cdot)\otimes\Psi)(P) \otimes \Psi(\mathscr O(x))) \\
	& \simeq (b(\lambda, \cdot)\otimes\Psi)(P|_{\Gamma_x}) \otimes \tr_x(\Psi(\mathscr O(x))) \simeq (b(\lambda, \cdot)\otimes\Psi)(P|_{\Gamma_x})
\end{align*}
from Remark \ref{rem-trace-map-projection-formula} and Remark \ref{rem-kummer-class-of-tautological-line-bundle}.
\end{proof}
\end{void}

\begin{void}
Given an $S$-point $x$ of $X$ ($S\in\Sch$), we have $\infty$-category $\derived_{\mathscr G, \zeta}(\Hec_{G, x})$ (\emph{cf.}~\S\ref{void-twisted-sheaves-hecke-stack}).

We shall use Proposition \ref{prop-gerbe-schubert-cell} to show that its objects are supported on (the base change of) Schubert cells $\Hec_{G, x}^{\lambda}$ corresponding to $\lambda \in \Lambda^{\sharp}$. (Recall that $\Lambda^{\sharp}$ denotes the kernel of $b$, \emph{cf.}~\S\ref{void-metaplectic-dual-group}.) More precisely, we have the following vanishing statement.
\end{void}

\begin{cor}
\label{cor-schubert-cell-reduction-of-support}
Let $x$ be an $S$-point of $X$ ($S\in\Sch$) and $\lambda \in \Lambda\setminus\Lambda^{\sharp}$. Then
$$
\derived_{\mathscr G, \zeta}(\Hec_{G, x}^{\lambda}) \simeq 0.
$$
\end{cor}

\begin{proof}
The vanishing of any $\mathscr A \in \derived_{\mathscr G, \zeta}(\Hec_{G, x}^{\lambda})$ may be checked on fibers, so we reduce to the case $S = \Spec \base$. The condition $\lambda \notin\Lambda^{\sharp}$ is equivalent to the nontriviality of the $A$-gerbe $b(\lambda, \cdot)\otimes\Psi$ over the classifying stack $\deloop(M^{\lambda})$ (\emph{cf.}~\S\ref{void-gerbe-schubert-cell-notation}).

By Proposition \ref{prop-gerbe-schubert-cell}, the restriction of $\mathscr G_{\Hec_G}$ to $\Hec_{G, x}^{\lambda}$ is \emph{non-canonically} isomorphic to the pullback of $b(\lambda, \cdot)\otimes\Psi$ along $\pi^{\lambda}$. Thus the result follows from Lemma \ref{lem-equivariance-vanishing}, applied to $\Spec\base$ and a finite type quotient of $L^+_x G \cap \varpi^{\lambda} L_x^+G \varpi^{-\lambda}$.
\end{proof}

\begin{void}
Next, let us study the case $\lambda \in \Lambda^{\sharp}$ more closely.
\end{void}

\begin{cor}
\label{cor-gerbe-schubert-cell-sharp-identification}
For each $\lambda \in \Lambda^{\sharp}$, there is a canonical isomorphism of $A$-gerbes
\begin{equation}
\label{eq-gerbe-schubert-cell-sharp-identification}
\mathscr G_{\Hec_G^{\lambda}} \simeq \mathscr G_{\varpi^{\lambda}}|_{\Hec_G^{\lambda}}.
\end{equation}
\end{cor}

\begin{proof}
The isomorphism \eqref{eq-gerbe-schubert-cell-sharp-identification} is supplied by \eqref{eq-gerbe-schubert-cell} and the trivialization of its right-hand-side coming from the vanishing of $b(\lambda, \cdot)$.
\end{proof}

\begin{void}
By Corollary \ref{cor-gerbe-schubert-cell-sharp-identification}, the $A$-gerbe $\mathscr G_{\Hec_G^{\lambda}}$ descends to $X$ for any $\lambda\in \Lambda^{\sharp}$. In particular, given $\lambda_1 \in W\lambda$, where $W$ is the Weyl group, we obtain an identification
\begin{equation}
\label{eq-gerbe-schubert-cell-weyl-equivariance}
\mathscr G_{\varpi^{\lambda_1}} \simeq \mathscr G_{\varpi^{\lambda}}
\end{equation}
from the fact that $\omega^{\lambda}$, $\omega^{\lambda_1}$ belong to the same $L^+G$-orbit.

On the other hand, the restriction of $\mu$ to $T^{\sharp}$ is $W$-equivariant by Lemma \ref{lem-etale-level-weyl-equivariance}. It follows that the $A$-gerbe $\mathscr G_{\Hec_T}$ is also $W$-equivariant. This yields another identification
\begin{equation}
\label{eq-gerbe-torus-cocharacter-weyl-equivariance}
\mathscr G_{\varpi^{\lambda_1}} \simeq \mathscr G_{\varpi^{\lambda}},
\end{equation}
which may \emph{a priori} differ from \eqref{eq-gerbe-schubert-cell-weyl-equivariance}. Let us show that this is not the case.
\end{void}

\begin{prop}
\label{prop-gerbe-schubert-cell-weyl-equivariance}
Given $\lambda, \lambda_1 \in \Lambda^{\sharp}$ in the same Weyl-orbit, the identifications \eqref{eq-gerbe-schubert-cell-weyl-equivariance} and \eqref{eq-gerbe-torus-cocharacter-weyl-equivariance} are canonically isomorphic.
\end{prop}

\begin{proof}
Recall the canonical $G$-equivariance structure on $\mu$ (\emph{cf.}~\S\ref{void-conjugation-equivariance}). It induces an $L^+G$-equivariance structure on $\mathscr G_{\Gr_G}$ as follows: Given $R$-points $g$ of $L^+G$ and $P^0 \overset{x}{\sim} P^1$ of $\Gr_G$ lying over the same $R$-point $x$ of $X$, the $G$-equivariance structure on $\mu$ yields an isomorphism in $\Gamma(D_x \text{ mod }\mathring D_x, \deloop^4 A(1))$:
\begin{equation}
\label{eq-conjugation-equivariance-grassmannian-isomorphism}
\mu(g(P^1)) - \mu(g(P^0)) \simeq \mu(P^1) - \mu(P^0),
\end{equation}
where $g(P^1)$, $g(P^0)$ are the twists of $P^1$, $P^0$ by $g$ (\emph{cf.}~Remark \ref{rem-conjugation-equivariance-trivialization}). The desired $L^+G$-equivariance structure on $\mathscr G_{\Gr_G}$ is obtained from \eqref{eq-conjugation-equivariance-grassmannian-isomorphism} via the trace map (\emph{cf.}~\S\ref{void-local-integration-map-construction}).

It suffices to prove that this $L^+G$-equivariance structure on $\mathscr G_{\Gr_G}$ agrees with the one coming from the descent data of $\mathscr G_{\Hec_G}$. Indeed, the former $L^+G$-equivariance structure induces \eqref{eq-gerbe-torus-cocharacter-weyl-equivariance} while the latter induces \eqref{eq-gerbe-schubert-cell-weyl-equivariance}.

Note that the $L^+G$-equivariance structure coming from the descent data of $\mathscr G_{\Hec_G}$ corresponds to the isomorphism
\begin{equation}
\label{eq-descent-data-grassmannian-isomorphism}
\mu(g(P^1)) - \mu(g(P^0)) \simeq \mu(P^1) - \mu(P^0)
\end{equation}
arising from the identification of $g(P^0) \overset{x}{\sim} g(P^1)$ and $P^0 \overset{x}{\sim} P^1$ as $R$-points of $\Hec_G$. This identification is in turn defined by the canonical isomorphisms of $G$-bundles $g(P^0) \simeq P^0$, $g(P^1) \simeq P^1$, \emph{cf.}~Remark \ref{rem-conjugation-equivariance-trivialization}. It follows from the \emph{loc.cit.}~that the difference between \eqref{eq-conjugation-equivariance-grassmannian-isomorphism} and \eqref{eq-descent-data-grassmannian-isomorphism} is given by the trivial section
$$
(\Omega \mu) (g) - (\Omega \mu)(g) \simeq 0
$$
of $\Gamma(D_x \text{ mod }\mathring D_x, \deloop^3A(1))$, so \eqref{eq-conjugation-equivariance-grassmannian-isomorphism} and \eqref{eq-descent-data-grassmannian-isomorphism} agree.
\end{proof}

\subsection{Properties of ULA sheaves}

\begin{void}
\label{void-dominant-sharp-cocharacters}
In this subsection, we establish a few useful properties of ULA sheaves over the local Hecke stack, by characterizing them in terms of Schubert cells. These are analogues of results in \cite[\S VI.6]{fargues2021geometrization}.

We remain in the same context as \S\ref{void-split-reductive-killing-pair-context}. Moreover, let us write $\Lambda^{\sharp, +} := \Lambda^{\sharp} \cap \Lambda^+$. This is the monoid parametrizing Schubert cells which support nonzero $(\mathscr G_{\Hec_G}, \zeta)$-twisted constructible complexes (\emph{cf.}~Corollary \ref{cor-schubert-cell-reduction-of-support}).
\end{void}

\begin{prop}
\label{prop-ula-stalkwise-characterization}
Let $x$ be an $S$-point of $X$ ($S \in \Sch$) and $\mathscr A$ be an object of $\derived_{\mathscr G, \zeta}(\Hec_{G, x})$. The following are equivalent:
\begin{enumerate}
	\item $\mathscr A$ is ULA relative to $S$;
	\item for each $\lambda \in \Lambda^{\sharp, +}$, the pullback of $\mathscr A$ along $\varpi^{\lambda} : S \rightarrow \Hec_{G, x}$ is locally constant.
\end{enumerate}
\end{prop}

\begin{proof}
For each $\lambda \in \Lambda^+$, we write $j^{\lambda} : \Hec_{G, x}^{\lambda} \hookrightarrow \Hec_{G, x}$ for the immersion of the Schubert cell. Recall that $\Hec_{G, x}^{\lambda}$ is the quotient of $S$ by an affine group scheme which is an extension of $M^{\lambda}$ by a pro-unipotent affine group scheme (\emph{cf.}~\S\ref{void-schubert-cell-notation}).

(2) $\Rightarrow$ (1). By Corollary \ref{cor-schubert-cell-reduction-of-support}, the pullback of $\mathscr A$ along $\varpi^{\lambda}$ vanishes for $\lambda \in \Lambda^+ \setminus \Lambda^{\sharp, +}$. Thus, if the pullback of $\mathscr A$ along $\varpi^{\lambda}$ is locally constant for $\lambda \in \Lambda^{\sharp, +}$, then $\mathscr A$ is an iterated extension of objects of the form $(j^{\lambda})_! \mathscr B$, where $\mathscr B \in \derived_{\mathscr G, \zeta}(\Hec_{G, x}^{\lambda})$ is locally constant.

We shall argue that $(j^{\lambda})_! \mathscr B$ is ULA relative to $S$. Indeed, by \cite[Corollary 3.9]{MR4630128}, we may assume that $x : S \rightarrow X$ is the spectrum of the henselian local ring of a $\base$-point $\bar x$ of $X$. By \eqref{eq-etale-level-fiber-sequence}, the pullback of $\mu$ along $S \rightarrow X$ then descends to $\Spec\base$. By choosing a uniformizer $\varpi$ at $x$, we realize $\Hec_{G, x}$ as the base change of $\Hec_{G, \bar x}$ along $S \rightarrow \Spec\base$, compatibly with the $A$-gerbe $\mathscr G_{\Hec_G}$. The assertion now follows from base change along
$$
\begin{tikzcd}[column sep = 1em]
	\Hec_{G, x}^{\lambda} \ar[r, "j^{\lambda}"]\ar[d] & \Hec_{G, x} \ar[d] \\ 
	\Hec_{G, \bar x}^{\lambda} \ar[r, "j_{\bar x}^{\lambda}"] & \Hec_{G, \bar x}
\end{tikzcd}
$$
and the fact that any object of $\derived_{\mathscr G, \zeta}(\Hec_{G, \bar x})$ is ULA relative to $\Spec\base$ (\emph{cf.}~Remark \ref{rem-ula-relative-to-geometric-point}).

(1) $\Rightarrow$ (2). Suppose that $\mathscr A$ is ULA relative to $S$. Let $\lambda$ be an element of $\Lambda^+$ whose corresponding Schubert cell is open in the support of $\mathscr A$. In particular, $\lambda \in \Lambda^{\sharp, +}$ (\emph{cf.}~Corollary \ref{cor-schubert-cell-reduction-of-support}). Then the pullback of $\mathscr A$ along $\varpi^{\lambda}$ is locally constant. This implies that $\mathscr B := (j^{\lambda})^*\mathscr A$ is locally constant as well, so $(j^{\lambda})_!\mathscr B$ is ULA relative to $S$ by the argument above.

We may now replace $\mathscr A$ by the cofiber of $(j^{\lambda})_!\mathscr B \rightarrow \mathscr A$ and proceed by induction.
\end{proof}

\begin{void}
Let $x$ be an $S$-point of $X$ ($S \in \Sch$). Given $\mathscr A, \mathscr B \in \derived_{\mathscr G, \zeta}(\Gr_{G, x})$, assumed ULA relative to $S$, the formation of their internal Hom
\begin{equation}
\label{eq-ula-internal-hom-base-change}
\SHom(\mathscr A, \mathscr B) \in \derived(\Gr_{G, x})
\end{equation}
commutes with arbitrary base change in $S \in \Sch$.

Indeed, writing $\mathbf D$ for the Verdier duality functor over $\Gr_{G, x}$ relative to $S$, we may express $\SHom(\mathscr A, \mathscr B)$ as $\mathbf D(\mathscr A \otimes \mathbf D\mathscr B)$ (\emph{cf.}~Remark \ref{rem-internal-hom-via-duality}). Furthermore, $\mathscr A \otimes \mathbf D\mathscr B$ is ULA relative to $S$ by Proposition \ref{prop-ula-stalkwise-characterization}, so the formation of its Verdier dual commutes with arbitrary base change in $S \in \Sch$ (\emph{cf.}~\S\ref{void-verdier-duality}).
\end{void}

\begin{void}
We shall use the base change property of \eqref{eq-ula-internal-hom-base-change} to relate the ``stalk" and ``fiber" of the category of ULA sheaves over the local Hecke stack.

To be more precise, let $\bar x$ be a $\base$-point of $X$ and $x : S \rightarrow X$ be the inverse limit of \'etale neighborhoods of $\bar x$. Pullback yields a functor
\begin{equation}
\label{eq-local-hecke-stack-derived-category-stalk}
	\derived_{\mathscr G, \zeta}(\Hec_{G, x}) \rightarrow \derived_{\mathscr G, \zeta}(\Hec_{G, \bar x}).
\end{equation}
\end{void}

\begin{prop}
\label{prop-local-hecke-stack-derived-category-stalk}
The functor \eqref{eq-local-hecke-stack-derived-category-stalk} restricts to an equivalence on the full subcategories consisting of ULA objects (relative to $S$, respectively $\base$).
\end{prop}

\begin{proof}
Given ULA objects $\mathscr A, \mathscr B \in \derived_{\mathscr G, \zeta}(\Gr_{G, x})$, we have an identification
\begin{equation}
\label{eq-affine-grassmannian-internal-hom-fiber}
(\pi_*\SHom(\mathscr A, \mathscr B))_{\bar x} \simeq \Hom(\mathscr A_{\bar x}, \mathscr B_{\bar x}),
\end{equation}
where $\pi : \Gr_{G, x} \rightarrow S$ is the structural map and $\mathscr A_{\bar x}, \mathscr B_{\bar x}$ are the pullbacks of $\mathscr A, \mathscr B$ to $\Gr_{G, \bar x}$. Indeed, \eqref{eq-affine-grassmannian-internal-hom-fiber} follows from the ind-properness of $\pi$ and the fact that \eqref{eq-ula-internal-hom-base-change} commutes with arbitrary base change.

From \eqref{eq-affine-grassmannian-internal-hom-fiber}, we deduce an isomorphism
\begin{equation}
\label{eq-affine-grassmannian-hom-fiber}
\Hom(\mathscr A, \mathscr B) \simeq \Hom(\mathscr A_{\bar x}, \mathscr B_{\bar x}),
\end{equation}
where the Homs are taken in $\derived_{\mathscr G, \zeta}(\Gr_{G, x})$, respectively $\derived_{\mathscr G, \zeta}(\Gr_{G, \bar x})$. Now, assuming that $\mathscr A$, $\mathscr B$ are endowed with $L_x^+G$-equivariance, we may compute their Homs in $\derived_{\mathscr G, \zeta}(\Hec_{G, x})$, respectively $\derived_{\mathscr G, \zeta}(\Hec_{G, \bar x})$ by taking $L_x^+G$-, respectively $L_{\bar x}^+G$-invariants on \eqref{eq-affine-grassmannian-hom-fiber}. They coincide by base change properties of the cohomology of $\deloop G$. This implies that \eqref{eq-local-hecke-stack-derived-category-stalk} is fully faithful on ULA objects.

To show that it is essentially surjective, we observe that $\derived_{\mathscr G, \zeta}(\Hec_{G, \bar x})$ is generated by $!$-extensions of constant sheaves from Schubert cells. These objects lift along \eqref{eq-local-hecke-stack-derived-category-stalk} to objects of the same kind, which are ULA by Proposition \ref{prop-ula-stalkwise-characterization}.
\end{proof}

\subsection{Generators}

\begin{void}
We fix $T \subset B \subset G$ as in \S\ref{void-split-reductive-killing-pair-context}. Let $S \in \Sch$ and $x$ be an $S$-point of $X$.

For each $\lambda \in \Lambda^{\sharp, +}$ (\emph{cf.}~the notation of \S\ref{void-dominant-sharp-cocharacters}), we denote by $\Lis_{\mathscr G_{\varpi^{\lambda}}, \zeta}(S)$ the $\coeff$-linear abelian category of $(\mathscr G_{\varpi^{\lambda}}, \zeta)$-twisted $\coeff$-local systems over $S$.

Pullback along the structural morphism $\pi : \Hec_{G, x}^{\lambda} \rightarrow S$, in view of the identification \eqref{eq-gerbe-schubert-cell-sharp-identification}, yields a functor
\begin{equation}
\label{eq-schubert-cell-pullback-functor}
\pi^* : \Lis_{\mathscr G_{\varpi^{\lambda}}, \zeta}(S) \rightarrow \derived_{\mathscr G, \zeta}(\Hec_{G, x}^{\lambda}).
\end{equation}
\end{void}

\begin{void}
\label{void-standard-costandard}
To make the target of \eqref{eq-schubert-cell-pullback-functor} perverse\footnote{As always, this means perverse relative to $S$ after pulling back to $\Gr_{G, x}^{\lambda}$, \emph{cf.}~\S\ref{void-satake-subcategory}.} relative to $S$ and to incorporate the ``correct" amount of Tate twists, we shall instead consider the functor
\begin{equation}
\label{eq-relevant-schubert-cell-constant-sheaves}
	\pi^*(\langle \check{\rho}, \lambda \rangle)[\langle 2\check{\rho}, \lambda\rangle] : \Lis_{\mathscr G_{\varpi^{\lambda}}, \zeta}(S) \rightarrow \derived_{\mathscr G, \zeta}(\Hec_{G, x}^{\lambda}).
\end{equation}
Since $\langle\check{\rho}, \lambda\rangle$ is a half-integer in general, the formation of \eqref{eq-relevant-schubert-cell-constant-sheaves} invokes $\coeff(\frac{1}{2})$.

For each $\mathscr E \in \Lis_{\mathscr G_{\varpi^{\lambda}}, \zeta}(S)$, we write $\widetilde{\mathscr E}$ for its image under \eqref{eq-relevant-schubert-cell-constant-sheaves}. Denote by $j^{\lambda} : \Hec_{G, x}^{\lambda} \rightarrow \Hec_{G, x}$ the locally closed immersion. We may form the following functors, called the \emph{standard}, respectively \emph{costandard functors}:
\begin{align*}
	\Delta^{\lambda} &: \Lis_{\mathscr G_{\varpi^{\lambda}}, \zeta}(S) \rightarrow \derived_{\mathscr G, \zeta}(\Hec_{G, x}),\quad \mathscr E \mapsto {}^pH^0 (j^{\lambda})_! \widetilde{\mathscr E} \\
	\nabla^{\lambda} &: \Lis_{\mathscr G_{\varpi^{\lambda}}, \zeta}(S) \rightarrow \derived_{\mathscr G, \zeta}(\Hec_{G, x}),\quad \mathscr E \mapsto {}^pH^0 (j^{\lambda})_* \widetilde{\mathscr E}
\end{align*}

Denote by $\mathbf D$ the Verdier duality functor over $\Gr_{G, x}$ relative to $S$ (\emph{cf.}~\S\ref{void-verdier-duality}). We have the following basic result concerning $\Delta^{\lambda}$ and $\nabla^{\lambda}$.
\end{void}

\begin{prop}
\label{prop-standard-costandard}
Given $\lambda \in \Lambda^{\sharp, +}$ and $\mathscr E \in \Lis_{\mathscr G_{\varpi^{\lambda}}, \zeta}(S)$, there holds
\begin{enumerate}
	\item the objects $\Delta^{\lambda}(\mathscr E)$, $\nabla^{\lambda}(\mathscr E)$ belong to $\Sat_{\mathscr G, \zeta}(\Hec_{G, x})$;
	\item the formations of $\Delta^{\lambda}(\mathscr E)$, $\nabla^{\lambda}(\mathscr E)$ commute with base change in $S \in \Sch$;
	\item there is a canonical isomorphism
	\begin{equation}
	\label{eq-standard-costandard-verdier-dual}
	\mathbf D\Delta^{\lambda}(\mathscr E) \simeq \nabla^{\lambda}(\mathscr E^{\vee}),\text{ for } \mathscr E^{\vee} := \SHom(\mathscr E, \coeff) \in \Lis_{-\mathscr G_{\varpi^{\lambda}}, \zeta}(S).
	\end{equation}
\end{enumerate}
\end{prop}

\begin{proof}
By Proposition \ref{prop-ula-stalkwise-characterization}, the complex $(j^{\lambda})_!\widetilde{\mathscr E} \in \derived_{\mathscr G, \zeta}(\Hec_{G, x})$ is ULA relative to $S$, for any $\mathscr E \in \Lis_{\mathscr G_{\varpi^{\lambda}}, \zeta}(S)$. We have an isomorphism
\begin{equation}
\label{eq-verdier-duality-schubert-cell-shriek-star-extensions}
\mathbf D(j^{\lambda})_! \simeq (j^{\lambda})_*\mathbf D,
\end{equation}
where by a slight abuse of notation, we also use $\mathbf D$ for the Verdier duality functor over $\Gr_{G, x}^{\lambda}$ relative to $S$. It follows from \eqref{eq-verdier-duality-schubert-cell-shriek-star-extensions} that $(j^{\lambda})_*\widetilde{\mathscr E}$ is also ULA relative to $S$. Moreover, the formation of both $(j^{\lambda})_!\widetilde{\mathscr E}$ and $(j^{\lambda})_*\widetilde{\mathscr E}$ commutes with base change in $S \in \Sch$.

Statement (1) now follows from the fact that truncation functor ${}^pH^0$ preserves universal local acyclicity (\emph{cf.}~Remark \ref{rem-ula-perverse-truncation}).

Statement (2) follows because ${}^pH^0$ commutes with base change in $S \in \Sch$.

To prove statement (3), we shall argue that the morphism
\begin{equation}
\label{eq-standard-truncation-morphism-verdier-dual}
\mathbf D\Delta^{\lambda}(\mathscr E) \rightarrow \mathbf D (j^{\lambda})_!\widetilde{\mathscr E},
\end{equation}
obtained by dualizing the degree-$0$ truncation $(j^{\lambda})_!\widetilde{\mathscr E} \rightarrow \Delta^{\lambda}(\mathscr E)$, is itself the degree-$0$ truncation of the connective complex $\mathbf D (j^{\lambda})_!\widetilde{\mathscr E}$.

This claim will yield the isomorphism \eqref{eq-standard-costandard-verdier-dual}, in view of the isomorphisms
\begin{align*}
	\mathbf D\Delta^{\lambda}(\mathscr E) & \simeq {}^pH^0 \mathbf D(j^{\lambda})_!\widetilde{\mathscr E} \\
	& \simeq {}^pH^0 (j^{\lambda})_*\mathbf D(\widetilde{\mathscr E}) \simeq \nabla^{\lambda}(\mathscr E^{\vee}),
\end{align*}
where we used \eqref{eq-verdier-duality-schubert-cell-shriek-star-extensions} and the fact that $\Gr_{G, x}^{\lambda} \rightarrow S$ is smooth of relative dimension $\langle 2\check{\rho}, \lambda\rangle$.

To prove the claim, we shall apply $\CT_B$ to \eqref{eq-standard-truncation-morphism-verdier-dual}. Realizing $\CT_B$ as a hyperbolic localization (\emph{cf.}~the proof of Proposition \ref{prop-constant-term-reflective-properties}), the image may be identified with
$$
\mathbf D \CT_{B^-}\Delta^{\lambda}(\mathscr E) \rightarrow \mathbf D \CT_{B^-} (j^{\lambda})_!\widetilde{\mathscr E},
$$
where $\mathbf D$ now stands for Verdier duality over $\Gr_{T, x}$ relative to $S$, \emph{i.e.}~the functor $\SHom(\cdot, \coeff)$. Since $\CT_{B^-}$ is $t$-exact (\emph{cf.}~Proposition \ref{prop-constant-term-reflective-properties}), we may identify $\CT_{B^-}\Delta^{\lambda}(\mathscr E)$ with the degree-$0$ trunction of the coconnective complex $\CT_{B^-}(j^{\lambda})_!\widetilde{\mathscr E}$. It follows that $\mathbf D \CT_{B^-}\Delta^{\lambda}(\mathscr E)$ is the degree-$0$ truncation of the connective complex $\mathbf D \CT_{B^-} (j^{\lambda})_!\widetilde{\mathscr E}$. The claim now follows from the $t$-exactness and conservativity of $\CT_B$ (\emph{cf.}~Proposition \ref{prop-constant-term-reflective-properties}).
\end{proof}

\begin{void}
\label{void-affine-grassmannian-ic-sheaf}
Given $\lambda \in \Lambda^{\sharp, +}$, we define the functor, called the \emph{intersection cohomology}
\begin{equation}
\label{eq-ic-functor}
\IC^{\lambda} : \Lis_{\mathscr G_{\varpi^{\lambda}}, \zeta}(S) \rightarrow \Sat_{\mathscr G, \zeta}(\Hec_{G, x}),
\end{equation}
which assigns to each $\mathscr E \in \Lis_{\mathscr G_{\varpi^{\lambda}}, \zeta}(S)$ the image of the natural morphism $\Delta^{\lambda}(\mathscr E) \rightarrow \nabla^{\lambda}(\mathscr E)$ in $\Sat_{\mathscr G, \zeta}(\Hec_{G, x})$ (\emph{cf.}~Proposition \ref{prop-standard-costandard}). Then the formation of $\IC^{\lambda}(\mathscr E)$ also commutes with base change in $S \in \Sch$. Summing \eqref{eq-ic-functor} yields a functor of $\coeff$-linear abelian categories
\begin{equation}
\label{eq-satake-category-generator-sum}
\bigoplus_{\lambda \in \Lambda^{\sharp, +}} \Lis_{\mathscr G_{\varpi^{\lambda}}, \zeta}(S) \rightarrow \Sat_{\mathscr G, \zeta}(\Hec_{G, x});
\end{equation}

By construction, we also have the following morphisms of functors:
\begin{equation}
\label{eq-standard-ic-costandard}
	\Delta^{\lambda} \twoheadrightarrow \IC^{\lambda} \hookrightarrow \nabla^{\lambda}.
\end{equation}
The remainder of this section is devoted to the proof of the following result.
\end{void}

\begin{thm}
\label{thm-semisimplicity}
The functor \eqref{eq-satake-category-generator-sum} is an equivalence.
\end{thm}

\begin{void}
Theorem \ref{thm-semisimplicity} has a number of pleasant consequences.
\end{void}

\begin{cor}
\label{cor-semisimplicity}
If $x$ is a geometric point of $X$, then $\Sat_{\mathscr G, \zeta}(\Hec_{G, x})$ is semisimple.
\end{cor}

\begin{proof}
When $S$ is the spectrum of an algebraically closed field, the category $\Lis_{\mathscr G_{\varpi^{\lambda}}, \zeta}(S)$ is \emph{non-canonically} equivalent to the category of finite-dimensional $\coeff$-vector spaces. Hence the assertion follows from Theorem \ref{thm-semisimplicity}.
\end{proof}

\begin{cor}
\label{cor-standard-ic-agreement}
The morphisms in \eqref{eq-standard-ic-costandard} are isomorphisms.
\end{cor}

\begin{proof}
Since the formation of $\Delta^{\lambda}$ and $\nabla^{\lambda}$ commutes with base change in $S$ (\emph{cf.}~Proposition \ref{prop-standard-costandard}), we reduce to the case $S = \Spec\base$. The same result implies that $\IC^{\lambda} \hookrightarrow \nabla^{\lambda}$ is the Verdier dual of $\Delta^{\lambda} \twoheadrightarrow \IC^{\lambda}$, so it suffices to prove that the latter is an isomorphism.

Let $\mathscr E$ be a $(\mathscr G_{\varpi^{\lambda}}, \zeta)$-twisted local system over $S$. The kernel $\mathscr A$ of the surjection $\Delta^{\lambda}(\mathscr E) \twoheadrightarrow \IC^{\lambda}(\mathscr E)$ is supported on the complement of $\Gr_{G, x}^{\lambda}$ in its closure $\overline{\Gr}{}_{G, x}^{\lambda}$. By Corollary \ref{cor-semisimplicity}, the injection $\mathscr A \hookrightarrow \Delta^{\lambda}(\mathscr E)$ admits a retraction $\Delta^{\lambda}(\mathscr E) \twoheadrightarrow \mathscr A$, but any such morphism vanishes.
\end{proof}

\subsection{Monodromic affine Hecke categories}
\label{sec-monodromic-affine-hecke-categories}

\begin{void}
The main input in Theorem \ref{thm-semisimplicity} is the parity vanishing statement for intersection cohomology sheaves on the affine flag variety. Contrary to the untwisted situation, the proof of this statement requires considering sheaves twisted by a \emph{family} of $A$-gerbes.

We fix $T\subset B \subset G$ as in \S\ref{void-split-reductive-killing-pair-context} and a $\base$-point $x$ of $X$. Denote by $W$ the Weyl group of $(G, T)$. It naturally acts on the cocharacter lattice $\Lambda$, and we may form the extended affine Weyl group $W^{\aff} := \Lambda \rtimes W$, whose elements may be represented by $\translate(\lambda) w$ ($\lambda \in \Lambda$, $w \in W$) where $\translate : \Lambda \subset W^{\aff}$ is the natural inclusion.

The group $W^{\aff}$ admits a length function $\ell : W^{\aff} \rightarrow \integers_{\ge 0}$, which may be expressed explicitly by Matsumoto's formula (\emph{cf.}~\cite[\S4.2]{MR3839695}).
\end{void}

\begin{void}[The $\infty$-category ${}_{\chi_1}(\derived_{\zeta}^{\aff})_{\chi_2}$]
Denote by $I \subset L_x^+G$ the preimage of $B^-$ under the projection map $L_x^+G \rightarrow G$, where $B^-\subset G$ is the Borel subgroup \emph{opposite} to $B$.

The \emph{Iwahori--Hecke stack} $\Hec_{G, x}^{\aff} := I \backslash L_xG /I$ admits the Bruhat stratification indexed by $W^{\aff}$. Furthermore, it is equipped with three projection maps
\begin{equation}
\label{eq-iwahori-hecke-stack-projections}
\begin{tikzcd}[column sep = 0.5em]
	& \Hec_{G, x}^{\aff} \ar[rr, "p"]\ar[dl, swap, "\overset{\leftarrow}{\pi}"]\ar[dr, "\overset{\rightarrow}{\pi}"] & & \Hec_{G, x} \\
	\deloop I & & \deloop I
\end{tikzcd}
\end{equation}
where $p$ is induced from the inclusion $I \subset L_x^+G$, and $\overset{\leftarrow}{\pi}$, $\overset{\rightarrow}{\pi}$ are defined by the left, respectively right $I$-action on $L_xG$. Denote by $\mathscr G^{\aff}$ the pullback of the $A$-gerbe $\mathscr G_{\Hec_G}$ along $p$.

Since the projection $I \twoheadrightarrow B^- \twoheadrightarrow T$ has a pro-unipotent kernel, we have a canonical isomorphism of abelian groups
$$
\Hom(\Lambda, A(-1)) \simeq \Maps_*(\deloop I, \deloop^2 A),
$$
so any character $\chi : \Lambda \rightarrow A(-1)$ defines an $A$-gerbe $\chi\otimes\Psi$ over $\deloop I$ (\emph{cf.}~\S\ref{void-gerbe-schubert-cell-notation}).

Given a pair of characters $\chi_1, \chi_2 \in \Hom(\Lambda, A(-1))$, we write ${}_{\chi_1}(\derived_{\zeta}^{\aff})_{\chi_2}$ for the $\infty$-category of $(\overset{\leftarrow}{\pi}{}^*(\chi_1\otimes\Psi) + \mathscr G^{\aff} - \overset{\rightarrow}{\pi}{}^*(\chi_2\otimes\Psi), \zeta)$-twisted constructible complexes over $\Hec_{G, x}^{\aff}$ (\emph{cf.}~\S\ref{void-twisted-sheaves-hecke-stack}), which we refer to as the \emph{monodromic affine Hecke category}.
\end{void}

\begin{void}[Convolution]
\label{void-monodromic-affine-hecke-convolution}
The monodromic affine Hecke categories have the following ``convolution product": Given a triple of characters $\chi_1, \chi_2, \chi_3 \in \Hom(\Lambda, A(-1))$, there is a functor
\begin{align}
\notag
{}_{\chi_1} (\derived_{\zeta}^{\aff})_{\chi_2} \times {}_{\chi_2}(\derived_{\zeta}^{\aff})_{\chi_3} &\rightarrow {}_{\chi_1} (\derived^{\aff}_{\zeta})_{\chi_3},\\
\label{eq-affine-hecke-convolution}
 \mathscr A,\mathscr B &\mapsto \mathscr A \circ \mathscr B := m_!(p_1^*\mathscr A \otimes p_2^*\mathscr B),
\end{align}
where $m$, $p_1$, $p_2$ refer to the three projections from $\Hec_{G, x}^{\aff, [2]}$ to $\Hec_{G, x}^{\aff}$ defined by the natural groupoid structure on $\Hec_{G, x}^{\aff}$ (\emph{cf.}~\S\ref{void-convolution-product}).

The formation of \eqref{eq-affine-hecke-convolution} relies on the following identification of $A$-gerbes over $\Hec_{G, x}^{\aff, [2]}$:
\begin{align*}
	(p_1)^*&(\overset{\leftarrow}{\pi}{}^*(\chi_1\otimes\Psi) + \mathscr G^{\aff} - \overset{\rightarrow}{\pi}{}^*(\chi_2\otimes\Psi)) + (p_2)^*(\overset{\leftarrow}{\pi}{}^*(\chi_2\otimes\Psi) + \mathscr G^{\aff} - \overset{\rightarrow}{\pi}{}^*(\chi_3\otimes\Psi)) \\
	&\simeq (p_1)^*\overset{\leftarrow}{\pi}{}^*(\chi_1\otimes\Psi) + (p_1)^*\mathscr G^{\aff} + (p_2)^*\mathscr G^{\aff} - (p_2)^*\overset{\rightarrow}{\pi}{}^*(\chi_3\otimes\Psi) \\
	&\simeq m^*(\overset{\leftarrow}{\pi}{}^*(\chi_1\otimes\Psi) + \mathscr G^{\aff} - \overset{\rightarrow}{\pi}{}^*(\chi_3\otimes\Psi)),
\end{align*}
where we used the identification $\overset{\rightarrow}{\pi} \circ p_1 = \overset{\leftarrow}{\pi} \circ p_2$ in the first isomorphism to cancel out the $(\chi_2\otimes\Psi)$-term, and the multiplicative structure on $\mathscr G^{\aff}$ (\emph{cf.}~Lemma \ref{lem-local-hecke-stack-gerbe-multiplicative}) with $\overset{\leftarrow}{\pi} \circ p_1 = \overset{\leftarrow}{\pi} \circ m$, $\overset{\rightarrow}{\pi}\circ p_2 = \overset{\rightarrow}{\pi} \circ m$ for the second isomorphism.

The convolution product \eqref{eq-affine-hecke-convolution} admits an associativity constraint in the evident sense, although we will not attempt to lift it to full homotopy coherence data.
\end{void}

\begin{void}[$W^{\aff}$-action on characters]
\label{void-affine-weyl-group-action-on-characters}
Observe that $W^{\aff}$ acts on the set $\Hom(\Lambda, A(-1))$: Given an element $v = \translate(\lambda)w \in W^{\aff}$ and $\chi \in \Hom(\Lambda, A(-1))$, we define
\begin{equation}
\label{eq-affine-weyl-group-action-on-characters}
v(\chi) := w(\chi) - b(\lambda, \cdot),
\end{equation}
where $w(\chi)$ is the contragredient action $w(\chi) := \chi(w^{-1}(\cdot))$.

The fact that the formula \eqref{eq-affine-weyl-group-action-on-characters} defines a $W^{\aff}$-action follows from the $W$-invariance of $b$, which is a consequence of \eqref{eq-strict-weyl-invariance}.
\end{void}

\begin{void}
\label{void-bruhat-cell}
For each $v = \translate(\lambda) w \in W^{\aff}$, we have the Bruhat cell
$$
j^v : \Hec_{G, x}^{\aff, v} \rightarrow \Hec_{G, x}^{\aff}.
$$
The corresponding Bruhat cell $\Fl_{G, x}^v$ in the affine flag variety $\Fl_{G, x} := L_xG/I$ has dimension $\ell(v)$. We shall consider the restrictions of the morphisms in \eqref{eq-iwahori-hecke-stack-projections} along $j^v$ without changing the notation. For a pair of characters $\chi_1, \chi_2 \in \Hom(\Lambda, A(-1))$, we write ${}_{\chi_1}(\derived_{\zeta}^{\aff, v})_{\chi_2}$ for the $\infty$-category of $(\overset{\leftarrow}{\pi}{}^*(\chi_1\otimes\Psi) + \mathscr G^{\aff} - \overset{\rightarrow}{\pi}{}^*(\chi_2\otimes\Psi), \zeta)$-twisted constructible complexes over $\Hec_{G, x}^{\aff, v}$.

By Proposition \ref{prop-gerbe-schubert-cell}, the restriction of the $A$-gerbe $\mathscr G^{\aff}$ along $j^v$ may be identified as
\begin{equation}
\label{eq-gerbe-one-dimensional-bruhat-cell}
(j^v)^*\mathscr G^{\aff} \simeq \overset{\leftarrow}{\pi}{}^*(b(\lambda, \cdot)\otimes\Psi) + (\mathscr G_{\varpi^{\lambda}}|_{\Hec_{G, x}^{\aff, v}}).
\end{equation}
This allows us to define pullback along $\pi : \Hec_{G, x}^{\aff, v} \rightarrow \Spec\base \simeq x$ as a functor
\begin{equation}
\label{eq-bruhat-cell-pullback-functor}
\pi^* : \Lis_{\mathscr G_{\varpi^{\lambda}}, \zeta}(x) \rightarrow {}_{-b(\lambda, \cdot)}(\derived_{\zeta}^{\aff, v})_0.
\end{equation}
\end{void}

\begin{void}
To proceed further, we will make additional choices. Let $\varpi$ be a uniformizer at $x$. It allows us to identify $L_x^+G$ (respectively $L_xG$) with $G\arc{\varpi}$ (respectively $G\loo{\varpi}$).

We furthermore choose a lift $\dot w \in G$ for each $w \in W$. With these choices, every $v = \translate(\lambda) w \in W^{\aff}$ admits a lift to $G\loo{\varpi}$, namely $\dot v := \varpi^{\lambda} \dot w$.

We shall use $\dot v$ to relate the left and right $I$-action on $I v I \subset G\loo{\varpi}$. Namely, consider the $(I \times I)$-action on $G\loo{\varpi}$ defined by $(b_1, b_2) \cdot g := b_1 g (b_2)^{-1}$. The stablizer of $\dot v$ is the subgroup $I_{\dot v} \subset I\times I$ consisting of points $(b_1, b_2)$ satisfying
$$
b_1 = \dot v b_2 \dot v^{-1}.
$$
In particular, the reduction $(\bar b_1, \bar b_2) \in T \times T$ of $(b_1, b_2)$ satisfies
$$
\bar b_1 = w(\bar b_2).
$$

We thus obtain a commutative diagram
\begin{equation}
\label{eq-bruhat-cell-left-right-torus-actions}
\begin{tikzcd}[column sep = 0.5em]
	\Hec_{G, x}^{\aff, v} \ar[rr, phantom, "\simeq"] & & \deloop(I_{\dot v}) \ar[dl, swap, "\overset{\leftarrow}{\pi}"] \ar[dr, "\overset{\rightarrow}{\pi}"] \\
	& \deloop T \ar[rr, "w^{-1}"] & & \deloop T
\end{tikzcd}
\end{equation}
where $\overset{\leftarrow}{\pi}$, $\overset{\rightarrow}{\pi}$ stand for the compositions of the morphisms in \eqref{eq-iwahori-hecke-stack-projections} with the projection onto $\deloop T$. It follows from \eqref{eq-bruhat-cell-left-right-torus-actions} that we have an isomorphism of $A$-gerbes over $\Hec_{G, x}^{\aff, v}$:
\begin{equation}
\label{eq-bruhat-cell-left-right-monodromy-identification}
\overset{\leftarrow}{\pi}{}^* (w(\chi)\otimes\Psi) \simeq \overset{\rightarrow}{\pi}{}^*(\chi\otimes\Psi).
\end{equation}
\end{void}

\begin{void}
\label{void-flag-variety-standard-costandard}
Using \eqref{eq-bruhat-cell-left-right-monodromy-identification}, we may view the pullback functor \eqref{eq-bruhat-cell-pullback-functor} as a functor
\begin{equation}
\label{eq-bruhat-cell-pullback-functor-monodromic}
\pi^* : \Lis_{\mathscr G_{\varpi^{\lambda}}, \zeta}(x) \rightarrow {}_{v(\chi)}(\derived_{\zeta}^{\aff, v})_{\chi}
\end{equation}
defined for any $v = \translate(\lambda) w \in W^{\aff}$ and $\chi \in \Hom(\Lambda, A(-1))$. Here, $v(\chi)$ refers to the $W^{\aff}$-action on $\Hom(\Lambda, A(-1))$ constructed in \S\ref{void-affine-weyl-group-action-on-characters}.

As in \S\ref{void-standard-costandard}, we shall incorporate shifts and Tate twists: For any $\mathscr E \in \Lis_{\mathscr G_{\varpi^{\lambda}}, \zeta}(x)$, we write $\widetilde{\mathscr E} := \pi^*\mathscr E(\frac{\ell(v)}{2})[ \ell(v)]$. Then we define the \emph{standard} and \emph{costandard functors}
\begin{align*}
	\Delta^v_{\chi} &: \Lis_{\mathscr G_{\varpi^{\lambda}}, \zeta}(x) \rightarrow {}_{v(\chi)}(\derived_{\zeta}^{\aff})_{\chi},\quad \mathscr E \mapsto (j^v)_! \widetilde{\mathscr E}; \\
	\nabla^v_{\chi} &: \Lis_{\mathscr G_{\varpi^{\lambda}}, \zeta}(x) \rightarrow {}_{v(\chi)}(\derived_{\zeta}^{\aff})_{\chi},\quad \mathscr E \mapsto (j^v)_* \widetilde{\mathscr E}.
\end{align*}

Since $j^v$ is affine, the images of $\Delta^v_{\chi}$ and $\nabla^v_{\chi}$ are perverse (after pulling back to the indscheme $\Fl_{G, x}$). Thus, we may also define the \emph{intersection cohomology functor}
\begin{equation}
\label{eq-affine-flag-variety-ic-functor}
\IC^v_{\chi} : \Lis_{\mathscr G_{\varpi^{\lambda}}, \zeta}(x) \rightarrow {}_{v(\chi)}(\derived_{\zeta}^{\aff})_{\chi}, 
\end{equation}
sending $\mathscr E$ to the image of the natural morphism $(j^v)_!\widetilde{\mathscr E} \rightarrow (j^v)_*\widetilde{\mathscr E}$.

The key statement about $\IC_{\chi}^v$, which we shall prove, is the ``parity vanishing" of its fiber cohomologies. It is an analogue of \cite[Proposition 3.12]{MR4108915} for \'etale levels.
\end{void}

\begin{prop}
\label{prop-parity-vanishing}
Fix $v \in W^{\aff}$ and $\chi \in \Hom(\Lambda, A(-1))$. The fibers of $\IC_{\chi}^v$ over $\Fl_{G, x}$ are concentrated in cohomological degrees $n$ such that $n + \ell(v)$ is even.
\end{prop}

\begin{rem}
For our application of Proposition \ref{prop-parity-vanishing} (\emph{cf.}~\S\ref{sec-parity-vanishing-application}), we will \emph{only} need the case $\chi = 0$. However, in order to prove Proposition \ref{prop-parity-vanishing} for $\chi = 0$, we must work with all possible choices of $\chi$. This feature is special to the twisted context.
\end{rem}

\subsection{Parity vanishing}
\label{sec-parity-vanishing}

\begin{void}
\label{void-parity-vanishing-context}
The goal of this subsection is to prove Proposition \ref{prop-parity-vanishing}, so we retain the notation of \S\ref{sec-monodromic-affine-hecke-categories}. We begin by describing the restriction of $\mathscr G^{\aff}$ to dimension-$1$ Bruhat varieties.

Let $S \subset W^{\aff}$ be the subset of length-$1$ elements. By Matsumoto's formula, an element of $S$ is one of the following kinds
\begin{enumerate}
	\item a finite simple reflection $s_{\check{\alpha}} \in W$ associated to $\check{\alpha} \in \check{\Delta}$, viewed as an element of $W^{\aff}$ via the natural inclusion $W \subset W^{\aff}$;
	\item a product $\translate(\theta) s_{\check{\theta}}$, where $\check{\theta}$ is a highest root (which we call an \emph{affine simple reflection}).
\end{enumerate}

For $s \in S$, write $P_s \subset G\loo{\varpi}$ for the corresponding subminimal parahoric group scheme. The quotient $P_s/I$ is isomorphic to $\mathbb P^1$, so isomorphism classes of $A$-gerbes over $P_s/I$ are in bijection with $A(-1)$, where $a \in A(-1)$ corresponds to the class of $a\otimes \Psi(\mathscr O_{\mathbb P^1}(1))$.
\end{void}

\begin{void}
Let us determine the restriction of $\mathscr G^{\aff}$ along the closed immersion
\begin{equation}
\label{eq-dimension-one-bruhat-variety}
i_s : P_s/I \hookrightarrow G\loo{\varpi}/I =: \Fl_{G, x}.
\end{equation}

The following result is an analogue of the computation of the degrees of line bundles on $\Gr_{G, x}$, for $G$ simple and simply connected, after pulling back to $P_s/I$ (\emph{cf.}~\cite[Theorem 7]{MR1961134}). It can be deduced from \emph{loc.cit.}, but let us offer a proof internal to $A$-gerbes.
\end{void}

\begin{lem}
\label{lem-gerbe-dimension-one-bruhat-variety}
The isomorphism class $a(s) \in A(-1)$ of $(i_s)^*\mathscr G^{\aff}$ is given by
$$
a(s)
=
\begin{cases}
	0 & \text{$s = s_{\check{\alpha}}$ for $\check{\alpha} \in \check{\Delta}$} \\
	Q(\theta) & \text{$s = \translate(\theta)s_{\check{\theta}}$ for $\check{\theta}$ a highest root}
\end{cases}
$$
\end{lem}

\begin{proof}
If $s = s_{\check{\alpha}}$ is a finite simple reflection, then $P_s \subset G\arc{\varpi}$ is the pre-image of the standard parabolic subgroup $P_{\check{\alpha}} \subset G$ corresponding to $\check{\alpha} \in \check{\Delta}$. This implies that \eqref{eq-dimension-one-bruhat-variety} factors as
\begin{equation}
\label{eq-dimension-one-bruhat-variety-finite-simple-reflection}
P_s/I \simeq P_{\check{\alpha}}/B \simeq G\arc{\varpi}/I \hookrightarrow G\loo{\varpi}/I.
\end{equation}
Since $\mathscr G^{\aff}$ is trivial over $G\arc{\varpi}/I$, we see that $(i_s)^*\mathscr G^{\aff}$ is trivial.

Suppose now that $s = \translate(\theta) s_{\check{\theta}}$ is an affine simple reflection. Let us identify the composition of \eqref{eq-dimension-one-bruhat-variety} with the projection onto $\Gr_{G, x}$:
\begin{equation}
\label{eq-affine-simple-reflection-line-in-grassmannian}
P_s/I \rightarrow \Gr_{G, x}.
\end{equation}
We first extend $\theta$ to a morphism $f_{\check{\theta}} : G_{\check{\theta}} \rightarrow G$, where $G_{\check{\theta}}$ is simply connected of semisimple rank $1$, covering the root subgroups $N_{\check{\theta}}$, $N_{-\check{\theta}}$ of $G$. The morphism \eqref{eq-affine-simple-reflection-line-in-grassmannian} factors through the morphism $\Gr_{G_{\check{\theta}}, x} \rightarrow \Gr_{G, x}$ induced from $f_{\check{\theta}}$, so we may replace $G$ by $G_{\check{\theta}}$ and assume that $G$ is simply connected with unique simple coroot $\theta$.

Denote by $\omega$ the fundamental coweight of $T_{\adjoint}$ and by $K_{\lambda} \subset L_xG$ the subgroup scheme $\varpi^{\lambda} L_x^+G \varpi^{-\lambda}$ for each $\lambda \in \Lambda_{\adjoint}$. Then $K_{\omega}\subset P_s$ and this inclusion induces an isomorphism
\begin{equation}
\label{eq-flag-variety-line-affine-grassmannian}
K_{\omega} / K_{\omega} \cap K_0 \simeq P_s/I,
\end{equation}
under which \eqref{eq-affine-simple-reflection-line-in-grassmannian} corresponds to the natural inclusion
\begin{equation}
\label{eq-affine-simple-reflection-line-in-grassmannian-explicit}
K_{\omega}/K_{\omega} \cap K_0 \subset \Gr_{G, x}.
\end{equation}

Let us form the group scheme $\widetilde G := G \rtimes T_{\adjoint}$.\footnote{To be more explicit, one may fix an isomorphism $G \simeq \SL_2$ and realize $\widetilde G$ as $\GL_2$.} Its center is the anti-diagonal copy of $T$, so the natural embedding $G \subset \widetilde G$ induces an isomorphism of their adjoint forms. The inclusion \eqref{eq-affine-simple-reflection-line-in-grassmannian-explicit} for $G$, $\widetilde G$ are related by the commutative diagram
\begin{equation}
\label{eq-maximal-compact-extension}
\begin{tikzcd}[column sep = 1em]
	K_{\omega} / K_{\omega} \cap K_0 \ar[r, phantom, "\subset"]\ar[d] & \Gr_{G, x} \ar[d] \\
	\widetilde K_{\omega} / \widetilde K_{\omega} \cap \widetilde K_0 \ar[r, phantom, "\subset"] & \Gr_{\widetilde G, x}
\end{tikzcd}
\end{equation}
Moreover, the lower inclusion in \eqref{eq-maximal-compact-extension} coincides with the $\varpi^{(0, \omega)}$-translate of the Schubert cell in $\Gr_{\widetilde G, x}$ containing $\varpi^{(0, -\omega)}$, with respect to the maximal torus $T\times T_{\adjoint}$ of $\widetilde G$.

Since $G$ is simply connected, we may extend $\mu$ to an \'etale level $\widetilde{\mu}$ of $\widetilde G$ (\emph{cf.}~\S\ref{void-whittaker-torsor}). Evaluating the symmetric form \eqref{eq-simply-connected-etale-level-extension-bilinear-form} of $\widetilde{\mu}$ at $(0, -\omega)$, we obtain the character
\begin{equation}
\label{eq-extended-adjoint-coweight-character}
\Lambda \oplus \Lambda_{\adjoint} \rightarrow A(-1)
\end{equation}
which annihilates the summand $\Lambda_{\adjoint}$ and sends $\theta \in \Lambda$ to $-Q(\theta)$. It now follows from Proposition \ref{prop-gerbe-schubert-cell} (applied to $\widetilde G$) that the restriction of $\mathscr G_{\Gr_G}$ along \eqref{eq-affine-simple-reflection-line-in-grassmannian-explicit} coincides with the pullback of $-Q(\theta)\otimes \Psi$ along the composition
\begin{align}
\notag
	K_{\omega} / K_{\omega} \cap K_0 &\rightarrow \widetilde K_{\omega} / \widetilde K_{\omega} \cap \widetilde K_0 \\
	& \simeq \widetilde K_0 / \widetilde K_0\cap \widetilde K_{-\omega} \rightarrow \deloop(\widetilde K_0 \cap \widetilde K_{-\omega}) \xrightarrow{\pi} \deloop(T \times T_{\adjoint}) \rightarrow \deloop T \simeq \deloop\mathbb G_m
\label{eq-affine-simple-reflection-line-bundle}
\end{align}
where $\pi := \pi^{(0, -\omega)}$ is the morphism defined in \S\ref{void-schubert-cell-notation}. The morphism \eqref{eq-affine-simple-reflection-line-bundle} classifies the tautological line bundle $\mathscr O_{\mathbb P^1}(-1)$, upon identifying its source with $\mathbb P^1$. It follows that the restriction of $\mathscr G_{\Gr_G}$ along \eqref{eq-affine-simple-reflection-line-in-grassmannian-explicit} is classified by $Q(\theta)$.
\end{proof}

\begin{void}
For each $\chi \in \Hom(\Lambda, A(-1))$, we define the subset
$$
S^{\circ}_{\chi} \subset S
$$
of \emph{integral elements} (with respect to $\chi$), which consists of
\begin{enumerate}
	\item finite simple reflections $s_{\check{\alpha}}$ satisfying $\chi(\alpha) = 0$;
	\item affine simple reflections $\translate(\theta)s_{\check{\theta}}$ satisfying $Q(\theta) + \chi(\theta) = 0$.
\end{enumerate}
\end{void}

\begin{rem}
\label{rem-integral-reflection-fixes-character}
With respect to the $W^{\aff}$-action on $\Hom(\Lambda, A(-1))$ (\emph{cf.}~\S\ref{void-affine-weyl-group-action-on-characters}), each element $s \in S^{\circ}_{\chi}$ fixes $\chi$. This is clear when $s$ is a finite simple reflection. When $s = \translate(\theta) s_{\check{\theta}}$ is an affine simple reflection, this follows from
\begin{align*}
	s(\chi) = s_{\check{\theta}}(\chi) - b(\theta, \cdot) &= \chi - \langle\check{\theta}, \cdot\rangle\chi(\theta) - b(\theta, \cdot) \\
	& = \chi - \langle\check{\theta}, \cdot\rangle(\chi(\theta) + Q(\theta)) = \chi,
\end{align*}
where we used \eqref{eq-strict-weyl-invariance} in the third equality.
\end{rem}

\begin{lem}
\label{lem-integral-reflection-gerbe-triviality}
Given $\chi \in \Hom(\Lambda, A(-1))$, an element $s \in S$ belongs to $S^{\circ}_{\chi}$ if and only if $(i_s)^*(\mathscr G^{\aff} - \overset{\rightarrow}{\pi}{}^*(\chi\otimes\Psi))$ is trivial.
\end{lem}

\begin{proof}
The $A$-gerbe $(i_s)^*\overset{\rightarrow}{\pi}{}^*(-\chi\otimes\Psi)$ is classified by an element $a_{\chi}(s) \in A(-1)$ (\emph{cf.}~\S\ref{void-parity-vanishing-context}). We claim that $a_{\chi}(s)$ is given by
\begin{equation}
\label{eq-rightward-pullback-gerbe-dimension-one-bruhat-variety}
a_{\chi}(s) =
\begin{cases}
	\chi(\alpha) & \text{$s = s_{\check{\alpha}}$ for $\check{\alpha} \in \check{\Delta}$} \\
	\chi(\theta) & \text{$s = \translate(\theta)s_{\check{\theta}}$ for $\check{\theta}$ a highest root}
\end{cases}
\end{equation}

Indeed, for $s = s_{\check{\alpha}}$ a finite simple reflection, this follows from the factorization \eqref{eq-dimension-one-bruhat-variety-finite-simple-reflection}. For $s = \translate(\theta)s_{\check{\theta}}$ an affine simple reflection, we argue as in the proof of Lemma \ref{lem-gerbe-dimension-one-bruhat-variety}, reducing to the fact that pulling back $\Psi^{-\chi}$ along \eqref{eq-flag-variety-line-affine-grassmannian} yields the $A$-gerbe classified by $\chi(\theta)$.

The assertion now follows by combining \eqref{eq-rightward-pullback-gerbe-dimension-one-bruhat-variety} with Lemma \ref{lem-gerbe-dimension-one-bruhat-variety}.
\end{proof}

\begin{void}
We now study the functor $\IC_{\chi}^s$ (\emph{cf.}~\S\ref{void-flag-variety-standard-costandard}) for $s \in S$.

Note that for $s \in S_{\chi}^{\circ}$, the $A$-gerbe $(i_s)^*(\mathscr G^{\aff} - \overset{\rightarrow}{\pi}{}^*(\chi\otimes\Psi))$ is canonically trivial: It is trivial by Lemma \ref{lem-integral-reflection-gerbe-triviality} and its fiber along $e : I/I \rightarrow P_s/I$ admits a canonical trivialization. In this case, $\IC_{\chi}^s$ may be understood as a functor
$$
\IC_{\chi}^s : \Lis(x) \rightarrow {}_{\chi} (\derived_{\zeta}^{\aff, s})_{\chi},
$$
as $s(\chi) = \chi$ by Remark \ref{rem-integral-reflection-fixes-character}. Write $\pi_s : P_s/I \rightarrow \Spec\base \simeq x$ for the structural map.
\end{void}

\begin{prop}
\label{prop-intersection-cohomology-length-one}
Given $\chi \in \Hom(\Lambda, A(-1))$ and $s\in S$, there holds
\begin{enumerate}
	\item if $s \in S^{\circ}_{\chi}$, then
	$$
	\IC^s_{\chi} \simeq (i_s)_*(\pi_s)^*(\frac{1}{2})[1].
	$$
	\item if $s \notin S^{\circ}_{\chi}$, then the canonical maps are isomorphisms:
	$$
	\Delta_{\chi}^s \simeq \IC_{\chi}^s \simeq \nabla_{\chi}^s.
	$$
\end{enumerate}
\end{prop}

\begin{proof}
By construction, $\IC_{\chi}^s$ is the intermediate extension of $((i_s)^*(\mathscr G^{\aff} - \overset{\rightarrow}{\pi}{}^*(\chi\otimes\Psi)), \zeta)$-twisted local systems along the inclusion $I s I/I \subset P_s/I$. Identifying this inclusion with $\mathbb A^1 \subset \mathbb P^1$, we deduce the result from Lemma \ref{lem-integral-reflection-gerbe-triviality}.
\end{proof}

\begin{void}[Proof of Proposition \ref{prop-parity-vanishing}]
Finally, we are ready to prove the parity vanishing of $\IC_{\chi}^v$. The proof will span \S\ref{void-parity-vanishing-equivalent-formulation}--\S\ref{void-parity-vanishing-case-analysis} below.

To simplify the notation, we write $\Delta_{\chi}^v(\coeff)$, $\IC_{\chi}^v(\coeff)$, $\nabla_{\chi}^v(\coeff)$ for the images of a rank-$1$ $(\mathscr G_{\varpi^{\lambda}}, \zeta)$-twisted local system over $x\simeq \Spec\base$ under $\Delta_{\chi}^v$, $\IC_{\chi}^v$, $\nabla_{\chi}^v$ (where $v = \translate(\lambda) w$). Since $\mathscr G_{\varpi^{\lambda}}$ is \emph{non-canonically} trivial, these objects are only well-defined up to \emph{non-unique} isomorphisms. The statement that we want to prove, \emph{i.e.}~the parity vanishing of fibers of $\IC_{\chi}^v(\coeff)$, only concerns its isomorphism class.
\end{void}

\begin{void}
\label{void-parity-vanishing-equivalent-formulation}
We begin by observing that an object $\mathscr A \in {}_{v(\chi)}(\derived^{\aff}_{\zeta})_{\chi}$ satisfies the property
\begin{equation}
\label{eq-parity-vanishing-property}
	H^n i_y^*\mathscr A \neq 0 \text{ for some } y \in \Fl_{G, x}(\base) \Rightarrow \ell(v) + n \text{ is even}
\end{equation}
if and only if $\mathscr A$ belongs to the full subcategory of ${}_{v(\chi)}(\derived^{\aff}_{\zeta})_{\chi}$ generated under extension by objects of the form
\begin{equation}
\label{eq-parity-vanishing-generator}
	\Delta_{\chi}^{v_1}(\coeff)[d]\text{ where } v_1(\chi) = v(\chi) \text{ and }\ell(v) - \ell(v_1) - d \text{ is even}.
\end{equation}

Indeed, any object of the form \eqref{eq-parity-vanishing-generator} satisfies \eqref{eq-parity-vanishing-property} and this property is preserved under extensions. Conversely, any $\mathscr A \in {}_{v(\chi)}(\derived^{\aff}_{\zeta})_{\chi}$ is supported on Bruhat cells $Iv_1 I/I$ satisfying $v_1(\chi) = v(\chi)$. If $\mathscr A$ satisfies \eqref{eq-parity-vanishing-property}, it belongs to the full subcategory generated by objects of the form \eqref{eq-parity-vanishing-generator} under extensions, by induction on support.
\end{void}

\begin{void}
We now prove that $\IC_{\chi}^v(\coeff)$ satisfies \eqref{eq-parity-vanishing-property} by induction on $\ell(v)$.

For the base step $\ell(v) = 0$, the object $\IC_{\chi}^v(\coeff)$ is a skyscraper, so $H^n i_y^*\IC_{\chi}^v(\coeff) \neq 0$ for some $y\in \Fl_{G, x}(\base)$ implies $n = 0$.

For the inductive step, we choose $s \in S$ such that $\ell(v) = \ell(vs) + 1$. Applying the decomposition theorem to the Demazure resolution (\emph{cf.}~\cite[\S6.2.4]{MR751966}), we see that $\IC_{\chi}^v(\coeff)$ is a direct summand of the convolution product (\emph{cf.}~\S\ref{void-monodromic-affine-hecke-convolution})
\begin{equation}
\label{eq-intersection-cohomology-summand}
\IC_{s(\chi)}^{vs}(\coeff) \circ \IC_{\chi}^s(\coeff),
\end{equation}
so it suffices to prove that the object \eqref{eq-intersection-cohomology-summand} satisfies \eqref{eq-parity-vanishing-property}.

Applying the induction hypothesis to $\IC_{s(\chi)}^{vs}(\coeff)$ and using the equivalence of \S\ref{void-parity-vanishing-equivalent-formulation}, we may replace it by an object of the form $\Delta_{s(\chi)}^{v_1}(\coeff)[d]$ where $v_1s(\chi) = v(\chi)$ and $\ell(vs) - \ell(v_1) - d$ is even. It then suffices to prove the following \emph{purity} statement:
\begin{equation}
\label{eq-convolution-parity-vanishing}
H^n i_y^*(\Delta_{s(\chi)}^{v_1}(\coeff) \circ \IC_{\chi}^s(\coeff)) \neq 0\text{ for some } y \in \Fl_{G, x}(\base) \Rightarrow \ell(v_1s) + n = 0
\end{equation}
for any $s \in S$ and $v_1 \in W^{\aff}$.
\end{void}

\begin{void}
\label{void-parity-vanishing-case-analysis}
Finally, we prove \eqref{eq-convolution-parity-vanishing} via a case-by-case analysis.
\begin{enumerate}
	\item \emph{$s\notin S_{\chi}^{\circ}$ and $\ell(v_1 s) = \ell(v_1) + 1$.} In this case, we have
	$$
	\Delta_{s(\chi)}^{v_1}(\coeff) \circ \IC_{\chi}^s(\coeff) \simeq \Delta_{s(\chi)}^{v_1}(\coeff) \circ \Delta_{\chi}^s(\coeff) \simeq \Delta_{\chi}^{v_1 s}(\coeff),
	$$
	where we used the isomorphism $\Delta_{\chi}^s(\coeff) \simeq \IC_{\chi}^s(\coeff)$ (\emph{cf.}~Proposition \ref{prop-intersection-cohomology-length-one}).
	
	\item \emph{$s \notin S_{\chi}^{\circ}$ and $\ell(v_1 s) = \ell(v_1) - 1$.} In this case, we have
	\begin{equation}
	\label{eq-standard-length-additive}
	\Delta_{s(\chi)}^{v_1}(\coeff) \simeq \Delta_{\chi}^{v_1s}(\coeff) \circ \Delta_{s(\chi)}^s(\coeff).
	\end{equation}
	
	On the other hand, Proposition \ref{prop-intersection-cohomology-length-one} yields $\IC_{\chi}^s(\coeff) \simeq \nabla_{\chi}^s(\coeff)$, and we have an isomorphism $\Delta^s_{s(\chi)}(\coeff) \circ \nabla_{\chi}^s(\coeff) \simeq \Delta_{\chi}^e(\coeff)$ which reduces to a calculation for $\SL_2$ (\emph{cf.}~\cite[Lemma 3.5]{MR4108915}). Putting these together gives
	\begin{align*}
		\Delta_{s(\chi)}^{v_1}(\coeff) \circ \IC_{\chi}^s(\coeff) & \simeq \Delta_{\chi}^{v_1 s}(\coeff) \circ \Delta_{s(\chi)}^s(\coeff) \circ \IC_{\chi}^s(\coeff) \\
		& \simeq \Delta_{\chi}^{v_1 s}(\coeff) \circ \Delta_{s(\chi)}^s(\coeff) \circ \nabla_{\chi}^s(\coeff) \simeq \Delta_{\chi}^{v_1 s}(\coeff).
	\end{align*}
	
	\item \emph{$s \in S^{\circ}_{\chi}$ and $\ell(v_1 s) = \ell(v_1) + 1$.} Consider the Cartesian square
	$$
	\begin{tikzcd}[column sep = 1em]
		G\loo{\varpi} \times^I P_s/I \ar[r, "m"]\ar[d, "p"] & G\loo{\varpi}/I \ar[d, "f"] \\
		G\loo{\varpi}/I \ar[r, "f"] & G\loo{\varpi}/P_s
	\end{tikzcd}
	$$
	where $f$ is induced from the inclusion $I \subset P_s$ and $p$, $m$ are the projection, respectively multiplication maps. The condition $s \in S^{\circ}_{\chi}$ allows us to identify $(\cdot)\circ \IC_{\chi}^s(\coeff)$ with the functor $f^*f_!(\frac{1}{2})[1]$ (\emph{cf.}~Proposition \ref{prop-intersection-cohomology-length-one}).
	
	On the other hand, the condition on $v_1$ implies that $I v_1 I \times^I P_s \simeq I v_1 P_s$ along the multiplication map. It induces an isomorphism $I v_1 I/I \simeq I v_1 P_s/P_s$, so
	$$
	\Delta_{s(\chi)}^{v_1}(\coeff) \circ \IC_{\chi}^s(\coeff) \simeq f^*f_! \Delta_{\chi}^{v_1}(\coeff)(\frac{1}{2})[1]
	$$
	is isomorphic to the $!$-extension along $I v_1 P_s/I \rightarrow G\loo{\varpi}/I$ of a constant sheaf placed in cohomological degree $-\ell(v_1 s)$.
	
	\item \emph{$s \in S^{\circ}_{\chi}$ and $\ell(v_1 s) = \ell(v_1) - 1$.} In this case, we again have \eqref{eq-standard-length-additive}. There is also an isomorphism $\Delta_{\chi}^s(\coeff) \circ \IC_{\chi}^s(\coeff) \simeq \IC_{\chi}^s(\coeff)(-\frac{1}{2})[-1]$ as $\IC_{\chi}^s(\coeff)$ is constant (\emph{cf.}~Proposition \ref{prop-intersection-cohomology-length-one}). Putting these together gives
	\begin{align*}
		\Delta_{s(\chi)}^{v_1}(\coeff) \circ \IC_{\chi}^s(\coeff) & \simeq \Delta_{\chi}^{v_1 s}(\coeff) \circ \Delta_{s(\chi)}^s(\coeff) \circ \IC_{\chi}^s(\coeff) \\
		& \simeq \Delta_{\chi}^{v_1 s}(\coeff) \circ \IC_{\chi}^s(\coeff)(-\frac{1}{2})[-1],
	\end{align*}
	whose fibers are concentrated in cohomological degree $-\ell(v_1s)$ by step (3).
\end{enumerate}

The proof of \eqref{eq-convolution-parity-vanishing}, hence of Proposition \ref{prop-parity-vanishing}, is now complete.\qed
\end{void}

\subsection{Applications to the Satake category}
\label{sec-parity-vanishing-application}

\begin{void}
The goal of this subsection is to deduce Theorem \ref{thm-semisimplicity} from Proposition \ref{prop-parity-vanishing}.

First, we note the analogue of Proposition \ref{prop-parity-vanishing} for the affine Grassmannian. We let $x$ be a $\base$-point of $X$ and consider the functor $\IC^{\lambda}$ associated to $\lambda \in \Lambda^{\sharp, +}$ (\emph{cf.}~\S\ref{void-affine-grassmannian-ic-sheaf}).
\end{void}

\begin{lem}
\label{lem-grassmannian-parity-vanishing}
Fix $\lambda \in \Lambda^{\sharp, +}$. The fibers of $\IC^{\lambda}$ over $\Gr_{G, x}$ are concentrated in cohomological degrees $n$ such that $n + \langle 2\check{\rho}, \lambda\rangle$ is even.
\end{lem}

\begin{proof}
The projection map $\pi : \Fl_{G, x} \rightarrow \Gr_{G, x}$ is smooth of relative dimension $\ell(w_0)$ and the pre-image of the Schubert stratification is refined by the Bruhat stratification (\emph{cf.}~\S\ref{void-bruhat-cell}). Thus, we have an isomorphism
$$
\pi^*\IC^{\lambda}(\frac{\ell(w_0)}{2})[\ell(w_0)] \simeq \IC^v_0,
$$
where $\IC^v_0$ is the functor \eqref{eq-affine-flag-variety-ic-functor} associated to $\chi = 0$ and $v \in W^{\aff}$ the unique element of $W \lambda W$ with length $\langle 2\check{\rho}, \lambda\rangle + \ell(w_0)$. The assertion now follows from Proposition \ref{prop-parity-vanishing}.
\end{proof}

\begin{void}
We shall now prove Theorem \ref{thm-semisimplicity}.

\begin{proof}[Proof of Theorem \ref{thm-semisimplicity}]
Recall that $\Sat_{\mathscr G, \zeta}(\Hec_{G, x})$ may be realized as a full subcategory of $\derived_{\mathscr G, \zeta}(\Gr_{G, x})$ (\emph{cf.}~Remark \ref{rem-satake-category-forgetful-to-grassmannian}).

Fix $\lambda \in \Lambda^{\sharp, +}$ and let $\mathscr E_1$, $\mathscr E_2$ be $(\mathscr G_{\varpi^{\lambda}}, \zeta)$-twisted local systems over $S$. Let us compute the space of maps from $\IC^{\lambda}(\mathscr E_1)$ to $\IC^{\lambda}(\mathscr E_2)[1]$ \emph{as objects of $\derived_{\mathscr G, \zeta}(\Gr_{G, x})$}.

Denote by $j : \Gr_{G, x}^{\lambda} \hookrightarrow \overline{\Gr}{}_{G, x}^{\lambda}$ the open immersion of the Schubert cell into its closure and $i$ its complement. Then $i^*\IC^{\lambda}(\mathscr E_1)$ is concentrated in perverse cohomological degrees $\le -1$, while $i^! \IC^{\lambda}(\mathscr E_2)[1]$ in degrees $\ge 0$. Applying excision to $\IC^{\lambda}(\mathscr E_1)$, we obtain
\begin{align}
\notag
\Maps_{\derived_{\mathscr G, \zeta}(\Gr_{G, x})}(\IC^{\lambda}(\mathscr E_1), \IC^{\lambda}(\mathscr E_2)[1]) &\simeq \Maps_{\derived_{\mathscr G, \zeta}(\Gr_{G, x}^{\lambda})}(\widetilde{\mathscr E}_1, \widetilde{\mathscr E}_2[1]) \\
\label{eq-ic-derived-orthogonality-same-cocharacter}
& \simeq \Maps_{\Lis_{\mathscr G_{\varpi^{\lambda}}, \zeta}(S)}(\mathscr E_1, \mathscr E_2[1]),
\end{align}
where the second isomorphism is due to the fact that the fibers of $\Gr_{G, x}^{\lambda} \rightarrow S$ have vanishing cohomology in degree $1$.

The isomorphism \eqref{eq-ic-derived-orthogonality-same-cocharacter} shows that the functor $\IC^{\lambda}$ is fully faithful and its essential image is closed under extensions.

Next, we fix $\lambda_1\neq \lambda_2 \in \Lambda^{\sharp, +}$ and $(\mathscr G_{\varpi^{\lambda_1}}, \zeta)$-, respectively $(\mathscr G_{\varpi^{\lambda_2}}, \zeta)$-twisted local systems $\mathscr E_1$, respectively $\mathscr E_2$ over $S$. We claim
\begin{equation}
\label{eq-ic-sheaf-orthogonality-distinct-cocharacters}
\Maps_{\derived_{\mathscr G, \zeta}(\Gr_{G, x})}(\IC^{\lambda_1}(\mathscr E_1), \IC^{\lambda_2}(\mathscr E_2)[1]) \simeq 0.
\end{equation}
If $S$ is a geometric point, this follows from Lemma \ref{lem-grassmannian-parity-vanishing} by the classical argument (\emph{cf.}~the proof of \cite[Proposition 1]{MR1826370}). The case for general $S$ follows from the base change property of Homs between ULA objects in $\derived_{\mathscr G, \zeta}(\Gr_{G, x})$ (\emph{cf.}~the proof of Proposition \ref{prop-local-hecke-stack-derived-category-stalk}).

The vanishing \eqref{eq-ic-sheaf-orthogonality-distinct-cocharacters} shows that the essential images of the functors $\IC^{\lambda_1}$ and $\IC^{\lambda_2}$ are mutually orthogonal and admit no nontrivial extensions.

The fully faithfulness of \eqref{eq-satake-category-generator-sum} follows from the fully faithfulness of $\IC^{\lambda}$, for each $\lambda \in \Lambda^{\sharp, +}$, and the orthogonality of the essential images of $\IC^{\lambda_1}$, $\IC^{\lambda_2}$, for $\lambda_1 \neq \lambda_2$. The essential surjectivity of \eqref{eq-satake-category-generator-sum} follows from the fact that $\Sat_{\mathscr G, \zeta}(\Hec_{G, x})$ is generated under extensions by the essential images of $\IC^{\lambda}$ for $\lambda \in \Lambda^{\sharp, +}$ (\emph{cf.}~Corollary \ref{cor-schubert-cell-reduction-of-support}), the closure of the essential image of $\IC^{\lambda}$ under extensions, and the absence of nontrivial extensions between the essential images of $\IC^{\lambda_1}$ and $\IC^{\lambda_2}$, for $\lambda_1 \neq\lambda_2$.
\end{proof}

\end{void}

\medskip

\section{Tori}
\label{sec-tori}

In this section, we fix a smooth curve $X$ over a field $\base$ and an $X$-torus $T$. Let $\coeff$ be a finite extension of $\rationals_{\ell}$, for a prime $\ell$ invertible in $\base$. Let $\zeta : A \subset \coeff^{\times}$ be a finite subgroup and $\mu$ be an $A$-valued \'etale level of $T$.

The goal of this section is to construct the geometric Satake equivalence (\emph{cf.}~Theorem \ref{thm-satake-equivalence}) for $(T, \mu)$, which will be supplied in \S\ref{void-satake-equivalence-tori}. The key ingredient is Proposition \ref{prop-sharp-torus-commutative-gerbe-identification}, which relates $\mathscr G_{\Gr_T}$ to the metaplectic dual datum $\nu$.

\subsection{Reduction to sharp tori}

\begin{void}
We begin by specializing the metaplectic dual data $(H, \nu)$ (\emph{cf.}~\S\ref{sec-metaplectic-dual-data}) to $(T, \mu)$, whose construction is substantially simpler than the general case.

Indeed, writing $\Lambda$ for the sheaf of cocharacters of $T$ and $\Lambda^{\sharp} \subset \Lambda$ for the subsheaf defined in \S\ref{void-metaplectic-dual-group}, we see that $\Lambda^{\sharp} \subset \Lambda$ corresponds to an isogeny of $X$-tori $T^{\sharp} \rightarrow T$ and $H$ is the Langlands dual of $T^{\sharp}$ (as a locally constant \'etale sheaf of $\integers$-tori).

On the other hand, the pullback $\mu^{\sharp}$ of the \'etale level $\mu$ to $T^{\sharp}$ acquires an $\mathbb E_{\infty}$-monoidal structure, so it induces an $\mathbb E_{\infty}$-monoidal morphism $\nu : \Lambda^{\sharp} \rightarrow \deloop^2_X A$ by applying the functor $\SMaps_*(\deloop_X\mathbb G_m, \cdot)$ (\emph{cf.}~\S\ref{void-sharp-torus-symmetric-monoidal-level}).
\end{void}

\begin{void}
For each finite set $I$, we write
\begin{equation}
\label{eq-torus-sharp-isogeny-grassmannian}
f_I : \Gr_{T^{\sharp}, I} \rightarrow \Gr_{T, I}
\end{equation}
for the induced morphism on affine Grassmannians (\emph{cf.}~Remark \ref{rem-satake-category-finite-set-functoriality}).

Since \eqref{eq-torus-sharp-isogeny-grassmannian} is ind-proper, the functor
\begin{equation}
\label{eq-torus-sharp-isogeny-grassmannian-pushforward}
(f_I)_! : \derived_{\mathscr G, \zeta}(\Gr_{T^{\sharp}, I}) \rightarrow \derived_{\mathscr G, \zeta}(\Gr_{T, I})
\end{equation}
preserves universal local acyclicity relative to $X^I$. Moreoever, the base change of \eqref{eq-torus-sharp-isogeny-grassmannian} to any geometric point of $X^I$ is a closed immersion, so $(f_I)_!$ is also perverse $t$-exact. Finally, $L^+_I T^{\sharp} \rightarrow L^+_I T$ is surjective in the \'etale topology. By Remark \ref{rem-satake-category-forgetful-to-grassmannian}, we see that \eqref{eq-torus-sharp-isogeny-grassmannian-pushforward} restricts to a functor on the Satake categories
\begin{equation}
\label{eq-torus-sharp-isogeny-satake-categories}
	(f_I)_! : \Sat_{\mathscr G, \zeta}(\Hec_{T^{\sharp}, I}) \rightarrow \Sat_{\mathscr G, \zeta}(\Hec_{T, I}).
\end{equation}

The functor \eqref{eq-torus-sharp-isogeny-satake-categories} is naturally symmetric monoidal, as its formation is compatible with the fusion product (\emph{cf.}~\S\ref{void-fusion-product}).
\end{void}

\begin{prop}
\label{prop-torus-sharp-isogeny-satake-categories}
The functor \eqref{eq-torus-sharp-isogeny-satake-categories} is an equivalence of categories.
\end{prop}

\begin{void}
\label{void-torus-grassmannian-canonical-sections}
Proposition \ref{cor-schubert-cell-reduction-of-support} is an immediate consequence of Corollary \ref{cor-schubert-cell-reduction-of-support} if $I = \{1\}$. We shall reduce the general case to this one using factorization.

Suppose that $T$ is split. For an $I$-tuple $\lambda^I = (\lambda^i)_{i\in I}$ of elements of $\Lambda$, we write $\varpi^{\lambda^I} : X^I \rightarrow \Gr_{T, I}$ for the closed immersion sending $x^I = (x^i)_{i\in I}$ to the modification of $T$-bundles
$$
\mathscr O \overset{x^I}{\sim} \mathscr O(\sum_{i\in I} \lambda^i \Gamma_{x^i}).
$$

Denote by $X^{I, \disj} \subset X^I$ the open subscheme consisting of pairwise disjoint points $x^i$ ($i\in I$) of $X$. Write $j : \Gr_{T, I}^{\disj} \rightarrow \Gr_{T, I}$ for its base change. Recall (\emph{cf.}~Proposition \ref{prop-satake-category-disjoint-pullback}) that the pullback functor is fully faithful:
\begin{equation}
\label{eq-satake-category-torus-disjoint-restriction}
j^* : \Sat_{\mathscr G, \zeta}(\Gr_{T, I}) \subset \Sat_{\mathscr G, \zeta}(\Gr_{T, I}^{\disj}),
\end{equation}
where $\Sat_{\mathscr G, \zeta}(\cdot)$, as usual, is the full subcategory of $\derived_{\mathscr G, \zeta}(\cdot)$ characterized by universal local acyclicity and perversity relative to $X^I$, respectively $X^{I, \disj}$.
\end{void}

\begin{lem}
\label{lem-satake-category-torus-extension-characterization}
An object $\mathscr A \in \Sat_{\mathscr G, \zeta}(\Gr_{T, I}^{\disj})$ belongs to the essential image of \eqref{eq-satake-category-torus-disjoint-restriction} if and only if $(\varpi^{\lambda^I})^*\mathscr A$ extends as a twisted local system along $X^{I, \disj} \subset X^I$ for every $\lambda^I \in \Lambda^I$.
\end{lem}

\begin{proof}
We note that the essential image of \eqref{eq-satake-category-torus-disjoint-restriction} is closed under direct summands. This is a twisted version of \cite[Theorem 6.8(ii)]{MR4630128} and follows from the same proof.

Suppose that $\mathscr A \simeq j^*\mathscr B$ for some $\mathscr B \in \Sat_{\mathscr G, \zeta}(\Gr_{T, I})$. Since each $(\varpi^{\lambda^I})_!(\varpi^{\lambda^I})^*\mathscr A$ is a direct summand of $\mathscr A$, it extends to $\Gr_{T, I}$ along \eqref{eq-satake-category-torus-disjoint-restriction}, which must then be an object of the form $(\varpi^{\lambda^I})_!\mathscr B^{\lambda^I}$, where $\mathscr B^{\lambda^I}$ is a twisted local system over $X^I$ extending $(\varpi^{\lambda^I})^*\mathscr A$.

Conversely, if each $(\varpi^{\lambda^I})^*\mathscr A$ extends to some twisted local system $\mathscr B^{\lambda^I}$ over $X^I$, then the sum $\bigoplus_{\lambda^I} (\varpi^{\lambda^I})_!\mathscr B^{\lambda^I}$ belongs to $\Sat_{\mathscr G, \zeta}(\Gr_{T, I})$ and extends $\mathscr A$.
\end{proof}

\begin{void}
We shall now prove Proposition \ref{prop-torus-sharp-isogeny-satake-categories}.

\begin{proof}[Proof of Proposition \ref{prop-torus-sharp-isogeny-satake-categories}]
The statement is of \'etale local nature over $X^I$, so we may assume that $T$ is split. Consider the commutative square
\begin{equation}
\label{eq-torus-sharp-isogeny-satake-disjoint}
\begin{tikzcd}[column sep = 1.5em]
	\Sat_{\mathscr G, \zeta}(\Hec_{T^{\sharp}, I}) \ar[r, phantom, "\subset"]\ar[d, "(f_I)_!"] & \Sat_{\mathscr G, \zeta}(\Hec_{T^{\sharp}, I}^{\disj}) \ar[d, "(f_I)_!"] \\
	\Sat_{\mathscr G, \zeta}(\Hec_{T, I}) \ar[r, phantom, "\subset"] & \Sat_{\mathscr G, \zeta}(\Hec_{T, I}^{\disj})
\end{tikzcd}
\end{equation}
where the horizontal embeddings are given by Proposition \ref{prop-satake-category-disjoint-pullback}.

By the factorization structure on $(\Hec_T, \mathscr G_{\Hec_T})$ (\emph{cf.}~Proposition \ref{prop-local-hecke-stack-gerbe-factorization}) and Corollary \ref{cor-schubert-cell-reduction-of-support}, the right vertical functor in \eqref{eq-torus-sharp-isogeny-satake-disjoint} is an equivalence. This implies that the left vertical functor is fully faithful. Applying the characterization of Lemma \ref{lem-satake-category-torus-extension-characterization} to both $T$ and $T^{\sharp}$, we see that it is also essentially surjective.
\end{proof}
\end{void}

\begin{void}
In \S\ref{sec-satake-equivalence-sharp-torus}, we shall prove the following statement.
\end{void}

\begin{prop}
\label{prop-satake-equivalence-sharp-torus}
For a finite set $I$, there is a canonical equivalence of symmetric monoidal $\coeff$-linear categories
\begin{equation}
\label{eq-satake-equivalence-sharp-torus}
\Sat_{\mathscr G, \zeta}(\Hec_{T^{\sharp}, I}) \simeq \Rep_{H^{\boxtimes I}, (\nu + \vartheta)^{\boxplus I}},
\end{equation}
which depends functorially on $I$ (\emph{cf.}~\S\ref{void-satake-equivalence-finite-set-functoriality}).
\end{prop}

\begin{void}[Theorem \ref{thm-satake-equivalence} for tori]
\label{void-satake-equivalence-tori}
Using Proposition \ref{prop-torus-sharp-isogeny-satake-categories} and Proposition \ref{prop-satake-equivalence-sharp-torus}, we shall construct the geometric Satake equivalence \eqref{eq-satake-equivalence} for $(T, \mu)$.

Namely, we set \eqref{eq-satake-equivalence} to be the composition of \eqref{eq-satake-equivalence-sharp-torus} with the inverse of \eqref{eq-torus-sharp-isogeny-satake-categories}:
$$
\begin{tikzcd}
	\Sat_{\mathscr G, \zeta}(\Hec_{T^{\sharp}, I}) \ar[r, phantom, "\simeq"] & \Rep_{H^{\boxtimes I}, (\nu + \vartheta)^{\boxplus I}} \\
	\Sat_{\mathscr G, \zeta}(\Hec_{T, I}) \ar[u, "\simeq"]
\end{tikzcd}
$$
\end{void}

\subsection{Equivalence for sharp tori}
\label{sec-satake-equivalence-sharp-torus}

\begin{void}
\label{void-satake-equivalence-sharp-torus-replacement}
The goal of this subsection is to prove Proposition \ref{prop-satake-equivalence-sharp-torus}. Note that $(T, \mu)$ and $(T^{\sharp}, \mu^{\sharp})$ share the same metaplectic dual data $(H, \nu)$, so we may replace $(T, \mu)$ by $(T^{\sharp}, \mu^{\sharp})$ and start with an \emph{$\mathbb E_{\infty}$-monoidal} morphism
\begin{equation}
\label{eq-symmetric-monoidal-etale-level}
\mu : \deloop_X T \rightarrow \deloop_X^4A(1).
\end{equation}

Recall that $\Gr_T$ is a factorization algebra in $\Stk$ (\emph{cf.}~Remark \ref{rem-affine-grassmannian-factorization}). Our first step is to lift $(\Gr_T, \mathscr G_{\Gr_T})$ to a factorization algebra in commutative multiplicative stacks endowed with a commutative multiplicative $A$-gerbe.
\end{void}

\begin{void}
Denote by $\Gp^{\com}(\Stk)$ the $\base$-presheaf assigning to $R$ the symmetric monoidal category of grouplike $\mathbb E_{\infty}$-monoids in \'etale $R$-stacks.

Recall that an $R$-point of $\Gr_T$ may be thought of as a morphism $D_{\underline x} \rightarrow \deloop_X T$ over $X$ (for $\underline x \in \Ran(R)$) with a trivialization over $\mathring D_{\underline x}$. The $\mathbb E_{\infty}$-monoid structure on $\deloop_X T$ thus lifts $\Gr_{T, \underline x}$ to an object of $\Gp^{\com}(\Stk)(R)$. Compatibility with the factorization structure then ensures that $\Gr_T$ lifts to a factorization algebra in $\Gp^{\com}(\Stk)$.

Denote by $\Gp^{\com}(\Stk)_{/\deloop^2 A}$ the $\base$-presheaf assigning to $R$ the symmetric monoidal category of objects $\mathscr Z \in \Gp^{\com}(\Stk)(R)$ endowed with an $\mathbb E_{\infty}$-monoidal morphism $\mathscr Z \rightarrow \deloop^2 A$.
\end{void}

\begin{prop}
\label{prop-sharp-torus-grassmannian-commutative-factorization}
The pair $(\Gr_T, \mathscr G_{\Gr_T})$ canonically lifts to a factorization algebra in $\Gp^{\com}(\Stk)_{/\deloop^2 A}$;
\end{prop}

\begin{proof}
This is analogous to the proof of Proposition \ref{prop-local-hecke-stack-gerbe-factorization}, the only difference being that the morphism \eqref{eq-hecke-groupoid-maps-to-linear-groupoid-ran-algebra} for $(T, \mu)$:
$$
\Hec_T(R) \rightarrow \Gamma(D_{(\cdot)}\text{ mod }\mathring D_{(\cdot)}, \deloop^4A(1))
$$
lifts to a morphism of $\Ran(R)$-algebras in $\Span(\Gp^{\com}(\Spc))$.
\end{proof}

\begin{rem}
The factorization structure supplied by Proposition \ref{prop-sharp-torus-grassmannian-commutative-factorization} may be compared with that of $(\Hec_T, \mathscr G_{\Hec_T})$ (\emph{cf.}~Proposition \ref{prop-local-hecke-stack-gerbe-factorization}) as follows.

Consider the morphism of $\Ran$-stacks
\begin{equation}
\label{eq-local-hecke-stack-torus-quotient-onto-grassmannian}
\Hec_T \rightarrow \Gr_T,\quad (P^0 \overset{\underline x}{\sim} P^1) \mapsto (\mathscr O \overset{\underline x}{\sim} P^1\otimes(P^0)^{-1}).
\end{equation}
which lifts to a morphism of factorization algebras in $\Span(\Gpd(\Stk))$, under the forgetful functor $\Gp^{\com}(\Stk) \rightarrow \Span(\Gpd(\Stk))$.

One can show that \eqref{eq-local-hecke-stack-torus-quotient-onto-grassmannian} canonically lifts to a morphism $(\Hec_T, \mathscr G_{\Hec_T}) \rightarrow (\Gr_T, \mathscr G_{\Gr_T})$ of factorization algebras in $\Span(\Gpd(\Stk))_{/\deloop^2 A}$. (This structure will not be used in the sequel, so we omit its construction.) To the contrary, the structural morphism $(\Gr_T, \mathscr G_{\Gr_T}) \rightarrow (\Hec_T, \mathscr G_{\Hec_T})$ only lifts to a morphism of factorization algebras in $\Span(\Stk)_{/\deloop^2 A}$, \emph{i.e.}~it is \emph{incompatible} with the multiplicative structures.
\end{rem}

\begin{void}
Given a finite set $I$, the system of morphisms $\varpi^{\lambda^I}$ (for $\lambda^I \in \Lambda^I$) of \S\ref{void-torus-grassmannian-canonical-sections} organizes into a morphism of \'etale sheaves of abelian groups over $X^I$
\begin{equation}
\label{eq-torus-grassmannian-canonical-sections-total}
	\varpi : \Lambda^{\boxplus I} \rightarrow \Gr_{T, I}.
\end{equation}

Indeed, \'etale locally over $X$, the $X$-torus $T$ is split and the \'etale sheaf $\Lambda^{\boxplus I}$ is represented by $\Lambda^I \times X^I$, for $\Lambda$ the abelian group of cocharacters of $T$. In this case, \eqref{eq-torus-grassmannian-canonical-sections-total} is the sum of $\varpi^{\lambda^I}$ over $\lambda^I \in \Lambda^I$. The general case follows from \'etale descent.

On the other hand, $\mathscr G_{\Gr_T}$ defines an $\mathbb E_{\infty}$-monoidal morphism by Proposition \ref{prop-sharp-torus-grassmannian-commutative-factorization}:
\begin{equation}
\label{eq-sharp-torus-grassmannian-commutative-morphism}
\mathscr G_{\Gr_{T, I}} : \Gr_{T, I} \rightarrow \deloop^2_{X^I} A.
\end{equation}
Recall the $\mathbb E_{\infty}$-monoidal morphism $(\nu + \vartheta)^{\boxplus I}$ (\emph{cf.}~\S\ref{void-twisted-representations-multiple-point}) appearing on the dual side of the Satake equivalence \eqref{eq-satake-equivalence-sharp-torus}. The following computation explains its geometric origin.
\end{void}

\begin{prop}
\label{prop-sharp-torus-commutative-gerbe-identification}
The composition of \eqref{eq-torus-grassmannian-canonical-sections-total} and \eqref{eq-sharp-torus-grassmannian-commutative-morphism} is canonically identified with the $\mathbb E_{\infty}$-monoidal morphism $(\nu + \vartheta)^{\boxplus I}$.
\end{prop}

\begin{proof}
Denote by $\nu^{\geom}_I$ the composition of \eqref{eq-torus-grassmannian-canonical-sections-total} and \eqref{eq-sharp-torus-grassmannian-commutative-morphism}. Both $\nu_I^{\geom}$ and $(\nu + \vartheta)^{\boxplus I}$ are defined by \'etale descent from the case of a split $X$-torus $T$, so we may assume that $T$ is split and view $\Lambda$ as an abelian group.

Let us identify $\nu_I^{\geom}$ with $(\nu_{\{1\}}^{\geom})^{\boxplus I}$ and thereby reduce to the case $I = \{1\}$. Indeed, given an $R$-point $x^I = (x^i)_{i\in I}$ of $X^I$ (with induced $R$-point $\underline x$ of $\Ran$), we have a commutative diagram of $\mathbb E_{\infty}$-monoids
\begin{equation}
\label{eq-sharp-torus-commutative-gerbe-reduction-to-singleton}
\begin{tikzcd}[column sep = 1em]
	\Lambda^{\oplus I} \ar[r] & \bigoplus_{i\in I} \Gamma(D_{x^i} \text{ mod }\mathring D_{x^i}, \deloop T) \ar[r, "\sum"]\ar[d, "\mu"] & \Gamma(D_{\underline x} \text{ mod }\mathring D_{\underline x}, \deloop T) \ar[d, "\mu"] \\
	&  \bigoplus_{i\in I} \Gamma(D_{x^i} \text{ mod }\mathring D_{x^i}, \deloop^4A(1)) \ar[r, "\sum"] & \Gamma(D_{\underline x} \text{ mod }\mathring D_{\underline x}, \deloop^4A(1)) \ar[r, "\tr_{\underline x}"] & \Gamma(\Spec R, \deloop^2 A)
\end{tikzcd}
\end{equation}
where the first map is the sum of the maps
$$
\Lambda \rightarrow \Gamma(D_{x^i} \text{ mod }\mathring D_{x^i}, \deloop T),\quad \lambda \mapsto (\mathscr O\overset{x^i}{\sim} \mathscr O(\lambda x^i)),
$$
and $\tr_{\underline x}$ is the (underlying $\mathbb E_{\infty}$-monoidal morphism of the) trace map \eqref{eq-trace-map-with-coefficients}. The upper circuit of \eqref{eq-sharp-torus-commutative-gerbe-reduction-to-singleton} is the pullback of $\nu_I^{\geom}$ along $x^I$, while the lower circuit is the pullback of $(\nu_{\{1\}}^{\geom})^{\boxplus I}$ along $x^I$ by additivity of the trace map.

We shall now assume $I = \{1\}$ and omit it from the notation. It remains to identify the $\mathbb E_{\infty}$-monoidal morphism $\nu^{\geom}$ with $\nu + \vartheta$. Fix an $R$-point $x$ of $X$. The pullback of $\nu^{\geom}$ along $x$ may be expressed as the composition
\begin{align}
\notag
	\Lambda \rightarrow \Maps_*(\deloop_X\mathbb G_m, \deloop_XT) & \xrightarrow{\mu} \Maps_*(\deloop_X\mathbb G_m, \deloop^4_XA(1)) \\
\label{eq-sharp-torus-commutative-gerbe-singleton-expression}
	& \xrightarrow{\mathscr O(x)^*} \Gamma(D_x \text{ mod }\mathring D_x, \deloop^4A(1)) \xrightarrow{\tr_x} \Gamma(\Spec R, \deloop^2 A),
\end{align}
where $\mathscr O(x)$ is viewed as a morphism $D_x \rightarrow \deloop_X\mathbb G_m$ trivialized over $\mathring D_x$. Note that all but the morphism $\mu$ in \eqref{eq-sharp-torus-commutative-gerbe-singleton-expression} come from $H\integers$-linear morphisms.

Let us invoke the decomposition of Remark \ref{rem-etale-level-fiber-sequence-multiplicative-group}:
\begin{equation}
\label{eq-multiplicative-group-etale-level-splitting}
\Maps_*(\deloop_X \mathbb G_m, \deloop_X^4A(1)) \simeq \Gamma(X, \deloop^2 A) \oplus \Gamma(X, A(-1)).
\end{equation}
Under this decomposition, the composition of the first two morphisms of \eqref{eq-sharp-torus-commutative-gerbe-singleton-expression} is the sum $\nu \oplus \epsilon$ (\emph{cf.}~\S\ref{void-sharp-torus-symmetric-monoidal-level}, \S\ref{void-theta-shift}). Let us identify the composition $\tr_x \circ \mathscr O(x)^*$ of the last two morphisms of \eqref{eq-sharp-torus-commutative-gerbe-singleton-expression}, individually for each summand in \eqref{eq-multiplicative-group-etale-level-splitting}.

\begin{enumerate}
	\item \emph{$\Gamma(X, \deloop^2 A)$-summand.} The inclusion $\Gamma(X, \deloop^2 A) \rightarrow \Maps_*(\deloop_X\mathbb G_m, \deloop_X^4A(1))$ is given by tensor product with the Kummer map $\Psi : \deloop_X\mathbb G_m \rightarrow \deloop^2_X\hat{\integers}(1)$.
	
	By the natural isomorphism (\emph{cf.}~Remark \ref{rem-trace-map-projection-formula}, Remark \ref{rem-kummer-class-of-tautological-line-bundle})
	\begin{align*}
	(\tr_x \circ \mathscr O(x)^*)(f\otimes \Psi) &\simeq \tr_x(f\otimes \Psi(\mathscr O(x))) \\
	& \simeq f\otimes \tr_x(\Psi(\mathscr O(x))) \simeq f
	\end{align*}
	linear in $f \in \Gamma(X, \deloop^2 A)$, we see that the restriction of $\tr_x \circ \mathscr O(x)^*$ to the $\Gamma(X, \deloop^2 A)$-summand is identified with the pullback along $x : \Spec R\rightarrow X$.
	
	\item \emph{$\Gamma(X, A(-1))$-summand.} The inclusion $\Gamma(X, A(-1)) \rightarrow \Maps_*(\deloop_X\mathbb G_m, \deloop^4_XA(1))$ is given by tensor product with $\Psi^{\otimes 2} : \deloop_X\mathbb G_m \rightarrow \deloop^4_X\hat{\integers}(2)$.
	
	By the natural isomorphism (\emph{cf.}~Remark \ref{rem-trace-map-projection-formula}, Remark \ref{rem-kummer-class-of-tautological-line-bundle})
	\begin{align*}
		(\tr_x \circ \mathscr O(x)^*)(a &\otimes \Psi^{\otimes 2}) \simeq \tr_x(a\otimes \Psi(\mathscr O(x))^{\otimes 2}) \\
		& \simeq a\otimes \Psi(\mathscr O(x))|_{\Gamma_x} \otimes \tr_x(\Psi(\mathscr O(x))) \simeq a \otimes \Psi(\mathscr O(x)|_{\Gamma_x})
	\end{align*}
	linear in $a \in \Gamma(X, A(-1))$, we see that the restriction of $\tr_x \circ \mathscr O(x)^*$ to the $\Gamma(X, A(-1))$ is the tensor product with $\Psi(\mathscr O(x)|_{\Gamma_x})$.
	
	Since $\mathscr O(x)|_{\Gamma_x}$ is canonically isomorphic to the pullback of $\Omega_X^{-1}$ along $x : \Spec R \rightarrow X$. This is the tensor product with $\Psi(\Omega_X^{-1})$, pulled back along $x$.
\end{enumerate}

In conclusion, the decomposition \eqref{eq-multiplicative-group-etale-level-splitting} exhibits \eqref{eq-sharp-torus-commutative-gerbe-singleton-expression} as the sum $\nu + \epsilon \otimes \Psi(\Omega_X^{-1})$, pulled back along $x : \Spec R \rightarrow X$. The latter is, by definition, the $\mathbb E_{\infty}$-monoidal morphism denoted by $\nu + \vartheta$ (\emph{cf.}~\S\ref{void-theta-shift}).
\end{proof}

\begin{void}
We shall now use Proposition \ref{prop-sharp-torus-grassmannian-commutative-factorization} and Proposition \ref{prop-sharp-torus-commutative-gerbe-identification} to construct \eqref{eq-satake-equivalence-sharp-torus}.

\begin{proof}[Proof of Proposition \ref{prop-satake-equivalence-sharp-torus}]
As explained in \S\ref{void-satake-equivalence-sharp-torus-replacement}, we may take as input an $\mathbb E_{\infty}$-monoidal morphism \eqref{eq-symmetric-monoidal-etale-level} and replace the left-hand-side of \eqref{eq-satake-equivalence-sharp-torus} by $\Sat_{\mathscr G, \zeta}(\Hec_{T, I})$. Furthermore, the desired equivalence \eqref{eq-satake-equivalence-sharp-torus} is of \'etale local nature over $X^I$, so we may assume that $T$ is split and view $\Lambda$ as an abelian group.

Since \eqref{eq-sharp-torus-grassmannian-commutative-morphism} is $\mathbb E_{\infty}$-monoidal, we may endow $\Sat_{\mathscr G, \zeta}(\Gr_{T, I})$ with a convolution \emph{symmetric} monoidal structure $\circ$. Pushforward along \eqref{eq-torus-grassmannian-canonical-sections-total} yields a symmetric monoidal functor
\begin{equation}
\label{eq-sharp-torus-canonical-section-pushforward}
\Rep_{H^{\boxtimes I}, (\nu + \vartheta)^{\boxplus I}} \simeq \bigoplus_{\lambda^I \in \Lambda^I} \Lis_{(\nu + \vartheta)^{\boxplus I}(\lambda^I)}(X^I) \xrightarrow{\varpi_!} \Sat_{\mathscr G, \zeta}(\Gr_{T, I}),
\end{equation}
where we used Proposition \ref{prop-sharp-torus-commutative-gerbe-identification} to identify the pullback of \eqref{eq-sharp-torus-grassmannian-commutative-morphism} with $(\nu + \vartheta)^{\boxplus I}$. By Lemma \ref{lem-satake-category-torus-extension-characterization}, the functor \eqref{eq-sharp-torus-canonical-section-pushforward} is an equivalence.

On the other hand, pullback along $\Gr_{T, I} \rightarrow \Hec_{T, I}$ yields an equivalence
\begin{equation}
\label{eq-sharp-torus-hecke-to-grassmannian-pullback}
\Sat_{\mathscr G, \zeta}(\Hec_{T, I}) \simeq \Sat_{\mathscr G, \zeta}(\Gr_{T, I})
\end{equation}
as $\mathscr G_{\Hec_{T, I}}$ descends along \eqref{eq-local-hecke-stack-torus-quotient-onto-grassmannian}. Thus, the construction of \S\ref{void-fusion-product} endows $\Sat_{\mathscr G, \zeta}(\Gr_{T, I})$ with a fusion symmetric monoidal structure $\star$. The equivalence \eqref{eq-sharp-torus-hecke-to-grassmannian-pullback} is symmetric monoidal when both sides are equipped with the fusion symmetric monoidal structure.

Finally, by Proposition \ref{prop-sharp-torus-grassmannian-commutative-factorization} and the proof of Proposition \ref{prop-convolution-equals-fusion}, the convolution symmetric monoidal structure lifts $(\Sat_{\mathscr G, \zeta}(\Gr_{T, I}), \star)$ to a commutative monoid in the (2-)category $\CAlg(\Cat)$. Thus we obtain equivalences of \emph{symmetric} monoidal categories
\begin{equation}
\label{eq-sharp-torus-convolution-equals-fusion}
(\Sat_{\mathscr G, \zeta}(\Gr_{T, I}), \circ) \simeq (\Sat_{\mathscr G, \zeta}(\Gr_{T, I}), \star).
\end{equation}

The desired equivalence \eqref{eq-satake-equivalence-sharp-torus} is the composition of \eqref{eq-sharp-torus-hecke-to-grassmannian-pullback}, \eqref{eq-sharp-torus-convolution-equals-fusion}, and \eqref{eq-sharp-torus-canonical-section-pushforward}. It is functorial in $I$ by construction.
\end{proof}
\end{void}

\begin{rem}
We emphasize that the equivalence \eqref{eq-sharp-torus-hecke-to-grassmannian-pullback} would be \emph{false} without the hypothesis that the \'etale level \eqref{eq-symmetric-monoidal-etale-level} is $\mathbb E_{\infty}$-monoidal.

Indeed, for general \'etale levels, objects on $\Sat_{\mathscr G, \zeta}(\Hec_{T, I})$ are supported on the image of the morphism \eqref{eq-torus-sharp-isogeny-grassmannian}, by Proposition \ref{prop-torus-sharp-isogeny-satake-categories}.
\end{rem}

\medskip

\section{The $\hat Z_H$-grading}
\label{sec-grading-by-virtual-connected-components}

We put ourselves in the context of \S\ref{void-satake-equivalence-statement-context} and assume in addition that $G$ is \emph{split} reductive. Under this assumption, the dual group $H$ is constant (\emph{cf.}~\S\ref{sec-metaplectic-dual-data}).

The goal of this section is to endow the Satake category for $I = \{1\}$ with a $\hat Z_H$-grading compatible with all of its essential structures. This $\hat Z_H$-grading will correspond to the natural $\hat Z_H$-grading on $\Rep_{H, \nu + \vartheta}$ under the Satake equivalence \eqref{eq-satake-equivalence}. However, our proof of the latter makes critical use of this $\hat Z_H$-grading, so it must be constructed independently.

\subsection{Central equivariance}

\begin{void}
\label{void-central-equivariance-context}
In this subsection, we work over one copy of the curve $X$, so we base change the local Hecke stack (\emph{cf.}~\S\ref{void-local-hecke-stack-definition}) along $X \rightarrow \Ran$ and denote the result by $\Hec_G$. The same convention applies to $\Gr_G$, $L^+G$, $LG$ and we abbreviate $\deloop_X$ to $\deloop$.

Consider the $\deloop Z$-action on $\deloop G$ defined by multiplication (\emph{cf.}~\S\ref{void-canonical-quadratic-structure}). It induces an action of $L^+(\deloop Z)$ on $\Hec_G$. Explicitly, the action of an $R$-point $P_Z$ of $L^+(\deloop Z)$, viewed as $Z$-bundle over $D_x$ (for $x$ an $R$-point of $X$), carries $P^0\overset{x}{\sim} P^1$ to the modification
$$
P^0\otimes P_Z \overset{x}{\sim} P^1 \otimes P_Z.
$$
\end{void}

\begin{void}
The inclusion $\Lambda^{\sharp} \subset \Lambda$ caries $\Lambda^{\sharp}_{\sconn}$ into $\Lambda_{\sconn}$, so it induces a map $\hat Z_H \rightarrow \pi_1 G$ (\emph{cf.}~Remark \ref{rem-metaplectic-dual-center-as-quotient}). This map is neither injective nor surjective in general.

Given $\theta \in \hat Z_H$, we shall denote by $[\theta]$ its image in $\pi_1 G$ and by $\Hec_G^{[\theta]} \subset \Hec_G$ the corresponding connected component. Recall the $A$-gerbe $\mathscr G_{\Hec_G}$ over $\Hec_G$ (\emph{cf.}~\S\ref{void-local-hecke-stack-gerbe}), whose restriction to $\Hec_G^{[\theta]}$ is denote by $\mathscr G_{\Hec_G^{[\theta]}}$.
\end{void}

\begin{prop}
\label{prop-gerbe-central-equivariance}
Given any $\theta \in \hat Z_H$ with image $[\theta] \in \pi_1 G$, the $A$-gerbe $\mathscr G_{\Hec_G^{[\theta]}}$ admits a canonical $L^+(\deloop Z)$-equivariance structure.
\end{prop}

\begin{proof}
Given $R$-points $P_Z$ of $L^+(\deloop Z)$, $(P^0\overset{x}{\sim}P^1)$ of $\Hec_G^{\theta}$ lying over the same $R$-point $x$ of $X$, we shall construct a natural isomorphism of $A$-gerbes over $\Spec R$:
\begin{equation}
\label{eq-hecke-stack-gerbe-central-equivariance}
\mathscr G_{\Hec_G}(P^0\otimes P_Z \overset{x}{\sim} P^1 \otimes P_Z) \simeq \mathscr G_{\Hec_G}(P^0 \overset{x}{\sim} P^1)
\end{equation}
along with cocycle data for multiple $R$-points of $L^+(\deloop Z)$.

By construction, the left-hand-side of \eqref{eq-hecke-stack-gerbe-central-equivariance} is the image of
$$
\mu(P^0\otimes P_Z) - \mu(P^1 \otimes P_Z) \in \Gamma(D_x\text{ mod }\mathring D_x, \deloop^4A(1))
$$
under the trace map (\emph{cf.}~\S\ref{void-local-integration-map-construction}). Using the canonical quadratic structure \eqref{eq-canonical-quadratic-structure}, we may identify this section as follows:
\begin{equation}
\label{eq-canonical-quadratic-structure-hecke-stack-point}
\mu(P^0\otimes P_Z) - \mu(P^1 \otimes P_Z) \simeq \mu(P^0) - \mu(P^1) + (b_2\otimes\Psi^{\otimes 2})(P^0_{\abelian} - P^1_{\abelian}, P_Z)
\end{equation}
where $P^0_{\abelian}$, $P^1_{\abelian}$ are the images of $P^0$, $P^1$ under $\deloop G \rightarrow \deloop G_{\abelian}$.

Consider the linear form $b_2([\theta], \cdot) : \Fib(\Lambda \rightarrow \Lambda_{\adjoint}) \rightarrow A(-1)$, the induced $\integers$-linear morphism $b_2([\theta], \cdot)\otimes\Psi : \deloop Z \rightarrow \deloop^2 A$, and its value $(b_2([\theta], \cdot)\otimes\Psi)(P_Z) \in \Gamma(D_x, \deloop^2 A)$. We shall construct a canonical isomorphism of $A$-gerbes over $\Spec R$ linear in $P_Z$:
\begin{equation}
\label{eq-canonical-quadratic-structure-pairing-trace}
\tr_x((b_2\otimes\Psi^{\otimes 2})(P^0_{\abelian} - P^1_{\abelian}, P_Z)) \simeq (b_2([\theta], \cdot)\otimes\Psi)(P_Z)|_{\Gamma_x}.
\end{equation}
This will give rise to \eqref{eq-hecke-stack-gerbe-central-equivariance} for the following reason: The lift $\theta \in \hat Z_H$ of $[\theta]$ trivializes $b_2([\theta], \cdot)$ as a linear form, because $b$ vanishes over $\Lambda^{\sharp} \otimes \Lambda$ and $b_1$ vanishes over $\Lambda^{\sharp}_{\sconn} \otimes \Lambda_{\adjoint}$ (\emph{cf.}~\S\ref{void-additional-bilinear-forms}, Remark \ref{rem-metaplectic-root-lattice-as-kernel}). This then induces a trivialization of the left-hand-side of \eqref{eq-canonical-quadratic-structure-pairing-trace}, so \eqref{eq-canonical-quadratic-structure-hecke-stack-point} gives rise to \eqref{eq-hecke-stack-gerbe-central-equivariance} under $\tr_x$. The cocycle data are induced from those of the canonical quadratic structure (\emph{cf.}~\S\ref{void-canonical-quadratic-structure}) and the linearity of \eqref{eq-canonical-quadratic-structure-pairing-trace} in $P_Z$.

To construct \eqref{eq-canonical-quadratic-structure-pairing-trace}, we note that for any $R$-point $(P^0 \overset{x}{\sim} P^1)$ of $\Hec_G$ and any $\integers$-linear morphism $f : \pi_1 G \rightarrow \deloop^2 A$, we have a canonical isomorphism in $\Gamma(\Spec R, \deloop^2 A)$:
\begin{equation}
\label{eq-local-hecke-stack-trace-degree-isomorphism}
\tr_x(\Psi^f(P^1_{\abelian} - P^0_{\abelian})) \simeq f(\deg(P^0 \overset{x}{\sim} P^1)),
\end{equation}
where $\Psi^f : \deloop G_{\abelian} \rightarrow \deloop^4A(1)$ is the tensor product $f$ with $\Psi$ and $\deg : \Hec_G \rightarrow \pi_1 G$ is the degree map classifying connected components. Indeed, \eqref{eq-local-hecke-stack-trace-degree-isomorphism} follows from the universal case with $f$ the reduction of $\pi_1 G$ along $\integers\rightarrow\hat{\integers}$, where it amounts to the compatibility between degree and \'etale Chern classes. The isomorphism \eqref{eq-local-hecke-stack-trace-degree-isomorphism} depends linearly on $f$.

The isomorphism \eqref{eq-canonical-quadratic-structure-pairing-trace} is obtained by specializing \eqref{eq-local-hecke-stack-trace-degree-isomorphism} to $f := (b_2\otimes\Psi)(\cdot, P_Z)$, where $b_2\otimes\Psi$ is the tensor product of $b_2$ with $\Psi$ along the second factor.
\end{proof}

\begin{rem}
\label{rem-gerbe-central-equivariance-as-adjoint-equivariance}
Proposition \ref{prop-gerbe-central-equivariance} implies that given $\theta \in \hat Z_H$, the restriction $\mathscr G_{\Gr_G^{[\theta]}}$ of $\mathscr G_{\Gr_G}$ to the connected component $\Gr_G^{[\theta]}$ admits a canonical $L^+G_{\adjoint}$-equivariance structure.

Indeed, this follows from the canonical map of \'etale $X$-stacks
$$
L^+G_{\adjoint}\backslash \Gr_G^{[\theta]} \rightarrow \Hec_G^{[\theta]} / L^+(\deloop Z)
$$
induced from the morphism $L^+(\deloop Z) \rightarrow \deloop L^+Z \simeq L^+G_{\adjoint} / L^+G$.
\end{rem}

\begin{void}
It is critical in Proposition \ref{prop-gerbe-central-equivariance} that we take as input $\theta \in \hat Z_H$ rather than an element of $\pi_1 G$. Namely, different lifts of the same element of $\pi_1 G$ to $\hat Z_H$ provide different $L^+(\deloop Z)$-equivariance structures.

Let us be more precise. Define a linear map
\begin{equation}
\label{eq-gerbe-central-equivariance-discrepancy-character}
\chi : \ker(\hat Z_H \rightarrow \pi_1 G) \rightarrow \SHom(\deloop Z, \deloop A) (\simeq \SHom(Z, A))
\end{equation}
as follows: Given $\lambda \in \Lambda^{\sharp} \cap \Lambda_{\sconn}$, the morphism $b_1(\lambda, \cdot)\otimes\Psi : T_{\adjoint} \rightarrow \deloop A$ is canonically trivialized along $T \rightarrow T_{\adjoint}$, so it factors through a map $\deloop Z \rightarrow \deloop A$. This map vanishes for $\lambda \in \Lambda^{\sharp}_{\sconn}$, so it induces \eqref{eq-gerbe-central-equivariance-discrepancy-character} after identifying the source with $(\Lambda^{\sharp} \cap \Lambda_{\sconn}) / \Lambda^{\sharp}_{\sconn}$.
\end{void}

\begin{lem}
\label{lem-gerbe-central-equivariance-discrepancy-character-injective}
The map \eqref{eq-gerbe-central-equivariance-discrepancy-character} is injective.
\end{lem}

\begin{proof}
If $\lambda \in \Lambda^{\sharp} \cap \Lambda_{\sconn}$ satisfies $\chi(\lambda) = 0$, then $b_1(\lambda, \cdot)\otimes\Psi : T_{\adjoint} \rightarrow \deloop A$ must also vanish because it is the composition of $T_{\adjoint} \rightarrow \deloop Z$ with $\chi(\lambda)$. This means that $b_1(\lambda, \cdot)$ is the zero form on $\Lambda_{\adjoint}$. We conclude because $\Lambda^{\sharp}_{\sconn} \subset \Lambda_{\sconn}$ is precisely the kernel of $b_1$ (\emph{cf.}~Remark \ref{rem-metaplectic-root-lattice-as-kernel}).
\end{proof}

\begin{void}
Given $\theta_1, \theta_2 \in \hat Z_H$ with equal image $[\theta]\in \pi_1 G$, Proposition \ref{prop-gerbe-central-equivariance} supplies two $L^+(\deloop Z)$-equivariance structures on $\mathscr G_{\Hec_G^{[\theta]}}$, or equivalently two $A$-gerbes $\mathscr G_{\Hec_G}^{\theta_1}$, $\mathscr G_{\Hec_G}^{\theta_2}$ over the quotient stack $\Hec_G^{[\theta]}/L^+(\deloop Z)$ whose pullbacks to $\Hec_G^{[\theta]}$ are identified with $\mathscr G_{\Hec_G}$.

The difference
\begin{equation}
\label{eq-gerbe-central-equivariance-difference}
\mathscr G_{\Hec_G}^{\theta_2} - \mathscr G_{\Hec_G}^{\theta_1}
\end{equation}
may thus be viewed as an $A$-gerbe over $\Hec_G^{[\theta]}/L^+(\deloop Z)$ neutralized over $\Hec_G^{[\theta]}$.
\end{void}

\begin{prop}
\label{prop-gerbe-central-equivariance-discrepancy}
Given $\theta_1, \theta_2 \in \hat Z_H$ with equal image $[\theta] \in \pi_1 G$, the $A$-gerbe \eqref{eq-gerbe-central-equivariance-difference} is canonically identified with the pullback of the deloop of $\chi(\theta_2 - \theta_1)$ along the map
$$
\Hec_G^{[\theta]} / L^+(\deloop Z) \rightarrow X / L^+(\deloop Z) \rightarrow \deloop^2 Z.
$$
\end{prop}

\begin{proof}
Let us interpret \eqref{eq-gerbe-central-equivariance-discrepancy-character} in terms of the form $b_2$ (\emph{cf.}~\S\ref{void-additional-bilinear-forms}). Namely, tensoring $b_2$ with $\Psi$ along its second factor, we obtain the pairing
$$
b_2\otimes\Psi : \pi_1 G \otimes \deloop Z \rightarrow \deloop^2 A
$$
which is canonically trivialized over $\hat Z_H \otimes \deloop Z$ (\emph{cf.}~\S\ref{void-additional-bilinear-forms}, Remark \ref{rem-metaplectic-root-lattice-as-kernel}). It then induces the pairing adjoint to \eqref{eq-gerbe-central-equivariance-discrepancy-character}:
\begin{equation}
\label{eq-gerbe-central-equivariance-discrepancy-adjoint}
\ker(\hat Z_H \rightarrow \pi_1 G) \otimes \deloop Z \rightarrow \deloop A.
\end{equation}

On the other hand, the choices $\theta_1, \theta_2 \in \hat Z_H$ intervene in the construction of the $L^+(\deloop Z)$-equivariance of $\mathscr G_{\Hec_G}$ over $\Hec_G^{[\theta]}$ by providing different trivializations of \eqref{eq-canonical-quadratic-structure-pairing-trace}, for any $R$-point $P_Z$ of $L^+(\deloop Z)$ lying over $x \in X(R)$. The difference in the trivializations is precisely the image of $(\theta_2 - \theta_1)\otimes P_Z|_{\Gamma_x}$ under \eqref{eq-gerbe-central-equivariance-discrepancy-adjoint}.
\end{proof}

\subsection{Compatibility lemmas}

\begin{void}
We keep the notation of \S\ref{void-central-equivariance-context}.

Given an element $\theta \in \hat Z_H$ with image $[\theta] \in \pi_1 G$, Proposition \ref{prop-gerbe-central-equivariance} allows us to descend $\mathscr G_{\Hec_G^{[\theta]}}$ to an $A$-gerbe $\mathscr G_{\Hec_G}^{\theta}$ along the quotient map
$$
\Hec_G^{[\theta]} \rightarrow \Hec_G^{[\theta]}/L^+(\deloop Z).
$$

In this subsection, we show that the association $\theta \mapsto \mathscr G_{\Hec_G}^{\theta}$ is compatible with the convolution product on $\Hec_G$, passage to Levis, and restrictions to Schubert cells.
\end{void}

\begin{void}[Convolution product]
Given $\theta_1,\cdots, \theta_n \in \hat Z_H$ with images $[\theta_1], \cdots, [\theta_n] \in \pi_1 G$, we denote by $\Hec_G^{[\theta_1], \cdots, [\theta_n]}$ the substack of $\Hec_G^{[n]}$ (\emph{cf.}~\S\ref{void-local-hecke-stack-convolution}) corresponding to modifications
\begin{equation}
\label{eq-hecke-stack-convolution-prescribed-degree}
P^0 \overset{x}{\sim} P^1 \overset{x}{\sim} \cdots \overset{x}{\sim} P^n
\end{equation}
where $P^{j-1}\overset{x}{\sim} P^j$ belongs to $\Hec_G^{[\theta_j]}$ for each $j = 1, \cdots, n$.

There is an $L^+\deloop Z$-action on $\Hec_G^{[\theta_1],\cdots, [\theta_n]}$, where $P_Z$ carries \eqref{eq-hecke-stack-convolution-prescribed-degree} to its termwise tensor product with $P_Z$. The structural morphisms \eqref{eq-convolution-diagram} thus induce morphisms
$$
\begin{tikzcd}[column sep = 1.5em]
	\Hec_G^{[\theta_1], \cdots, [\theta_n]} / L^+\deloop Z \ar[r, "\prod p_j"]\ar[d, "m"] & \prod_{j = 1}^n \Hec_G^{[\theta_j]}/ L^+\deloop Z \\
	\Hec_G^{[\theta_1 + \cdots + \theta_n]} / L^+\deloop Z
\end{tikzcd}
$$
\end{void}

\begin{lem}
\label{lem-central-equivariance-convolution-compatibility}
The isomorphism \eqref{eq-gerbe-multiplicativity-isomorphism} over $\Hec_G^{[\theta_1], \cdots, [\theta_n]}$ canonically extends to an isomorphism of $A$-gerbes over $\Hec_G^{[\theta_1], \cdots, [\theta_n]} / L^+\deloop Z$:
$$
m^* \mathscr G_{\Hec_G}^{\theta_1 + \cdots + \theta_n} \simeq p_1^*\mathscr G_{\Hec_G}^{\theta_1} + \cdots + p_n^* \mathscr G_{\Hec_G}^{\theta_n}.
$$
\end{lem}

\begin{proof}
By construction, \eqref{eq-gerbe-multiplicativity-isomorphism} arises as the image of the isomorphism
\begin{equation}
\label{eq-gerbe-convolution-central-equivariance}
\mu(P^n) - \mu(P^0) \simeq \sum_{j = 1}^n (\mu(P^j) - \mu(P^{j-1}))
\end{equation}
in $\Gamma(D_x \text{ mod }\mathring D_x, \deloop^4A(1))$ under the trace map (\emph{cf.}~the proof of Lemma \ref{lem-local-hecke-stack-gerbe-multiplicative}). By inspecting the proof of Proposition \ref{prop-gerbe-central-equivariance}, it suffices to identify the trivializations of
$$
b_2([\theta_1 + \cdots + \theta_n], \cdot)\otimes\Psi \simeq \sum_{j = 1}^n b_2([\theta_j], \cdot)\otimes\Psi
$$
supplied by $\theta_1 + \cdots + \theta_n \in \hat Z_H$, respectively the product of those supplied by $\theta_1,\cdots, \theta_n \in \hat Z_H$. This holds because $b_2$ is trivialized over $\hat Z_H \otimes \Fib(\Lambda \rightarrow \Lambda_{\adjoint})$ as a \emph{bilinear} form.
\end{proof}

\begin{void}[Passage to Levis]
\label{void-central-equivariance-passage-to-levis}
Let $P \subset G$ be a parabolic subgroup with Levi quotient $P \twoheadrightarrow M$. We have natural maps $Z \rightarrow Z_M$ and $\pi_1 M \rightarrow \pi_1 G$. Moreover, the pullback of $\mu$ to $P$ canonically descends to an \'etale level $\mu_M$ for $M$ and we obtain the dual group $H_M$ for $(M, \mu_M)$. The map $\pi_1 M \rightarrow \pi_1 G$ lifts to a map $\hat Z_{H_M} \rightarrow \hat Z_H$.

Let $\lambda$ be an element of $\hat Z_{H_M}$. Its images in $\pi_1 M$, $\hat Z_H$, and $\pi_1 G$ are denoted by $[\lambda]$, $\theta$, and $[\theta]$, respectively. Consider the following quotients of the structural morphisms:
$$
\begin{tikzcd}[column sep = -2.5em]
	& \Hec_P^{[\lambda]} / L^+\deloop Z \ar[dl, swap, "q"]\ar[dr, "p"] \\
	\Hec_G^{[\theta]} / L^+\deloop Z & & \Hec_M^{[\lambda]} / L^+\deloop Z_M
\end{tikzcd}
$$
where $\Hec_P^{[\lambda]}$ denotes the pre-image of $\Hec_M^{[\lambda]}$ under $\Hec_P \rightarrow \Hec_M$.
\end{void}

\begin{lem}
\label{lem-central-equivariance-levi-compatibility}
The isomorphism \eqref{eq-constant-term-gerbe-compatibility} over $\Hec_P^{[\lambda]}$ canonically extends to an isomorphism of $A$-gerbes over $\Hec_P^{[\lambda]} / L^+\deloop Z$:
$$
q^* \mathscr G_{\Hec_G}^{\theta} \simeq p^* \mathscr G_{\Hec_M}^{\lambda}.
$$
\end{lem}

\begin{proof}
By inspecting the proof of Proposition \ref{prop-gerbe-central-equivariance}, this follows from the compatibility of the canonical quadratic structure \eqref{eq-canonical-quadratic-structure} with passage to Levis.

Namely, we consider the $\deloop Z$-action on $\deloop P$:
$$
a : \deloop Z \times \deloop P \rightarrow \deloop P,
$$
with projection maps $p_1, p_2$ from $\deloop Z\times \deloop P$ onto $\deloop Z$, respectively $\deloop P$. Then the canonical quadratic structures \eqref{eq-canonical-quadratic-structure} for $(G, \mu)$ and $(M, \mu_M)$ restrict to the same isomorphism
$$
a^* \mu_P \simeq (p_1)^*\mu_P + (p_2)^*\mu_Z + b_2\otimes\Psi^{\otimes 2},
$$
where $\mu_P$, $\mu_Z$ denote the restrictions of $\mu$ to $P$, respectively $Z$.
\end{proof}

\begin{void}[Schubert cells]
Let us now fix a maximal torus and a Borel subgroup $T \subset B \subset G$.

Given $\lambda \in \Lambda^{\sharp}$, the restriction of $\mathscr G_{\Hec_G}$ to the Schubert cell $\Hec_G^{\lambda}$ canonically descends to the $A$-gerbe $\mathscr G_{\varpi^{\lambda}}$ over $X$ (\emph{cf.}~Corollary \ref{cor-gerbe-schubert-cell-sharp-identification}). On the other hand, $\lambda$ has class $\theta \in \hat Z_H$, which gives rise to the $A$-gerbe $\mathscr G_{\Hec_G}^{\theta}$.

Consider the structural morphisms
$$
\begin{tikzcd}[column sep = 1em]
	\Hec_G^{\lambda} / L^+\deloop Z \ar[r, "j"]\ar[d, "\pi"] & \Hec_G^{[\theta]} / L^+\deloop Z \\
	X
\end{tikzcd}
$$
\end{void}

\begin{lem}
\label{lem-central-equivariance-schubert-cell-compatibility}
The isomorphism \eqref{eq-gerbe-schubert-cell-sharp-identification} over $\Hec_G^{\lambda}$ canonically extends to an isomorphism of $A$-gerbes over $\Hec_G^{\lambda} / L^+\deloop Z$:
\begin{equation}
\label{eq-central-equivariance-schubert-cell-compatibility}
j^*\mathscr G_{\Hec_G}^{\theta} \simeq \pi^*\mathscr G_{\varpi^{\lambda}}.
\end{equation}
\end{lem}

\begin{proof}
Note that $A$-gerbes over $\Hec_G^{\lambda} / L^+\deloop Z$ canonically descend to $\deloop(M^{\lambda}/Z)$ along the quotient of \eqref{eq-bialynicki-birula-map}. Moreover, $T_{\adjoint}$ may be thought of as a maximal torus of $M^{\lambda}/Z$, so $A$-gerbes over $\deloop(M^{\lambda}/Z)$ are uniquely determined after pulling back to $\deloop T_{\adjoint}$.

It therefore suffices to construct \eqref{eq-central-equivariance-schubert-cell-compatibility} after pulling back to $\Hec_T^{\lambda} / L^+\deloop Z$. The pullback of $\mathscr G_{\Hec_G}^{\theta}$ to $\Hec_T^{\lambda} / L^+\deloop Z$ extends to the $A$-gerbe $\mathscr G_{\Hec_T}^{\lambda}$ over $\Hec_T^{\lambda} / L^+\deloop T$ (\emph{cf.}~Lemma \ref{lem-central-equivariance-levi-compatibility}). It thus remains to compare $\mathscr G_{\Hec_T}^{\lambda}$ with the descent of $\mathscr G_{\Hec_T^{\lambda}}$ to $X$ supplied by \eqref{eq-gerbe-schubert-cell-sharp-identification} for the torus $T$. In the case of tori, the equivariance structure constructed in Proposition \ref{prop-gerbe-central-equivariance} is identical to one constructed in Corollary \ref{cor-gerbe-schubert-cell-sharp-identification}.
\end{proof}

\subsection{Construction of the $\hat Z_H$-grading}
\label{sec-virtual-connected-components}

\begin{void}
Recall that for any $S$-point $x$ of $X$ (with $S \in \Sch$), we have constructed the Satake category $\Sat_{\mathscr G, \zeta}(\Hec_{G, x})$ (\emph{cf.}~\S\ref{void-satake-subcategory}), which is an $\coeff$-linear category equipped with the convolution monoidal structure (\emph{cf.}~\S\ref{void-satake-category-general-base-monoidal}).

The goal of this subsection is to equip $\Sat_{\mathscr G, \zeta}(\Hec_{G, x})$ with a $\hat Z_H$-grading compatible with the monoidal structure and constant term functors (\emph{cf.}~\S\ref{void-constant-term-satake-category}).
\end{void}

\begin{void}
For each $\theta \in \hat Z_H$ with image $[\theta] \in \pi_1 G$, we shall consider the quotient map
\begin{equation}
\label{eq-hecke-stack-adjoint-quotient}
\Hec_{G, x}^{[\theta]} \rightarrow L^+_xG_{\adjoint} \backslash \Gr_{G, x}^{[\theta]}
\end{equation}
and view (the pulback of) $\mathscr G_{\Hec_G}^{\theta}$ as an $L^+_xG_{\adjoint}$-equivariant $A$-gerbe over $\Gr_{G, x}$ (\emph{cf.}~Remark \ref{rem-gerbe-central-equivariance-as-adjoint-equivariance}). Consider the full subcategory $\Sat_{\mathscr G, \zeta}(\Gr_{G, x}^{[\theta]}) \subset \derived_{\mathscr G, \zeta}(\Gr_{G, x}^{[\theta]})$ of perverse ULA sheaves relative to $S$ and the $\coeff$-linear category
$$
\Sat_{\mathscr G, \zeta}^{\theta}(\Hec_{G, x}) := \Sat_{\mathscr G, \zeta}(\Gr_{G, x}^{[\theta]})^{L_x^+G_{\adjoint}}
$$
of its $L^+_xG_{\adjoint}$-equivariant objects defined with respect to $\mathscr G_{\Hec_G}^{\theta}$.

Pullback along \eqref{eq-hecke-stack-adjoint-quotient} yields a functor
\begin{equation}
\label{eq-hecke-stack-adjoint-quotient-pullback}
\Sat_{\mathscr G, \zeta}^{\theta}(\Hec_{G, x}) \rightarrow \Sat_{\mathscr G, \zeta}(\Hec_{G, x}^{[\theta]}),
\end{equation}
whose target is naturally a full subcategory of $\Sat_{\mathscr G, \zeta}(\Hec_{G, x})$. Furthermore, \eqref{eq-hecke-stack-adjoint-quotient-pullback} is fully faithful and its essential image is closed under extensions (\emph{cf.}~Remark \ref{rem-satake-category-forgetful-to-grassmannian}).
\end{void}

\begin{prop}
\label{prop-virtual-connected-components}
The sum of \eqref{eq-hecke-stack-adjoint-quotient-pullback} over $\theta \in \hat Z_H$ yields an equivalence
\begin{equation}
\label{eq-virtual-connected-components}
	\bigoplus_{\theta \in \hat Z_H} \Sat_{\mathscr G, \zeta}^{\theta}(\Hec_{G, x}) \simeq \Sat_{\mathscr G, \zeta}(\Hec_{G, x}).
\end{equation}
\end{prop}

\begin{proof}
We first prove that the essential images of distinct summands in \eqref{eq-virtual-connected-components} are \emph{derived} orthogonal in $\derived_{\mathscr G, \zeta}(\Gr_{G, x})$. More precisely, given objects $\mathscr A_1, \mathscr A_2 \in \Sat_{\mathscr G, \zeta}(\Hec_{G, x})$ coming from summands corresponding to $\theta_1 \neq \theta_2 \in \hat Z_H$, we shall prove
\begin{equation}
\label{eq-gerbe-central-equivariance-derived-orthogonality}
\Hom(\mathscr A_1, \mathscr A_2) \simeq 0,
\end{equation}
where $\Hom(\cdot, \cdot)$ is the derived Hom functor of $\derived_{\mathscr G, \zeta}(\Gr_{G, x})$.

The statement is clear if $\theta_1$, $\theta_2$ have distinct images in $\pi_1 G$, as then $\mathscr A_1$, $\mathscr A_2$ would be supported on distinct connected components of $\Hec_{G, x}$, so let us assume that $\theta_1$, $\theta_2$ have equal image $[\theta] \in \pi_1 G$. Note that $\Hom(\mathscr A_1, \mathscr A_2)$ is the global section of $\SHom(\mathscr A_1, \mathscr A_2)$, which admits an $L^+_xG_{\adjoint}$-equivariance structure with respect to
\begin{equation}
\label{eq-gerbe-central-equivariance-discrepancy-grassmannian}
\mathscr G_{\Hec_G}^{\theta_2} - \mathscr G_{\Hec_G}^{\theta_1}.
\end{equation}
By Proposition \ref{prop-gerbe-central-equivariance-discrepancy}, \eqref{eq-gerbe-central-equivariance-discrepancy-grassmannian} coincides with the pullback of the $A$-gerbe over $\deloop G_{\adjoint}$ classified by the character $b_1(\theta_2 - \theta_1, \cdot) : \pi_1 G_{\adjoint} \rightarrow A(-1)$ along the projection
$$
L^+_x G_{\adjoint} \backslash \Gr_{G, x}^{[\theta]} \rightarrow \deloop L^+_x G_{\adjoint} \rightarrow \deloop G_{\adjoint}.
$$
The assumption $\theta_1 \neq \theta_2$ implies that $b_1(\theta_2 - \theta_1, \cdot)$ is nontrivial (\emph{cf.}~Remark \ref{rem-metaplectic-root-lattice-as-kernel}). The global section of $\SHom(\mathscr A_1, \mathscr A_2)$ thus vanishes by Proposition \ref{prop-vanishing-by-equivariance}.

The vanishing \eqref{eq-gerbe-central-equivariance-derived-orthogonality} implies that the essential images of distinct summands in \eqref{eq-virtual-connected-components} are orthogonal and have no nontrivial extensions in $\Sat_{\mathscr G, \zeta}(\Hec_{G, x})$. It remains to prove that they generate $\Sat_{\mathscr G, \zeta}(\Hec_{G, x})$ under extensions. For this, it suffices to prove that the standard functor $\Delta^{\lambda}$ for $\lambda \in \Lambda^{\sharp, +}$ (\emph{cf.}~\S\ref{void-standard-costandard}) factors through \eqref{eq-hecke-stack-adjoint-quotient-pullback}, for $\theta$ the class of $\lambda$.

Consider the commutative square
\begin{equation}
\label{eq-standard-functor-adjoint-equivariance}
\begin{tikzcd}[column sep = 1.5em]
	\Hec_{G, x}^{\lambda} \ar[r, phantom, "\subset"]\ar[d] & \Hec_{G, x}^{[\theta]} \ar[d] \\
	L_x^+G_{\adjoint} \backslash \Gr_{G, x}^{\lambda} \ar[r, phantom, "\subset"] & L_x^+G_{\adjoint} \backslash \Gr_{G, x}^{[\theta]}
\end{tikzcd}
\end{equation}
defined by inclusions of the Schubert cell. Using Lemma \ref{lem-central-equivariance-schubert-cell-compatibility}, we may factor the standard functor $\Delta^{\lambda}$ through extension along $\Gr_{G, x}^{\lambda} \subset \Gr_{G, x}^{[\theta]}$ as \emph{$L_x^+G_{\adjoint}$-equivariant} perverse sheaves, followed by the functor \eqref{eq-hecke-stack-adjoint-quotient-pullback}.
\end{proof}

\begin{cor}
\label{cor-virtual-connected-components-convolution-compatibility}
The $\hat Z_H$-grading \eqref{eq-virtual-connected-components} on $\Sat_{\mathscr G, \zeta}(\Hec_{G, x})$ is compatible with the convolution monoidal structure, \emph{i.e.}~
\begin{enumerate}
	\item the monoidal unit is homogeneous of degree $0$;
	\item given $\mathscr A_1, \mathscr A_2 \in \Sat_{\mathscr G, \zeta}(\Hec_{G, x})$ which are homogeneous of degrees $\theta_1, \theta_2$, the convolution product $\mathscr A_1 \circ \mathscr A_2$ is homogeneous of degree $\theta_1 + \theta_2$.
\end{enumerate}
\end{cor}

\begin{proof}
Statement (1) is clear.

Statement (2) follows from Lemma \ref{lem-central-equivariance-convolution-compatibility} and the definition of the convolution product (\emph{cf.}~\S\ref{void-convolution-product}).
\end{proof}

\begin{void}
Next, we formulate the compatibility between the $\hat Z_H$-grading \eqref{eq-virtual-connected-components} with constant term functors. We fix a parabolic subgroup $P\subset G$ with Levi quotient $P \twoheadrightarrow M$.

Recall the constant term functor $\CT_P$ on the Satake categories (\emph{cf.}~\S\ref{void-constant-term-satake-category}) and the natural map $\hat Z_{H_M} \rightarrow \hat Z_H$ (\emph{cf.}~\S\ref{void-central-equivariance-passage-to-levis}).
\end{void}

\begin{cor}
\label{cor-virtual-connected-components-constant-term-compatibility}
Given an object $\mathscr A \in \Sat_{\mathscr G, \zeta}(\Hec_{G, x})$ homogeneous of degree $\theta \in \hat Z_H$, its image $\CT_P(\mathscr A) \in \Sat_{\mathscr G, \zeta}(\Hec_{M, x})$ belongs to the sum of $\lambda$-graded components, where $\lambda \in \hat Z_{H_M}$ lies over $\theta$.
\end{cor}

\begin{proof}
Let $\lambda \in \hat Z_{H_M}$ be an element whose image in $\hat Z_H$ is $\theta_1 \neq \theta$. It suffices to prove that $\CT_P(\mathscr A)$ is orthogonal to any object $\mathscr B \in \Sat_{\mathscr G, \zeta}(\Hec_{M, x})$ homogeneous of degree $\lambda$. Since $\CT_P(\mathscr A)$ is supported on the connected components of $\Hec_{M, x}$ corresponding to elements of $\pi_1 M$ lying over $[\theta] \in \pi_1 G$, we may assume that $[\theta_1] = [\theta]$ in this proof.

We shall prove the vanishing of the derived Hom
\begin{equation}
\label{eq-constant-term-central-equivariance-orthogonality}
\Hom(\CT_P^{[\lambda]}(\mathscr A), \mathscr B) \simeq 0
\end{equation}
in $\derived_{\mathscr G, \zeta}(\Gr_{M, x}^{[\lambda]})$, where $\CT_P^{[\lambda]}(\mathscr A)$ denotes the restriction of $\CT_P(\mathscr A)$ to the connected component corresponding to $[\lambda] \in \pi_1 M$.

Note that $\CT_P^{[\lambda]}(\mathscr A)$ and $\mathscr B$ admit $L^+_x(M/Z)$-equivariance structures with respect to \emph{different} $L_x^+(M/Z)$-equivariance structures on $\mathscr G_{\Gr_{M, x}}$, induced from $\mathscr G_{\Hec_G}^{\theta}$, respectively $\mathscr G_{\Hec_G}^{\theta_1}$ (\emph{cf.}~Lemma \ref{lem-central-equivariance-levi-compatibility}). The difference $\mathscr G_{\Hec_G}^{\theta_1} - \mathscr G_{\Hec_G}^{\theta}$ is pulled back from the $A$-gerbe over $\deloop G_{\adjoint}$ classified by the character $b_1(\theta - \theta_1, \cdot) : \pi_1 G_{\adjoint} \rightarrow A(-1)$ (\emph{cf.}~Proposition \ref{prop-gerbe-central-equivariance-discrepancy}). Since $\theta \neq \theta_1$, this character is nonzero, so it defines a nontrivial $A$-gerbe over $\deloop(M/Z)$. The vanishing \eqref{eq-constant-term-central-equivariance-orthogonality} thus follows from Proposition \ref{prop-vanishing-by-equivariance} as in the proof of Proposition \ref{prop-virtual-connected-components}.
\end{proof}

\begin{rem}
It follows from the proof of Proposition \ref{prop-virtual-connected-components} that \eqref{eq-virtual-connected-components} is a \emph{coarsening} of the $\Lambda^{\sharp, +}$-grading provided by Theorem \ref{thm-semisimplicity}. More precisely, for each $\lambda \in \Lambda^{\sharp, +}$, the $\lambda$-graded component of Theorem \ref{thm-semisimplicity} is contained in the $\theta$-graded component of \eqref{eq-virtual-connected-components} for $\theta \in \hat Z_H$ the class of $\lambda$. However, there are \emph{no} analogues of Corollary \ref{cor-virtual-connected-components-convolution-compatibility} and Corollary \ref{cor-virtual-connected-components-constant-term-compatibility} for the $\Lambda^{\sharp, +}$-grading. If one \emph{defines} the $\hat Z_H$-grading by coarsening the $\Lambda^{\sharp, +}$-grading, there also seems to be no convenient way of proving these statements.

We find it more profitable to view \eqref{eq-virtual-connected-components} as a \emph{refinement} of the $\pi_1 G$-grading given by connected components of $\Hec_{G, x}$. Unlike Theorem \ref{thm-semisimplicity}, the construction of \eqref{eq-virtual-connected-components} can be extended to the derived category (replacing $L^+_xG_{\adjoint}$-equivariant objects by $L^+_xG_{\adjoint}$-monodromic objects), as well as to more general coefficients.
\end{rem}

\medskip

\section{Tannakian reconstruction}

In this section, we complete the construction of the geometric Satake equivalence (\emph{cf.}~Theorem \ref{thm-satake-equivalence}). We shall begin by lifting the Satake category ${}^+\Sat_{\mathscr G, \zeta}(\Hec_{G, I})$, for split $G$, to a $(\hat Z_H)^{\oplus I}$-graded \'etale sheaf of symmetric monoidal $\coeff$-linear categories over $X^I$ in the sense of \S\ref{void-graded-sheaf-symmetric-monoidal-categories}. This will allow us to construct an ``untwisted" Satake category, for which it is possible to define a fiber functor and apply the Tannakian formalism.

We remain in the context of \S\ref{void-satake-equivalence-statement-context} throughout this subsection.

\subsection{The fiber functor}
\label{sec-fiber-functor}

\begin{void}
We begin with a mild generalization of some constructions of \S\ref{sec-tori} and \S\ref{sec-virtual-connected-components}, adapted to a sheaf-theoretic version of the Satake category.

Let us return to the context of \S\ref{sec-satake-category}. Given a finite set $I$, the construction of the symmetric monoidal $\coeff$-linear category ${}^+\Sat_{\mathscr G, \zeta}(\Hec_{G, I})$ extends to the case where $X^I$ is replaced by an \'etale $X^I$-scheme $U$: The only change is that the fusion product is defined using the base change of $X^{\bigsqcup_n I, \disj} \subset X^{\bigsqcup_n I}$ to $U$ (\emph{cf.}~\S\ref{void-fusion-product}).

This construction is functorial in $U$, defining an \'etale sheaf of symmetric monoidal $\coeff$-linear categories ${}^+\SSat_{\mathscr G, \zeta}(\Hec_{G, I})$ over $X^I$.
\end{void}

\begin{void}[The case for tori]
For $G = T$ a torus, the construction of the geometric Satake equivalence for $(T, \mu)$ (\emph{cf.}~\S\ref{void-satake-equivalence-tori}) yields an equivalence of \'etale sheaves of symmetric monoidal $\coeff$-linear categories over $X^I$:
\begin{equation}
\label{eq-torus-satake-equivalence-sheaf}
{}^+\SSat_{\mathscr G, \zeta}(\Hec_{T, I}) \simeq \SRep_{T_H^{\boxtimes I}, (\nu + \vartheta)^{\boxplus I}},
\end{equation}
by the \'etale local nature of the construction.

Here, the normalization of the commutativity constraint has no effect (\emph{cf.}~\S\ref{void-normalized-satake-category}).
\end{void}

\begin{void}
Next, we assume that $G$ is split. Our goal is to generalize the $\hat Z_H$-grading \eqref{eq-virtual-connected-components} to ${}^+\SSat_{\mathscr G, \zeta}(\Hec_{G, I})$ for any finite set $I$.

The pair $(\Hec_G, \mathscr G_{\Hec_G})$ factorizes over the pairwise disjoint locus $X^{I, \disj} \subset X^I$ (\emph{cf.}~Proposition \ref{prop-local-hecke-stack-gerbe-factorization}), so Proposition \ref{prop-gerbe-central-equivariance} supplies an $L_I^+(\deloop Z)$-equivariance on the pullback of $\mathscr G_{\Hec_G}$ to the component
$$
\Hec_{G, I}^{\disj, [\theta^I]} := (\prod_{i\in I} \Hec_{G, \{i\}}^{[\theta^i]}) \times_{X^I} X^{I, \disj}
$$
associated to the each tuple $\theta^I \in (\hat Z_H)^{\oplus I}$ with image $[\theta^I] \in (\pi_1G)^{\oplus I}$. As in Remark \ref{rem-gerbe-central-equivariance-as-adjoint-equivariance}, this leads to an $L_I^+G_{\adjoint}$-equivariance on the pullback of $\mathscr G_{\Gr_G}$ to the corresponding component $\Gr_{G, I}^{\disj, [\theta^I]}$ of the affine Grassmannian over $X^{I, \disj}$.
\end{void}

\begin{void}
For each $\theta^I \in (\hat Z_H)^{\oplus I}$, we define a full subsheaf
\begin{equation}
\label{eq-virtual-connected-component-finite-set-sheaf}
\SSat_{\mathscr G, \zeta}^{\theta^I}(\Hec_{G, I}) \subset {}^+\SSat_{\mathscr G, \zeta}(\Hec_{G, I})
\end{equation}
as follows: Over each \'etale $X^I$-scheme $U$, it consists of objects of ${}^+\Sat_{\mathscr G, \zeta}(\Hec_{G, U})$ whose restriction to the pairwise disjoint locus comes from $L^+_I G_{\adjoint}$-equivariant objects of the category $\Sat_{\mathscr G, \zeta}(\Gr_{G, U}^{\disj, [\theta^I]})$, where the $L^+_IG_{\adjoint}$-equivariance on the restriction of $\mathscr G_{\Gr_G}$ is supplied by $\theta^I$, and the subscript $U$ means base change along $U \rightarrow X^I$.

The following results are straightforward adaptations of those of \S\ref{sec-virtual-connected-components}.
\end{void}

\begin{prop}
The functors \eqref{eq-virtual-connected-component-finite-set-sheaf} induce an equivalence
\begin{equation}
\label{eq-virtual-connected-component-grading-sheaf}
	\bigoplus_{\theta^I \in (\hat Z_H)^{\oplus I}} \SSat_{\mathscr G, \zeta}^{\theta^I}(\Hec_{G, I}) \simeq {}^+\SSat_{\mathscr G, \zeta}(\Hec_{G, I}).
\end{equation}
\end{prop}

\begin{proof}
This follows from Proposition \ref{prop-virtual-connected-components} together with the fact the restriction functor realizes $\Sat_{\mathscr G, \zeta}(\Gr_{G, I} \times_{X^I} U)$ as a full subcategory of $\Sat_{\mathscr G, \zeta}(\Gr_{G, I}^{\disj} \times_{X^{I, \disj}} U^{\disj})$ closed under direct summands (\emph{cf.}~\cite[Theorem 6.8]{MR4630128}).
\end{proof}

\begin{cor}
\label{cor-virtual-connected-components-convolution-sheaf}
The $(\hat Z_H)^{\oplus I}$-grading \eqref{eq-virtual-connected-component-grading-sheaf} on $\SSat_{\mathscr G, \zeta}(\Hec_{G, I})$ is compatible with the convolution monoidal structure, \emph{i.e.}~
\begin{enumerate}
	\item the monoidal unit is homogeneous of degree $0$;
	\item given $\mathscr A_1, \mathscr A_2 \in \SSat_{\mathscr G, \zeta}(\Hec_{G, I})$ which are homogeneous of degrees $\theta_1^I, \theta_2^I$, the convolution product $\mathscr A_1 \circ \mathscr A_2$ is homogeneous of degree $\theta_1^I + \theta_2^I$.
\end{enumerate}
\end{cor}

\begin{proof}
This follows from Corollary \ref{cor-virtual-connected-components-convolution-compatibility} and the fact that the convolution product factorizes over $X^I_{\disj}$ (\emph{cf.}~the proof of Proposition \ref{prop-convolution-equals-fusion}).
\end{proof}

\begin{cor}
\label{cor-virtual-connected-components-constant-term-sheaf}
The $(\hat Z_H)^{\oplus I}$-grading \eqref{eq-virtual-connected-component-grading-sheaf} on $\SSat_{\mathscr G, \zeta}(\Hec_{G, I})$ is compatible with constant term functors, \emph{i.e.}~given a parabolic subgroup $P \subset G$ with Levi quotient $P \twoheadrightarrow M$, the functor $\CT_P$ carries the $\theta^I$-graded component into the sum of the $\lambda^I$-graded components where $\lambda^I \in (\hat Z_{H_M})^{\oplus I}$ lies over $\theta^I$.
\end{cor}

\begin{proof}
This follows from Corollary \ref{cor-virtual-connected-components-constant-term-compatibility} and the fact that $\CT_P$ factorizes over $X^I_{\disj}$.
\end{proof}

\begin{void}
It follows from Corollary \ref{cor-virtual-connected-components-convolution-sheaf} that the assignment of the $\theta^I$-graded component \eqref{eq-virtual-connected-component-finite-set-sheaf} to each $\theta^I \in (\hat Z_H)^{\oplus I}$ is lax monoidal. By equipping it with the commutativity constraint inherited from ${}^+\SSat_{\mathscr G, \zeta}(\Hec_{G, I})$, we lift the latter to a $(\hat Z_H)^{\oplus I}$-graded \'etale sheaf of \emph{symmetric} monoidal $\coeff$-linear categories in the sense of \S\ref{void-graded-sheaf-symmetric-monoidal-categories}.

In particular, we may use any $\mathbb E_{\infty}$-monoidal morphism $(\hat Z_H)^{\oplus I} \rightarrow \deloop^2\coeff^{\times}$ to form the symmetric monoidal twist of ${}^+\SSat_{\mathscr G, \zeta}(\Hec_{G, I})$ (\emph{cf.}~\S\ref{void-symmetric-monoidal-twist}). We apply this construction to the \emph{opposite} of \eqref{eq-metaplectic-dual-morphism-external-sum} to obtain a sheaf of symmetric monoidal $\coeff$-linear categories
\begin{equation}
\label{eq-satake-category-untwisted-sheaf}
{}^+\SSat_{\mathscr G, \zeta}(\Hec_{G, I})_{-(\nu + \vartheta)^{\boxplus I}},
\end{equation}
which inherits a $(\hat Z_H)^{\oplus I}$-grading. We refer to its global section ${}^+\Sat_{\mathscr G, \zeta}(\Hec_{G, I})_{-(\nu + \vartheta)^{\boxplus I}}$ as the \emph{untwisted Satake category} and to \eqref{eq-satake-category-untwisted-sheaf} as the \emph{sheaf of untwisted Satake categories}.
\end{void}

\begin{void}
\label{void-constant-term-functor-untwisted}
Let $B \subset G$ be a Borel subgroup. Using the compatibility between the $(\hat Z_H)^{\oplus I}$-grading \eqref{eq-virtual-connected-component-grading-sheaf} with constant term functors (\emph{cf.}~Corollary \ref{cor-virtual-connected-components-constant-term-sheaf}), we obtain a constant term functor between the sheaves of untwisted Satake categories for $G$ and $T$:
\begin{equation}
\label{eq-constant-term-functor-untwisted-sheaf}
{}^+\SSat_{\mathscr G, \zeta}(\Hec_{G, I})_{-(\nu + \vartheta)^{\boxplus I}} \rightarrow {}^+\SSat_{\mathscr G, \zeta}(\Hec_{T, I})_{-(\nu + \vartheta)^{\boxplus I}}
\end{equation}

We shall argue that \eqref{eq-constant-term-functor-untwisted-sheaf} is ``independent of the choice of $B$." To make sense of this statement, we need to justify why the target of \eqref{eq-constant-term-functor-untwisted-sheaf} is indepedent of $B$. (Note that the induced \'etale level $\mu_T$ of $T$ \emph{depends} on $B$.) By Proposition \ref{prop-torus-sharp-isogeny-satake-categories}, we may replace the target of \eqref{eq-constant-term-functor-untwisted-sheaf} by ${}^+\SSat_{\mathscr G, \zeta}(\Hec_{T^{\sharp}, I})_{-(\nu + \vartheta)^{\boxplus I}}$, but the induced \'etale level $\mu_{T^{\sharp}}$ of $T^{\sharp}$ is canonically independent of $B$ (\emph{cf.}~Remark \ref{rem-commutative-cover-independence-of-borel}).

It thus suffices to show that the constant term functor $\CT_B$ is independent of the choice of $B$ upon replacing its target by ${}^+\Sat_{\mathscr G, \zeta}(\Hec_{T^{\sharp}, I})$, as we do now.
\end{void}

\begin{lem}
\label{lem-constant-term-independence-of-borel}
The symmetric monoidal functor
\begin{equation}
\label{eq-constant-term-functor-sharp-valued}
\CT_B : {}^+\Sat_{\mathscr G, \zeta}(\Hec_{G, I}) \rightarrow {}^+\Sat_{\mathscr G, \zeta}(\Hec_{T^{\sharp}, I})
\end{equation}
is canonically independent of the Borel subgroup $B \subset G$.
\end{lem}

\begin{proof}
Recall that the target of \eqref{eq-constant-term-functor-sharp-valued} embeds fully faithfully into ${}^+\Sat_{\mathscr G, \zeta}(\Hec_{T^{\sharp}, I} \times_{X^I} X^{I, \disj})$, where $X^{I, \disj} \subset X^I$ is the pairwise disjoint locus (\emph{cf.}~Proposition \ref{prop-satake-category-disjoint-pullback}(1)). By compatibility of $\CT_B$ with factorization, we reduce to $I = \{1\}$ and suppress it from the notation. It remains to prove that the functor
\begin{equation}
\label{eq-constant-term-functor-sharp-valued-singleton}
\CT_B : \Sat_{\mathscr G, \zeta}(\Hec_G) \rightarrow \Sat_{\mathscr G, \zeta}(\Hec_{T^{\sharp}})
\end{equation}
is independent of $B$.

To do so, we construct \eqref{eq-constant-term-functor-sharp-valued-singleton} using the universal Borel subgroup of $G$ over $\Torel$ and verify that the result descends along $\Torel \rightarrow X$ (\emph{cf.}~Remark \ref{rem-commutative-cover-independence-of-borel}). However, since the geometric fibers of $\Torel$ are simply connected, every $\coeff$-local system over $\Torel$ descends to $X$.
\end{proof}

\begin{void}[The fiber functor]
\label{void-fiber-functor-definition}
Define the \emph{fiber functor} to be the composition
\begin{align}
\notag
	\omega^I : {}^+\SSat_{\mathscr G, \zeta}(\Hec_{G, I})_{-(\nu + \vartheta)^{\boxplus I}} &\xrightarrow{\eqref{eq-constant-term-functor-untwisted-sheaf}} {}^+\SSat_{\mathscr G, \zeta}(\Hec_{T, I})_{-(\nu + \vartheta)^{\boxplus I}} \\
\label{eq-fiber-functor-finite-set}
	& \simeq \SRep_{T_H^{\boxtimes I}} \rightarrow \SLis(X^I),
\end{align}
where the equivalence is the untwist of \eqref{eq-torus-satake-equivalence-sheaf} and the last functor is the forgetful functor, for $T_H \subset H$ the canonical maximal torus.

The basic properties of \eqref{eq-fiber-functor-finite-set} are summarized as follows.
\end{void}

\begin{prop}
\label{prop-fiber-functor-properties}
The functor \eqref{eq-fiber-functor-finite-set} is symmetric monoidal, exact, conservative, and $\SLis(X^I)$-linear.
\end{prop}

\begin{proof}
The symmetric monoidal structure on \eqref{eq-fiber-functor-finite-set} comes from the symmetric monoidal structure on $\CT_B$ (\emph{cf.}~Lemma \ref{lem-constant-term-symmetric-monoidal}) and the functoriality of symmetric monoidal twists. The exactness and conservativity of \eqref{eq-fiber-functor-finite-set} are consequences of Proposition \ref{prop-constant-term-reflective-properties}.

The $\SLis(X^I)$-linearity is immediate, as all constructions we perform respect the $\SLis(X^I)$-module structure of the Satake category.
\end{proof}

\subsection{The Tannaka dual $\check G$}

\begin{void}
We continue to assume that $G$ is split.

For any finite set $I$ and any \'etale $X^I$-scheme $U$, the value of \eqref{eq-fiber-functor-finite-set} at $U$ yields a symmetric monoidal, exact, conservative $\Lis(U)$-linear functor
\begin{equation}
\label{eq-fiber-functor-value-etale-cover}
\omega^I : {}^+\Sat_{\mathscr G, \zeta}(\Hec_{G, U})_{-(\nu + \vartheta)^{\boxplus I}} \rightarrow \Lis(U),
\end{equation}
where $\Hec_{G, U}$ is the base change of $\Hec_{G, I}$ to $U$ (\emph{cf.}~Proposition \ref{eq-fiber-functor-finite-set}).

We wish to apply the relative Tannakian formalism (\emph{cf.}~\cite[Proposition VI.10.2]{fargues2021geometrization}) to \eqref{eq-fiber-functor-value-etale-cover} and identify the resulting bi-algebra object $\mathscr A_U \in \Ind\Lis(U)$ with (the restriction of) the external tensor product of $\mathscr A_X$, associated to $I = \{1\}$ and $U = X$. For this purpose, we need to filter the source of \eqref{eq-fiber-functor-value-etale-cover} by full subcategories on which it is co-representable.
\end{void}

\begin{void}[Weight functors]
\label{void-weight-functors}
Let us begin with the case $I = \{1\}$.

For each $\lambda \in \Lambda^{\sharp}$, we consider the post-composition of \eqref{eq-constant-term-functor-untwisted-sheaf} with the restriction functor along $\varpi^{\lambda} : X \rightarrow \Hec_T$ (\emph{cf.}~\S\ref{void-schubert-cell-notation}):
\begin{equation}
\label{eq-weight-functor-satake-category}
\omega^{\lambda} : \SSat_{\mathscr G, \zeta}(\Hec_{G, \{1\}})_{-(\nu + \vartheta)} \rightarrow \SLis(X),
\end{equation}
where we used the canonical identification of $(\varpi^{\lambda})^*\mathscr G_{\Hec_T}$ with $(\nu + \vartheta)(\lambda)$ to identify the target with $\SLis(X)$ (\emph{cf.}~Proposition \ref{prop-sharp-torus-commutative-gerbe-identification}).

Note that the functor \eqref{eq-weight-functor-satake-category} vanishes on all but one component of the source with respect to the $\hat Z_H$-grading, namely the one corresponding to the image $\theta \in \hat Z_H$ of $\lambda$ (\emph{cf.}~Corollary \ref{cor-virtual-connected-components-constant-term-compatibility}). It may thus be viewed as a functor
\begin{equation}
\label{eq-weight-functor-satake-category-homogeneous}
\omega^{\lambda} : \SSat_{\mathscr G, \zeta}^{\theta}(\Hec_{G, \{1\}})_{-(\nu + \vartheta)(\lambda)} \rightarrow \SLis(X).
\end{equation}

By construction, \eqref{eq-weight-functor-satake-category-homogeneous} is the $-(\nu + \vartheta)(\lambda)$-twist of the $\lambda$-graded component of the constant term functor $\CT_B$. We refer to \eqref{eq-weight-functor-satake-category-homogeneous} as the \emph{$\lambda$-weight functor}.
\end{void}

\begin{void}
\label{void-fiber-functor-corepresetable-singleton}
Given an $L^+_{\{1\}} G$-stable closed subscheme $Z$ of $\Gr_{G, \{1\}}$, we may consider the full subsheaf of categories
\begin{equation}
\label{eq-arc-group-stable-closed-subscheme-satake-subcategory}
\SSat_{\mathscr G, \zeta}(L^+_{\{1\}}G \backslash Z)_{-(\nu + \vartheta)} \subset {}^+\SSat_{\mathscr G, \zeta}(\Hec_{G, \{1\}})_{-(\nu + \vartheta)}
\end{equation}
consisting of objects supported on (the base change of) $Z$.

Restricting the functor \eqref{eq-fiber-functor-value-etale-cover} (for $I = \{1\}$) along \eqref{eq-arc-group-stable-closed-subscheme-satake-subcategory}, we obtain a functor
\begin{equation}
\label{eq-fiber-functor-arc-group-stable-restriction}
	\omega^{\{1\}} : \SSat_{\mathscr G, \zeta}(L^+_{\{1\}}G \backslash Z)_{-(\nu + \vartheta)} \rightarrow \SLis(X).
\end{equation}
\end{void}

\begin{lem}
\label{lem-fiber-functor-corepresentable-singleton}
The functor \eqref{eq-fiber-functor-arc-group-stable-restriction} is co-representable (over any \'etale $X$-scheme).
\end{lem}

\begin{proof}
It suffices to prove that for each $\lambda \in \Lambda^{\sharp}$, the restriction of the $\lambda$-weight functor \eqref{eq-weight-functor-satake-category-homogeneous} to objects supported on $Z$ is co-representable. More precisely, writing $\theta \in \hat Z_H$ for the image of $\lambda$, we need to construct an object
\begin{equation}
\label{eq-weight-functor-corepresenting-object}
	\mathscr P_Z^{\lambda} \in \Sat^{\theta}_{\mathscr G, \zeta}(L^+_{\{1\}} G \backslash Z)_{-(\nu + \vartheta)},
\end{equation}
together with an isomorphism of $\coeff$-local systems over $X$:
\begin{equation}
\label{eq-weight-functor-corepresenting-isomorphism}
\pi_*\SHom(\mathscr P_Z^{\lambda}, \mathscr A) \simeq \omega^{\lambda}(\mathscr A)
\end{equation}
natural in $\mathscr A \in \Sat^{\theta}_{\mathscr G, \zeta}(L^+_{\{1\}} G \backslash Z)_{-(\nu + \vartheta)}$. Here,  $\pi : L^+_{\{1\}}G \backslash Z \rightarrow X$ denotes the projection, and the formation of \eqref{eq-weight-functor-corepresenting-object}, \eqref{eq-weight-functor-corepresenting-isomorphism} must commute with \'etale base change along $X$.

Using hyperbolic localization (\emph{cf.}~the proof of Proposition \ref{prop-constant-term-reflective-properties}), we may present the restriction of \eqref{eq-weight-functor-satake-category-homogeneous} to objects supported on $Z$ as the functor $(p^-)_*(q^-)^!(\frac{d_{G, T}}{2})[d_{G, T}]$, formed with respect to the morphisms
$$
\begin{tikzcd}[column sep = 0em]
	& Z \cap \Gr_{B, \{1\}}^{-, \lambda} \ar[dl, swap, "q^-"]\ar[dr, "p^-"] & \\
	Z & & X
\end{tikzcd}
$$
and the trivialization of $\mathscr G_{\Gr_G}- (\nu + \vartheta)(\theta)$ over $Z \cap \Gr_{B, \{1\}}^{-, \lambda}$ specified in \S\ref{void-weight-functors}.

On the other hand, the functor of forgetting the $L^+_{\{1\}}G_{\adjoint}$-equivariance
$$
\derived_{\mathscr G, \zeta}(L^+_{\{1\}} G_{\adjoint} \backslash Z) \rightarrow \derived_{\mathscr G, \zeta}(Z),
$$
where $\mathscr G_{\Gr_G}$ is endowed with the $L^+_{\{1\}}G_{\adjoint}$-equivariance structure corresponding to $\theta$, admits a left adjoint $\Av_!^{\theta}$ (\emph{cf.}~\cite[\S A.2]{MR3839695}). Since the formation of $\Av_!^{\theta}$ is of \'etale local nature over $X$, we may twist it by $-(\nu + \vartheta)(\theta)$ and still write $\Av_!^{\theta}$ for the resulting functor. By adjunction, we have a natural isomorphism
$$
\omega^{\lambda}(\mathscr A) \simeq \pi_*\SHom({}^pH^0 \Av^{\theta}_! (q^-)_! \coeff(-\frac{d_{G, T}}{2})[-d_{G, T}], \mathscr A).
$$

It follows that the following $L^+_{\{1\}}G_{\adjoint}$-equivariant perverse sheaf over $Z$:
$$
\mathscr P^{\lambda}_Z := {}^pH^0 \Av^{\theta}_! (q^-)_! \coeff(-\frac{d_{G, T}}{2})[-d_{G, T}]
$$
co-represents the functor $\omega^{\lambda}$. We need to show that it belongs to the Satake category.

The only remaning property to check is that $\mathscr P_Z^{\lambda}$ is ULA relative to $X$. To prove this, we induct on the number of $L^+_{\{1\}}G$-orbits in $Z$. Let $j : \Gr_{G, \{1\}}^{\lambda'} \hookrightarrow Z$ be an open Schubert cell, with complementary (reduced) closed subscheme $i : Z' \hookrightarrow Z$. Then $j^*\mathscr P_Z^{\lambda}$ is locally constant perverse, hence of the form $\widetilde{\mathscr E}$ for some $\coeff$-local system $\mathscr E$ over $X$ (\emph{cf.}~notation of \S\ref{void-standard-costandard}). We claim that the co-unit of the adjunction
\begin{equation}
\label{eq-weight-functor-corepresenting-object-open-cell}
{}^pH^0 j_! \widetilde{\mathscr E} \rightarrow \mathscr P_Z^{\lambda}
\end{equation}
is \emph{injective}: This holds because of Corollary \ref{cor-standard-ic-agreement}, which exhibits $\mathscr P^{\lambda}_Z \rightarrow {}^pH^0 j_*\widetilde{\mathscr E}$ as a retraction of \eqref{eq-weight-functor-corepresenting-object-open-cell}. By the Cousin triangle, the cokernel of \eqref{eq-weight-functor-corepresenting-object-open-cell} is identified with $i_* {}^pH^0 i^*\mathscr P_Z^{\lambda}$, but we have an isomorphism
$$
{}^pH^0 i^*\mathscr P_Z^{\lambda} \simeq \mathscr P_{Z'}^{\lambda}
$$
because they co-represent the same functor.

It now follows from the ULA property of $\Delta^{\lambda'}(\mathscr E)$ (\emph{cf.}~Proposition \ref{prop-standard-costandard}) and the induction hypothesis that $\mathscr P_Z^{\lambda}$ is an extension of ULA objects, hence ULA.
\end{proof}

\begin{void}
\label{void-satake-category-support-filtration-finite-set}
In the context of \S\ref{void-fiber-functor-corepresetable-singleton}, we write
$$
\mathscr P_Z \in \Sat_{\mathscr G, \zeta}(L^+_{\{1\}}G \backslash Z)_{-(\nu + \vartheta)}
$$
for the object co-representing \eqref{eq-fiber-functor-arc-group-stable-restriction}, which exists by Lemma \ref{lem-fiber-functor-corepresentable-singleton}. (More precisely, for any \'etale $X$-scheme $U$, the pullback $\mathscr P_Z |_U$ co-represents the section of \eqref{eq-fiber-functor-arc-group-stable-restriction} over $U$.)

More generally, given a finite set $I$ and an $I$-tuple $\{Z_i\}_{i\in I}$ of $L^+_{\{1\}}G$-stable closed subschemes of $\Gr_{G, \{1\}}$, we may consider the full subsheaf of categories
\begin{equation}
\label{eq-satake-category-support-filtration}
\mathscr C_{\{Z_i\}} \subset {}^+\SSat_{\mathscr G, \zeta}(\Hec_G)_{-(\nu + \vartheta)^{\boxplus I}}
\end{equation}
consisting of objects whose restrictions along the pairwise disjoint locus $X^{I, \disj} \subset X^I$ are supported on (the base change of) $\prod_{i\in I} Z_i$ under the factorization isomorphism.

The $I$-tuples $\{Z_i\}_{i\in I}$ form a poset under inclusion and \eqref{eq-satake-category-support-filtration} exhibits the target as the filtered colimit of $\mathscr C_{\{Z_i\}}$ over $\{Z_i\}_{i\in I}$.
\end{void}

\begin{lem}
\label{lem-fiber-functor-corepresentable}
The restriction of \eqref{eq-fiber-functor-value-etale-cover} to $\mathscr C_{\{Z_i\}}$ is co-represented by the external fusion product of $\mathscr P_{Z_i}$ over $i\in I$.
\end{lem}

\begin{proof}
Since restriction to the pairwise disjoint locus is fully faithful (\emph{cf.}~Proposition \ref{prop-satake-category-disjoint-pullback}), this reduces to the case where $I$ is a singleton.
\end{proof}

\begin{void}
For any finite set $I$ and \'etale $X^I$-scheme $U$, we apply the relative Tannakian reconstruction \cite[Proposition VI.10.2]{fargues2021geometrization} to the fiber functor \eqref{eq-fiber-functor-value-etale-cover}, where the source is endowed with the filtration \eqref{eq-satake-category-support-filtration} on which \eqref{eq-fiber-functor-value-etale-cover} is co-representable (\emph{cf.}~Lemma \ref{lem-fiber-functor-corepresentable}).

This yields a bi-algebra object $\mathscr A_U$ of $\Ind\Lis(U)$ and lifts \eqref{eq-fiber-functor-value-etale-cover} to a symmetric monoidal, $\Lis(U)$-linear equivalence natural in $U$:
\begin{equation}
\label{eq-tannaka-reconstruction-equivalence}
	{}^+\Sat_{\mathscr G, \zeta}(\Hec_{G, U})_{-(\nu + \vartheta)^{\boxplus I}} \simeq \mathscr A_U\Comod(\Lis(U)).
\end{equation}

Denote by $\mathscr A_X$ the result of this construction for $I = \{1\}$ and $U = X$.
\end{void}

\begin{prop}
\label{prop-tannaka-group-reduction-to-singleton}
For any \'etale $X^I$-scheme $U$, the bi-algebra $\mathscr A_U$ is canonically identified with the pullback along $U \rightarrow X^I$ of the external tensor product $\boxtimes_{i\in I} \mathscr A_X$.
\end{prop}

\begin{proof}
By \cite[Proposition VI.10.2]{fargues2021geometrization} and Lemma \ref{lem-fiber-functor-corepresentable}, there is a canonical isomorphism of bi-algebras
$$
\mathscr A_U \simeq \colim_{\{Z_i\}} \omega^I(\mathscr P_{\{Z_i\}}|_U)^{\vee},
$$
where each $\mathscr P_{\{Z_i\}}$ is the external fusion product of $\mathscr P_{Z_i}$ over $i\in I$. Since $\omega^I$ carries external fusion products to external tensor products, the result follows.
\end{proof}

\begin{void}
\label{void-tannaka-dual-monoid}
Finally, we may define $\check G$ to be the Tannaka dual of $\mathscr A_X$, viewed as a locally constant \'etale sheaf of affine monoid $\coeff$-schemes over $X$.
\end{void}

\begin{cor}
\label{cor-tannaka-dual-equivalence}
For any finite set $I$, there is a canonical equivalence of \'etale sheaves of symmetric monoidal $\coeff$-linear categories over $X^I$:
\begin{equation}
\label{eq-tannaka-dual-equivalence}
	{}^+\SSat_{\mathscr G, \zeta}(\Hec_{G, I})_{-(\nu + \vartheta)^{\boxplus I}} \simeq \SRep_{\check G^{\boxtimes I}}.
\end{equation}
\end{cor}

\begin{proof}
This is a restatement of the equivalence \eqref{eq-tannaka-reconstruction-equivalence}, where $U$ ranges over all \'etale $X^I$-schemes, and we use Proposition \ref{prop-tannaka-group-reduction-to-singleton} to identify $\mathscr A_U\Comod(\Lis(U))$ with $\SRep_{\check G^{\boxtimes I}}(U)$.
\end{proof}

\subsection{Discrepancy $A$-torsors}

\begin{void}
We assume that $G$ is split and choose a Borel subgroup $B \subset G$ and a splitting of the quotient $B \twoheadrightarrow T$. We shall also specialize to $I = \{1\}$ and omit $\{1\}$ from the notation (\emph{i.e.}~we write $\Gr_G$ instead of $\Gr_{G, \{1\}}$).

The goal of this subsection is to give an explicit description of the weight functor \eqref{eq-weight-functor-satake-category-homogeneous} applied to generators of the untwisted Satake category. For this purpose, we shall define certain $A$-torsors over Mirkovi\'c--Vilonen cycles which measure the discrepancy between two trivializations of $\mathscr G_{\Gr_G}$.
\end{void}

\begin{void}
Fix $\lambda_1 \in \Lambda^{\sharp, +}$, $\lambda\in\Lambda^{\sharp}$ with equal image $\theta\in\hat Z_H$.

Denote by $S^{-, \lambda}$ the $LN^-$-orbit of $\varpi^{\lambda}$ in $\Gr_G$, where $N^- \subset B^-$ is the unipotent radical of the Borel opposite to $B$. Consider the \emph{Mirkovi\'c--Vilonen cycle} $\Gr_G^{\lambda_1} \cap S^{-, \lambda}$ in $\Gr_G$. It comes equipped with structural morphisms
$$
\begin{tikzcd}[column sep = 1em]
	\Gr_G^{\lambda_1} \ar[r, phantom, "\supset"]\ar[d] & \Gr_G^{\lambda_1} \cap S^{-, \lambda} \ar[r, phantom, "\subset"] & S^{-, \lambda} \ar[d] \\
	X \ar[rr, phantom, "\simeq"] & & X
\end{tikzcd}
$$

The restriction of $\mathscr G_{\Gr_G}$ to $\Gr_G^{\lambda_1}$ is identified with the pullback of $\mathscr G_{\varpi^{\lambda_1}}$ (\emph{cf.}~Corollary \ref{cor-gerbe-schubert-cell-sharp-identification}). The restriction of $\mathscr G_{\Gr_G}$ to $S^{-, \lambda}$ is identified with the pullback of $\mathscr G_{\varpi^{\lambda}}$ by descent along $\Gr_{B^-} \rightarrow \Gr_T$. Meanwhile, $\mathscr G_{\varpi^{\lambda}}$, $\mathscr G_{\varpi^{\lambda_1}}$ may be identified with the values of $\nu + \vartheta$ at $\lambda$, respectively $\lambda_1$ by Proposition \ref{prop-sharp-torus-commutative-gerbe-identification}, which coincide because $\lambda$, $\lambda_1$ map to the same element $\theta \in \hat Z_H$ and $\nu + \vartheta$ factors through $\hat Z_H$ (\emph{cf.}~\S\ref{void-whittaker-torsor}).

The composition of these identifications
\begin{align}
\notag
(\nu + \vartheta)(\theta)|_{\Gr_G^{\lambda_1} \cap S^{-, \lambda}} \simeq \mathscr G_{\varpi^{\lambda_1}}|_{\Gr_G^{\lambda_1} \cap S^{-, \lambda}} &\simeq \mathscr G_{\Gr_G}|_{\Gr_G^{\lambda_1} \cap S^{-, \lambda}} \\
\label{eq-discrepancy-torsor-definition}
& \simeq \mathscr G_{\varpi^{\lambda}} |_{\Gr_G^{\lambda_1} \cap S^{-, \lambda}} \simeq (\nu + \vartheta)(\theta)|_{\Gr_G^{\lambda_1} \cap S^{-, \lambda}}
\end{align}
thus defines an $A$-torsor $\tau^{\lambda_1, \lambda}$ over $\Gr_G^{\lambda_1} \cap S^{-, \lambda}$. We refer to $\tau^{\lambda_1, \lambda}$ as the \emph{discrepancy $A$-torsor} over the Mirkovi\'c--Vilonen cycle $\Gr_G^{\lambda_1} \cap S^{-, \lambda}$.

The following lemma gives us some qualitative control of $\tau^{\lambda_1, \lambda}$ in special cases. It is a variant of the calculation performed in \cite[\S4.4]{MR2684259}.
\end{void}

\begin{lem}
\label{lem-discrepancy-torsor-constancy}
Let $\lambda_1 \in \Lambda^{\sharp, +}$ and $\lambda \in \Lambda^{\sharp}$ be two elements with equal image in $\hat Z_H$. Then $\tau^{\lambda_1, \lambda}$ descends to $X$ under either of the following assumptions:
\begin{enumerate}
	\item $\lambda$ belongs to the Weyl-orbit of $\lambda_1$;
	\item $G$ has semisimple rank $1$.
\end{enumerate}
\end{lem}

\begin{proof}
Under (1), the morphism $\Gr_G^{\lambda_1} \cap S^{-, \lambda} \rightarrow X$ is a fibration in affine spaces. Indeed, it is identified with the $L^+N^-$-orbit of $\varpi^{\lambda}$ (\emph{cf.}~\cite[proof of Theorem 3.2]{MR2342692}). This implies that any $A$-torsor descends to $X$.

Let us now assume (2). For $\Gr_G^{\lambda_1} \cap S^{-, \lambda}$ to be nonempty and $\lambda_1$, $\lambda$ to have the same image in $\hat Z_H$, the cocharacter $\lambda$ must belong to the string
$$
\{\lambda_1, \lambda_1 - \alpha^{\sharp}, \cdots, s_{\check{\alpha}}(\lambda_1)\},
$$
where $\alpha$ is the unique simple coroot of $G$. We shall further assume $\lambda \neq \lambda_1, s_{\check{\alpha}}(\lambda_1)$, or else we fall back to case (1).

Choose an isomorphism between $\mathbb G_a$ and the root subgroup $N_{-\check{\alpha}}$ of $G$. We filter $L\mathbb G_a$ by the group subschemes $L^{\ge n}\mathbb G_a := \varpi^n L^+\mathbb G_a$ for $n \in \integers$. Then the $L\mathbb G_a$-action on $\varpi^{\lambda}$ has stabilizer $L^{\ge - \langle\check{\alpha}, \lambda\rangle}\mathbb G_a$. Furthermore, $u \in L\mathbb G_a$ satisfies $u \varpi^{\lambda} \in \Gr_G^{\lambda_1}$ if and only if
\begin{equation}
\label{eq-semiinfinite-orbit-inequality}
u \in L^{\ge - \langle\check{\alpha}, \lambda\rangle - d}\mathbb G_a \setminus L^{\ge - \langle\check{\alpha}, \lambda\rangle - d + 1}\mathbb G_a,
\end{equation}
where $0 < d < \langle\check{\alpha}, \lambda_1\rangle$ is the integer defined by $\lambda_1 - \lambda = d\cdot\alpha$. (Thus, $d$ coincides with the dimension of the Mirkovi\'c--Vilonen cycle $\Gr_G^{\lambda_1} \cap S^{-, \lambda}$.)

Projecting $u\varpi^{\lambda}$ onto the leading coefficient of $u$ defines a morphism
\begin{equation}
\label{eq-mirkovic-vilonen-cycle-morphism-to-canonical-bundle}
	\Gr_G^{\lambda_1} \cap S^{-, \lambda} \rightarrow \Omega_X^{-\langle\check{\alpha}, \lambda\rangle - d} \setminus 0,
\end{equation}
where $\Omega_X^{-\langle\check{\alpha}, \lambda\rangle - d}$ denotes the total space of the corresponding power of the canonical line bundle $\Omega_X$, and $0$ denotes its zero section. The morphism \eqref{eq-mirkovic-vilonen-cycle-morphism-to-canonical-bundle} intertwines the $L_+T_{\adjoint}$-action on the source with the scaling $\mathbb G_m$-action on the target, with respect to
\begin{equation}
\label{eq-adjoint-torus-simple-root-quotient}
L^+T_{\adjoint} \simeq L^+\mathbb G_m \twoheadrightarrow \mathbb G_m,
\end{equation}
where the isomorphism is supplied by $\check{\alpha}$ and the second morphism is the projection.

We shall now endow $\mathscr G_{\Gr_G^{[\theta]}}$ with the $L^+G_{\adjoint}$-equivariance determined by the image $\theta \in \hat Z_H$ of both $\lambda_1$ and $\lambda$ (\emph{cf.}~Proposition \ref{prop-gerbe-central-equivariance}, Remark \ref{rem-gerbe-central-equivariance-as-adjoint-equivariance}). This in particular restricts to an $L^+T_{\adjoint}$-equivariance on $\mathscr G_{\Gr_G^{[\theta]}}$ such that all isomorphisms in \eqref{eq-discrepancy-torsor-definition} are $L^+T_{\adjoint}$-equivariant (\emph{cf.}~Lemma \ref{lem-central-equivariance-levi-compatibility}, Lemma \ref{lem-central-equivariance-schubert-cell-compatibility}). Therefore, $\tau^{\lambda_1, \lambda}$ is an $L^+T_{\adjoint}$-equivariant $A$-torsor over $\Gr_G^{\lambda_1} \cap S^{-, \lambda}$. Since \eqref{eq-mirkovic-vilonen-cycle-morphism-to-canonical-bundle} is a fibration in affine spaces and \eqref{eq-adjoint-torus-simple-root-quotient} has a pro-unipotent kernel, all such $A$-torsors canonically descend to
$$
(\Omega_X^{-\langle\check{\alpha}, \lambda\rangle - d} \setminus 0) / \mathbb G_m \simeq X.
$$

This concludes the proof.
\end{proof}

\begin{rem}
Lemma \ref{lem-discrepancy-torsor-constancy} may suggest that $\tau^{\lambda_1, \lambda}$ descends to $X$ in general. However, this is \emph{not} true.

More precisely, when $G$ has semisimple rank $\ge 2$, it can happen that the restrictions of $\tau^{\lambda_1, \lambda}$ to \emph{certain} irreducible components of $\Gr_G^{\lambda_1} \cap S^{-, \lambda}$ descend to $X$, while its restrictions to other irreducible components do not. We refer the reader to \cite[\S7]{MR4303733}, where a similar phenomenon for line bundles over Zastava spaces is studied.
\end{rem}

\begin{void}[Trivialization of $\tau^{\alpha^{\sharp}, 0}$]
Let us now assume that $G$ has a unique simple root $\check{\alpha}$.

According to Lemma \ref{lem-discrepancy-torsor-constancy}, the discrepancy $A$-torsor $\tau^{\alpha^{\sharp}, 0}$ canonically descends to $X$. We shall strengthen this result by constructing a \emph{canonical} trivialization of $\tau^{\alpha^{\sharp}, 0}$. It makes essential use of the trivalization of $\nu(\alpha^{\sharp})$ defined in \S\ref{void-whittaker-torsor}.
\end{void}

\begin{lem}
\label{lem-whittaker-torsor-canonical-trivialization}
If $G$ has a unique simple root $\check{\alpha}$, then $\tau^{\alpha^{\sharp}, 0}$ is canonically trivial.
\end{lem}

\begin{proof}
It suffices to show that \eqref{eq-discrepancy-torsor-definition} is the identity automorphism for $\lambda_1 = \alpha^{\sharp}$ and $\lambda = 0$. Equivalently, we need to prove that the following diagram commutes
\begin{equation}
\label{eq-whittaker-torsor-canonical-trivialization-circuit}
\begin{tikzcd}[column sep = 1em]
	\mathscr G_{\varpi^{\alpha^{\sharp}}}|_{\Gr_G^{\alpha^{\sharp}} \cap S^{-, 0}} \ar[r, phantom, "\simeq"]\ar[d, "\simeq"] & \mathscr G_{\Gr_G} |_{\Gr_G^{\alpha^{\sharp}} \cap S^{-, 0}} \ar[r, phantom, "\simeq"] & \mathscr G_{\varpi^0}|_{\Gr_G^{\alpha^{\sharp}} \cap S^{-, 0}} \ar[d, "\simeq"] \\
	(\nu + \vartheta)(\alpha^{\sharp})|_{\Gr_G^{\alpha^{\sharp}} \cap S^{-, 0}} \ar[rr, phantom, "\simeq"] & & (\nu + \vartheta)(0)|_{\Gr_G^{\alpha^{\sharp}} \cap S^{-, 0}}
\end{tikzcd}
\end{equation}
where the vertical isomorphisms are given by Proposition \ref{prop-sharp-torus-commutative-gerbe-identification}, and the lower horizontal isomorphism is induced from the canonical trivialization of $\nu(\alpha^{\sharp})$ (\emph{cf.}~\S\ref{void-whittaker-torsor}). The statement immediately reduces to the case where $G$ is simply connected.

By the criterion \eqref{eq-semiinfinite-orbit-inequality} for the containment $u\varpi^0 \in \Gr_G^{\alpha^{\sharp}}$, we have
$$
\Gr_G^{\alpha^{\sharp}} \cap S^{-, 0} \subset K_{d \omega} \cdot \varpi^0
$$
where for any $\lambda \in \Lambda_{\adjoint}$, we write $K_{\lambda}$ for the subgroup scheme $\varpi^{\lambda} L^+G \varpi^{-\lambda}$ of $LG$, and here $\omega$ is the fundamental coweight of $\Lambda_{\adjoint}$ and $d := \ord Q(\alpha)$. Furthermore, the orbit $K_{d \omega} \cdot \varpi^0$ contains both $\varpi^0$ and $\varpi^{\alpha^{\sharp}}$. It thus suffices to prove the following statements
\begin{enumerate}
	\item the restriction of $\mathscr G_{\Gr_G}$ to $K_{d\omega} \cdot \varpi^0$ descends to $X$;
	\item the isomorphism $\mathscr G_{\varpi^{\alpha^{\sharp}}} \simeq \mathscr G_{\varpi^0}$ arising from (1) fits into the commutative square
	\begin{equation}
	\label{eq-whittaker-torsor-canonical-trivialization-skew-orbit}
	\begin{tikzcd}[column sep = 1.5em]
		\mathscr G_{\varpi^{\alpha^{\sharp}}} \ar[r, phantom, "\simeq"]\ar[d, "\simeq"] & \mathscr G_{\varpi^0} \ar[d, "\simeq"] \\
		(\nu + \vartheta)(\alpha^{\sharp}) \ar[r, phantom, "\simeq"] & (\nu + \vartheta)(0)
	\end{tikzcd}
	\end{equation}
	where the vertical isomorphisms are again given by Proposition \ref{prop-sharp-torus-commutative-gerbe-identification} and the lower horizontal isomorphism is induced from the canonical trivialization of $\nu(\alpha^{\sharp})$.
\end{enumerate}

We shall prove both statements using the embedding $G \subset \widetilde G := G \rtimes T_{\adjoint}$ and the canonical extension of $\mu$ to an \'etale level $\widetilde{\mu}$ of $\widetilde G$ (\emph{cf.}~the proof of Lemma \ref{lem-gerbe-dimension-one-bruhat-variety}). The latter allows us to extend $\mathscr G_{\Gr_G}$ to an $A$-gerbe $\mathscr G_{\Gr_{\widetilde G}}$. Identifying the maximal torus of $\widetilde G$ with $T \times T_{\adjoint}$ and its adjoint group with $G_{\adjoint}$, we may also consider the subgroup schemes
$$
\widetilde K_{\lambda} := \varpi^{\lambda} L^+\widetilde G \varpi^{-\lambda} \simeq \varpi^{(0, \lambda)} L^+\widetilde G \varpi^{-(0, \lambda)}
$$
of $L\widetilde G$ for each $\lambda \in \Lambda_{\adjoint}$. The inclusion $\Gr_G \rightarrow \Gr_{\widetilde G}$ carries $K_{d\omega} \cdot \varpi^0$ into $\widetilde K_{d\omega}\cdot \varpi^0$, which is in turn the $\varpi^{(0, d\omega)}$-translate of the Schubert cell in $\Gr_{\widetilde G}$ containing $\varpi^{(0, -d\omega)}$. Since the symmetric form \eqref{eq-simply-connected-etale-level-extension-bilinear-form} of $\widetilde{\mu}$ evaluates to the trivial character at $(0, -d\omega)$, Corollary \ref{cor-gerbe-schubert-cell-sharp-identification} (applied to $\widetilde G$) implies that the restriction of $\mathscr G_{\Gr_{\widetilde G}}$ to this Schubert cell descends to $X$. Furthermore, the induced isomorphism
\begin{equation}
\label{eq-extended-level-schubert-cell-weyl-equivariance}
\mathscr G_{\varpi^{s_{\check{\alpha}}(0, -d\omega)}} \simeq \mathscr G_{\varpi^{(0, -d\omega)}},
\end{equation}
for $s_{\check{\alpha}}$ the simple reflection, coincides with the isomorphism arising from Weyl-equivariance (\emph{cf.}~Proposition \ref{prop-gerbe-schubert-cell-weyl-equivariance}).

Therefore, the restriction of $\mathscr G_{\Gr_G}$ to $K_{d\omega} \cdot \varpi^0$ descends to $X$, and the induced isomorphism $\mathscr G_{\varpi^{\alpha^{\sharp}}} \simeq \mathscr G_{\varpi^0}$ is the difference between \eqref{eq-extended-level-schubert-cell-weyl-equivariance} and the identity automorphism of $\mathscr G_{\varpi^{(0, -d\omega)}}$:
$$
\mathscr G_{\varpi^{s_{\check{\alpha}}(0, -d\omega)}} - \mathscr G_{\varpi^{(0, -d\omega)}} \simeq \mathscr G_{\varpi^{(0, -d\omega)}} - \mathscr G_{\varpi^{(0, -d\omega)}}.
$$
This isomorphism renders \eqref{eq-whittaker-torsor-canonical-trivialization-skew-orbit} commutative, by comparing with the construction of the trivialization of $\nu(\alpha^{\sharp})$ (\emph{cf.}~\S\ref{void-whittaker-torsor}).
\end{proof}

\begin{void}[Application to weight functors]
We shall use discrepancy $A$-torsors to give an explicit description of the image of standard objects under weight functors (\emph{cf.}~\S\ref{void-weight-functors}).

More precisely, consider the standard functor $\Delta^{\lambda_1}$ associated to $\lambda_1 \in \Lambda^{\sharp, +}$ (\emph{cf.}~\S\ref{void-standard-costandard}). Upon untwisting the Satake category, we may regard it as a functor
\begin{equation}
\label{eq-untwisted-standard-functor}
	\Delta^{\lambda_1} : \SLis(X) \rightarrow {}^+\SSat^{\theta}_{\mathscr G, \zeta}(\Hec_G)_{-(\nu + \vartheta)(\lambda_1)},
\end{equation}
where $\theta \in \hat Z_H$ is the image of $\lambda_1$. Here, we invoked the fact that the image of $\Delta^{\lambda_1}$ is contained in the $\theta$-graded component of the Satake category (\emph{cf.}~the proof of Proposition \ref{prop-virtual-connected-components}).

Consider now the composition of \eqref{eq-untwisted-standard-functor} and \eqref{eq-weight-functor-satake-category-homogeneous} for any $\lambda \in \Lambda^{\sharp}$ with image $\theta$ in $\hat Z_H$:
\begin{equation}
\label{eq-weight-of-standard-object}
	\omega^{\lambda} \circ \Delta^{\lambda_1} : \SLis(X) \rightarrow \SLis(X),
\end{equation}
where we identify $(\nu + \vartheta)(\lambda_1)$ with $(\nu + \vartheta)(\lambda)$. On the other hand, we have the discrepancy $A$-torsor $\tau^{\lambda_1, \lambda}$ over the Mirkovi\'c--Vilonen cycle $\Gr_G^{\lambda_1} \cap S^{-, \lambda}$. Along the character $\zeta : A \subset \coeff^{\times}$, it defines a rank-$1$ $\coeff$-local system $\tau_{\zeta}^{\lambda_1, \lambda}$.
\end{void}

\begin{lem}
\label{lem-weight-functor-discrepancy-torsor-cohomology}
The functor \eqref{eq-weight-of-standard-object} is identified with tensor product with
\begin{equation}
\label{eq-discrepancy-torsor-cohomology}
H^{2d} (p^-)_! \tau_{\zeta}^{\lambda_1, \lambda}(d),
\end{equation}
where $p^- : \Gr_G^{\lambda_1} \cap S^{-, \lambda} \rightarrow X$ is the structural map and $d$ is its relative dimension.
\end{lem}

\begin{proof}
We define the weight functor $\omega^{\lambda}$ using the constant term functor $\CT_{B^-}$ associated to $B^- \subset G$ (\emph{cf.}~Lemma \ref{lem-constant-term-independence-of-borel}). Since $\CT_{B^-}$ is perverse $t$-exact (\emph{cf.}~Proposition \ref{prop-constant-term-reflective-properties}), we have an isomorphism for any $(\mathscr G_{\varpi^{\lambda_1}}, \zeta)$-twisted local system $\mathscr E$:
\begin{equation}
\label{eq-constant-term-standard-top-cohomology}
\CT_{B^-}(\Delta^{\lambda_1}(\mathscr E)) \simeq {}^pH^0(\CT_{B^-}((j^{\lambda_1})_!\widetilde{\mathscr E})),
\end{equation}
where the notation $j^{\lambda_1}$, $\widetilde{\mathscr E}$ are as in \S\ref{void-standard-costandard}. The right-hand-side coincides with tensor product with \eqref{eq-discrepancy-torsor-cohomology}, by base change and the projection formula.
\end{proof}

\begin{prop}
\label{prop-weight-functor-nondegeneracy}
Let $\lambda_1 \in \Lambda^{\sharp, +}$ and $\lambda \in \Lambda^{\sharp}$ be two elements with equal image in $\hat Z_H$. The functor $\omega^{\lambda} \circ \Delta^{\lambda_1}$ is an equivalence under either of the following assumptions:
\begin{enumerate}
	\item $\lambda$ belongs to the Weyl-orbit of $\lambda_1$;
	\item $G$ has semisimple rank $1$.
\end{enumerate}
\end{prop}

\begin{proof}
This follows from Lemma \ref{lem-discrepancy-torsor-constancy}, Lemma \ref{lem-weight-functor-discrepancy-torsor-cohomology}, together with the fact that the structural map $p^- : \Gr_G^{\lambda_1} \cap S^{-, \lambda} \rightarrow X$ has irreducible fibers under assumption (1) or (2).
\end{proof}

\begin{prop}
\label{prop-metaplectic-simple-root-weight-zero-trivialization}
If $G$ has a unique simple root $\check{\alpha}$, there is a canonical isomorphism
\begin{equation}
\label{eq-metaplectic-simple-root-weight-zero-trivialization}
\omega^0 \circ \Delta^{\alpha^{\sharp}} \simeq \id_{\SLis(X)}.
\end{equation}
\end{prop}

\begin{proof}
This follows from Lemma \ref{lem-whittaker-torsor-canonical-trivialization} and Lemma \ref{lem-weight-functor-discrepancy-torsor-cohomology}.
\end{proof}

\begin{rem}
\label{rem-metaplectic-simple-root-weight-zero-trivialization-canonical}
Proposition \ref{prop-metaplectic-simple-root-weight-zero-trivialization} is stated with a given Borel subgroup $B\subset G$ and a section of $B \twoheadrightarrow T$. It is, however, indepedent of these choices.

First, let us argue that the left-hand-side of \eqref{eq-metaplectic-simple-root-weight-zero-trivialization} is independent of these choices. We have seen that this is the case for the functor $\omega^0$ (\emph{cf.}~\S\ref{void-constant-term-functor-untwisted}). As for $\Delta^{\alpha^{\sharp}}$, note that the Schubert cell $\Gr_G^{\alpha^{\sharp}} \subset \Gr_G$ is independent of $T \subset B \subset G$, as different choices are related by $G$-conjugation. It remains to show that the identification of $A$-gerbes
$$
\mathscr G_{\Gr_G^{\alpha^{\sharp}}} \simeq (\nu + \vartheta)(\alpha^{\sharp})|_{\Gr_G^{\alpha^{\sharp}}}
$$
is independent of $T \subset B \subset G$. This follows from comparing the $L^+G$-equivariance on $\mathscr G_{\Gr_G^{\alpha^{\sharp}}}$ with the canonical $G$-equivariance on $\mu$ (\emph{cf.}~the proof of Proposition \ref{prop-gerbe-schubert-cell-weyl-equivariance}).

Next, let us show that the isomorphism \eqref{eq-metaplectic-simple-root-weight-zero-trivialization} is independent of $T \subset B \subset G$. Unwinding the definitions, this amounts to showing that the trivialization of the descent of $\tau^{\alpha^{\sharp}, 0}$ to $X$, constructed in Lemma \ref{lem-whittaker-torsor-canonical-trivialization}, is independent of these choices. For this, we may perform the same construction using the universal Borel and maximal torus over $\Torel$ and verify that the resulting trivialization of $\tau^{\alpha^{\sharp}, 0}$ descends along $\Torel \rightarrow X$ (\emph{cf.}~Remark \ref{rem-commutative-cover-independence-of-borel}). However, since the geometric fibers of $\Torel$ are simply connected, pullback along $\Torel \rightarrow X$ induces an equivalence on the categories of $A$-torsors.
\end{rem}

\subsection{The isomorphism $\check G \simeq H$}

\begin{void}
We remain in the setting where $I = \{1\}$ and omit it from the notation. Recall the \'etale sheaf $\check G$ of affine monoid $\coeff$-schemes over $X$.

The goal of this subsection is to construct an isomorphism
\begin{equation}
\label{eq-tannaka-group-identification}
	\check G \simeq H.
\end{equation}
This isomorphism, combined with Corollary \ref{cor-tannaka-dual-equivalence}, will allow us to complete the construction of the Satake equivalence (\emph{cf.}~\S\ref{void-satake-equivalence-construction}).
\end{void}

\begin{void}
Let us first note a ``pointwise" version of \eqref{eq-tannaka-dual-equivalence}. Indeed, given a geometric point $\bar x$ of $X$, the stalk of $\SSat_{\mathscr G, \zeta}(\Hec_G)$ at $\bar x$ is canonically identified with $\Sat_{\mathscr G, \zeta}(\Hec_{G, \bar x})$ via the pullback functor (\emph{cf.}~Proposition \ref{prop-local-hecke-stack-derived-category-stalk}).

From this fact and \eqref{eq-tannaka-dual-equivalence}, we deduce the equivalence
\begin{equation}
\label{eq-tannaka-dual-equivalence-pointwise}
	{}^+\Sat_{\mathscr G, \zeta}(\Hec_{G, \bar x})_{-(\nu + \vartheta)} \simeq \Rep_{\check G_{\bar x}},
\end{equation}
where $\check G_{\bar x}$ is the fiber of $\check G$ at $\bar x$ and the symmetric monoidal structure on the left-hand-side is induced from that of ${}^+\SSat_{\mathscr G, \zeta}(\Hec_G)_{-(\nu + \vartheta)}$.

We use \eqref{eq-tannaka-dual-equivalence-pointwise} to translate properties of the pointwise Satake category to those of $\check G$.
\end{void}

\begin{lem}
\label{lem-tannaka-dual-reductive}
Sections of $\check G$ are reductive group $\coeff$-schemes.
\end{lem}

\begin{proof}
It suffices to prove that the geometric fibers of $\check G$ are reductive group $\coeff$-schemes.

Let $\bar x$ be a geometric point of $X$. The Satake category ${}^+\Sat_{\mathscr G, \zeta}(\Hec_{G, \bar x})$ is rigid (\emph{cf.}~Proposition \ref{prop-convolution-monoidal-dual}), showing that $\check G_{\bar x}$ is grouplike.

To prove that $\check G_{\bar x}$ is of finite type, we need to find a finite set of objects of $\Rep_{\check G_{\bar x}}$ which generate it under sums, monoidal products, and subquotients (\emph{cf.}~\cite[Proposition 2.20(b)]{Deligne1982}). For this, we may choose a set $\lambda_1, \cdots, \lambda_n$ of generators of $\Lambda^{\sharp, +}$ as a monoid. The images of $\IC^{\lambda_1}, \cdots, \IC^{\lambda_n}$ generate ${}^+\Sat_{\mathscr G, \zeta}(\Hec_{G, \bar x})$ and are homogeneous under the $\hat Z_H$-grading (\emph{cf.}~the proof of Proposition \ref{prop-virtual-connected-components}). They thus define generators of $\Rep_{\check G_{\bar x}}$ upon untwisting.

Finally, ${}^+\Sat_{\mathscr G, \zeta}(\Hec_{G, \bar x})$ is semisimple (\emph{cf.}~Corollary \ref{cor-semisimplicity}), showing that $\check G_{\bar x}$ is linearly reductive. To prove that $\check G_{\bar x}$ is reductive, it remains to show that it is connected, but this holds by the criterion of \cite[Corollary 2.22]{Deligne1982}.
\end{proof}

\begin{void}[Maximal torus]
\label{void-tannaka-dual-maximal-torus}
By construction of \eqref{eq-tannaka-dual-equivalence}, the forgetful functor $\Rep_{\check G} \rightarrow \Lis(X)$ factors through the forgetful functor $\Rep_{T_H}$ (\emph{cf.}~\S\ref{void-fiber-functor-definition}).

Taking Tannaka duals, we obtain a morphism
\begin{equation}
\label{eq-tannaka-group-maximal-torus}
T_H \rightarrow \check G
\end{equation}
of locally constant \'etale sheaves of group $\coeff$-schemes.
\end{void}

\begin{lem}
\label{lem-tannaka-group-maximal-torus}
Sections of \eqref{eq-tannaka-group-maximal-torus} are maximal tori.
\end{lem}

\begin{proof}
It again suffices to prove this assertion for geometric fibers.

Fix a geometric point $\bar x$ of $X$. By Proposition \ref{prop-weight-functor-nondegeneracy}(1), to each $\lambda \in \Lambda^{\sharp}$, viewed as a character of $T_{H, \bar x}$, there exists an object $V \in \Rep_{\check G_{\bar x}}$ whose image in $\Rep_{T_{H, \bar x}}$ has a $1$-dimensional $\lambda$-weight space. Namely, this $V$ may be chosen as the image of $\coeff$ under $\Delta^{\lambda_1}$, where we take $\lambda_1$ to be the dominant representative of $\lambda$. It follows that $T_{H, \bar x} \rightarrow \check G_{\bar x}$ is a closed immersion.

To prove that $T_{H, \bar x}$ is a maximal torus of $G_{\bar x}$, one may now proceed as in the classical argument (\emph{cf.}~\cite[\S9.1]{MR3839695}).
\end{proof}

\begin{void}[Borel subgroup]
\label{void-tannaka-dual-borel-subgroup}
Recall that the dual group $H$ of $(G, \mu)$, being pinned, has a canonical Borel subgroup $B_H$ containing the maximal torus $T_H$ (\emph{cf.}~\S\ref{sec-metaplectic-dual-data}).

Let $\check B \subset \check G$ be the subsheaf of Borel subgroups containing $T_H$ which renders the character $2\check{\rho}$ dominant. It is uniquely characterized by the property that $T_H$ has the same set of dominant weights with respect to $\check B \subset \check G$ and $B_H \subset H$ (\emph{cf.}~\cite[Lemma 9.5]{MR3839695}).
\end{void}

\begin{prop}
\label{prop-tannaka-group-root-data}
The root data of $T_H \subset \check B \subset \check G$ and $T_H \subset B_H \subset H$ coincide.
\end{prop}

\begin{proof}
We again pass to fibers at a geometric point $\bar x$ of $X$ and view $\Lambda^{\sharp}$ as the character lattice of $T_{H, \bar x}$. Furthermore, we fix a maximal torus and a Borel subgroup $T \subset B \subset G$.

Let us first assume that $G$ has a unique simple root $\check{\alpha}$. In this case, $H$ admits a unique simple root $\alpha^{\sharp} = \ord Q(\alpha) \cdot \alpha$. By Corollary \ref{cor-virtual-connected-components-constant-term-compatibility} and Proposition \ref{prop-weight-functor-nondegeneracy}(2) (and the agreement of standard and intersection cohomology functors, \emph{cf.}~Corollary \ref{cor-standard-ic-agreement}), the weights occurring in a simple $\check G_{\bar x}$-representation of highest weight $\lambda_1$ are precisely elements of the set
$$
\{\lambda_1, \lambda_1 - \alpha^{\sharp}, \cdots, s_{\check{\alpha}}(\lambda_1)\}.
$$
This shows that $\alpha^{\sharp}$ is the unique simple root of $\check G_{\bar x}$. Since $\check B_{\bar x}$ and $B_{H, \bar x}$ determine the same subset of dominant elements of $\Lambda^{\sharp}$, the coroot of $\check G_{\bar x}$ associated to $\alpha^{\sharp}$ must be proportional to $\check{\alpha}$, as they annihilate the same hyperplane in $\Lambda^{\sharp} \otimes \rationals$. The proportion is $\ord(Q(\alpha))^{-1}$, because $\langle\check{\alpha}^{\sharp}, \alpha^{\sharp}\rangle = 2$.

The general case reduces to the above case by compatibility with constant term functors. More precisely, each simple root $\check{\alpha}$ of $G$ determines a subminimal parabolic $P$ with Levi quotient $P \twoheadrightarrow M$. By construction of \eqref{eq-tannaka-dual-equivalence} and the composition property of constant term functors (\emph{cf.}~Remark \ref{rem-constant-term-composition}), the forgetful functor $\Rep_{\check G} \rightarrow \Lis(X)$ factors through $\Rep_{\check M}$, where $\check M$ is the Tannaka dual of the untwisted Satake category for $M$, equipped with the induced \'etale level. This yields a morphism
\begin{equation}
\label{eq-tannaka-group-levi-subgroup}
\check M \rightarrow \check G
\end{equation}
of locally constant \'etale sheaves of group $\coeff$-schemes.

Note that \eqref{eq-tannaka-group-levi-subgroup} is a closed immersion: This follows from the same argument as in the proof of Lemma \ref{lem-tannaka-group-maximal-torus} and boils down to the fact that the $L^+_{\bar x}N_P$-orbit of $\varpi^{\lambda}$ in $\Gr_{G, \bar x}$ (for $N_P$ the nilradical of $P$) is an affine space.

Thus, $\alpha^{\sharp}$ is a root of $\check G$, so the set of roots of $H$ is contained in the set of roots of $\check G$, preserving positivity. This containment is an equality because the weights of $\Delta^{\lambda_1}(\coeff) \in \Rep_{\check G_{\bar x}}$, for any $\lambda_1 \in \Lambda^{\sharp, +}$, belong to $\lambda_1 - \integers^{\ge 0} \Delta^{\sharp}$ (\emph{cf.}~Corollary \ref{cor-virtual-connected-components-constant-term-compatibility}). The identification of coroots follows, as \eqref{eq-tannaka-group-levi-subgroup} carries root subgroups of $\check M_{\bar x}$ to those of $\check G_{\bar x}$.
\end{proof}

\begin{void}[Pinning]
\label{void-tannaka-group-pinning}
Next, we shall construct a pinning on $\check G$ under the assumption that $G$ is split. Namely, we shall construct a canonical isomorphism
\begin{equation}
\label{eq-tannaka-group-pinning}
\mathbb G_a \simeq N_{\check{\alpha}^{\sharp}},
\end{equation}
where $N_{\check{\alpha}^{\sharp}}$ is the root subgroup of $\check G$ corresponding to $\check{\alpha}^{\sharp}$.

Indeed, consider the object $\Delta^{\alpha^{\sharp}}(\coeff) \in \Rep_{\check G}$. By Proposition \ref{prop-metaplectic-simple-root-weight-zero-trivialization} and Remark \ref{rem-metaplectic-simple-root-weight-zero-trivialization-canonical}, the weight space $\omega^0 \circ \Delta^{\alpha^{\sharp}}(\coeff)$ is canonically isomorphic to the constant local system $\coeff$. Since $\omega^{\alpha^{\sharp}} \circ \Delta^{\alpha^{\sharp}}(\coeff)$ is tautologically isomorphic to $\coeff$, we have distinguished sections
\begin{align*}
	v_{\alpha^{\sharp}} : \coeff & \simeq \omega^{\alpha^{\sharp}} \circ \Delta^{\alpha^{\sharp}}(\coeff), \\
	v_0 : \coeff & \simeq \omega^0 \circ \Delta^{\alpha^{\sharp}}(\coeff).
\end{align*}

Under the $\check G$-action on the underlying local system of $\Delta^{\alpha^{\sharp}}(\coeff)$, the subgroup $N_{\check{\alpha}^{\sharp}}$ acts on the subspace spanned by $\{v_{\alpha^{\sharp}}, v_0\}$ in a strictly upper triangular manner. We obtain \eqref{eq-tannaka-group-pinning} as the unique isomorphism under which $a \in \mathbb G_a$ acts by $a\cdot v_0 = v_0 - 2 a v_{\alpha}$.
\end{void}

\begin{rem}
Our definition of the pinning on $\check G$ ensures that when $G$ has semisimple rank $1$, under the identification $(\check G)_{\adjoint} \simeq \PGL_2$ induced from this pinning, the object $\Delta^{\alpha^{\sharp}}(\coeff) \in \Rep_{\check G}$ is canonically the pullback of the adjoint representation of $\PGL_2$ along $\check G \rightarrow (\check G)_{\adjoint}$.
\end{rem}

\begin{void}[Construction of \eqref{eq-tannaka-group-identification}]
When $G$ is split, we have endowed $\check G$ with the canonical structure of a constant sheaf of pinned split reductive group $\coeff$-schemes whose root data coincide with those of $H$ (\emph{cf.}~Proposition \ref{prop-tannaka-group-root-data}). We obtain \eqref{eq-tannaka-group-identification} as the unique pinning-preserving isomorphism.

For general $G$, we choose an \'etale cover of $X$ splitting $G$, over which we have a pinning-preserving isomorphism $\check G \simeq H$. To see that this isomorphism intertwines the descent data of $\check G$ and $H$, we may pass to the corresponding isomorphism of their root data, which is induced from the identity automorphism of $T_H$ (as a maximal torus of both $\check G$ and $H$). We thus obtain the desired isomorphism \eqref{eq-tannaka-group-identification} by \'etale descent.
\end{void}

\begin{void}
\label{void-satake-equivalence-construction}
Finally, we complete the construction of the geometric Satake equivalence.

\begin{proof}[Proof of Theorem \ref{thm-satake-equivalence}]
Combining \eqref{eq-tannaka-dual-equivalence} and \eqref{eq-tannaka-group-identification}, we obtain an equivalence
\begin{equation}
\label{eq-satake-equivalence-local-twisted}
{}^+\SSat_{\mathscr G, \zeta}(\Hec_{G, I})_{-(\nu + \vartheta)^{\boxplus I}} \simeq \SRep_{H^{\boxtimes I}}
\end{equation}
of \'etale sheaves of $(\hat Z_H)^{\oplus I}$-graded symmetric monoidal $\coeff$-linear categories over $X^I$ for any finite set $I$. Twisting both sides of \eqref{eq-satake-equivalence-local-twisted} by $(\nu + \vartheta)^{\boxplus I}$ (\emph{cf.}~\S\ref{void-symmetric-monoidal-twist}), we obtain an equivalence of \'etale sheaves of symmetric monoidal $\coeff$-linear categories
\begin{equation}
\label{eq-satake-equivalence-local}
{}^+\SSat_{\mathscr G, \zeta}(\Hec_{G, I}) \simeq \SRep_{H^{\boxtimes I}, (\nu + \vartheta)^{\boxplus I}}.
\end{equation}
The desired equivalence \eqref{eq-satake-equivalence} arises by taking sections of \eqref{eq-satake-equivalence-local} over $X^I$.
\end{proof}
\end{void}

\newpage

\part{Global function fields}

\section{Preparation}
\label{sec-global-function-fields-preparation}

In this section, we begin our treatment of global function fields. We let $\base$ be a field and $X$ be a smooth, proper, geometrically connected curve over $\base$. Let $D \subset X$ be a $\base$-finite closed subscheme (the ``ramification divisor") and write $\mathring X := X\setminus D$.

Let $A$ be a finite abelian group whose order is invertible in $\base$. Let $G$ be a smooth affine group $X$-scheme equipped with an $A$-valued \'etale level \emph{over $\mathring X$}, \emph{i.e.}~a morphism of pointed \'etale stacks $\mu : \deloop_{\mathring X} G \rightarrow \deloop^4_{\mathring X} A(1)$.

We begin by defining the ``global $A$-gerbe" in \S\ref{sec-global-gerbe}. Then we explain a trace-of-Frobenius construction for gerbes (\emph{cf.}~\S\ref{sec-trace-of-frobenius}) and use it to define $\zeta$-genuine automorphic forms (\emph{cf.}~\S\ref{sec-automorphic-forms}). These constructions are mild generalizations of standard ones in the literature (\emph{cf.}~\cite{MR2265675}), where we treat ramification and general $A$-gerbes. In \S\ref{sec-meta-galois-twist}, we prove a novel result, Theorem \ref{thm-meta-galois-twist}, which gives a geometric interpretation of Weissman's meta-Galois group in the function field setting.

\subsection{The global $A$-gerbe}
\label{sec-global-gerbe}

\begin{void}
Denote by $\Bun_{G, D}$ the moduli stack of $G$-bundles over $X$ rigidified along $D$, \emph{i.e.}~the prestack assigning to a $\base$-algebra $R$ the groupoid of pairs $(P, \phi)$, where $P$ is a $G$-bundles over $X_R := \Spec R\times X$ and $\phi$ is a trivialization of the restriction of $P$ to $D_R$.

By \cite[Proposition 1]{MR2640041}, $\Bun_{G, \emptyset}$ is representable by a smooth algebraic $\base$-stack. The same then holds for $\Bun_{G, D}$, since it is a torsor over $\Bun_{G, \emptyset}$ with structure group the Weil restriction of $G$ along $D \rightarrow \Spec \base$.
\end{void}

\begin{void}
By convention, we shall write $\Ran$ for the Ran space \emph{of $\mathring X$} (\emph{cf.}~\S\ref{void-ran-space-definition}).

Denote by $\Hec(\Bun_{G, D})$ the global Hecke stack, \emph{i.e.}~an $R$-point of $\Hec(\Bun_{G, D})$ consists of an $R$-point $\underline x$ of $\Ran$, a pair of $R$-points $(P^0, \phi^0)$, $(P^1, \phi^1)$ of $\Bun_{G, D}$ together with an isomorphism over $X_R\setminus \Gamma_{\underline x}$. This makes sense because $D_R$ is contained in $X_R \setminus \Gamma_{\underline x}$.

We shall represent an $R$-point of $\Hec(\Bun_{G, D})$ as:
$$
(P^0, \phi^0) \overset{\underline x}{\sim} (P^1, \phi^1)
$$
and refer to it as a \emph{modification}, in parallel with \S\ref{void-local-hecke-stack-definition}.
\end{void}

\begin{void}
\label{void-global-iterated-hecke-stack}
More generally, we may define an ``outer convolution diagram" version of the global Hecke stack as follows.

Let $n\ge 0$ be an integer. Denote by $\Hec^{[n]}(\Bun_{G, D})$ the moduli stack parametrizing $n$ points $\underline x{}_1, \cdots, \underline x{}_n$ of $\Ran$, together with a chain of modifications
\begin{equation}
\label{eq-global-iterated-hecke-stack-modification}
(P^0, \phi^0) \overset{\underline x{}_1}{\sim} (P^1, \phi^1) \overset{\underline x{}_2}{\sim} \cdots \overset{\underline x{}_n}{\sim} (P^n, \phi^n).
\end{equation}

By construction, we have structural morphisms
\begin{equation}
\label{eq-global-iterated-hecke-stack-structural-morphisms}
\begin{tikzcd}[column sep = 1em]
	\Hec^{[n]}(\Bun_{G, D}) \ar[r, "\prod r_j"]\ar[d, "{(p_0, p_n)}"] & \prod_{j = 1}^n \Hec_G \\
	\Bun_{G, D} \times \Bun_{G, D}
\end{tikzcd}
\end{equation}
where $p_0$, $p_n$ send \eqref{eq-global-iterated-hecke-stack-modification} to $(P^0, \phi^0)$, respectively $(P^n, \phi^n)$, and each $r_j$ ($j = 1, \cdots, n$) sends it to the restriction of the modification $P^{j-1} \overset{\underline x{}_j}{\sim} P^j$ to $D_{\underline x{}_j}$.
\end{void}

\begin{void}
\label{void-global-integration-map}
In parallel with \S\ref{void-local-hecke-stack-gerbe}, we shall construct a morphism of spaces
\begin{equation}
\label{eq-global-integration-map}
\int_X : \Maps_*(\deloop_{\mathring X} G, \deloop_{\mathring X}^4A(1)) \rightarrow \Maps(\Bun_{G, D}, \deloop^2 A).
\end{equation}

Then we will set $\mathscr G_{\Bun_{G, D}}$ to be the image of $\mu$ under \eqref{eq-global-integration-map}. We shall refer to $\mathscr G_{\Bun_{G, D}}$ as the \emph{global $A$-gerbe} associated to the \'etale level $\mu$.
\end{void}

\begin{rem}
Note that $A$-gerbes over $\Bun_{G, D}$ are canonically equivalent to those over $\Bun_{G, D_{\red}}$, where $D_{\red} \subset D$ is the reduced subscheme.

One may view this fact as saying that the our covers are ``tamely ramified". We shall see another manifestation of this fact in \S\ref{void-covering-group-construction}.
\end{rem}

\begin{void}[The global trace map]
Given a $\base$-algebra $R$, we shall construct a morphism in the pro-category of $H\integers$-module spectra
\begin{equation}
\label{eq-global-trace-map}
\tr_X : \Gamma(X_R, \hat{\integers}(1)[2]) \rightarrow \Gamma(\Spec R, \hat{\integers}).
\end{equation}

Indeed, write $\pi : X_R \rightarrow \Spec R$ for the projection morphism. Since $\pi$ is smooth of relative dimension $1$, we may identify $\hat{\integers}(1)[2]$ with $\pi^!\hat{\integers}$. Since $\pi$ is also proper, we may use the co-unit of the adjunction to form the composition
\begin{align*}
\Gamma(X_R, \hat{\integers}(1)[2]) & \simeq \Gamma(X_R, \pi^!\hat{\integers}) \\
& \simeq \Gamma(\Spec R, \pi_*\pi^!\hat{\integers}) \rightarrow \Gamma(\Spec R, \hat{\integers}),
\end{align*}
which is our definition of \eqref{eq-global-trace-map}.
\end{void}

\begin{rem}
In parallel with Remark \ref{rem-trace-map-projection-formula}, we comment on the behavior of \eqref{eq-global-trace-map} with respect to change of coefficients: Given any complex $\mathscr A$ of torsion \'etale sheaves of invertible order, the same construction yields a morphism
\begin{equation}
\label{eq-global-trace-map-with-coefficients}
\tr_X : \Gamma(X_R, \mathscr A(1)[2]) \rightarrow \Gamma(\Spec R, \mathscr A),
\end{equation}
which is related to \eqref{eq-global-trace-map} by the commutative diagram
$$
\begin{tikzcd}[column sep = 1em]
	\Gamma(\Spec R, \mathscr A) \otimes \Gamma(X_R, \hat{\integers}(1)[2]) \ar[r, "\otimes"]\ar[d, "\id\otimes\tr_X"] & \Gamma(X_R, \mathscr A(1)[2]) \ar[d, "\tr_X"] \\
	\Gamma(\Spec R, \mathscr A) \otimes \Gamma(\Spec R, \hat{\integers}) \ar[r, "\otimes"] & \Gamma(\Spec R, \mathscr A)
\end{tikzcd}
$$
\end{rem}

\begin{void}[Construction of \eqref{eq-global-integration-map}]
\label{void-global-integration-map-construction}
Using the cohomological interpretation of $\deloop^n A$ (\emph{cf.}~\S\ref{void-bar-construction}), it suffices to construct a morphism of $H\integers$-module spectra
\begin{equation}
\label{eq-global-integration-map-linear}
\Gamma(\deloop_{\mathring X} G \text{ mod }\mathring X, A(1)[4]) \rightarrow \Gamma(\Bun_{G, D}, A[2]).
\end{equation}

Since the closed immersion $D \subset X$ is cohomologically pure, the associated Cousin sequence yields a fiber sequence of $H\integers$-module spectra
\begin{equation}
\label{eq-ramification-divisor-cousin-sequence}
	\Gamma(\deloop_D G \text{ mod }D, A[2]) \rightarrow \Gamma(\deloop_X G \text{ mod }X, A(1)[4]) \rightarrow \Gamma(\deloop_{\mathring X} G \text{ mod }\mathring X, A(1)[4]).
\end{equation}

To construct \eqref{eq-global-integration-map-linear}, it thus suffices to construct a morphism
\begin{equation}
\label{eq-global-integration-map-linear-unramified}
\Gamma(\deloop_X G, A(1)[4]) \rightarrow \Gamma(\Bun_{G, D}, A[2])
\end{equation}
and trivialize its restriction to $\Gamma(\deloop_D G\text{ mod }D, A[2])$.

Let $P$ denote the universal $G$-bundle over $X\times \Bun_{G, D}$, which we view as a morphism $X \times \Bun_{G, D} \rightarrow \deloop_X G$. We define \eqref{eq-global-integration-map-linear-unramified} as the composition $\tr_X \circ P^*$, where $\tr_X$ is the global trace map \eqref{eq-global-trace-map-with-coefficients} (for $\mathscr A = A[2]$). To trivialize $\tr_X \circ P^*$ over $\Gamma(\deloop_D G \text{ mod }D, A[2])$, we use the commutative square of $H\integers$-module spectra
$$
\begin{tikzcd}[column sep = 1em]
	\Gamma(\deloop_D G, A[2]) \ar[d, "(P|_D)^*"]\ar[r] & \Gamma(\deloop_XG, A(1)[4]) \ar[d, "P^*"] \\
	\Gamma(D\times \Bun_{G, D}, A[2]) \ar[r] & \Gamma(X \times \Bun_{G, D}, A(1)[4])
\end{tikzcd}
$$
where the horizontal arrows are given by cohomological purity. Note that $(P|_D)^*$ is trivial over $\Gamma(\deloop_DG \text{ mod }D, A[2])$ since $P|_D$ factors through $D$ via $\phi$, as desired.

This concludes the construction of \eqref{eq-global-integration-map}, hence the global $A$-gerbe $\mathscr G_{\Bun_{G, D}}$.
\end{void}

\begin{rem}
For $D = \emptyset$, the morphism \eqref{eq-global-integration-map} is simply the morphism of spaces induced from \eqref{eq-global-integration-map-linear-unramified}. In particular, it factors through $\Maps(\deloop_X G, \deloop^4_XA(1))$, \emph{i.e.}~the rigidification of $\mu$ at the neutral point is irrelevant for the construction.

On the other hand, the construction of $\mathscr G_{\Bun_{G, D}}$ shows that it is canonically trivial along the neutral $\base$-point of $\Bun_{G, D}$, \emph{i.e.}~the pair $(P, \phi)$ where $P$ is the trivial $G$-bundle over $X$ and $\phi$ is the identity.
\end{rem}

\begin{void}
Next, we relate $\mathscr G_{\Bun_{G, D}}$ to the local $A$-gerbe $\mathscr G_{\Hec_G}$ (\emph{cf.}~\S\ref{void-local-hecke-stack-gerbe}) via the global Hecke action diagram \eqref{eq-global-iterated-hecke-stack-structural-morphisms} for a fixed integer $n\ge 0$.
\end{void}

\begin{lem}
\label{lem-global-gerbe-local-compatibility}
There is a canonical isomorphism of $A$-gerbes over $\Hec^{[n]}(\Bun_{G, D})$:
\begin{equation}
\label{eq-global-gerbe-local-compatibility}
(p_n)^*\mathscr G_{\Bun_{G, D}} - (p_0)^*\mathscr G_{\Bun_{G, D}} \simeq \sum_{j = 1}^n (r_j)^*\mathscr G_{\Hec_G}.
\end{equation}
\end{lem}

\begin{proof}
For $n = 0$, \eqref{eq-global-gerbe-local-compatibility} is the tautological trivialization of $(p_0)^*\mathscr G_{\Bun_{G, D}} - (p_0)^*\mathscr G_{\Bun_{G, D}}$.

For $n \ge 1$, we may reduce to the case $n = 1$ by rewriting the left-hand-side of \eqref{eq-global-gerbe-local-compatibility} as a telescopic sum. Then we need to prove the following: Given $G$-bundles $P^0$, $P^1$ over $X_R$ with an isomorphism over $X_R \setminus \Gamma_{\underline x}$, where $\underline x$ is an $R$-point of $\Ran$, the image of
\begin{align*}
\mu(P^1) - \mu(P^0) &\in \Gamma(\mathring X_R \text{ mod }\mathring X_R \setminus \Gamma_{\underline x}, \deloop^4A(1)) \\
& \simeq \Gamma(X_R \text{ mod } X_R \setminus \Gamma_{\underline x}, \deloop^4A(1))
\end{align*}
under the global trace map \eqref{eq-global-trace-map-with-coefficients} coincides with the image of
$$
\mu(P^1 |_{D_{\underline x}}) - \mu(P^0 |_{D_{\underline x}}) \in \Gamma(D_{\underline x}\text{ mod }\mathring D_{\underline x}, \deloop^4 A(1))
$$
under the local trace map \eqref{eq-trace-map-with-coefficients}.

This follows by identifying the morphism of $H\integers$-module spectra
\begin{align}
\notag
	\Gamma(X_R \text{ mod }X_R\setminus \Gamma_{\underline x}, A(1)[4]) & \xrightarrow{\simeq} \Gamma(D_{\underline x}\text{ mod }\mathring D_{\underline x}, A(1)[4]) \\
\label{eq-restriction-followed-by-local-trace}
	& \xrightarrow{\tr_{\underline x}} \Gamma(\Spec R, A[2])
\end{align}
with the restriction of the global trace map, where the isomorphism is given by formal base change (\emph{cf.}~\S\ref{void-trace-map}). Note, however, that the composition \eqref{eq-restriction-followed-by-local-trace} is defined by adjunction for the proper morphism $x_R \rightarrow \Spec R$, so the desired identification follows from naturality with respect to the morphism $x_R \rightarrow X_R$ of proper $R$-schemes.
\end{proof}

\begin{void}[Variant for $\infty D$]
\label{void-infinite-level}
Let us note a variant of the construction of $\mathscr G_{\Bun_{G, D}}$ at ``infinite level", where it is enough to start with a group scheme over $\mathring X$.

More precisely, we let $G$ be a smooth affine group $\mathring X$-scheme. Denote by $\Bun_{G, \infty D}$ the moduli stack of pairs $(P, \phi)$, where $P$ is a $G$-bundle over $\mathring X_R := \mathring X \times \Spec R$ and $\phi$ is a trivialization of $P$ over $\hat D_R \setminus D_R$, for $\hat D_R$ the formal completion of $X_R$ along $D_R$. (Here, $R$ is a test $\base$-algebra.)

There is a morphism of spaces
\begin{equation}
\label{eq-infinite-level-global-trace-map}
	\int_X : \Maps(\deloop_{\mathring X}G, \deloop^4_{\mathring X} A(1)) \rightarrow \Maps(\Bun_{G, \infty D}, \deloop^2 A),
\end{equation}
constructed in analogy with \eqref{eq-global-integration-map}: We use the fact that any section of $\deloop^4A(1)$ over $\mathring X_R$ trivialized over $\hat D_R \setminus D_R$ admits a canonical extension to $X_R$ and then apply the global trace map \eqref{eq-global-trace-map-with-coefficients} (for $\mathscr A := A[2]$). Given an \'etale level $\mu$ of $G$, we thus obtain an $A$-gerbe $\mathscr G_{\Bun_{G, \infty D}}$ over $\Bun_{G, \infty D}$ as its image under \eqref{eq-infinite-level-global-trace-map}.
\end{void}

\begin{rem}
\label{rem-infinite-level-structural-morphism-to-finite-level}
Let $G$ be a smooth affine group $X$-scheme. By the Beauville--Laszlo lemma (\emph{cf.}~\cite{MR1320381}), an $R$-point of $\Bun_{G, \infty D}$ is equivalent to a $G$-bundle over $X_R$ trivialized over $\hat D_R$. In particular, we have a canonical morphism
$$
\Bun_{G, \infty D} \rightarrow \Bun_{G, D}
$$
under which $\mathscr G_{\Bun_{G, D}}$ pulls back to $\mathscr G_{\Bun_{G, \infty D}}$.
\end{rem}

\subsection{Interlude: trace of Frobenius}
\label{sec-trace-of-frobenius}

\begin{void}
We now suppose that $\base$ is a finite field of cardinality $q$. For a $\base$-prestack $\mathscr X$, we write $\Fr_{\mathscr X} : \mathscr X \rightarrow \mathscr X$ for the $q$th power absolute Frobenius endomorphism, \emph{i.e.}~it sends $x \in \mathscr X(R)$ to $x\circ \Fr_R$, where $\Fr_R$ is induced from the ring map $R\rightarrow R$, $a\mapsto a^q$.

We shall describe a ``trace-of-Frobenius" construction for $A$-gerbes, which takes as input a $\base$-prestack $\mathscr X$ equipped with an $A$-gerbe $\mathscr G$, and gives as output an $A$-torsor $\Tr(\Fr \mid\mathscr G)$ over the fixed-point prestack $\mathscr X^{\Fr}$, \emph{i.e.}~the fiber product
\begin{equation}
\label{eq-fixed-point-prestack}
\begin{tikzcd}[column sep = 2.5em]
	\mathscr X^{\Fr} \ar[r]\ar[d] & \mathscr X \ar[d, "\Delta"] \\
	\mathscr X \ar[r, "{(\id, \Fr_{\mathscr X})}"] & \mathscr X \times \mathscr X
\end{tikzcd}
\end{equation}
\end{void}

\begin{rem}
\label{rem-rational-points-embedding-in-frobenius-fixed-locus}
Let $\mathscr X$ be an algebraic $\base$-stack locally of finite type. Regarding the groupoid $\mathscr X(\base)$ as a discrete $\base$-stack, we obtain a natural morphism of $\base$-stacks
\begin{equation}
\label{eq-rational-points-embedding-in-frobenius-fixed-locus}
\mathscr X(\base) \rightarrow \mathscr X^{\Fr}.
\end{equation}

By \cite[Lemma 3.3]{MR2061225}, $\mathscr X^{\Fr}$ is an \'etale Deligne--Mumford $\base$-stack and the morphism \eqref{eq-rational-points-embedding-in-frobenius-fixed-locus} is an open and closed immersion; it is an isomorphism if $\Delta : \mathscr X \rightarrow \mathscr X \times \mathscr X$ has connected geometric fibers. Furthermore, by the proof of \cite[Proposition 2.16(c)]{MR2061225} (\emph{cf.}~\cite[\S12.3.2]{MR3787407}), \eqref{eq-rational-points-embedding-in-frobenius-fixed-locus} is an isomorphism for $\mathscr X := \Bun_{G, D}$, where $G$ is a smooth group $X$-scheme with connected geometric fibers, which is reductive over $\mathring X$.
\end{rem}

\begin{void}
Recall that for any $\base$-scheme $S$, the endofunctor $\Fr_S^*$ on the \'etale site of $S$ is naturally isomorphic to the identity (\emph{cf.}~\cite[03SN]{stacks-project}).

Let us explicitly describe the value of this natural isomorphism at an $A$-gerbe $\mathscr G$:
\begin{equation}
\label{eq-baffling-isomorphism}
\Fr_S^*\mathscr G \simeq \mathscr G.
\end{equation}
For any \'etale morphism $f : S_1 \rightarrow S$, the groupoid $\Fr_S^*(\mathscr G)(S_1)$ is the filtered colimit of $\mathscr G(U)$ over \'etale morphisms $u : U \rightarrow S$ through which $\Fr_S\circ f$ factors. This index category has an initial object, namely $(U, u) = (S_1, f)$ with the factorization $\Fr_S \circ f = f\circ \Fr_{S_1}$. The colimit is thus identified with $\mathscr G(S_1)$.

Since the isomorphism \eqref{eq-baffling-isomorphism} is natural in $S$, we obtain a natural isomorphism
\begin{equation}
\label{eq-baffling-isomorphism-prestack}
	\Fr_{\mathscr X}^*\mathscr G \simeq \mathscr G
\end{equation}
for any $\base$-prestack $\mathscr X$.
\end{void}

\begin{void}[Construction of $\Tr(\Fr\mid\mathscr G)$]
\label{void-trace-of-frobenius-gerbe-construction}
Let $\mathscr X$ be a $\base$-prestack equipped with an $A$-gerbe $\mathscr G$. Since $\Fr_{\mathscr X}$ restricts to the identity map on $\mathscr X^{\Fr}$, we obtain an isomorphism
\begin{equation}
\label{eq-tautological-isomorphism-fixed-point}
\Fr_{\mathscr X}^*(\mathscr G|_{\mathscr X^{\Fr}}) \simeq \mathscr G|_{\mathscr X^{\Fr}}.
\end{equation}

The $A$-gerbe $\Tr(\Fr\mid\mathscr G)$ is defined as the difference between the restriction of \eqref{eq-baffling-isomorphism-prestack} to $\mathscr X^{\Fr}$ and \eqref{eq-tautological-isomorphism-fixed-point}, \emph{i.e.}~it corresponds to the composition
\begin{equation}
\label{eq-gerbe-trace-of-frobenius-definition}
\mathscr G|_{\mathscr X^{\Fr}} \xrightarrow{\eqref{eq-tautological-isomorphism-fixed-point}} \Fr_{\mathscr X}^*(\mathscr G|_{\mathscr X^{\Fr}}) \simeq (\Fr_{\mathscr X}^*\mathscr G)|_{\mathscr X^{\Fr}} \xrightarrow{\eqref{eq-baffling-isomorphism}} \mathscr G|_{\mathscr X^{\Fr}}.
\end{equation}
\end{void}

\begin{void}
\label{void-trace-of-frobenius-genuine-functions}
We now fix an algebraic $\base$-stack $\mathscr X$ locally of finite type such that $\mathscr X(\base) \rightarrow \mathscr X^{\Fr}$ is an isomorphism (\emph{cf.}~Remark \ref{rem-rational-points-embedding-in-frobenius-fixed-locus}). Let $\mathscr G$ be an $A$-gerbe over $\mathscr X$.

Taking $\base$-points of $\Tr(\Fr\mid\mathscr G)$, we obtain a morphism of groupoids
\begin{equation}
\label{eq-trace-of-frobenius-set-theoretic-torsor}
\widetilde{\mathscr X} := \Tr(\Fr\mid\mathscr G)(\base) \rightarrow \mathscr X^{\Fr}(\base) \simeq \mathscr X(\base).
\end{equation}
Note that \eqref{eq-trace-of-frobenius-set-theoretic-torsor} is surjective: By functoriality of the construction, it suffices to show that for any $\base$-point $x$ of $\mathscr X$, the $A$-torsor $\Tr(\Fr\mid\mathscr G|_x)$ over $x$ is trivial, but this holds because $\mathscr G|_x$ is trivial, as $H^2_{\etale}(\Spec\base, A) \simeq 0$. It follows that \eqref{eq-trace-of-frobenius-set-theoretic-torsor} is a set-theoretic $A$-torsor.
\end{void}

\begin{rem}
\label{rem-trace-of-frobenius-via-fundamental-group}
Let us give a more concrete description of \eqref{eq-trace-of-frobenius-set-theoretic-torsor}. Given $x\in\mathscr X(\base)$, we choose an algebraic closure $\bar{\base}$ of $\base$ and a trivialization $\bar g$ of $\mathscr G$ over $\bar x := \Spec \bar{\base}$. The resulting fundamental group of $\mathscr G|_x$ fits into a short exact sequence (\emph{cf.}~\cite[Theorem 19.6]{MR3802418})
\begin{equation}
\label{eq-trace-of-frobenius-via-fundamental-group}
1 \rightarrow A \rightarrow \pi_1^{\etale}(\mathscr G|_x, \bar g) \rightarrow \pi_1^{\etale}(x, \bar x) \rightarrow 1.
\end{equation}

The fiber of \eqref{eq-trace-of-frobenius-set-theoretic-torsor} at $x$, as an $A$-torsor, is identified with the preimage of the \emph{geometric} Frobenius element under the surjection of \eqref{eq-trace-of-frobenius-via-fundamental-group}. To see this, it suffices to observe that the composition \eqref{eq-gerbe-trace-of-frobenius-definition} (for $\mathscr X := \Spec\base$) is \emph{inverse} to the Frobenius-pullback on $\bar{\base}$-points.
\end{rem}

\begin{rem}
\label{rem-gerbe-norm}
Let $\base_1$ be a finite extension of $\base$ of cardinality $q_1$ and $\mathscr X_1$ be a $\base_1$-prestack. Write $\mathscr X$ for the Weil restriction of $\mathscr X_1$ along $\Spec \base_1 \rightarrow \Spec\base$. There is a \emph{norm} map:
\begin{equation}
\label{eq-gerbe-norm-map}
	\Nm : \Maps(\mathscr X_1, \deloop^2 A) \rightarrow \Maps(\mathscr X, \deloop^2 A).
\end{equation}

Indeed, by adjunction, it suffices to construct a natural map of $H\integers$-module spectra
\begin{equation}
\label{eq-gerbe-norm-map-adjunction}
\Gamma(\Spec R\otimes_{\base} \base_1, A) \rightarrow \Gamma(\Spec R, A)
\end{equation}
for any $\base$-algebra $R$. This is given by the trace map on \'etale cochains, using the fact that $\Spec\base_1 \rightarrow \Spec\base$ is finite \'etale.

The norm map \eqref{eq-gerbe-norm-map} is compatible with the trace-of-Frobenius construction in the following sense: For $(\mathscr X_1, \mathscr G_1)$ as in \S\ref{void-trace-of-frobenius-genuine-functions}, the $A$-torsor $\Tr(\Fr \mid \mathscr G)(\base_1)$ over $\mathscr X_1(\base_1)$, formed with respect to the $q_1$th power Frobenius, coincides with $\Tr(\Fr \mid \Nm \mathscr G)$ over $\mathscr X(\base)$, formed with respect to the $q$th power Frobenius, under the natural identification $\mathscr X(\base_1) \simeq \mathscr X(\base)$.
\end{rem}

\begin{void}
\label{void-constructing-genuine-functions}
Given the coefficient data of \S\ref{void-satake-category-coefficients}, we write $\Fun_{\zeta}(\widetilde{\mathscr X}, \coeff)$ for the $\coeff$-vector space of \emph{$\zeta$-genuine} functions on $\widetilde{\mathscr X}$, \emph{i.e.}~functions $f : \widetilde{\mathscr X} \rightarrow \coeff$ satisfying
$$
f(x \cdot a) = f(x) \cdot \zeta(a)
$$
for each $x \in \widetilde{\mathscr X}$ and $a\in A$. We denote by $\Fun_{c, \zeta}(\widetilde{\mathscr X}, \coeff)$ the subspace of compactly supported $\zeta$-genuine functions on $\widetilde{\mathscr X}$.

Denote by $\mathscr L_{\Tr(\Fr \mid \mathscr G), \zeta}$ the $1$-dimensional $\coeff$-local system over $\mathscr X^{\Fr}$ induced from $\Tr(\Fr \mid \mathscr G)$ along $\zeta$. Since $\mathscr X^{\Fr}$ is identified with $\mathscr X(\base)$, we have a canonical isomorphism
\begin{equation}
\label{eq-genuine-functions-as-cohomology}
	\Fun_{\zeta}(\widetilde{\mathscr X}, \coeff) \simeq H^0(\mathscr X^{\Fr}, (\mathscr L_{\Tr(\Fr \mid \mathscr G), \zeta})^{\otimes -1}).
\end{equation}
Here, the inverse appears because sections of $(\mathscr L_{\tau, \zeta})^{\otimes -1}$, for a trivial $A$-torsor $\tau$ over $\Spec \base$, are canonically identified with $\zeta$-genuine functions $\tau(\base) \rightarrow \coeff$.

Moreover, any $\mathscr A \in \derived_{\mathscr G, \zeta}(\mathscr X)$ (\emph{cf.}~\S\ref{void-twisted-sheaves}) induces a $\zeta$-genuine function $\Tr(\Fr \mid \mathscr A)$ on $\widetilde{\mathscr X}$ as follows: For each $x \in \mathscr X(\base)$, a choice of an algebraic closure $\bar{\base}$ of $\base$ and a trivialization $\bar g$ of $\mathscr G$ over $\bar x := \Spec\bar{\base}$ defines a $\zeta$-genuine representation of $\pi_1^{\etale}(\mathscr G|_x, \bar g)$ on the constructible complex of $\coeff$-vector spaces $\mathscr A_{\bar x}$, whose trace defines $\Tr(\Fr \mid \mathscr A)$ over $x$ (\emph{cf.}~Remark \ref{rem-trace-of-frobenius-via-fundamental-group}).
\end{void}

\subsection{Automorphic forms}
\label{sec-automorphic-forms}

\begin{void}
\label{void-automorphic-form-context}
We continue to assume $\base$ to be a finite field of cardinality $q$ and invoke the coefficient data of \S\ref{void-satake-category-coefficients}. Instead of $\coeff$, we shall work with an algebraic closure $\overline{\rationals}_{\ell}$ of $\rationals_{\ell}$ and define ``$\zeta$-genuine automorphic forms" using the trace-of-Frobenius construction of $A$-gerbes.

Assume that $G \rightarrow X$ has connected geometric fibers and is reductive over $\mathring X$. Then $(\Bun_{G, D})^{\Fr}$ coincides with the discrete stack $\Bun_{G, D}(\base)$ (\emph{cf.}~Remark \ref{rem-rational-points-embedding-in-frobenius-fixed-locus}). Applying the trace-of-Frobenius construction to the global $A$-gerbe $\mathscr G_{\Bun_{G, D}}$, we obtain an $A$-torsor
\begin{equation}
\label{eq-global-set-theoretic-torsor}
\widetilde{\Bun}_{G, D} := \Tr(\Fr \mid \mathscr G_{\Bun_{G, D}})(\base)
\end{equation}
over $\Bun_{G, D}(\base)$. We shall refer to $\zeta$-genuine functions over $\widetilde{\Bun}_{G, D}$ as \emph{$\zeta$-genuine automorphic forms} (with ramification divisor $D$). They form the $\overline{\rationals}_{\ell}$-vector space
\begin{equation}
\label{eq-genuine-automorphic-forms}
	\Fun_{\zeta}(\widetilde{\Bun}_{G, D}, \overline{\rationals}_{\ell}).
\end{equation}
\end{void}

\begin{rem}
\label{rem-infinite-level-set-theoretic-torsor}
Let us note a variant of $\widetilde{\Bun}_{G, D}$ at ``infinite level", for any smooth affine group $\mathring X$-scheme $G$ (\emph{cf.}~\S\ref{void-infinite-level}).

Namely, restricting $\Tr(\Fr \mid \mathscr G_{\Bun_{G, \infty D}})(\base)$ along $\Bun_{G, \infty D}(\base) \subset (\Bun_{G, \infty D})^{\Fr}(\base)$, we obtain a set-theoretic $\base$-torsor $\widetilde{\Bun}_{G, \infty D}$ over $\Bun_{G, \infty D}(\base)$.
\end{rem}

\begin{void}
The $A$-torsor \eqref{eq-global-set-theoretic-torsor} may be expressed in an ad\`elic form, as follows.

Denote by $F$ the field of fractions of $X$. For each closed point $x\in X$, we write $F_x$ for the local field and $\mathscr O_x$ for its ring of integers. Denote by $\mathbb A_F$ (respectively, $\mathbb O_F$) the topological ring of (respectively, integral) ad\`eles of $F$. By Weil uniformization, we have a fully faithful embedding of groupoids
\begin{equation}
\label{eq-weil-uniformization-map}
	G(F) \backslash G(\mathbb A_F) / K_D \subset \Bun_{G, D}(\base),
\end{equation}
where $K_D \subset G(\mathbb O_F)$ is the kernel of $G(\mathbb O_F) \rightarrow G(\mathscr O_D)$.

Let us describe the pullback of $\widetilde{\Bun}_{G, D}$ along \eqref{eq-weil-uniformization-map} in terms of a cover of $G(\mathbb A_F)$ equipped with canonical sections over $G(F)$ and over $K_D$. To do so, we need a mild generalization of the construction of $\mathscr G_{\Hec_G}$ (\emph{cf.}~\S\ref{void-local-hecke-stack-gerbe}). Indeed, we shall construct a morphism
\begin{equation}
\label{eq-local-multiplicative-gerbe-rational-form}
	\int_{\mathring D_x} : \Maps_*(\deloop_{F_x} G, \deloop^4_{F_x} A(1)) \rightarrow \Maps_{\mathbb E_1}(L_x G, \deloop^2 A),
\end{equation}
where $\deloop_{F_x}$ denotes the deloop functor on the big \'etale site of $\Spec F_x$, and $L_x G$ is the loop group associated to the $\base_x$-point $x : \Spec \base_x \rightarrow X$ (\emph{cf.}~\S\ref{void-local-hecke-stack-definition}).
\end{void}

\begin{void}[Construction of \eqref{eq-local-multiplicative-gerbe-rational-form}]
For any $\base_x$-algebra $R$, we shall construct a morphism
\begin{equation}
\label{eq-local-multiplicative-gerbe-rational-form-value}
\Maps_*(\deloop_{F_x} G, \deloop^4_{F_x} A(1)) \rightarrow \Maps_{\mathbb E_1}(G(\mathring D_{x_R}), \Gamma(\Spec R, \deloop^2 A))
\end{equation}
natural in $R$. Here, $x_R$ denotes the $R$-point of $X$ defined by $x$.

By rewriting the left-hand-side of \eqref{eq-local-multiplicative-gerbe-rational-form-value} as $\Maps_{\mathbb E_1}(G_{F_x}, \deloop^3_{F_x} A(1))$, we shall obtain \eqref{eq-local-multiplicative-gerbe-rational-form-value} as the composition of the evaluation over $\mathring D_{x_R}$, followed by the morphism of $\mathbb E_1$-monoids induced from some morphism of $H\integers$-module spectra
\begin{equation}
\label{eq-local-trace-parametrized-rational-point}
\Gamma(\mathring D_{x_R}, A(1)[3]) \rightarrow \Gamma(\Spec R, A[2]).
\end{equation}

The morphism \eqref{eq-local-trace-parametrized-rational-point} is the composition of the second map in the fiber sequence
\begin{equation}
\label{eq-parametrized-rational-point-cousin-sequence}
\Gamma(D_{x_R}, A(1)[3]) \rightarrow \Gamma(\mathring D_{x_R}, A(1)[3]) \rightarrow \Gamma(D_{x_R} \text{ mod }\mathring D_{x_R}, A(1)[4])
\end{equation}
with the local trace map \eqref{eq-trace-map-with-coefficients} (for $\mathscr A := A[2]$).
\end{void}

\begin{rem}
\label{rem-rational-point-local-gerbe-unramified}
Given a pointed morphism $\deloop_{F_x} G \rightarrow \deloop^4_{F_x} A(1)$ equipped with an extension to $\Spec \mathscr O_x$, its image under \eqref{eq-local-multiplicative-gerbe-rational-form} admits a canonical trivialization over $L_x^+G$. Indeed, this follows from the fiber sequence \eqref{eq-parametrized-rational-point-cousin-sequence}.

In this situation, we obtain an $A$-gerbe over $\Hec_{G, x} \simeq L_x^+G \backslash L_x G/ L_x^+G$. This coincides with the (pullback of the) $A$-gerbe $\mathscr G_{\Hec_G}$ constructed in \S\ref{void-local-hecke-stack-gerbe}.
\end{rem}

\begin{void}
\label{void-covering-group-construction}
For each closed point $x \in X$, we may restrict the \'etale level $\mu$ to $\Spec F_x$ and apply \eqref{eq-local-multiplicative-gerbe-rational-form} to obtain a monoidal $A$-gerbe $\mathscr G_{L_x G}$ over $L_x G$.

Taking trace-of-Frobenius of $\mathscr G_{L_x G}$ (with respect to the $q_x$th power Frobenius, for $q_x$ the cardinality of $\base_x$), we obtain a central extension
\begin{equation}
\label{eq-local-covering-group}
1 \rightarrow A \rightarrow \widetilde G_x \rightarrow G(F_x) \rightarrow 1
\end{equation}
canonically split over $G(X\setminus x)$ by the sum-of-residue formula.

Since the principal congruence subgroup $L_x^{\ge 1} G \subset L_x^+G$ is pro-unipotent, any monoidal $A$-gerbe over $L_x^{\ge 1}G$ is canonically trivial. This implies that \eqref{eq-local-covering-group} admits a canonical section over $\ker(G(\mathscr O_x) \rightarrow G(\base_x))$. On the other hand, if $x \in \mathring X$, then this section extends to $G(\mathscr O_x)$ by Remark \ref{rem-rational-point-local-gerbe-unramified}.

Repeating the construction of \cite[\S2.2]{zhao2022metaplectic}, we obtain from \eqref{eq-local-covering-group} a central extension of the ad\`elic group
\begin{equation}
\label{eq-global-covering-group}
1 \rightarrow A \rightarrow \widetilde G \rightarrow G(\mathbb A_F) \rightarrow 1
\end{equation}
equipped with canonical sections over $G(F)$ and over $K_D$. In particular, it gives rise to an $A$-torsor $G(F) \backslash \widetilde G / K_D$ over $G(F) \backslash G(\mathbb A_F) / K_D$.
\end{void}

\begin{rem}
\label{rem-covering-group-tate-duality}
Note that \cite[\S2]{zhao2022metaplectic} supplies an \emph{a priori} different construction of the local covering group \eqref{eq-local-covering-group}, and consequently another global covering group \eqref{eq-global-covering-group}. However, these constructions coincide for an appropriate normalization of the Tate duality isomorphism $H^2_{\etale}(\Spec F_x, A(1)) \simeq A$. Namely, we exhibit it as the composition
\begin{align*}
H^2_{\etale}(\Spec F_x, A(1)) &\simeq H^1_{\etale}(\Spec\base_x, A) \\
& \simeq \Hom(\Gal(\bar{\base}_x / \base_x), A) \simeq A,
\end{align*}
where the last map is evaluation on the \emph{geometric} Frobenius.
\end{rem}

\begin{lem}
\label{lem-global-gerbe-adelic-description}
The $A$-torsor $G(F)\backslash \widetilde G / K_D$ over $G(F) \backslash G(\mathbb A_F) / K_D$ is canonically identified with the restriction of $\widetilde{\Bun}_{G, D}$ along \eqref{eq-weil-uniformization-map}.
\end{lem}

\begin{proof}
For each $x\in X$, consider the Weil restriction $\res_{\base}(L_x G)$ of the loop group $L_x G$ along $\Spec \base_x \rightarrow \Spec\base$. Beauville--Laszlo gluing yields a map
\begin{equation}
\label{eq-loop-group-weil-restriction-to-moduli-of-bundles}
\res_{\base}(L_x G) \rightarrow \Bun_{G, D},
\end{equation}
which recovers the map $G(F_x) \rightarrow \Bun_{G, D}(\base)$ in Weil uniformization upon taking $\base$-points.

By Remark \ref{rem-gerbe-norm}, it suffices to identify the pullback of $\mathscr G_{\Bun_{G, D}}$ along \eqref{eq-loop-group-weil-restriction-to-moduli-of-bundles} with the norm of the $A$-gerbe $\mathscr G_{L_x G}$, compatibly with the canonical sections. Unwinding the definitions, this amounts to identifying the composition
\begin{align*}
	\Gamma(X_R \text{ mod }X_R \setminus x_R, A(1)[4]) & \xrightarrow{\simeq} \Gamma(D_{x_R} \text{ mod }\mathring D_{x_R}, A(1)[4]) \\
	& \xrightarrow{\tr_{x_R}} \Gamma(\Spec R\otimes_{\base} \base_x, A[2]) \xrightarrow{\eqref{eq-gerbe-norm-map-adjunction}} \Gamma(\Spec R, A[2])
\end{align*}
with the restriction of the global trace map \eqref{eq-global-trace-map-with-coefficients} (for $\mathscr A = A[2]$), where $x_R$ denotes the $(R\otimes_{\base} \base_x)$-point of $X$ defined by $x$. As in the proof of Lemma \ref{lem-global-gerbe-local-compatibility}, this follows from naturality with respect to the morphism of proper $R$-schemes $x_R \rightarrow X_R$. We omit verifying the compatibilities with canonical sections.
\end{proof}

\begin{rem}
\label{rem-shafarevich-set}
The essential image of the Weil uniformization map \eqref{eq-weil-uniformization-map} consists of pairs $(P, \phi)$, where $P$ is trivial over the generic point $\Spec F$ of $X$. To describe the remaining points of $\Bun_{G, D}(\base)$, we need to consider inner forms of $G_F$.

Define the set
$$
\ker^1(F, G) := \ker(H^1_{\etale}(\Spec F, G) \rightarrow \prod_{x\in X} H^1_{\etale}(\Spec F_x, G)),
$$
where the product is taken over closed points of $X$. By \cite[Remark 12.2]{MR3787407}, $\ker^1(F, G)$ is finite and, after choosing a $G$-bundle $P_{\alpha}$ over $\Spec F$ trivialized over $\Spec F_x$ ($x\in X$) representing each $\alpha \in \ker^1(F, G)$, the collection of Weil uniformization maps for all $G_{\alpha} := \Aut(P_{\alpha})$ defines an equivalence of groupoids
\begin{equation}
\label{eq-weil-uniformization-with-inner-forms}
	\bigsqcup_{\alpha \in \ker^1(F, G)} G_{\alpha}(F) \backslash G_{\alpha}(\mathbb A_F) / K_D \simeq \Bun_{G, D}(\base).
\end{equation}
Here, we implicitly extended $P_{\alpha}$ (hence $G_{\alpha}$) to $X$ using its trivializations over $\Spec F_x$.

On the other hand, the \'etale level $\mu$ of $G$ induces an \'etale level $\mu_{\alpha}$ of $G_{\alpha}$, thanks to its $G$-equivariance structure (\emph{cf.}~\S\ref{void-conjugation-equivariance}). Repeating the construction of \eqref{eq-global-covering-group} with $(G_{\alpha}, \mu_{\alpha})$ instead of $(G, \mu)$, we obtain a central extension $\widetilde G_{\alpha}$ of $G_{\alpha}(\mathbb A_F)$ by $A$, and by Lemma \ref{lem-global-gerbe-adelic-description}, an identification of $A$-torsors over \eqref{eq-weil-uniformization-with-inner-forms}:
\begin{equation}
\label{eq-weil-uniformization-with-inner-forms-torsors}
	\bigsqcup_{\alpha \in \ker^1(F, G)} G_{\alpha}(F) \backslash \widetilde G_{\alpha} / K_D \simeq \widetilde{\Bun}_{G, D}.
\end{equation}
\end{rem}

\begin{void}[The ``sharp center"]
\label{void-sharp-center}
In order to obtain some finiteness, we need to restrict the action of the connected component $Z^{\circ}$ of the center $Z$ of $G$ over $\mathring X$. In the twisted setting, the role of the connected center is played by a torus isogenic to $Z^{\circ}$.

To wit, denote by $\Lambda_{Z^{\circ}}$ the sheaf of cocharacters of $Z^{\circ}$ and set
$$
\Lambda_{Z^{\sharp}} := \Lambda_{Z^{\circ}} \cap \Lambda^{\sharp} \subset \Lambda,
$$
where $\Lambda^{\sharp}$ is the kernel of the symmetric form $b$ (\emph{cf.}~\S\ref{void-metaplectic-dual-group}). Then the inclusion $\Lambda_{Z^{\sharp}} \subset \Lambda_{Z^{\circ}}$ corresponds to an isogeny of tori $Z^{\sharp} \rightarrow Z^{\circ}$. Since the \'etale level $\mu$ is $\mathbb E_{\infty}$-monoidal over $Z^{\sharp}$ (\emph{cf.}~\cite[Proposition 4.6.2]{zhao2022metaplectic}), the associated cover $\widetilde Z^{\sharp}$ of $Z^{\sharp}(\mathbb A_F)$ is commutative. Moreover, the ``infinite level" construction (\emph{cf.}~Remark \ref{rem-infinite-level-set-theoretic-torsor}) applies to $Z^{\sharp}$ and yields a symmetric monoidal extension
$$
A \rightarrow \widetilde{\Bun}_{Z^{\sharp}, \infty D} \rightarrow \Bun_{Z^{\sharp}, \infty D}(\base).
$$

There is a natural $\Bun_{Z^{\sharp}, \infty D}$-action on $\Bun_{G, D}$ defined by the natural $\deloop Z$-action and Beauville--Laszlo gluing. More precisely, an $R$-point $(P_{Z^{\sharp}}, \phi_{Z^{\sharp}})$ of $\Bun_{Z^{\sharp}, \infty D}$ carries $(P, \phi) \in \Bun_{G, D}$ to the pair $(P', \phi')$, where $P'$ is the gluing of $P \otimes P_{Z^{\sharp}}$ (defined over $X_R \setminus D_R$) with $P|_{\hat D_R}$ along the isomorphism
$$
\id_P \otimes \phi_{Z^{\sharp}} : P\otimes P_{Z^{\sharp}} |_{\hat D_R \setminus D_R} \simeq P |_{\hat D_R \setminus D_R},
$$
and $\phi'$ is the trivialization of $P'|_{D_R}$ corresponding to $\phi$.

Since $\mu$ is $\deloop Z^{\sharp}$-equivariant (\emph{cf.}~\S\ref{void-canonical-quadratic-structure}), the $\Bun_{Z^{\sharp}, \infty D}(\base)$-action on $\Bun_{G, D}$ lifts to an action of $\widetilde{\Bun}_{Z^{\sharp}, \infty D}$-action on $\widetilde{\Bun}_{G, D}$, compatibly with the natural $A$-action.
\end{void}

\begin{void}
\label{void-connected-center-lattice}
We shall fix a cocompact lattice $\Xi \subset Z^{\sharp}(F) \backslash \widetilde Z^{\sharp}$ which maps isomorphically onto its image in $\Bun_{Z^{\sharp}, \infty D}(\base)$.

Using the $\widetilde{\Bun}_{Z^{\sharp}, \infty D}$-action on $\widetilde{\Bun}_{G, D}$, we may consider $\Xi$-invariant $\zeta$-genuine automorphic forms, \emph{i.e.}~the $\overline{\rationals}_{\ell}$-vector space
\begin{equation}
\label{eq-genuine-forms-center-invariant}
\Fun_{\zeta}(\widetilde{\Bun}_{G, D}/\Xi, \overline{\rationals}_{\ell}).
\end{equation}

Next, we shall define the subspace of \eqref{eq-genuine-forms-center-invariant} consisting of ``cusp forms".
\end{void}

\begin{void}[Cusp forms]
Note that the central extension \eqref{eq-global-covering-group} splits canonically over $N(\mathbb A_F) \subset G(\mathbb A_F)$, where $N \subset G$ is any unipotent subgroup.

We shall call a $\zeta$-genuine function $f : G(F)\backslash \widetilde G \rightarrow \overline{\rationals}_{\ell}$ \emph{cuspidal} if for any $\tilde g \in \widetilde G$ and the unipotent radical $N$ of any proper parabolic subgroup of $G_F$, there holds
$$
\int_{N(F) \backslash N(\mathbb A_F)} f(n\tilde g) dn = 0,
$$
with respect to some Haar measure $dn$ over $N(\mathbb A_F)$.

We likewise have the notion of cuspidality of $\zeta$-genuine functions over $G_{\alpha}(F)\backslash \widetilde G_{\alpha}$ for each $\alpha \in \ker^1(F, G)$ (\emph{cf.}~Remark \ref{rem-shafarevich-set}). An element $f$ of \eqref{eq-genuine-forms-center-invariant} is called a \emph{cusp form} if its restriction to each $G_{\alpha}(F) \backslash \widetilde G_{\alpha}$ satisfies the cuspidality condition.

In other words, the space of $\zeta$-genuine cusp forms is given by
\begin{equation}
\label{eq-genuine-cusp-forms}
\Fun_{\cusp, \zeta}(\widetilde{\Bun}_{G, D}/\Xi, \overline{\rationals}_{\ell}) \simeq \bigoplus_{\alpha \in \ker^1(F, G)} \Fun_{\cusp, \zeta}(G_{\alpha}(F) \backslash \widetilde G_{\alpha}/K_D \Xi, \overline{\rationals}_{\ell}).
\end{equation}
\end{void}

\begin{lem}
\label{lem-cusp-form-finite-dimensional}
The $\overline{\rationals}_{\ell}$-vector space \eqref{eq-genuine-cusp-forms} is finite-dimensional.
\end{lem}

\begin{proof}
Since $\ker^1(F, G)$ is finite, it suffices to show that each summand in \eqref{eq-genuine-cusp-forms} is finite-dimensional. However, by the proof of \cite[Corollary 1.2.3]{MR563090}, cuspidal $\zeta$-genuine functions over $G_{\alpha}(F) \backslash \widetilde G_{\alpha}$ are uniformly supported on a compact subset.
\end{proof}

\begin{rem}
Instead of choosing $\Xi$ in the definition of $\zeta$-genuine cusp forms \eqref{eq-genuine-cusp-forms}, we may fix a $\zeta$-genuine character $\omega : \widetilde{\Bun}_{Z^{\sharp}, \infty D} \rightarrow \overline{\rationals}{}^{\times}_{\ell}$ and consider the $\overline{\rationals}_{\ell}$-vector space
\begin{equation}
\label{eq-cusp-form-fixed-central-character}
\Fun_{\cusp, \zeta}^{\omega}(\widetilde{\Bun}_{G, D}, \overline{\rationals}_{\ell})
\end{equation}
of $\zeta$-genuine cusp forms which are $\widetilde{\Bun}_{Z^{\sharp}, \infty D}$-equivariant against $\omega$.

An analogue of Lemma \ref{lem-cusp-form-finite-dimensional} asserts that \eqref{eq-cusp-form-fixed-central-character} is finite-dimensional. In fact, all our considerations concerning \eqref{eq-genuine-cusp-forms} apply to \eqref{eq-cusp-form-fixed-central-character} with minor modifications.
\end{rem}

\begin{void}[The spherical Hecke algebra]
\label{void-hecke-operators}
As the final topic of this subsection, we define the spherical Hecke algebra and its action on the space of $\zeta$-genuine automorphic forms.

Let $x$ be a closed point of $\mathring X$. Recall that the \'etale level $\mu$ defines a multiplicative $A$-gerbe $\mathscr G_{\Hec_{G, x}}$ over the local Hecke stack $\Hec_{G, x}$ (\emph{cf.}~\S\ref{void-local-hecke-stack-gerbe}). Its trace-of-Frobenius yields the groupoid $G(\mathscr O_x) \backslash \widetilde G_x /G(\mathscr O_x)$ over $\Hec_{G, x}(\base)$, where $\widetilde G_x$ is the cover of $G(F_x)$ associated to $\mu$ (\emph{cf.}~Remark \ref{rem-rational-point-local-gerbe-unramified}).\footnote{By Remark \ref{rem-hecke-stack-ind-presentation}, $\Hec_{G, x}$ admits an ind-presentation by limits of algebraic stacks $\mathscr X$ such that $\mathscr X(\base)$ is canonically isomorphic to $\mathscr X^{\Fr}$ (\emph{cf.}~Remark \ref{rem-rational-points-embedding-in-frobenius-fixed-locus}).} The induced multiplicative structure on $G(\mathscr O_x) \backslash \widetilde G_x / G(\mathscr O_x)$ coincides with the one defined by the group structure on $\widetilde G_x$ its splitting over $G(\mathscr O_x)$.

In particular, the $\overline{\rationals}_{\ell}$-vector space
\begin{equation}
\label{eq-local-hecke-algebra}
\Fun_{c, \zeta}(G(\mathscr O_x) \backslash \widetilde G_x/G(\mathscr O_x), \overline{\rationals}_{\ell})
\end{equation}
of compactly supported $\zeta$-genuine functions on $G(\mathscr O_x) \backslash \widetilde G_x/G(\mathscr O_x)$ admits an algebra structure by convolution. We shall refer to \eqref{eq-local-hecke-algebra} as the \emph{spherical Hecke algebra} at $x$.
\end{void}

\begin{void}
\label{void-local-ell-group}
We shall obtain distinguished elements of \eqref{eq-local-hecke-algebra} using the geometric Satake equivalence \eqref{eq-satake-equivalence-ell-group} (for $I = \{1\}$). We fix a square root $q^{1/2} \in \overline{\rationals}_{\ell}$ which determines a half-integral Tate twist $\overline{\rationals}_{\ell}(\frac{1}{2})$.

Indeed, denote by ${}^LH_x$ the \emph{unramified local $L$-group} at $x$, \emph{i.e.}~the base change of the $L$-group \eqref{eq-ell-group} along $\pi_1^{\etale}(x, \bar x) \rightarrow \pi_1^{\etale}(X, \bar x)$ (where we assume that $\bar x$ lies over $x$). The stalk of \eqref{eq-satake-equivalence-ell-group} at $x$ (\emph{cf.}~Proposition \ref{prop-local-hecke-stack-derived-category-stalk}) yields a \emph{monoidal} equivalence
\begin{equation}
\label{eq-satake-equivalence-finite-field}
\Sat_{\mathscr G, \zeta}(\Hec_{G, x}) \simeq \Rep({}^LH_x).
\end{equation}
Thus, to each $V \in \Rep({}^LH_x)$ corresponds an object $\mathscr S_V$ of $\Sat_{\mathscr G, \zeta}(\Hec_{G, x})$.

Taking Grothendieck groups and applying trace-of-Frobenius of twisted constructible complexes (\emph{cf.}~\S\ref{void-constructing-genuine-functions}), we obtain a map of $\overline{\rationals}_{\ell}$-algebras
\begin{align}
\notag
K_0(\Rep({}^LH_x)) &\simeq K_0(\Sat_{\mathscr G, \zeta}(\Hec_{G, x})) \\
\label{eq-distinguished-hecke-operators}
& \rightarrow \Fun_{c, \zeta}(G(\mathscr O_x) \backslash \widetilde G_x / G(\mathscr O_x), \overline{\rationals}_{\ell}) \\
\notag
V & \mapsto h_{V, x} := \Tr(\Fr \mid \mathscr S_V).
\end{align}
\end{void}

\begin{rem}
Denote by ${}^LH_{x}^1$ the fiber of ${}^LH_x \rightarrow \pi_1^{\etale}(x, \bar x)$ at the geometric Frobenius element $\Fr_x$. It admits the natual structure of an affine variety over $\overline{\rationals}_{\ell}$ and an $H$-action by conjugation. Consider the $\overline{\rationals}_{\ell}$-algebra of $H$-invariant algebraic functions on it:
$$
\Gamma({}^LH_x^1 /\!/ H, \mathscr O) := \Gamma({}^LH_x^1, \mathscr O)^H.
$$
It receives a map from $K_0(\Rep({}^LH_x))$, mapping $V$ to the character of its restriction to ${}^LH_x^1$.

By \cite[\S 5.6]{MR3752460}, \eqref{eq-distinguished-hecke-operators} factors through an isomorphism of $\overline{\rationals}_{\ell}$-algebras
\begin{equation}
\label{eq-classical-satake-isomorphism}
	\Gamma({}^LH_x^1 /\!/ H, \mathscr O) \simeq \Fun_{c, \zeta}(G(\mathscr O_x) \backslash \widetilde G_x/G(\mathscr O_x), \overline{\rationals}_{\ell}),
\end{equation}
which is our version of the \emph{classical} Satake isomorphism.

We note that a version of the classical Satake isomorphism for covering groups (over any nonarchimedean local field) has been obtained by McNamara (\emph{cf.}~\cite{MR2963537}). In the remainder of this article, we will not use the isomorphism \eqref{eq-classical-satake-isomorphism}.
\end{rem}

\begin{void}[Hecke action]
The spherical Hecke algebra \eqref{eq-local-hecke-algebra} acts on the $\overline{\rationals}_{\ell}$-vector space of $\zeta$-genuine automorphic forms \eqref{eq-genuine-automorphic-forms}.

Indeed, by Lemma \ref{lem-global-gerbe-local-compatibility} (for $n = 1$), the trace-of-Frobenius of $(p_1)^*\mathscr G_{\Bun_{G, D}}$ over the global Hecke stack at $x$ may be identified with the groupoid
$$
\widetilde{\Bun}_{G, D + \infty x} \times^{G(\mathscr O_x) \times A} (G(\mathscr O_x) \backslash \widetilde G_x / G(\mathscr O_x)),
$$
where $\widetilde{\Bun}_{G, D + \infty x}$ is the pullback of $\widetilde{\Bun}_{G, D}$, and the superscript $G(\mathscr O_x) \times A$ means quotient by the anti-diagonal action. It admits structural morphisms $p_0$, $p_1$, $r$:
\begin{equation}
\label{eq-global-hecke-action-correspondence}
\begin{tikzcd}[column sep = 1em]
	\widetilde{\Bun}_{G, D + \infty x} \times^{G(\mathscr O_x) \times A} (G(\mathscr O_x) \backslash \widetilde G_x / G(\mathscr O_x)) \ar[d, "{(p_0, p_1)}"] \ar[r, "r"] & G(\mathscr O_x) \backslash \widetilde G_x / G(\mathscr O_x) \\
	\widetilde{\Bun}_{G, D} \times \widetilde{\Bun}_{G, D}
\end{tikzcd}
\end{equation}
where $p_0$, $r$ are projections onto the first, respectively second factor, and $p_1$ is the multiplication map. Integral transform along \eqref{eq-global-hecke-action-correspondence} thus defines an action
\begin{align}
\label{eq-global-hecke-action}
\Fun_{\zeta}(\widetilde{\Bun}_{G, D}, \overline{\rationals}_{\ell}) &\circlearrowleft \Fun_{c, \zeta}(G(\mathscr O_x) \backslash \widetilde G_x / G(\mathscr O_x), \overline{\rationals}_{\ell}) \\
\notag
f \star h &:= (p_1)_! (p_0^*f \cdot r^* h),
\end{align}
where upper-$*$ indicates pullback and lower-$*$ indicates summation along fibers.

The action \eqref{eq-global-hecke-action} preserves the subspace of $\Xi$-invariant automorphic forms \eqref{eq-genuine-forms-center-invariant}, as it commutes with the $\widetilde{\Bun}_{Z^{\sharp}, \infty D}$-action. It also preserves the subspace of cusp forms \eqref{eq-genuine-cusp-forms}, by compatibility with the constant term functors (\emph{cf.}~\S\ref{sec-constant-term-functors}).
\end{void}

\subsection{Meta-Galois group \emph{vs.}~$\vartheta$-characteristics}
\label{sec-meta-galois-twist}

\begin{void}
As the final topic before we tackle the spectral decomposition of \eqref{eq-genuine-cusp-forms}, we shall prove a result that ``explains" the appearance of Weissman's meta-Galois twist in the formation of the $L$-group (\emph{cf.}~\cite[\S4]{MR3802418}).

This relies on results of \S\ref{sec-satake-equivalence-sharp-torus}, together with a categorical analogue of geometric class field theory, which we shall presently explain.
\end{void}

\begin{void}
\label{void-artin-reciprocity-map}
We fix a geometric point $\bar{\eta}$ lying over the generic point $\eta := \Spec F$ of $X$. We also assume that $D$ is nonempty, so $\Bun_{\mathbb G_m, \infty D}$ is a $\base$-scheme.

The Artin reciprocity map is an isomorphism of topological abelian groups:
\begin{equation}
\label{eq-artin-reciprocity-map}
\Art : \pi_1(\mathring X, \bar{\eta})^{\abelian} \xrightarrow{\simeq} \Bun_{\mathbb G_m, \infty D}(\base)^{\profin},
\end{equation}
where the target denotes the profinite completion of $\Bun_{\mathbb G_m, \infty D}(\base)$. The isomorphism \eqref{eq-artin-reciprocity-map} is normalized so that the geometric Frobenius element $\varphi_x \in \Gal(\bar{\base}_x / \base_x)$, for each closed point $x \in \mathring X$, maps to $\mathscr O(x)$.

On the other hand, consider the Abel--Jacobi morphism:
\begin{equation}
\label{eq-abel-jacobi-map}
\AJ : \mathring X \rightarrow \Bun_{\mathbb G_m, \infty D}
\end{equation}
sending an $R$-point $x$ to $\mathscr O(x)$, equipped with its canonical trivialization over $\hat D_R$. The Abel--Jacobi map \eqref{eq-abel-jacobi-map} recovers \eqref{eq-artin-reciprocity-map} via trace-of-Frobenius as follows: Given a multiplicative $A$-torsor $\tau$ over $\Bun_{\mathbb G_m, \infty D}$, defining a character $\Tr(\Fr \mid \tau) : \Bun_{\mathbb G_m, \infty D}(\base) \rightarrow A$ by the classical trace-of-Frobenius construction, the character $\Tr(\Fr \mid \tau) \circ \Art$ coincides with the character of $\pi_1(\mathring X, \bar{\eta})$ associated to the $A$-torsor $\AJ^*(\tau)$ over $\mathring X$.
\end{void}

\begin{void}
We shall prove a categorical analogue of the aforementioned compatibility.

Denote by $\CExt(\pi_1^{\etale}(\mathring X, \bar{\eta}), A)$ the groupoid of central extensions of $\pi_1^{\etale}(\mathring X, \bar{\eta})$ by $A$, so the map sending an $A$-gerbe $\mathscr G$ over $\mathring X$ to the \'etale fundamental group of $\mathscr G - (\mathscr G|_{\bar{\eta}})$ (\emph{cf.}~\cite[Theorem 19.6]{MR3802418}), where $\mathscr G|_{\bar{\eta}}$ is viewed as a constant $A$-gerbe over $\mathring X$, defines a map
\begin{equation}
\label{eq-gerbe-fundamental-group-rigidified}
\Maps(\mathring X, \deloop^2 A) \rightarrow \CExt(\pi_1(\mathring X, \bar{\eta}), A).
\end{equation}
\end{void}

\begin{prop}
The following diagram canonically commutes:
\begin{equation}
\label{eq-higher-class-field-theory}
\begin{tikzcd}
	\Maps_{\mathbb E_{\infty}}(\Bun_{\mathbb G_m, \infty D}, \deloop^2 A) \ar[r, "\AJ^*"]\ar[d, "\Tr(\Fr\mid\cdot)(\base)"] & \Maps(\mathring X, \deloop^2 A) \ar[d, "\eqref{eq-gerbe-fundamental-group-rigidified}"] \\
	\Maps_{\mathbb E_{\infty}}(\Bun_{\mathbb G_m, \infty D}(\base), \deloop A) \ar[r, "\Art^*"] & \CExt(\pi_1(\mathring X, \bar{\eta}), A)
\end{tikzcd}
\end{equation}
\end{prop}

\begin{proof}
Recall that the containment
$$
\Maps_{\integers}(\Bun_{\mathbb G_m, \infty D}, \deloop^2 A) \subset \Maps_{\mathbb E_{\infty}}(\Bun_{\mathbb G_m, \infty D}, \deloop^2 A)
$$
admits a retract, sending an $\mathbb E_{\infty}$-monoidal morphism to a $\integers$-linear morphism with the same underlying $\mathbb E_1$-monoidal morphism (\emph{cf.}~\S\ref{void-symmetric-monoidal-dual-data-fiber-sequence}).

Note that both $\Tr(\Fr\mid\cdot)(\base)$ and $\AJ^*$ factor through this retract: This is obvious for $\AJ^*$, while for $\Tr(\Fr\mid\cdot)(\base)$, this holds as any symmetric monoidal morphism $\Bun_{\mathbb G_m, \infty D}(\base) \rightarrow \deloop A$ is uniquely determined by its underlying monoidal morphism. We thus reduce to proving the commutativity of \eqref{eq-higher-class-field-theory} over the subspace $\Maps_{\integers}(\Bun_{\mathbb G_m, \infty D}, \deloop^2 A)$.

Over the neutral connected component of $\Maps_{\integers}(\Bun_{\mathbb G_m, \infty D}, \deloop^2 A)$, we may establish the commutativity of \eqref{eq-higher-class-field-theory} after taking loop spaces, where it reduces to the classical compatibility between $\Art$ and $\AJ$ (\emph{cf.}~\S\ref{void-artin-reciprocity-map}).

By \'etale descent, it remains to prove that any $\integers$-linear morphism $\Bun_{\mathbb G_m, \infty D} \rightarrow \deloop^2 A$ is trivial over a finite extension of $\base$, \emph{i.e.}~
$$
\SExt^2(\Bun_{\mathbb G_m, \infty D}, A) \simeq 0,
$$
where $\SExt$ denotes the internal Ext-group for \'etale sheaves over $\Spec \base$.

Replacing $\base$ by a finite extension if necessary, we may assume that the reduced subscheme of $D$ is a finite (nonempty) collection of $\base$-points $x^I = (x^i)_{i\in I}$ of $X$. We choose an element $i_0 \in I$ and fit $\Bun_{\mathbb G_m, \infty D}$ into a system of three short exact sequences:
$$
\begin{tikzcd}[column sep = -0.5em]
	& \prod_{i\in I} \ker(L^+_{x^i}\mathbb G_m \rightarrow \mathbb G_m) \ar[d, hookrightarrow] \\
	& \Bun_{\mathbb G_m, \infty D} \ar[d, twoheadrightarrow] & \Pic^0\ar[d, hookrightarrow] \\
	\prod_{i\neq i_0}\mathbb G_m \ar[r, hookrightarrow] & \Bun_{\mathbb G_m, x^I} \ar[r, twoheadrightarrow] & \Pic \ar[d, twoheadrightarrow] \\
	& & \integers
\end{tikzcd}
$$
where $\Pic \simeq \Bun_{\mathbb G_m, x^{i_0}}$ is the Picard scheme of $X$ and $\Pic^0$ is its neutral component. It thus suffices to prove $\SExt^2(M, A) \simeq 0$, for $M$ a pro-unipotent group scheme, $\mathbb G_m$, an abelian variety, or the constant sheaf $\integers$.

The pro-unipotent case and the case $M = \integers$ are clear. For $M = \mathbb G_m$ or an abelian variety, we may replace $A$ by $\mu_n$, for an integer $n\ge 1$ invertible in $\base$. Morphisms $M \rightarrow \mu_n[2]$ of complexes are equivalent to morphisms $M_{n\tors} \rightarrow \mathbb G_m[1]$, where $M_{n\tors} \subset M$ is the subgroup $\base$-scheme of $n$-torsion elements. However, $\SExt^1(M_{n\tors}, \mathbb G_m) \simeq 0$ because $M_{n\tors}$ is finite (locally) free (\emph{cf.}~\cite[Expos\'e VIII, Proposition 3.3.1]{SGA7-I}).
\end{proof}

\begin{void}
\label{void-hilbert-central-extension}
Let us now assume that $\characteristic \base \neq 2$ and $A$ is the subgroup $\{\pm 1\}$ of $\coeff^{\times}$.

Consider the \'etale level $\mu : \deloop\mathbb G_m \rightarrow \deloop^4 \{\pm 1\}^{\otimes 2}$ defined by self-tensor product of $\Psi \text{ mod }2$, \emph{i.e.}~the morphism $\deloop\mathbb G_m \rightarrow \deloop^2\{\pm 1\}$ given by the double cover of $\mathbb G_m$. The topological cover of $\mathbb A_F^{\times}$ induced from $\mu$ is the central extension
\begin{equation}
\label{eq-hilbert-central-extension}
1 \rightarrow \{\pm 1\} \rightarrow \widetilde{\mathbb A}_F^{\times} \rightarrow \mathbb A_F^{\times} \rightarrow 1
\end{equation}
defined by the cocycle
\begin{align*}
\mathbb A_F^{\times} \otimes \mathbb A_F^{\times} &\rightarrow \{\pm 1\} \\
a\otimes b &\mapsto \prod_{x\in X} \{a, b\}_x,
\end{align*}
where $\{\cdot, \cdot\}_x$ is the quadratic Hilbert symbol at $x\in X$. Indeed, by Remark \ref{rem-covering-group-tate-duality}, the central extension \eqref{eq-hilbert-central-extension} may be obtained from the construction of \cite[\S2.2]{zhao2022metaplectic}, so this identification follows from \cite[Proposition 2.3.12]{zhao2022metaplectic}.

The central extension \eqref{eq-hilbert-central-extension} is equipped with canonical splittings over $F^{\times}$ and $\mathbb O_F^{\times}$. It induces, in particular, a central extension of $\Bun_{\mathbb G_m, \infty D}(\base)$, whose pullback along $\Art$ is the \emph{meta-Galois group} of $\mathring X$ (\emph{cf.}~\cite[\S4]{MR3802418}):
\begin{equation}
\label{eq-meta-galois-group}
	1 \rightarrow \{\pm 1\} \rightarrow \widetilde{\pi}{}^{\etale}_1(\mathring X, \bar{\eta}) \rightarrow \pi_1^{\etale}(\mathring X, \bar{\eta}) \rightarrow 1.
\end{equation}

We shall identify \eqref{eq-meta-galois-group} as the \'etale fundamental group of the $\{\pm 1\}$-gerbe $\Psi(\Omega_{\mathring X}^{-1})\text{ mod }2$ of $\vartheta$-characteristics of $\mathring X$ (\emph{cf.}~\S\ref{void-theta-shift}).
\end{void}

\begin{thm}
\label{thm-meta-galois-twist}
The topological central extension \eqref{eq-meta-galois-group} is canonically identified with the image of $\Psi(\Omega_{\mathring X}^{-1})\text{ mod }2$ under \eqref{eq-gerbe-fundamental-group-rigidified}.
\end{thm}

\begin{proof}
First, we observe that the following diagram commutes
\begin{equation}
\label{eq-trace-of-frobenius-symmetric-monoidal-compatibility}
\begin{tikzcd}[column sep = 1em]
	\Maps_{\mathbb E_{\infty}}(\deloop_{\mathring X}\mathbb G_m, \deloop_{\mathring X}^4\{\pm 1\}^{\otimes 2}) \ar[r, "\int_X"]\ar[d] & \Maps_{\mathbb E_{\infty}}(\Bun_{\mathbb G_m, \infty D}, \deloop^2 \{\pm 1\}) \ar[d, "\Tr(\Fr\mid\cdot)(\base)"] \\
	\Maps_{\mathbb E_{\infty}}(F^{\times} \backslash \mathbb A_F^{\times}/ K_{\infty D}, \deloop \{\pm 1\}) \ar[r, phantom, "\simeq"] & \Maps_{\mathbb E_{\infty}}(\Bun_{\mathbb G_m, \infty D}(\base), \deloop \{\pm 1\})
\end{tikzcd}
\end{equation}
where the top horizontal arrow is given by the $\mathbb E_{\infty}$-monoidal version of \eqref{eq-infinite-level-global-trace-map} and the left vertical arrow is the construction of ad\`elic covers from \'etale levels. The commutativity of \eqref{eq-trace-of-frobenius-symmetric-monoidal-compatibility} is an $\mathbb E_{\infty}$-monoidal version of Lemma \ref{lem-global-gerbe-adelic-description} and follows from its proof.

By construction, \eqref{eq-meta-galois-group} is the image of $\mu$ under the lower circuit of \eqref{eq-trace-of-frobenius-symmetric-monoidal-compatibility}, followed by pullback along Artin reciprocity \ref{eq-artin-reciprocity-map}. By the commutativity of \eqref{eq-trace-of-frobenius-symmetric-monoidal-compatibility} and \eqref{eq-higher-class-field-theory}, we may also express it as the image of $\mu$ under the composition
\begin{align*}
	\Maps_{\mathbb E_{\infty}}(\deloop_{\mathring X}\mathbb G_m, \deloop^4_{\mathring X}\{\pm 1\}^{\otimes 2}) & \xrightarrow{\int_X} \Maps_{\mathbb E_{\infty}}(\Bun_{\mathbb G_m, \infty D}, \deloop^2 \{\pm 1\}) \\
	& \xrightarrow{\AJ^*} \Maps(\mathring X, \deloop^2 \{\pm 1\}) \xrightarrow{\eqref{eq-gerbe-fundamental-group-rigidified}} \CExt(\pi_1^{\etale}(\mathring X, \bar{\eta}), \{\pm 1\}).
\end{align*}

It therefore remains to identify the image of $\mu$ under $\AJ^* \circ \int_X$ with $\Psi(\Omega_{\mathring X}^{-1})\text{ mod }2$. However, since $\AJ$ factors through the map $\varpi^1 : \mathring X \rightarrow \Gr_{\mathbb G_m, \{1\}}$ (\emph{cf.}~\S\ref{void-torus-grassmannian-canonical-sections} for $\Lambda \simeq \integers$), this identification is a special case of Proposition \ref{prop-sharp-torus-commutative-gerbe-identification}.
\end{proof}

\begin{rem}
It follows from Theorem \ref{thm-meta-galois-twist} that the meta-Galois group for function fields (global, local, and local integral) is \emph{non-canonically} split, and is functorial with respect to finite separable extensions. These facts have been established by Weissman by different means (\emph{cf.}~\cite[\S4.2, \S4.4]{MR3802418}).

Moreover, for $\base$ of characteristic $2$, Weissman stipulates the meta-Galois group to be the split extension of $\pi_1^{\etale}(\mathring X, \bar{\eta})$ by $\{\pm 1\}$ (\emph{cf.}~\cite[\S4.1]{MR3802418}). This appears to align well with the fact that $\mathring X$ admits a canonical $\vartheta$-characteristic in characteristic $2$ (\emph{cf.}~\cite{MR286136}).
\end{rem}

\medskip

\section{Spectral decomposition}
\label{sec-spectral-decomposition}

The goal of this section is to construct the spectral decomposition of $\zeta$-genuine cusp forms, following V.~Lafforgue's method (\emph{cf.}~\cite{MR3787407}). All arguments of this section are minor variants of those of \emph{op.cit.}~and are already sketched in \cite[\S14]{MR3787407}. Thus, one may view this section as explaining how to put our results in contact with Lafforgue's work.

Let $\base$ be a finite field of cardinality $q$. Let $X$ be a smooth, proper, geometrically connected $\base$-curve with a nonempty $\base$-finite closed subscheme $D \subset X$. Write $\mathring X := X \setminus D$ and $\Ran$ for the Ran space of $\mathring X$ (\emph{cf.}~\S\ref{void-ran-space-definition}). Let $\ell$ be a prime not dividing $q$ and $\overline{\rationals}_{\ell}$ be an algebraic closure of $\rationals_{\ell}$ with a chosen square root $q^{1/2}$. Let $A$ be a finite subgroup of $\overline{\rationals}_{\ell}$ and write $\zeta : A \hookrightarrow \overline{\rationals}_{\ell}$ for the inclusion. Let $G$ be a smooth affine group $X$-scheme which is reductive over $\mathring X$. Let $\mu$ be an $A$-valued \'etale level of the base change of $G$ to $\mathring X$. Let $\Xi$ be as in \S\ref{void-connected-center-lattice}.

\subsection{Statement}

\begin{void}[$L$-parameters]
We fix an algebraic closure $\bar F$ of the field of fractions $F$ of $X$ and write $\bar{\eta} := \Spec \bar F$ for the corresponding geometric point of $X$. The L-group \eqref{eq-ell-group}, applied to $\mathring X$, yields a short exact sequence of topological groups
\begin{equation}
\label{eq-ell-group-finite-form}
1 \rightarrow H_{\bar{\eta}}(\overline{\rationals}_{\ell}) \rightarrow {}^L H_{\mathring X} \rightarrow \pi_1^{\etale}(\mathring X, \bar{\eta}) \rightarrow 1.
\end{equation}

By an \emph{$L$-parameter}, we mean an $H_{\bar{\eta}}(\overline{\rationals}_{\ell})$-conjugacy class of continuous sections
\begin{equation}
\label{eq-ell-parameter}
\sigma : \pi_1^{\etale}(\mathring X, \bar{\eta}) \rightarrow {}^LH_{\mathring X}
\end{equation}
of the surjection in \eqref{eq-ell-group-finite-form}.
\end{void}

\begin{void}
\label{void-ell-parameter-semisimplicity}
By construction, \eqref{eq-ell-group-finite-form} is the pullback of an extension ${}^LH_{\Gamma}$ of some finite quotient $\pi_1^{\etale}(\mathring X, \bar{\eta}) \twoheadrightarrow \Gamma$ by $H_{\bar{\eta}}(\overline{\rationals}_{\ell})$. Thus, we may view ${}^LH_{\Gamma}$ as (the $\overline{\rationals}_{\ell}$-points of) an algebraic group over $\overline{\rationals}_{\ell}$ whose neutral component is reductive.

We may use this observation to define ``semisimplicity" of an $L$-parameter. Namely, a continuous section \eqref{eq-ell-parameter} is \emph{semisimple} if the Zariski closure of its image in some ${}^LH_{\Gamma}$ has reductive neutral component. An $L$-parameter $[\sigma]$ is \emph{semisimple} if some (equivalently, any) representative of $[\sigma]$ is semisimple.
\end{void}

\begin{rem}
\label{rem-ell-parameter-rationality}
Since \eqref{eq-ell-group-finite-form} is induced from an extension ${}^LH_{\mathring X, \coeff}$ of $\pi_1^{\etale}(\mathring X, \bar{\eta})$ by $H_{\bar{\eta}}(\coeff)$ for some finite extension $\coeff$ of $\rationals_{\ell}$, it also makes sense to say that an $L$-parameter is ``defined" over an intermediate field between $\coeff$ and $\overline{\rationals}_{\ell}$.

To be more precise, an $L$-parameter is \emph{defined} over $\coeff_1$, for an intermediate field $\coeff \subset \coeff_1 \subset \overline{\rationals}_{\ell}$, if a representative of it factors through a continuous homomorphism $\pi_1^{\etale}(\mathring X, \bar{\eta}) \rightarrow {}^LH_{\mathring X, \coeff_1}$, where ${}^LH_{\mathring X, \coeff_1}$ is the extension of $\pi_1^{\etale}(\mathring X, \bar{\eta})$ by $H_{\bar{\eta}}(\coeff_1)$ induced from ${}^LH_{\mathring X, \coeff}$ along the inclusion $H_{\bar{\eta}}(\coeff) \subset H_{\bar{\eta}}(\coeff_1)$.
\end{rem}

\begin{void}
\label{void-local-ell-parameter}
For each closed point $x \in \mathring X$, we choose a $\bar{\base}$-point $\bar x$ lying over $x$ to form $\pi_1^{\etale}(x, \bar x)$. A choice of an \'etale path $\bar x \simeq \bar{\eta}$ in $\mathring X$ yields a continuous homomorphism $\pi_1^{\etale}(x, \bar x) \rightarrow \pi_1^{\etale}(\mathring X, \bar{\eta})$, along which we may pull back \eqref{eq-ell-group-finite-form} to obtain the unramified local $L$-group (\emph{cf.}~\S\ref{void-local-ell-group})
\begin{equation}
\label{eq-unramified-local-ell-group}
1 \rightarrow H_{\bar{\eta}}(\overline{\rationals}_{\ell}) \rightarrow {}^LH_x \rightarrow \pi_1^{\etale}(x, \bar x) \rightarrow 1.
\end{equation}

Thus, given an $L$-parameter $[\sigma]$, we obtain a $H_{\bar{\eta}}(\overline{\rationals}_{\ell})$-conjugacy class $[\sigma_x]$ of continuous sections $\sigma_x : \pi_1^{\etale}(x, \bar x) \rightarrow {}^LH_x$ of the surjection in \eqref{eq-unramified-local-ell-group}.

Writing $\Fr_x \in \pi_1^{\etale}(x, \bar x)$ for the geometric Frobenius element, we obtain an $H_{\bar{\eta}}(\overline{\rationals}_{\ell})$-conjugacy class $[\sigma_x](\Fr_x)$ in ${}^LH_x$, which we may regard as the (unramified) local $L$-parameter associated to $\sigma$ at $x$.
\end{void}

\begin{void}
In \S\ref{sec-automorphic-forms}, we have defined the $\overline{\rationals}_{\ell}$-vector space of $\zeta$-genuine cusp forms \eqref{eq-genuine-cusp-forms} equipped with the action of the spherical Hecke algebra at each $x \in \mathring X$ (\emph{cf.}~\S\ref{void-hecke-operators}).

Our main result is a deomposition of \eqref{eq-genuine-cusp-forms} indexed by semisimple $L$-parameters, where each summand is a simultaneous eigenspace for all the Hecke operators.
\end{void}

\begin{thm}
\label{thm-spectral-decomposition}
There is a canonical isomorphism
\begin{equation}
\label{eq-spectral-decomposition}
\Fun_{\cusp, \zeta}(\widetilde{\Bun}_{G, D} / \Xi, \overline{\rationals}_{\ell}) \simeq \bigoplus_{[\sigma]} \mathbf H_{D, [\sigma]}
\end{equation}
indexed by a set of $L$-parameters $[\sigma]$, such that:
\begin{enumerate}
	\item each $[\sigma]$ is semisimple (\emph{cf.}~\S\ref{void-ell-parameter-semisimplicity});
	\item each $[\sigma]$ is defined over some finite extension of $\rationals_{\ell}$ (\emph{cf.}~Remark \ref{rem-ell-parameter-rationality});
	\item each $\mathbf H_{D, [\sigma]}$ is a Hecke eigenspace: For any $x \in \mathring X$ and $V \in \Rep({}^LH_x)$, the operator $h_{V, x}$ acts on $\mathbf H_{D, [\sigma]}$ by the scalar $\Tr([\sigma_x](\Fr_x) \mid V)$ (\emph{cf.}~\S\ref{void-local-ell-parameter}).
\end{enumerate}
\end{thm}

\begin{void}
The spectral decomposition \eqref{eq-spectral-decomposition} is compatible with inclusions of nonempty $\base$-finite closed subschemes $D \subset D_1$ of $X$, \emph{i.e.}~the restriction map
\begin{equation}
\label{eq-ramification-divisor-enlargement}
\Fun_{\cusp, \zeta}(\widetilde{\Bun}_{G, D}/\Xi, \overline{\rationals}_{\ell}) \subset \Fun_{\cusp, \zeta}(\widetilde{\Bun}_{G, D_1}/\Xi, \overline{\rationals}_{\ell})
\end{equation}
carries $\mathbf H_{D, [\sigma]}$ into $\mathbf H_{D_1, [\sigma]}$, where we view $[\sigma]$ as an $L$-parameter over $X \setminus D_1$ by restriction along $X \setminus D_1 \subset X \setminus D = \mathring X$.

Assuming that $G$ splits, we may take the union of \eqref{eq-spectral-decomposition} over the poset of nonempty $\base$-finite closed subschemes $D$ to obtain the ``generic form" of spectral decomposition
\begin{equation}
\label{eq-spectral-decomposition-generic-form}
	\Fun_{\cusp, \zeta}(G(F) \backslash \widetilde G/\Xi, \overline{\rationals}_{\ell}) \simeq \bigoplus_{[\sigma]} \mathbf H_{[\sigma]},
\end{equation}
where $[\sigma]$ ranges over $H_{\bar{\eta}}(\overline{\rationals}_{\ell})$-conjugacy classes of continuous sections $\Gal(\bar F/F) \rightarrow {}^LH_{\eta}$. Here, ${}^LH_{\eta}$ denotes the pullback of \eqref{eq-ell-parameter} along $\Gal(\bar F/F) \simeq \pi_1^{\etale}(\eta, \bar{\eta}) \rightarrow \pi_1^{\etale}(\mathring X, \bar{\eta})$.
\end{void}

\begin{rem}
If one is only interested in the generic form \eqref{eq-spectral-decomposition-generic-form} of the spectral decomposition, it suffices to replace our input $(G, \mu)$ by a reductive group $F$-scheme $G$ equipped with an $A$-valued \'etale level $\mu$. Indeed, one extends $G$ to a smooth group $X$-scheme as in \cite[\S12.1]{MR3787407} and $\mu$ likewise extends to some open subscheme of $X$ by \cite[Lemma 2.2.5]{zhao2022metaplectic}.

However, we note that an extension of $\mu$ across a closed point $x$ of $X$ is \emph{additional data}. Thus, in order to state the compatibility between \eqref{eq-spectral-decomposition-generic-form} and the Satake isomorphism at a point $x \in X$ (\emph{i.e.}~the analogue of statement (3) of Theorem \ref{thm-spectral-decomposition}), we must assume that $G$ is reductive over $\mathscr O_x$ and \emph{choose} an extension of $\mu$ along $\Spec F_x \subset \Spec \mathscr O_x$.
\end{rem}

\begin{rem}
In the statement of Theorem \ref{thm-spectral-decomposition}, we assumed $D \neq \emptyset$. However, it also applies when $D = \emptyset$ and $X$ is \emph{not} isomorphic to $\mathbb P^1$. Indeed, this is because we only use the property that $\mathring X$ is an algebraic $K(\pi, 1)$-space (\emph{cf.}~\S\ref{void-ell-group}), which holds as long as it is not isomorphic to $\mathbb P^1$.
\end{rem}

\subsection{Cohomology of Shtukas}
\label{sec-cohomology-of-shtukas}

\begin{void}
The remainder of this article is devoted to the proof of Theorem \ref{thm-spectral-decomposition}. The main construction involved is the moduli space of Shtukas, which dates back to Drinfeld's pioneering works \cite{MR918745, MR902291, MR936697}.

In this subsection, we shall define the cohomology of Shtukas in the twisted setting. This relies on the twisted geometric Satake equivalence (\emph{cf.}~Theorem \ref{thm-satake-equivalence}) and the compatibility between the local and global $A$-gerbes (\emph{cf.}~Lemma \ref{lem-global-gerbe-local-compatibility}).
\end{void}

\begin{void}[Moduli of Shtukas]
\label{void-moduli-of-shtukas}
For any integer $n\ge 0$, we define the moduli stack of Shtukas to be the fiber product of $\base$-stacks
\begin{equation}
\label{eq-shtuka-definition}
\begin{tikzcd}
	\Sht_{G, D}^{[n]} \ar[r]\ar[d] & \Hec^{[n]}(\Bun_{G, D}) \ar[d, "{(p_0, p_n)}"] \\
	\Bun_{G, D} \ar[r, "{(\id, \Fr)}"] & \Bun_{G, D} \times \Bun_{G, D}
\end{tikzcd}
\end{equation}
where the structural maps $p_0$, $p_n$ are as in \S\ref{void-global-iterated-hecke-stack}.

In other words, an $R$-point of $\Sht_{G, D}^{[n]}$ consists of $R$-points $\underline x{}_1, \cdots, \underline x{}_n$ of $\Ran$, a chain of $n$ modifications along with an isomorphism
\begin{equation}
\label{eq-shtuka-point-as-modifications}
(P^0, \phi^0) \overset{\underline x{}_1}{\sim} (P^1, \phi^1) \overset{\underline x{}_2}{\sim} \cdots \overset{\underline x{}_n}{\sim} (P^n, \phi^n) \simeq ({}^{\tau} P^0, {}^{\tau} \phi^n),
\end{equation}
where the superscript $\tau$ means pullback along the endomorphism $\id_X \times \Fr_R$ of $X_R$.

The $R$-points $\underline x{}_1, \cdots, \underline x{}_n$ of $\Ran^n$ are called the \emph{paws} of the Shtuka \eqref{eq-shtuka-point-as-modifications}. By \cite[Proposition 2.16]{MR2061225} (\emph{cf.}~\cite[Lemma 12.19]{MR3787407} for nonsplit $G$), $\Sht_{G, D}^{[n]}$ is a relative ind-Deligne--Mumford stack over $\Ran^n$.
\end{void}

\begin{rem}
\label{rem-shtuka-no-leg}
For $n = 0$, the Cartesian square \eqref{eq-shtuka-definition} reduces to the square \eqref{eq-fixed-point-prestack} defining the Frobenius-fix point stack of $\Bun_{G, D}$. Thanks to Remark \ref{rem-rational-points-embedding-in-frobenius-fixed-locus}, the latter is identified with the discrete stack $\Bun_{G, D}(\base)$, so we have a canonical isomorphism
\begin{equation}
\label{eq-shtuka-no-leg}
	\Sht_{G, D}^{[0]} \simeq \Bun_{G, D}(\base).
\end{equation}
\end{rem}

\begin{void}
\label{void-shtuka-gerbe-trivialization}
For each $j = 1, \cdots, n$, we have a morphism
\begin{equation}
\label{eq-shtuka-projection-to-local-hecke-stack}
r_j : \Sht_{G, D}^{[n]} \rightarrow \Hec_G
\end{equation}
sending \eqref{eq-shtuka-point-as-modifications} to the restriction of $P^{n - 1} \overset{\underline x{}_n}{\sim} P^n$ to $D_{\underline x{}_n}$.

Recall the local $A$-gerbe $\mathscr G_{\Hec_G}$ over $\Hec_G$ (\emph{cf.}~\S\ref{void-local-hecke-stack-gerbe}). We shall need a canonical trivialization of the sum of pullbacks:
\begin{equation}
\label{eq-shtuka-pullback-gerbe}
\mathscr G_{\Sht_{G, D}^{[n]}} := \sum_{j = 1}^n (r_j)^* \mathscr G_{\Hec_G}.
\end{equation}

Indeed, by Lemma \ref{lem-global-gerbe-local-compatibility}, the $A$-gerbe \eqref{eq-shtuka-pullback-gerbe} is canonically identified with the pullback of $(p_n)^*\mathscr G_{\Bun_{G, D}} - (p_0)^*\mathscr G_{\Bun_{G, D}}$ along the top horizontal arrow of \eqref{eq-shtuka-definition}. Using the commutativity of \eqref{eq-shtuka-definition}, we may rewrite it as the pullback of
\begin{equation}
\label{eq-shtuka-gerbe-as-frobenius-difference}
(\Fr_{\Bun_{G, D}})^*\mathscr G_{\Bun_{G, D}} - \mathscr G_{\Bun_{G, D}}
\end{equation}
along the left vertical arrow. However, the $A$-gerbe \eqref{eq-shtuka-gerbe-as-frobenius-difference} admits a canonical trivialization by \eqref{eq-baffling-isomorphism-prestack}. This yields the desired trivialization of \eqref{eq-shtuka-pullback-gerbe}.
\end{void}

\begin{void}
Given a family of finite sets $I_1, \cdots, I_n$, we write $\Sht_{G, D}^{I_1, \cdots, I_n}$ for the base change of $\Sht_{G, D}^{[n]}$ along the product of the maps $\mathring X^{I_j} \rightarrow \Ran$ over $j = 1, \cdots, n$, with respect to the paws.

We shall construct a functor
\begin{equation}
\label{eq-satake-coefficients}
\Rep(\prod_{j = 1}^n {}^LH_{\mathring X^{I_j}}) \rightarrow \derived(\Sht_{G, D}^{I_1, \cdots, I_n}),
\end{equation}
where ${}^LH_{\mathring X^{I_j}}$ is the L-group defined in \S\ref{void-ell-group-multiple-point} (for $j = 1, \cdots, n$). Here, $\prod_{j = 1}^n {}^LH_{\mathring X^{I_j}}$ is viewed as an extension of $\prod_{j = 1}^n \pi_1^{\etale}(\mathring X^{I_j}, \bar{\eta})$ by $H_{\bar{\eta}}^I(\overline{\rationals}_{\ell})$ (with $I := I_1 \sqcup \cdots \sqcup I_n$) and the source of \eqref{eq-satake-coefficients} is the category of its continuous representations on finite-dimensional $\overline{\rationals}_{\ell}$-vector spaces which are algebraic over $H^I_{\bar{\eta}}(\overline{\rationals}_{\ell})$.
\end{void}

\begin{void}[Construction of \eqref{eq-satake-coefficients}]
\label{void-shtuka-coefficient-construction}
Recall the outer convolution diagram $\widetilde{\Hec}{}_G^{I_1, \cdots, I_n}$ together with the structural morphisms $\tilde p_1, \cdots, \tilde p_n$ (\emph{cf.}~\S\ref{void-outer-convolution-diagram-construction}). Denote by $\Sat_{\mathscr G, \zeta}(\widetilde{\Hec}{}^{I_1, \cdots, I_n}_G)$ the category of $(\sum_{j = 1}^n \tilde p{}_j^*\mathscr G_{\Hec_G}, \zeta)$-twisted perverse ULA sheaves over $\widetilde{\Hec}{}_G^{I_1, \cdots, I_n}$ relative to $\mathring X^I$. (As in \S\ref{void-satake-subcategory}, perversity is defined with respect to the pullback to the corresponding version of the affine Grassmannian $\widetilde{\Gr}{}_G^{I_1, \cdots, I_n}$, parametrizing modifications \eqref{eq-outer-convolution-point} where $P^0$ is trivial.) There is a morphism
\begin{equation}
\label{eq-shtuka-maps-to-convolution-hecke-stack}
	r^{I_1, \cdots, I_n} : \Sht_{G, D}^{I_1, \cdots, I_n} \rightarrow \widetilde{\Hec}{}_G^{I_1, \cdots, I_n}
\end{equation}
sending \eqref{eq-shtuka-point-as-modifications} to its restriction over $D_{\underline x{}_1\cup\cdots\cup\underline x{}_n}$ (and forgetting the last isomorphism).

Consider the functor
\begin{align}
\label{eq-satake-functor-out-of-product}
	\prod_{j = 1}^n \Rep({}^LH_{\mathring X^{I_j}}) &\rightarrow \Sat_{\mathscr G, \zeta}(\widetilde{\Hec}{}_G^{I_1, \cdots, I_n}) \\
\notag
	(V_1, \cdots, V_n) & \mapsto \bigotimes_{j = 1}^n \tilde p_j^* \mathscr S_{I_j, V_j},
\end{align}
where $\mathscr S_{I_j, V_j} \in \Sat_{\mathscr G, \zeta}(\Hec_{G, I_j})$ denotes the image of $V_j$ under the equivalence \eqref{eq-satake-equivalence-ell-group}. Since \eqref{eq-satake-functor-out-of-product} is $\overline{\rationals}_{\ell}$-multilinear and right exact in each factor, it factors through the tensor product of the Tannakian categories $\Rep({}^LH_{\mathring X^{I_j}})$, which is identified with $\Rep(\prod_{j = 1}^n {}^LH_{\mathring X^{I_j}})$ (\emph{cf.}~\cite[\S 5.18]{MR1106898}). This yields a functor
\begin{equation}
\label{eq-satake-functor-iterated-hecke-stack}
\Rep(\prod_{j = 1}^n {}^LH_{\mathring X^{I_j}}) \rightarrow \Sat_{\mathscr G, \zeta}(\widetilde{\Hec}{}_G^{I_1, \cdots, I_n}).
\end{equation}

The desired morphism \eqref{eq-satake-coefficients} is the composition of \eqref{eq-satake-functor-iterated-hecke-stack} with the pullback along \eqref{eq-shtuka-maps-to-convolution-hecke-stack}, which is an \emph{untwisted} complex by the canonical trivialization of the $A$-gerbe \eqref{eq-shtuka-pullback-gerbe}.
\end{void}

\begin{rem}
\label{rem-shtuka-coefficients-perverse-ula}
Let us express $\widetilde{\Hec}{}_G^{I_1, \cdots, I_n}$ as the quotient of $\widetilde{\Gr}{}_G^{I_1, \cdots, I_n}$ by $L^+_I G$, which acts by modifying the trivialization of $P^0$. It acts on any closed subscheme $Z_{\alpha}$ of $\widetilde{\Gr}{}_G^{I_1, \cdots, I_n}$ through some finite type quotient $L_I^+G \rightarrow H_{\beta}$ (\emph{cf.}~Remark \ref{rem-hecke-stack-ind-presentation}).

The pullback of $\Sht_{G, D}^{I_1, \cdots, I_n}$ to the stack $L^+_IG \backslash Z_{\alpha}$ along \eqref{eq-shtuka-maps-to-convolution-hecke-stack} is smooth over $H_{\beta} \backslash Z_{\alpha}$ of the same relative dimension as $H_{\beta}$ over $\mathring X^I$ (\emph{cf.}~\cite[Proposition 2.8]{MR3787407}). This implies that the image of \eqref{eq-satake-coefficients} consists of perverse ULA sheaves relative to $\mathring X^I$.
\end{rem}

\begin{void}[$\Bun_{Z^{\sharp}, \infty D}(\base)$-equivariance]
\label{void-shtuka-coefficients-sharp-center-equivariance}
Recall the $\Bun_{Z^{\sharp}, \infty D}$-action on $\Bun_{G, D}$ (\emph{cf.}~\S\ref{void-sharp-center}). This induces a $\Bun_{Z^{\sharp}, \infty D}$-action on $\Hec^{[n]}(\Bun_{G, D})$, by acting on each term in a chain of modifications \eqref{eq-global-iterated-hecke-stack-modification}. This action lifts to a $\Bun_{Z^{\sharp}, \infty D}(\base)$-action on $\Sht_{G, D}^{[n]}$.

We shall argue that \eqref{eq-satake-coefficients} factors through the category of $\Bun_{Z^{\sharp}, \infty D}(\base)$-equivariant objects in $\derived(\Sht_{G, D}^{I_1, \cdots, I_n})$. Indeed, by the equivariance structure on \eqref{eq-shtuka-maps-to-convolution-hecke-stack}:
$$
\begin{tikzcd}
	\Bun_{Z^{\sharp}, \infty D}(\base) \ar[r, phantom, "\circlearrowright"]\ar[d] & \Sht_{G, D}^{I_1, \cdots, I_n} \ar[d, "\eqref{eq-shtuka-maps-to-convolution-hecke-stack}"] \\
	L^+_I \deloop Z^{\sharp} \ar[r, phantom, "\circlearrowright"] & \widetilde{\Hec}{}_G^{I_1, \cdots, I_n}
\end{tikzcd}
$$
where the left vertical arrow is the restriction to the appropriate formal disk, this follows from the observation below.
\end{void}

\begin{lem}
The functor forgetting $L_I^+\deloop Z^{\sharp}$-equivariance is an equivalence:
$$
\Sat_{\mathscr G, \zeta}(\widetilde{\Hec}{}_G^{I_1, \cdots, I_n})^{L_I^+\deloop Z^{\sharp}} \simeq \Sat_{\mathscr G, \zeta}(\widetilde{\Hec}{}_G^{I_1, \cdots, I_n}).
$$
\end{lem}

\begin{proof}
The forgetful functor is fully faithful, as $L^+_I \deloop Z^{\sharp}$ may be written as an inverse limit of connected smooth algebraic $\mathring X^I$-stacks. It remains to prove that it is essentially surjective.

Let $\mathscr A$ be an object of $\Sat_{\mathscr G, \zeta}(\widetilde{\Hec}{}_G^{I_1, \cdots, I_n})$. Denote by $a$ and $p$ the action, respectively projection map from $L_I^+\deloop Z^{\sharp} \times \widetilde{\Hec}{}_G^{I_1, \cdots, I_n}$ to $\widetilde{\Hec}{}_G^{I_1, \cdots, I_n}$. It suffices to find an isomorphism
\begin{equation}
\label{eq-sharp-center-equivariance-structure}
a^*\mathscr A \simeq p^*\mathscr A
\end{equation}
extending the natural isomorphism over the neutral section $e \times \widetilde{\Hec}{}_G^{I_1, \cdots, I_n}$, as such an isomorphism is unique if it exists.

By universal local acyclicity, it suffices to construct \eqref{eq-sharp-center-equivariance-structure} over the pointwise disjoint locus in $\mathring X^I$ (\emph{cf.}~\cite[Theorem 6.8]{MR4630128}). Using the factorization structure, we may further reduce to the case where $I$ is a singleton, where \eqref{eq-sharp-center-equivariance-structure} exists by Proposition \ref{prop-virtual-connected-components}.
\end{proof}

\begin{void}
\label{void-cohomology-of-shtukas-construction}
Let us now invoke the lattice $\Xi$. Consider the structural morphism
\begin{equation}
\label{eq-shtuka-structural-morphism-lattice-quotient}
\nu_{\Xi}^{I_1, \cdots, I_n} : \Sht_{G, D}^{I_1, \cdots, I_n} / \Xi \rightarrow \mathring X^I,
\end{equation}
which is a relative ind-Deligne--Mumford stack (\emph{cf.}~\S\ref{void-moduli-of-shtukas}).

On the other hand, the functor \eqref{eq-satake-coefficients} factors through the category of $\Xi$-equivariant objects (\emph{cf.}~\S\ref{void-shtuka-coefficients-sharp-center-equivariance}), defining a functor
\begin{align}
\label{eq-satake-coefficients-lattice-equivariant}
	\Rep(\prod_{j = 1}^n {}^LH_{\mathring X^{I_j}}) &\rightarrow \derived(\Sht_{G, D}^{I_1, \cdots, I_n} / \Xi), \\
	\notag
	V & \mapsto \mathscr F_{I_1, \cdots, I_n, V}.
\end{align}

For any object $V\in \Rep(\prod_{j = 1}^n {}^LH_{\mathring X^{I_j}})$, we shall write
\begin{equation}
\label{eq-cohomology-of-shtukas}
\mathscr H_{I_1, \cdots, I_n, V} := (\nu_{\Xi}^{I_1, \cdots, I_n})_! \mathscr F_{I_1, \cdots, I_n, V} \in \Ind \derived(\mathring X^I)
\end{equation}
and refer to it as the \emph{cohomology of Shtukas} with coefficients in (the perverse sheaf associated to) the representation $V$.
\end{void}

\begin{rem}
\label{rem-cohomology-of-shtukas-external-fusion}
There is a canonical isomorphism in $\Ind \derived(\mathring X^I)$:
\begin{equation}
\label{eq-cohomology-of-shtukas-external-fusion}
\mathscr H_{I_1, \cdots, I_n, V} \simeq \mathscr H_{I, V},
\end{equation}
where the right-hand-side is the cohomology of Shtukas with coefficients in the restriction of $V$ along the natural map ${}^LH_{\mathring X^I} \rightarrow \prod_{j = 1}^n {}^LH_{\mathring X^{I_j}}$.

Indeed, \eqref{eq-cohomology-of-shtukas-external-fusion} follows from the compatibility between \eqref{eq-satake-functor-iterated-hecke-stack} and (a mild generalization of) the outer convolution product (\emph{cf.}~\S\ref{void-outer-convolution-product}, for $n$ possibly distinct finite sets $I_1, \cdots, I_n$.)
\end{rem}

\begin{void}[The case $n = 0$]
\label{void-shtuka-no-leg-automorphic-forms}
For $n = 0$, the source of \eqref{eq-satake-coefficients-lattice-equivariant} is the category of finite-dimensional $\overline{\rationals}_{\ell}$-vector spaces, while the right-hand-side is $\derived(\Bun_{G, D}(\base) / \Xi)$ (\emph{cf.}~Remark \ref{rem-shtuka-no-leg}). Let us identify the image $\mathscr S_{\mathbf 1}$ of the $1$-dimensional vector space $\mathbf 1 := \overline{\rationals}_{\ell}$.

By construction, $\mathscr S_{\mathbf 1}$ is obtained as follows: Consider the constant $1$-dimensional local system $\overline{\rationals}_{\ell}$ over $\Bun_{G, D}$, regarded as twisted by the pullback of $(-\mathscr G_{\Bun_{G, D}}) \boxplus \mathscr G_{\Bun_{G, D}}$ along the diagonal of $\Bun_{G, D}$. The pullback of $\overline{\rationals}_{\ell}$ to $\Bun_{G, D}(\base)$ may be viewed as an \emph{untwisted} complex by trivializing the pullback of $(-\mathscr G_{\Bun_{G, D}}) \boxplus \mathscr G_{\Bun_{G, D}}$ along
$$
(\id, \Fr_{\Bun_{G, D}}) : \Bun_{G, D} \rightarrow \Bun_{G, D} \times \Bun_{G, D}
$$
via \eqref{eq-baffling-isomorphism-prestack}. This complex is the $1$-dimesional local system $(\mathscr L_{\Tr(\Fr \mid \mathscr G), \zeta})^{\otimes -1}$ (\emph{cf.}~\S\ref{void-trace-of-frobenius-gerbe-construction}).

Thus, \eqref{eq-genuine-functions-as-cohomology} (applied to quasi-compact substacks of $\Bun_{G, D}$) yields an isomorphism
\begin{equation}
\label{eq-cohomology-of-shtuka-no-leg}
	\mathscr H_{\emptyset, \overline{\rationals}_{\ell}} \simeq \Fun_{c, \zeta}(\widetilde{\Bun}_{G, D} / \Xi, \overline{\rationals}_{\ell}).
\end{equation}
\end{void}

\subsection{Hecke action}

\begin{void}
In this subsection, we construct an action of the spherical Hecke algebra (\emph{cf.}~\S\ref{void-hecke-operators}) on the cohomology of Shtukas (\emph{cf.}~\S\ref{void-cohomology-of-shtukas-construction})

More precisely, for each closed point $x \in \mathring X$ and element $h \in \Fun_{c, \zeta}(G(\mathscr O_x) \backslash \widetilde G_x / \mathscr O_x)$, we shall construct an endomorphism $T(h)$ of the \emph{restriction} of $\mathscr H_{I_1, \cdots, I_n, V}$ to $(\mathring X \setminus x)^I$, where $I_1, \cdots, I_n$ is any family of finite sets with $I := I_1 \sqcup \cdots\sqcup I_n$ and $V \in \Rep(\prod_{j = 1}^n {}^LH_{\mathring X^{I_j}})$.
\end{void}

\begin{void}
Theis construction will makes use of the Hecke correspondence for $\Sht_{G, D}^{[n]}$.

Namely, write $\Hec(\Sht_{G, D}^{[n]})$ for the prestack parametrizing two $R$-points
\begin{align}
\label{eq-shtuka-hecke-first-row}
	(P^0, \phi^0) &\overset{\underline x{}_1}{\sim} \cdots \overset{\underline x{}_{n - 1}}{\sim} (P^{n - 1}, \phi^{n - 1}) \overset{\underline x{}_n}{\sim} ({}^{\tau} P^0, {}^{\tau}\phi^0) \\
\label{eq-shtuka-hecke-second-row}
	(Q^0, \psi^0) &\overset{\underline x{}_1}{\sim} \cdots \overset{\underline x{}_{n - 1}}{\sim} (Q^{n - 1}, \psi^{n-1}) \overset{\underline x{}_n}{\sim} ({}^{\tau} Q^0, {}^{\tau} \psi^0)
\end{align}
of $\Sht_{G, D}^{[n]}$ over the same $R$-point $(\underline x{}_1, \cdots, \underline x{}_n)$ of $\Ran^n$, equipped with an isomorphism off the graph of another $R$-point $\underline x$ of $\Ran$ disjoint from $\underline x{}_1, \cdots, \underline x{}_n$.

Sending an $R$-point of $\Hec(\Sht_{G, D}^{[n]})$ to \eqref{eq-shtuka-hecke-first-row}, respectively \eqref{eq-shtuka-hecke-second-row}, as well as the restriction of $P^0 \overset{\underline x}{\sim} Q^0$ to $D_{\underline x}$ define morphisms $p$, $q$, and $r$ in the following diagram:
\begin{equation}
\label{eq-shtuka-hecke-correspondence}
\begin{tikzcd}[column sep = 0em]
	\Hec(\Sht_{G, D}^{[n]}) \ar[d, "{(p, q)}"]\ar[r, "r"] & \Hec_G(\base) \\
	\Sht_{G, D}^{[n]} \times_{\Ran^n} \Sht_{G, D}^{[n]}
\end{tikzcd}
\end{equation}
where we have identified $\Hec_G(\base)$ with the Frobenius-fixed point stack of $\Hec_G$.
\end{void}

\begin{void}
We shall formulate a compatibility statement for $A$-gerbes with respect to \eqref{eq-shtuka-hecke-correspondence}.

Recall the $A$-gerbe $\mathscr G_{\Sht_{G, D}^{[n]}}$ over $\Sht_{G, D}^{[n]}$ (\emph{cf.}~\S\ref{void-shtuka-gerbe-trivialization}). Since the post-compositions of $p$, $q$ with each $r_j$ ($j = 1, \cdots, n$) coincide, we obtain an isomorphism of $A$-gerbes
\begin{equation}
\label{eq-gerbe-shtuka-hecke-isomorphism}
p^* \mathscr G_{\Sht_{G, D}^{[n]}} \simeq q^* \mathscr G_{\Sht_{G, D}^{[n]}}.
\end{equation}

On the other hand, \emph{loc.cit.}~supplies a canonical trivialization of $\mathscr G_{\Sht_{G, D}^{[n]}}$, whose pullback along $p$ and $q$ do \emph{not} coincide under the isomorphism \eqref{eq-gerbe-shtuka-hecke-isomorphism}.
\end{void}

\begin{lem}
\label{lem-gerbe-shtuka-hecke-compatibility}
With respect to the pullbacks of the canonical trivialization of $\mathscr G_{\Sht_{G, D}^{[n]}}$ under $p$ and $q$, the isomorphism \eqref{eq-gerbe-shtuka-hecke-isomorphism} is given by the $A$-torsor $r^*(-\Tr(\Fr \mid \mathscr G_{\Hec_G}))$.
\end{lem}

\begin{proof}
The value of \eqref{eq-gerbe-shtuka-hecke-isomorphism} at an $R$-point of $\Hec(\Sht_{G, D}^{[n]})$ is the sum of isomorphisms
\begin{equation}
\label{eq-shtuka-hecke-gerbe-identification-summand}
\mathscr G_{\Hec_G}(P^{j - 1} \overset{\underline x{}_j}{\sim} P^j) \simeq \mathscr G_{\Hec_G}(Q^{j - 1} \overset{\underline x{}_j}{\sim} Q^j),
\end{equation}
supplied by the identifications $(P^{j - 1} \overset{\underline x{}_j}{\sim} P^j) \simeq (Q^{j - 1} \overset{\underline x{}_j}{\sim} Q^j)$ over $D_{\underline x{}_j}$, over $j = 1, \cdots, n$.

On the other hand, we may use Lemma \ref{lem-global-gerbe-local-compatibility} to rewrite \eqref{eq-shtuka-hecke-gerbe-identification-summand} as an isomorphism
\begin{equation}
\label{eq-shtuka-hecke-gerbe-identification-summand-difference-expression}
\mathscr G_{\Bun_{G, D}}(P^j, \phi^j) - \mathscr G_{\Bun_{G, D}}(P^{j - 1}, \phi^{j - 1}) \simeq \mathscr G_{\Bun_{G, D}}(Q^j, \psi^j) - \mathscr G_{\Bun_{G, D}}(Q^{j - 1}, \psi^{j-1}).
\end{equation}
Since \eqref{eq-shtuka-hecke-gerbe-identification-summand-difference-expression} is obtained from an identification of Hecke modifications along $\underline x{}_j$, it admits the following description: Taking the difference of the isomorphisms from Lemma \ref{lem-global-gerbe-local-compatibility}
\begin{align*}
	\mathscr G_{\Bun_{G, D}}(P^j, \phi^j) + \mathscr G_{\Hec_G}(P^j \overset{\underline x}{\sim} Q^j) & \simeq \mathscr G_{\Bun_{G, D}}(Q^j, \psi^j), \\
	\mathscr G_{\Bun_{G, D}}(P^{j - 1}, \phi^{j - 1}) + \mathscr G_{\Hec_G}(P^{j - 1} \overset{\underline x}{\sim} Q^{j - 1}) & \simeq \mathscr G_{\Bun_{G, D}}(Q^{j - 1}, \psi^{j-1}),
\end{align*}
and applying the identification in $\Hec_{G, x}$
$$
(P^{j - 1} \overset{\underline x}{\sim} Q^{j - 1}) \simeq (P^j \overset{\underline x}{\sim} Q^j),
$$
supplied by the Hecke modifications along $\underline x{}_j$, we again arrive at \eqref{eq-shtuka-hecke-gerbe-identification-summand-difference-expression}.

Using this description, we see that the sum of \eqref{eq-shtuka-hecke-gerbe-identification-summand} over $j = 1, \cdots, n$ yields the difference of the two isomorphisms from Lemma \ref{lem-global-gerbe-local-compatibility}
\begin{align}
\label{eq-shtuka-hecke-gerbe-identification-terminal}
	\mathscr G_{\Bun_{G, D}}({}^{\tau }P^0, {}^{\tau}\phi^0) + \mathscr G_{\Hec_G}({}^{\tau} P^0 \overset{\underline x}{\sim} {}^{\tau} Q^0) & \simeq \mathscr G_{\Bun_{G, D}}({}^{\tau} Q^0, {}^{\tau}\psi^0), \\
\label{eq-shtuka-hecke-gerbe-identification-initial}
	\mathscr G_{\Bun_{G, D}}(P^0, \phi^0) + \mathscr G_{\Hec_G}(P^0 \overset{\underline x}{\sim} Q^0) & \simeq \mathscr G_{\Bun_{G, D}}(Q^0, \psi^0)
\end{align}
under the composition of identifications in $\Hec_{G, x}$:
\begin{equation}
\label{eq-shtuka-hecke-as-hecke-shtuka}
(P^0 \overset{\underline x}{\sim} Q^0) \simeq (P^1 \overset{\underline x}{\sim} Q^1) \simeq \cdots \simeq ({}^{\tau} P^0 \overset{\underline x}{\sim} {}^{\tau} Q^0)
\end{equation}
supplied by the chains of Hecke modifications along $\underline x{}_1, \cdots, \underline x{}_n$.

Identifying the values of $\mathscr G_{\Bun_{G, D}}$ at $({}^{\tau} P^0, {}^{\tau}\phi^0)$ and $(P^0, \phi^0)$, respectively $({}^{\tau} Q^0, {}^{\tau}\psi^0)$ and $(Q^0, \psi^0)$, using Frobenius-equivariance \eqref{eq-baffling-isomorphism-prestack}, we may write the difference of \eqref{eq-shtuka-hecke-gerbe-identification-terminal} from \eqref{eq-shtuka-hecke-gerbe-identification-initial} as an automorphism of the trivial $A$-gerbe, \emph{i.e.}~an $A$-torsor.

The above discussion allows us to identify this $A$-torsor with the difference of the Frobenius-equivariance isomorphism
$$
\mathscr G_{\Hec_G}(P^0 \overset{\underline x}{\sim} Q^0) \simeq \mathscr G_{\Hec_G}({}^{\tau}P^0 \overset{\underline x}{\sim} {}^{\tau}Q^0)
$$
from the isomorphism induced from the identification \eqref{eq-shtuka-hecke-as-hecke-shtuka}. By construction, this is inverse to the $A$-torsor $\Tr(\Fr \mid \mathscr G_{\Hec_G})$ over $\Hec_G(\base)$ (\emph{cf.}~\S\ref{void-trace-of-frobenius-gerbe-construction}).
\end{proof}

\begin{void}[Construction of Hecke action]
\label{void-hecke-operator-cohomology-of-shtuka-construction}
Let $I_1, \cdots, I_n$ be a family of finite sets, with $I := I_1\sqcup \cdots \sqcup I_n$. Fix a closed point $x$ of $\mathring X$. Consider the base change of \eqref{eq-shtuka-hecke-correspondence}, where the paws of the Shtukas are restricted to $(\mathring X \setminus x)^I$ and the copy of local Hecke stack $\Hec_G(\base)$ is restricted to $x$:
\begin{equation}
\label{eq-shtuka-hecke-correspondence-multiset}
\begin{tikzcd}[column sep = -0.5em]
	\Hec_x(\Sht_{G, D}^{I_1, \cdots, I_n}) \ar[r, "r"]\ar[d, "{(p, q)}"] & \Hec_{G, x}(\base) \\
	\Sht_{G, D}^{I_1, \cdots, I_n} \times_{(\mathring X\setminus x)^I} \Sht_{G, D}^{I_1, \cdots, I_n}
\end{tikzcd}
\end{equation}

Given $V \in \Rep(\prod_{j = 1}^n {}^LH_{\mathring X^{I_j}})$, we have defined the $\Xi$-equivariant perverse sheaf $\mathscr F_{I_1, \cdots, I_n, V}$ over $\Sht_{G, D}^{I_1, \cdots, I_n}$ (\emph{cf.}~\S\ref{void-cohomology-of-shtukas-construction}). Using Lemma \ref{lem-gerbe-shtuka-hecke-compatibility}, we obtain an isomorphism of $\Xi$-equivariant perverse sheaves
\begin{equation}
\label{eq-hecke-correspondence-shtuka-coefficients}
p^* \mathscr F_{I_1, \cdots, I_n, V} \otimes r^*(\mathscr L_{\Tr(\Fr \mid \mathscr G_{\Hec_G})})^{\otimes -1} \simeq q^* \mathscr F_{I_1, \cdots, I_n, V}.
\end{equation}

By expressing $q$ as a union of \'etale morphisms, we obtain from \eqref{eq-hecke-correspondence-shtuka-coefficients} a morphism
$$
q_! (p^* \mathscr F_{I_1, \cdots, I_n, V} \otimes r^*(\mathscr L_{\Tr(\Fr \mid \mathscr G_{\Hec_G})})^{\otimes -1}) \rightarrow \mathscr F_{I_1, \cdots, I_n, V}
$$
and thus by taking compactly supported cohomology along $\nu^{I_1, \cdots, I_n}_{\Xi}$ (\emph{cf.}~\S\ref{void-cohomology-of-shtukas-construction}) and using the isomorphism \eqref{eq-genuine-functions-as-cohomology}, we obtain a morphism in $\Ind \derived((\mathring X \setminus x)^I)$:
$$
\mathscr H_{I_1, \cdots, I_n, V} |_{(\mathring X\setminus x)^I} \otimes \Fun_{c, \zeta}(G(\mathscr O_x) \backslash \widetilde G_x / G(\mathscr O_x)) \rightarrow \mathscr H_{I_1, \cdots, I_n, V}|_{(\mathring X\setminus x)^I},
$$
and thus a map
\begin{equation}
\label{eq-hecke-action-cohomology-of-shtukas}
T : \Fun_{c, \zeta}(G(\mathscr O_x) \backslash \widetilde G_x / G(\mathscr O_x)) \rightarrow \End (\mathscr H_{I_1, \cdots, I_n, V}|_{(\mathring X\setminus x)^I}).
\end{equation}
\end{void}

\begin{rem}
By working with the convolution local Hecke stack instead of $\Hec_{G, x}$, we also obtain an isomorphism of endomorphisms
\begin{equation}
\label{eq-hecke-action-cohomology-of-shtukas-multiplicative}
T(h_1 \star h_2) \simeq T(h_1) \circ T(h_2),
\end{equation}
for any $h_1, h_2 \in \Fun_{c, \zeta}(G(\mathscr O_x) \backslash \widetilde G_x / G(\mathscr O_x))$. Furthermore, $T$ carries the unit for the convolution product to the identity endomorphism.

Note that we certainly expect \eqref{eq-hecke-action-cohomology-of-shtukas} to lift to a morphism of $\mathbb E_1$-monoids, but we will not supply the construction. (We only need that $T$ induces a morphism of associative algebras after taking $H^0$ of $\mathscr H_{I_1, \cdots, I_n, V}|_{(\mathring X\setminus x)^I}$, which is guaranteed by \eqref{eq-hecke-action-cohomology-of-shtukas-multiplicative}.)
\end{rem}

\subsection{Partial Frobenius-equivariance}

\begin{void}
\label{void-partial-frobenius}
Given $n$ finite sets $I_1, \cdots, I_n$ ($n\ge 0$) with $I := I_1 \sqcup \cdots \sqcup I_n$, we may consider the \emph{partial Frobenius} endomorphism $F_{I_j} := \Fr_{X^{I_j}} \times \id_{X^{I\setminus I_j}}$ of $X^I$ (for $j = 1, \cdots, n$). Note that the composition $F_{I_1} \circ \cdots \circ F_{I_n}$ equals $\Fr_{X^I}$. Let $V$ be an object of $\Rep(\prod_{j = 1}^n {}^LH_{\mathring X^{I_j}})$.

In this subsection, we shall equip the corresponding cohomology of Shtukas \eqref{eq-cohomology-of-shtukas} with an $F_{I_1}$-equivariance structure, in the sense of an isomorphism
\begin{equation}
\label{eq-partial-frobenius-equivariance}
	(F_{I_1})^* \mathscr H_{I_1, \cdots, I_n, V} \simeq \mathscr H_{I_2, \cdots, I_n, I_1, V},
\end{equation}
such that its $n$-fold composition
\begin{align}
\notag
	(\Fr_{\mathring X^I})^*\mathscr H_{I_1, \cdots, I_n, V} & \simeq (F_{I_n})^* \circ \cdots \circ (F_{I_1})^* \mathscr H_{I_1, \cdots, I_n, V} \\
\label{eq-partial-frobenius-equivariance-composition}
	& \simeq (F_{I_n})^* \circ \cdots \circ (\Fr_{I_2})^* \mathscr H_{I_2, \cdots, I_1, V} \simeq \cdots \simeq \mathscr H_{I_1, \cdots, I_n, V}
\end{align}
is induced from the isomorphism $(\Fr_{\mathring X^I})^* \simeq \id$ on the \'etale site of $\mathring X^I$.
\end{void}

\begin{void}
\label{void-shtuka-partial-frobenius}
To construct \eqref{eq-partial-frobenius-equivariance}, we consider the morphism
$$
F_{I_1} : \Sht_{G, D}^{I_1, \cdots, I_n} \rightarrow \Sht_{G, D}^{I_2, \cdots, I_1},
$$
lifting the corresponding partial Frobenius endomorphism of $\mathring X^I$, sending an $R$-point
\begin{equation}
\label{eq-shtuka-point-as-modification}
(P^0, \phi^0) \overset{\underline x{}_1}{\sim} (P^1, \phi^1) \overset{\underline x{}_2}{\sim} \cdots \overset{\underline x{}_n}{\sim} ({}^{\tau} P^0, {}^{\tau} \phi^0)
\end{equation}
to the $R$-point
$$
(P^1, \phi^1) \overset{\underline x{}_2}{\sim} \cdots \overset{\underline x{}_n}{\sim} ({}^{\tau}P^0, {}^{\tau}\phi^0) \overset{{}^{\tau}\underline x{}_1}{\sim} ({}^{\tau} P^1, {}^{\tau}\phi^1).
$$

By definition, we have a commutative square
\begin{equation}
\label{eq-shtuka-frobenius-restriction-local-hecke-stack}
\begin{tikzcd}[column sep = 1.5em]
	\Sht_{G, D}^{I_1, \cdots, I_n} \ar[d, "F_{I_1}"]\ar[r, "\prod r_j"] & \prod_{j = 1}^n \Hec_{G, I_j} \ar[d, "\Fr_{\Hec_{G, I_1}} \times \id"] \\
	\Sht_{G, D}^{I_2, \cdots, I_1} \ar[r, "\prod r_j"] & \prod_{j = 1}^n \Hec_{G, I_j}
\end{tikzcd}
\end{equation}
where the horizontal morphisms are restrictions to the formal disks.

It remains to construct, for an $n$-tuple $V_j \in \Rep({}^L H_{\mathring X^{I_j}})$, an isomorphism
\begin{equation}
\label{eq-shtuka-coefficient-partial-frobenius-equivariance}
(F_{I_1})^* \bigotimes_{j = 1}^n (r_j)^* \mathscr S_{I_j, V_j} \simeq \bigotimes_{j = 1}^n (r_j)^* \mathscr S_{I_j, V_j}
\end{equation}
of (untwisted) perverse sheaves over $\Sht_{G, D}^{I_1, \cdots, I_n}$, compatibly with the $\Xi$-equivariance structure (\emph{cf.}~\S\ref{void-shtuka-coefficient-construction}). Indeed, we then obtain \eqref{eq-partial-frobenius-equivariance} for $V := \bigboxtimes_{j = 1}^n V_j$ by taking compactly supported direct image of \eqref{eq-shtuka-coefficient-partial-frobenius-equivariance} along $\nu_{\Xi}^{I_1, \cdots, I_n}$ (using the remark on base change in \cite[\S4.3]{MR3787407}) and the case for general $V$ by the universal property.
\end{void}

\begin{void}[Construction of \eqref{eq-shtuka-coefficient-partial-frobenius-equivariance}]
\label{void-shtuka-coefficient-partial-frobenius-equivariance-construction}
Since $(\Fr_{\Hec_{G, I_1}})^*$ acts as the identity on the \'etale site of $\Hec_{G, I_1}$, we have a canonical isomorphism
\begin{equation}
\label{eq-local-hecke-stack-frobenius-equivariance}
(\Fr_{\Hec_{G, I_1}})^* \mathscr S_{I_1, V_1} \simeq \mathscr S_{I_1, V_1}
\end{equation}
with respect to the canonical isomorphism \eqref{eq-baffling-isomorphism-prestack} for the local $A$-gerbe over $\Hec_{G, I_1}$.

The isomorphism \eqref{eq-shtuka-coefficient-partial-frobenius-equivariance} is obtained as the tensor product of the pullback of \eqref{eq-local-hecke-stack-frobenius-equivariance} along $r_1$ with the identity on $(r_j)^*\mathscr S_{I_j, V_j}$ for $j = 2,\cdots, n$, using the commutativity of \eqref{eq-shtuka-frobenius-restriction-local-hecke-stack}.

For this procedure to yield an isomorphism of \emph{untwisted} perverse sheaves, the following compatibility of the trivializations of $A$-gerbes is required: Along the isomorphisms
\begin{align}
\notag
	(F_{I_1})^*\sum_{j = 1}^n (r_j)^* \mathscr G_{\Hec_{G, I_j}} &\simeq (r_1)^* (\Fr_{\Hec_{G, I_1}})^*\mathscr G_{\Hec_{G, I_1}} + \sum_{j = 2}^n (r_j)^* \mathscr G_{\Hec_{G, I_j}} \\
\label{eq-partial-frobenius-gerbe-compatibility}
	&\simeq \sum_{j = 1}^n (r_j)^* \mathscr G_{\Hec_{G, I_j}}
\end{align}
supplied by \eqref{eq-shtuka-frobenius-restriction-local-hecke-stack} and \eqref{eq-baffling-isomorphism-prestack} (for $\mathscr G_{\Hec_{G, I_1}}$), the trivialization of $\sum_{j = 1}^n (r_j)^* \mathscr G_{\Hec_{G, I_j}}$ constructed in \S\ref{void-shtuka-coefficient-construction} corresponds to its own pullback by $F_{I_1}$. However, the two sides of \eqref{eq-partial-frobenius-gerbe-compatibility} evaluate at an $R$-point \eqref{eq-shtuka-point-as-modification} of $\Sht_{G, D}^{I_1, \cdots, I_n}$ to the same telescopic sum (\emph{cf.}~the proof of Lemma \ref{lem-global-gerbe-local-compatibility}), where both trivializations are induced from \eqref{eq-baffling-isomorphism-prestack} for $\mathscr G_{\Bun_{G, D}}$.

The fact that the $n$-fold composition \eqref{eq-partial-frobenius-equivariance-composition} is the canonical $\Fr_{X^I}$-equivariance structure on $\mathscr H_{I_1, \cdots, I_n, V}$ follows immediately from the construction.
\end{void}

\begin{rem}
\label{rem-partial-frobenius-independence-of-partition}
Using \eqref{eq-cohomology-of-shtukas-external-fusion}, we may write \eqref{eq-partial-frobenius-equivariance} as an isomorphism
\begin{equation}
\label{eq-partial-frobenius-equivariance-independent}
	\varphi_{I_1} : (F_{I_1})^*\mathscr H_{I, V} \simeq \mathscr H_{I, V}.
\end{equation}

By the compatibility between the construction of \eqref{eq-partial-frobenius-equivariance} with outer convolution product (along $I_2, \cdots, I_n$), the isomorphism \eqref{eq-partial-frobenius-equivariance-independent} depends only on $I_1$ as a subset of $I$, and not on the partition of $I \setminus I_1$ into $I_2, \cdots, I_n$.

Likewise, for two disjoint subsets $I_1$, $I_2$ of $I$, there is a commutative diagram
$$
\begin{tikzcd}
	(F_{I_2})^*(F_{I_1})^*\mathscr H_{I, V} \ar[r, "(F_{I_2})^*\varphi_{I_1}"]\ar[d, "\simeq"] & (F_{I_2})^*\mathscr H_{I, V} \ar[d, "\varphi_{I_2}"] \\
	(F_{I_1\sqcup I_2})^* \mathscr H_{I, V} \ar[r, "\varphi_{I_1\sqcup I_2}"] & \mathscr H_{I, V}
\end{tikzcd}
$$
\end{rem}

\begin{void}
\label{void-constructible-complex-partial-frobenius-equivariance}
Let us now specialize to the case where $I_1, \cdots, I_n$ are all singletons, so $I \simeq \{1,\cdots, n\}$. We write $F_i$ for the partial Frobenius endomorphism $F_{\{i\}}$ of $X^I$.

We write $(\Ind \derived(\mathring X^I))^{F_1, \cdots, F_n}$ for the category of objects $\mathscr A \in \Ind\derived(\mathring X^I)$ equipped with $n$ isomorphisms $(F_i)^*\mathscr A \simeq \mathscr A$, which pairwise commute (in the evident sense) and such that their composite
\begin{align*}
(\Fr_{\mathring X^I})^*\mathscr A & \simeq (F_n)^* \circ \cdots \circ (F_1)^* \mathscr A \\
& \simeq (F_n)^*\circ \cdots \circ (F_2)^*\mathscr A \simeq \cdots \simeq \mathscr A
\end{align*}
is the canonical $\Fr_{\mathring X^I}$-equivariance on $\mathscr A$.

The construction of \S\ref{void-partial-frobenius} (\emph{cf.}~Remark \ref{rem-partial-frobenius-independence-of-partition}) thus furnishes a functor:
\begin{align}
\label{eq-cohomology-of-shtukas-partial-frobenius-equivariance}
	\Rep(({}^L H_{\mathring X})^I) &\rightarrow (\Ind \derived(\mathring X^I))^{F_1, \cdots, F_n} \\
	\notag
	V & \mapsto \mathscr H_{I, V} \text{ with } \varphi_{\{1\}}, \cdots, \varphi_{\{n\}}.
\end{align}
Furthermore, \eqref{eq-cohomology-of-shtukas-partial-frobenius-equivariance} is natural in $I$ (\emph{cf.}~\cite[Proposition 4.12, Proposition 4.14]{MR3787407}). The isomorphisms $\varphi_{\{1\}}, \cdots, \varphi_{\{n\}}$ are the \emph{partial Frobenius-equivariance} structures on $\mathscr H_{I, V}$.
\end{void}

\begin{rem}
It follows immediately from the construction that the partial Frobenius-equivariance structures on $\mathscr H_{I, V}$ commute with Hecke operators (\emph{cf.}~\S\ref{void-hecke-operator-cohomology-of-shtuka-construction}).
\end{rem}

\begin{rem}
\label{rem-cohomology-of-shtukas-local-form}
Given a subset $J \subset I$ and a $J$-tuple of closed points $\{x_j\}_{j\in J}$ of $\mathring X$, we have the restriction functor
\begin{equation}
\label{eq-constructible-sheaves-curves-restrictions-to-closed-points}
	(\Ind \derived(\mathring X^I))^{F_1, \cdots, F_n} \rightarrow (\Ind \derived(\mathring X^{I\setminus J} \times x^J))^{F_1, \cdots, F_n},
\end{equation}
where the target denotes the category of objects of $\Ind \derived(\mathring X^{I\setminus J} \times x^J)$ equipped with partial Frobenius-equivariance defined as in \S\ref{void-constructible-complex-partial-frobenius-equivariance}. The composition of \eqref{eq-cohomology-of-shtukas-partial-frobenius-equivariance} and \eqref{eq-constructible-sheaves-curves-restrictions-to-closed-points} canonically factors through a functor
\begin{equation}
\label{eq-cohomology-of-shtukas-local-form}
	\Rep(({}^LH_{\mathring X})^{I\setminus J} \times \prod_{j\in J} {}^LH_{x_j}) \rightarrow (\Ind \derived(\mathring X^{I\setminus J} \times x^J))^{F_1, \cdots, F_n},
\end{equation}
where each ${}^LH_{x_j}$ is the unramified local $L$-group at $x_j$ (\emph{cf.}~\S\ref{void-local-ell-group}).

To see this, we repeat the construction of \eqref{eq-cohomology-of-shtukas-partial-frobenius-equivariance}, with the paws of the Shtukas restricted to $\mathring X^{I \setminus J} \times x^J$, and use the local form of the Satake equivalence (\emph{cf.}~\S\ref{void-local-ell-group}).

Note that the family of functors \eqref{eq-cohomology-of-shtukas-local-form} is natural with respect to maps of inclusions of finite sets $(J_1 \subset I_1) \rightarrow (J_2 \subset I_2)$.
\end{rem}

\subsection{Excursion operators}

\begin{void}
\label{void-verlinde-loop-operator-context}
In this subsection, we summarize Lafforgue's construction of the spectral decomposition of cusp forms, taking as input the family of functors \eqref{eq-cohomology-of-shtukas-partial-frobenius-equivariance}. This corresponds to the material \cite[\S6-11 and \S12.3]{MR3787407}.

We shall start with the construction of the operators $S_{V, x}$ (\emph{cf.}~\cite[\S12.3.3]{MR3787407}). Namely, for any finite set $I$ and $W \in \Rep(({}^LH_{\mathring X})^I)$, we shall construct an endomorphism on the cohomology of Shtukas
\begin{equation}
\label{eq-verlinde-loop-operator}
	S_{V, x} : \mathscr H_{I, W} \rightarrow \mathscr H_{I, W}
\end{equation}
in $\Ind \derived(\mathring X^I)$, associated to any closed point $x\in \mathring X$ and any $V \in \Rep({}^LH_x)$, where ${}^LH_x$ is the unramified local $L$-group at $x$ (\emph{cf.}~\S\ref{void-local-ell-group}).
\end{void}

\begin{void}[Construction of \eqref{eq-verlinde-loop-operator}]
By naturality of \eqref{eq-cohomology-of-shtukas-partial-frobenius-equivariance} with respect to the inclusion of finite sets $I \rightarrow I \sqcup \{1\}$, we obtain an isomorphism (``insertion of vacuum"):
\begin{equation}
\label{eq-insertion-of-vacuum}
\mathscr H_{I, W} |_{\mathring X^I \times \mathring X} \simeq \mathscr H_{I \sqcup \{1\}, W \boxtimes \mathbf 1}
\end{equation}
where $\mathbf 1$ denotes the trivial representation of ${}^LH_{\mathring X}$.

Let us restrict \eqref{eq-insertion-of-vacuum} along the inclusion $\mathring X^I \times x \rightarrow \mathring X^I \times\mathring X$. We may express $\mathscr H_{I \sqcup\{1\}, W \boxtimes \mathbf 1} |_{\mathring X^I \times x}$ as the image of $W \boxtimes \mathbf 1 \in \Rep(({}^LH_{\mathring X})^I \times {}^LH_x)$ under \eqref{eq-cohomology-of-shtukas-local-form}. Applying the naturality of the latter with respect to the map of inclusions of finite sets
$$
\begin{tikzcd}[column sep = 1.5em]
\{1, 2\} \ar[r, phantom, "\subset"]\ar[d, "f"] & I \sqcup \{1, 2\} \ar[d, "\id_I \sqcup f"] \\
\{1\} \ar[r, phantom, "\subset"] & I \sqcup \{1\}
\end{tikzcd}
$$
with $f(1) = f(2) = 1$, we may identify $\mathscr H_{I \sqcup\{1\}, W \boxtimes \mathbf 1} |_{\mathring X^I \times x}$ with $\mathscr H_{I \sqcup \{1, 2\}, W \boxtimes \mathbf 1} |_{\mathring X^I \times \Delta(x)}$, where $\mathbf 1$ now means the trivial representation of $({}^LH_x)^{\{1, 2\}}$. In other words, the restriction of \eqref{eq-insertion-of-vacuum} to $\mathring X^I \times x$ yields an isomorphism
\begin{equation}
\label{eq-insertion-of-vacuum-diagonal}
\mathscr H_{I, W} |_{\mathring X^I \times x} \simeq \mathscr H_{I \sqcup \{1, 2\}, W \boxtimes \mathbf 1}|_{\mathring X^I \times \Delta(x)}.
\end{equation}

Consider now the composition
\begin{align}
\notag
	\mathscr H_{I, W} |_{\mathring X^I \times x} &\simeq \mathscr H_{I \sqcup \{1, 2\}, W \boxtimes \mathbf 1}|_{\mathring X^I \times \Delta(x)} \xrightarrow{\unit} \mathscr H_{I \sqcup \{1, 2\}, W \boxtimes V \boxtimes V^*}|_{\mathring X^I \times \Delta(x)} \\
\label{eq-verlinde-loop-oprator-construction}
	& \xrightarrow{\varphi_{\{1\}}^{\deg x}} \mathscr H_{I \sqcup \{1, 2\}, W \boxtimes V \boxtimes V^*}|_{\mathring X^I \times \Delta(x)} \xrightarrow{\counit} \mathscr H_{I \sqcup \{1, 2\}, W \boxtimes \mathbf 1}|_{\mathring X^I \times \Delta(x)} \simeq \mathscr H_{I, W} |_{\mathring X^I \times x},
\end{align}
where $\unit$, $\counit$ denote the morphisms induced by functoriality of \eqref{eq-cohomology-of-shtukas-local-form} with respect to the unit $\mathbf 1 \rightarrow V\boxtimes V^*$, respectively counit $V \boxtimes V^* \rightarrow \mathbf 1$, the first and last isomorphisms are \eqref{eq-insertion-of-vacuum-diagonal}, and the endomorphism $\varphi_{\{1\}}^{\deg x}$ is the $(\deg x)$-fold iteration of the partial Frobenius equivariance structure along the factor corresponding to $\{1\}$. (Here, $\deg x := [\base_x : \base]$, where $\base_x$ is the residue field of $x$.)

The composition \eqref{eq-verlinde-loop-oprator-construction} respects partial Frobenius equivariance along the $x$-factor, so it descends to an endomorphism of $\mathscr H_{I, W}$, giving the desired operator $S_{V, x}$.
\end{void}

\begin{rem}
\label{rem-verlinde-loop-operator-local-form}
In the construction of \eqref{eq-verlinde-loop-operator}, we may also replace $W$ by a representation of a product of global and local $L$-groups, as follows: Given a subset $J \subset I$, a $J$-tuple of closed points $\{x_j\}_{j\in J}$ of $\mathring X$, and $W \in \Rep(({}^LH_{\mathring X})^{I\setminus J} \times \prod_{j \in J}{}^LH_{x_j})$, we have an endomorphism
\begin{equation}
\label{eq-verlinde-loop-operator-local-form}
S_{V, x} : \mathscr H_{I, W} \rightarrow \mathscr H_{I, W}
\end{equation}
in $\Ind \derived(\mathring X^{I \setminus J} \times x^J)$, associated to any $x \in \mathring X$ and $V \in \Rep({}^LH_x)$.

The construction of \eqref{eq-verlinde-loop-operator-local-form} is completely parallel to that of \eqref{eq-verlinde-loop-operator}, except that we use the local form of the cohomology of Shtukas \eqref{eq-cohomology-of-shtukas-local-form} also along $x^J$.
\end{rem}

\begin{void}
As in \cite[Proposition 6.2]{MR3787407}, the endomorphism $S_{V, x}$ extends the Hecke operator $h_{V, x}$ on the restriction of $\mathscr H_{I, W}$ to $(\mathring X \setminus x)^I$ (\emph{cf.}~\S\ref{void-hecke-operators}, \S\ref{void-hecke-operator-cohomology-of-shtuka-construction}).
\end{void}

\begin{prop}[``$S = T$"]
\label{prop-s=t}
Given a finite set $I$ and $W \in \Rep(({}^LH_{\mathring X})^I)$, as well as a closed point $x \in \mathring X$ and $V \in \Rep({}^LH_x)$, there is a canonical isomorphism
\begin{equation}
\label{eq-s=t}
S_{V, x} \simeq T(h_{V, x})
\end{equation}
of endomorphisms of $\mathscr H_{I, W} |_{(\mathring X \setminus x)^I}$.
\end{prop}

\begin{void}[Cohomological correspondences]
\label{void-cohomological-correspondences}
In the untwisted setting (\emph{cf.}~\cite[\S 6.3]{MR3787407}), the isomorphism \eqref{eq-s=t} is obtained from an identification of ``cohomological correspondences" from the coefficient sheaf $\mathscr F_{I, W}$ over $\Sht_{G, D}^I$ to itself (\emph{cf.}~\S\ref{void-cohomology-of-shtukas-construction}). The twisted setting is similar, except that we need an appropriate notion of ``cohomological correspondences."

Let $(\mathscr X_1, \mathscr G_1)$, $(\mathscr X_2, \mathscr G_2)$ be algebraic $\base$-stacks equipped with $A$-gerbes. Given objects $\mathscr A_1 \in \derived_{\mathscr G_1, \zeta}(\mathscr X_1)$, $\mathscr A_2 \in \derived_{\mathscr G_2, \zeta}(\mathscr X_2)$, a \emph{cohomological correspondence}
$$
(\mathscr X_1, \mathscr G_1, \mathscr A_1) \rightarrow (\mathscr X_2, \mathscr G_2, \mathscr A_2)
$$
consists of a span of the following data
\begin{enumerate}
	\item a span of algebraic $\base$-stacks $\mathscr X_1 \xleftarrow{\overset{\leftarrow}{\pi}} \mathscr X \xrightarrow{\overset{\rightarrow}{\pi}} \mathscr X_2$, where $\overset{\rightarrow}{\pi}$ is schematic, separated, and of finite type (so $\overset{\rightarrow}{\pi}{}^!$ is well-defined, \emph{cf.}~Remark \ref{rem-etale-local-notion-twisted-complex});
	\item an isomorphism of $A$-gerbes $\gamma : \overset{\leftarrow}{\pi}{}^*\mathscr G_1 \simeq \overset{\rightarrow}{\pi}{}^*\mathscr G_2$;
	\item a morphism in $\derived_{\mathscr G, \zeta}(\mathscr X)$, where $\mathscr G$ is the $A$-gerbe $\overset{\leftarrow}{\pi}{}^*\mathscr G_1 \simeq \overset{\rightarrow}{\pi}{}^*\mathscr G_2$ identified via $\gamma$:
	$$
	\overset{\leftarrow}{\pi}{}^* \mathscr A_1 \rightarrow \overset{\rightarrow}{\pi}{}^! \mathscr A_2.
	$$
\end{enumerate}

Compositions of cohomological correspondences are defined as in \cite[\S4.1]{MR3787407}, the only modification being that one composes the pullbacks of the isomorphisms of the $A$-gerbes over the fiber product of algebraic $\base$-stacks.
\end{void}

\begin{eg}
Let $Y$ be a separated $\base$-scheme of finite type equipped with an $A$-gerbe $\mathscr G$. For any $\mathscr A \in \derived_{\mathscr G, \zeta}(Y)$, we have composable cohomological correspondences
\begin{align}
\notag
(\Spec \base, \mathbf 0, \overline{\rationals}{}_{\ell}) & \rightarrow (Y \times Y, \mathscr G \boxplus (-\mathscr G), \mathscr A \boxtimes \mathbf D\mathscr A)\\
\label{eq-diagonal-lowering-correspondence}
	& \rightarrow (Y \times Y, \mathscr G \boxplus (-\mathscr G), \mathscr A \boxtimes \mathbf D\mathscr A) \rightarrow (\Spec \base, \mathbf 0, \overline{\rationals}{}_{\ell}).
\end{align}

Here, $\mathbf 0$ denotes the trivial $A$-gerbe over $\Spec \base$ and $\mathbf D$ is the Verdier duality over $Y$ relative to $\base$ (\emph{cf.}~\S\ref{void-verdier-duality}). The spans supporting them are represented in the following diagram
\begin{equation}
\label{eq-frobenius-fixed-point-as-composition-of-spans}
\begin{tikzcd}[column sep = 1em]
 & Y\ar[dl]\ar[dr, "\Delta"] & & Y \times Y \ar[dl, swap, "{(\Fr_Y, \id)}"] \ar[dr, equal] & & Y \ar[dl, swap, "\Delta"]\ar[dr]  \\
 \Spec \base & & Y \times Y & & Y \times Y & & \Spec \base
\end{tikzcd}
\end{equation}
The identifications of $A$-gerbes are given by the canonical trivialization of $\Delta^*(\mathscr G \boxplus(-\mathscr G))$, as well as the isomorphism $(\Fr_Y, \id)^*(\mathscr G \boxplus (-\mathscr G)) \simeq \mathscr G \boxplus (-\mathscr G)$ supplied by \eqref{eq-tautological-isomorphism-fixed-point}. The morphisms on constructible complexes are given by the unit of Verdier duality, the identity on $\mathscr A \boxtimes \mathbf D \mathscr A$, and the co-unit of Verdier duality, respectively.

Composing the spans \eqref{eq-frobenius-fixed-point-as-composition-of-spans}, we obtain the span $Y^{\Fr}$ from $\Spec \base$ to itself. The composition of isomorphisms of $A$-gerbes yields an automorphism of the trivial $A$-gerbe over $Y^{\Fr}$: By construction, \emph{this is the $A$-torsor $-\Tr(\Fr \mid \mathscr G)$} (\emph{cf.}~\S\ref{void-trace-of-frobenius-gerbe-construction}). With respect to this automorphism, the composition of the morphisms of constructible complexes yields a $\zeta$-genuine function on $\Tr(\Fr \mid \mathscr G)(\base)$, which equals $\Tr(\Fr \mid \mathscr A)$ (\emph{cf.}~\S\ref{void-constructing-genuine-functions}).
\end{eg}

\begin{void}
Let us return to the context of Proposition \ref{prop-s=t}. Our goal is to reformulate the isomorphism \eqref{eq-s=t} in terms of cohomological correspondences.

We begin with the endomorphism $T(h_{V, x})$. Recall the stack $\Hec_x(\Sht_{G, D}^I)$ of Hecke modifcations of $\Sht_{G, D}^I$ at $x$ (\emph{cf.}~\S\ref{void-hecke-operator-cohomology-of-shtuka-construction}). Viewed as a span, $\Hec_x(\Sht_{G, D}^I)$ supports a cohomological correspondence
$$
(\Sht_{G, D}^I, \mathbf 0, \mathscr F_{I, W}) \rightarrow (\Sht_{G, D}^I, \mathbf 0, \mathscr F_{I, W}),
$$
where the isomorphism of the $A$-gerbes is given by $r^*(-\Tr(\Fr \mid \mathscr G_{\Hec}))$ and the morphism of constructible complexes is given by \eqref{eq-hecke-correspondence-shtuka-coefficients}, which we may rewrite as
\begin{equation}
\label{eq-hecke-operation-as-cohomological-correspondence}
p^* \mathscr F_{I, W} \otimes r^*(\mathscr L_{\Tr(\Fr \mid \mathscr G_{\Hec_G}), \zeta})^{\otimes -1} \rightarrow q^! \mathscr F_{I, W}.
\end{equation}

The element $h_{V, x}$ may be viewed as a global section of $(\mathscr L_{\Tr(\Fr \mid \mathscr G_{\Hec_G}), \zeta})^{\otimes -1}$ over $\Hec_{G, x}(\base)$, so by evaluating \eqref{eq-hecke-operation-as-cohomological-correspondence} on $h_{V, x}$ along its second factor, we obtain a morphism
\begin{equation}
\label{eq-hecke-operator-as-cohomological-correspondence}
\widetilde T(h_{V, x}) : p^* \mathscr F_{I, W} \rightarrow q^! \mathscr F_{I, W}.
\end{equation}

The morphism \eqref{eq-hecke-operator-as-cohomological-correspondence} gives rise to $T(h_{V, x})$ by adjunction and compactly supported cohomology along $\nu_{\Xi}^{I_1, \cdots, I_n}$.
\end{void}

\begin{void}
As for $S_{V, x}$, we shall consider the spans defined in \cite[\S 6.3]{MR3787407}:
\begin{equation}
\label{eq-verlinde-operator-composition-of-spans}
\begin{tikzcd}[column sep = -0.5em]
	& \mathscr Y_{\sharp} \ar[dl, swap, "p_{\sharp}"] \ar[dr, hookrightarrow, "i_{\sharp}"] & & \Sht_{G, D}^{\{1\}, \{2\}, I} |_{\mathring X^I \times \Delta(x)} \ar[dl, swap, "F_{\{1\}}^{\deg x}"] \ar[dr, equal] & & \mathscr Y_{\flat} \ar[dl, hookrightarrow, swap, "i_{\flat}"]\ar[dr, "p_{\flat}"]  \\
	\Sht_{G, D}^I \times x & & \Sht_{G, D}^{I, \{2\}, \{1\}} |_{\mathring X^I \times \Delta(x)} & & \Sht_{G, D}^{\{1\}, \{2\}, I} |_{\mathring X^I \times \Delta(x)} & & \Sht_{G, D}^I \times x
\end{tikzcd}
\end{equation}

Here, $\mathscr Y_{\sharp}$ is the closed substack of $\Sht_{G, D}^{I, \{2\}, \{1\}} |_{\mathring X^I \times \Delta(x)}$ consisting of $R$-points
$$
(P^0, \phi^0) \overset{x^I}{\sim} (P^1, \phi^1) \overset{x}{\sim} (P^2, \phi^2) \overset{x}{\sim} ({}^{\tau}P^0, {}^{\tau}\phi^0)
$$
where the composition of the two last modifications extends to an isomorphism $(P^1, \phi^1) \simeq ({}^{\tau}P^0, {}^{\tau}\phi^0)$ of $R$-points of $\Bun_{G, D}$. The closed substack $\mathscr Y_{\flat}$ of $\Sht_{G, D}^{\{1\}, \{2\}, I}|_{\mathring X^I \times \Delta(x)}$ is similarly defined. The morphism $F_{\{1\}}^{\deg x}$ is the $(\deg x)$-fold iteration of the lift of the partial Frobenius endomorphism along $\{1\}$ (\emph{cf.}~\S\ref{void-shtuka-partial-frobenius}).\footnote{The formation of $F_{\{1\}}^{\deg x}$ requires a re-ordering of $\{1\}, \{2\}, I$; see \cite[\S 6.3(c)]{MR3787407} for details. }

According to \cite[Lemme 6.10]{MR3787407}, the composition of the spans \eqref{eq-verlinde-operator-composition-of-spans} is canonically isomorphic to $\Hec_x(\Sht_{G, D}^I) \times x$, \emph{i.e.}~we have a canonical isomorphism
\begin{equation}
\label{eq-s=t-support-identification}
	\mathscr Y_{\sharp} \underset{\Sht_{G, D}^{I, \{2\}, \{1\}}|_{\mathring X^I \times \Delta x}}{\times} \Sht_{G, D}^{\{1\}, \{2\}, I}|_{\mathring X^I \times \Delta(x)} \underset{\Sht_{G, D}^{\{1\}, \{2\}, I}|_{\mathring X^I \times \Delta(x)}}{\times} \mathscr Y_{\flat} \simeq \Hec_x(\Sht_{G, D}^I) \times x,
\end{equation}
where the structural morphisms $p$, $q$ correspond to the projection onto $\mathscr Y_{\sharp}$ followed by $p_{\sharp}$, respectively the projection onto $\mathscr Y_{\flat}$, followed by $p_{\flat}$.
\end{void}

\begin{void}
Let us now invoke the $A$-gerbes $\mathscr G_{\Sht_{G, D}}$ over the various moduli stacks of Shtukas, temporarily ignoring their canonical trivializations (\emph{cf.}~\S\ref{void-shtuka-gerbe-trivialization}). Note that we have a canonical identification of $A$-gerbes
\begin{equation}
\label{eq-gerbe-shtuka-unit-identification}
(p_{\sharp})^* \mathscr G_{\Sht_{G, D}^I} \simeq (i_{\sharp})^* \mathscr G_{\Sht_{G, D}^{I, \{2\}, \{1\}}}
\end{equation}
arising from the multiplicative structure on $\mathscr G_{\Hec_G}$ (\emph{cf.}~Lemma \ref{lem-local-hecke-stack-gerbe-multiplicative}).

We shall extend \eqref{eq-gerbe-shtuka-unit-identification} to a cohomological correspondence supported on $\mathscr Y_{\sharp}$:
\begin{equation}
\label{eq-unit-as-cohomological-correspondence}
(\Sht_{G, D}^I \times x, \mathscr G_{\Sht_{G, D}} \boxplus \mathbf 0, \mathscr F_{I, W}) \rightarrow (\Sht_{G, D}^{I, \{2\}, \{1\}}|_{\mathring X^I \times \Delta(x)}, \mathscr G_{\Sht_{G, D}}, \mathscr F_{I, \{2\}, \{1\}, W \boxtimes V^* \boxtimes V}).
\end{equation}
Indeed, what remains to be constructed is a morphism of constructible complexes over $\mathscr Y_{\sharp}$ twisted by the $A$-gerbes \eqref{eq-gerbe-shtuka-unit-identification}:
\begin{equation}
\label{eq-unit-as-cohomological-correspondence-morphism-of-constructible-complexes}
(p_{\sharp})^*\mathscr F_{I, W} \rightarrow (i_{\sharp})^! \mathscr F_{I, \{2\}, \{1\}, W \boxtimes V^* \boxtimes V}.
\end{equation}

The morphism \eqref{eq-unit-as-cohomological-correspondence-morphism-of-constructible-complexes} is the pullback of the unit morphism in the Satake category (\emph{cf.}~Proposition \ref{prop-convolution-monoidal-dual}). To be more precise, we consider the Cartesian square
\begin{equation}
\label{eq-unit-cohomological-correspondence-cartesian-diagram}
\begin{tikzcd}
	\mathscr Y_{\sharp} \ar[r, hookrightarrow, "i_{\sharp}"]\ar[d, "{(r^I, r^{\{2\}})}"] & \Sht_{G, D}^{I, \{2\}, \{1\}}|_{\mathring X^I \times \Delta(x)} \ar[d, "{(r^I, r^{\{2\},\{1\}})}"] \\
	\Hec_{G, I} \times \Hec_{G, x} \ar[r, hookrightarrow, "{(\id, \delta)}"] & \Hec_{G, I} \times \Hec_{G, x}^{\{2\}, \{1\}}
\end{tikzcd}
\end{equation}
where the vertical arrows are the products of restriction maps (\emph{cf.}~\S\ref{void-shtuka-coefficient-construction}) and $\delta$ sends $P^0 \overset{x}{\sim} P^1$ to the concatenation with its own inverse $P^0 \overset{x}{\sim} P^1 \overset{x}{\sim} P^0$. The unit morphism of the Satake category \eqref{eq-convolution-dual-unit} (applied to the objects $\mathscr S_V, \mathscr S_{V^*} \in \Sat_{\mathscr G, \zeta}(\Hec_{G, x})$ corresponding to $V$ and $V^*$) yields a morphism in $\derived(\Hec_{G, x})$:
\begin{equation}
\label{eq-unit-satake-category-as-cohomological-correspondence}
\overline{\rationals}_{\ell} \rightarrow \delta^!(p_{\{2\}}^* \mathscr S_{V^*} \otimes p_{\{1\}}^* \mathscr S_V),
\end{equation}
where $p_{\{2\}}$, $p_{\{1\}}$ are the projections onto the first, respectively the second modification. The desired morphism \eqref{eq-unit-as-cohomological-correspondence-morphism-of-constructible-complexes} is the pullback along $(r^I, r^{\{2\}})$ of the external tensor product of \eqref{eq-unit-satake-category-as-cohomological-correspondence} with the identity on $\mathscr S_{I, W} \in \Sat_{\mathscr G, \zeta}(\Hec_{G, I})$. Here, we use the base change property for $*$- and $!$-pullbacks along \eqref{eq-unit-cohomological-correspondence-cartesian-diagram}, which holds by smoothness after passing to an ind-pro-presentation of the local Hecke stacks (\emph{cf.}~Remark \ref{rem-shtuka-coefficients-perverse-ula}).

Let us now invoke the canonical trivialization of $\mathscr G_{\Sht_{G, D}}$ (\emph{cf.}~\S\ref{void-shtuka-gerbe-trivialization}). It is compatible with \eqref{eq-gerbe-shtuka-unit-identification} in the following sense: After trivializing both $\mathscr G_{\Sht_{G, D}^I}$ and $\mathscr G_{\Sht_{G, D}^{I, \{2\}, \{1\}}}$, the isomorphism \eqref{eq-gerbe-shtuka-unit-identification} becomes the \emph{identity} automorphism of the trivial $A$-gerbe. In particular, we may view \eqref{eq-unit-as-cohomological-correspondence-morphism-of-constructible-complexes} as a morphism of \emph{untwisted} constructible complexes over $\mathscr Y_{\sharp}$, and thus \eqref{eq-unit-as-cohomological-correspondence} as a cohomological correspondence in the usual sense:
\begin{equation}
\label{eq-creation-cohomological-correspondence}
(\Sht_{G, D}^I\times x, \mathbf 0, \mathscr F_{I, W}) \rightarrow (\Sht_{G, D}^{I, \{2\}, \{1\}}|_{\mathring X^I \times \Delta(x)}, \mathbf 0, \mathscr F_{I, \{2\}, \{1\}, W \boxtimes V^* \boxtimes V}),
\end{equation}
which gives rise to the morphism labeled ``unit" in \eqref{eq-verlinde-loop-oprator-construction}.

Dually, we have the cohomological correspondence supported on $\mathscr Y_{\flat}$:
\begin{equation}
\label{eq-annihilation-cohomological-correspondence}
(\Sht_{G, D}^{\{1\}, \{2\}, I}|_{\mathring X^I \times \Delta(x)}, \mathbf 0, \mathscr F_{\{1\}, \{2\}, I, V \boxtimes V^* \boxtimes W}) \rightarrow (\Sht_{G, D}^I \times x, \mathbf 0, \mathscr F_{I, W}),
\end{equation}
which gives rise to teh morphism ``counit" in \eqref{eq-verlinde-loop-oprator-construction}.
\end{void}

\begin{void}
Finally, we have the cohomological correspondence from the target of \eqref{eq-creation-cohomological-correspondence} to the source of \eqref{eq-annihilation-cohomological-correspondence}, supported on $\Sht_{G, D}^{\{1\}, \{2\}, I}|_{\mathring X^I \times \Delta(x)}$:
\begin{equation}
\label{eq-frobenius-equivariance-cohomological-correspondence}
(\Sht_{G, D}^{I, \{2\}, \{1\}}|_{\mathring X^I \times \Delta(x)}, \mathbf 0, \mathscr F_{I, \{2\}, \{1\}, W \boxtimes V^* \boxtimes V}) \rightarrow (\Sht_{G, D}^{\{1\}, \{2\}, I}|_{\mathring X^I \times \Delta(x)}, \mathbf 0, \mathscr F_{\{1\}, \{2\}, I, V \boxtimes V^* \boxtimes W})
\end{equation}
defined by the $F_{\{1\}}^*$-equivariance of $\mathscr F_{\{1\}, \{2\}, I, V \boxtimes V^* \boxtimes W}$ (\emph{cf.}~\S\ref{void-shtuka-partial-frobenius}).

The composition of the three cohomological correspondences \eqref{eq-creation-cohomological-correspondence}, \eqref{eq-frobenius-equivariance-cohomological-correspondence}, \eqref{eq-annihilation-cohomological-correspondence} then corresponds to a morphism of (untwisted) constructible complexes over $\Hec_x(\Sht_{G, D}^I)$:
\begin{equation}
\label{eq-verlinde-operator-cohomological-correspondence}
	\widetilde S_{V, x} : p^* \mathscr F_{I, W} \rightarrow q^! \mathscr F_{I, W},
\end{equation}
which gives rise to $S_{V, x}$ by adjunction and compactly supported cohomology.
\end{void}

\begin{rem}
\label{rem-s=t-gerbe-identifications}
Let us point out an aspect of the construction of \eqref{eq-verlinde-operator-cohomological-correspondence} which may cause confusion. Indeed, \eqref{eq-verlinde-operator-cohomological-correspondence} invokes the composition of three identifications of $A$-gerbes, after pulling back to $\Hec_x(\Sht_{G, D}^I) \times x$ along \eqref{eq-s=t-support-identification}:
\begin{align*}
	(p_{\sharp})^*\mathscr G_{\Sht_{G, D}^I} & \simeq (i_{\sharp})^* \mathscr G_{\Sht_{G, D}^{I, \{2\}, \{1\}}}, \\
	(F_{\{1\}}^{\deg x})^* \mathscr G_{\Sht_{G, D}^{I, \{2\}, \{1\}}} & \simeq \mathscr G_{\Sht_{G, D}^{\{1\}, \{2\}, I}}, \\
	(i_{\flat})^*\mathscr G_{\Sht_{G, D}^{\{1\}, \{2\}, I}} & \simeq (p_{\flat})^* \mathscr G_{\Sht_{G, D}^I},
\end{align*}
where the first isomorphism is \eqref{eq-gerbe-shtuka-unit-identification}, the third isomorphism is similarly defined, and the second isomorphism comes from the $(F_{\{1\}}^{\deg x})^*$-equivariance (\emph{cf.}~\S\ref{void-shtuka-coefficient-partial-frobenius-equivariance-construction}).

Their composition yields an isomorphism of $A$-gerbes
\begin{equation}
\label{eq-verlinde-operator-composition-gerbe-isomorphism}
p^* \mathscr G_{\Sht_{G, D}^I} \simeq q^* \mathscr G_{\Sht_{G, D}^I}.
\end{equation}
This isomorphism is \emph{different} from the canonical isomorphism \eqref{eq-gerbe-shtuka-hecke-isomorphism} defined over $(\mathring X \setminus x)^I \times x$. Indeed, under the canonical trivialization of $\mathscr G_{\Sht_{G, D}^I}$ (\emph{cf.}~\S\ref{void-shtuka-gerbe-trivialization}), \eqref{eq-verlinde-operator-composition-gerbe-isomorphism} corresponds to the identity automorphism of the trivial $A$-gerbe while \eqref{eq-gerbe-shtuka-hecke-isomorphism} corresponds to the multiplication by the pullback of $-\Tr(\Fr \mid \mathscr G_{\Hec_{G, x}})$ (\emph{cf.}~Lemma \ref{lem-gerbe-shtuka-hecke-compatibility}).
\end{rem}

\begin{void}
With the above constructions in place, we may prove Proposition \ref{prop-s=t} by repeating the proof of \cite[Proposition 6.2]{MR3787407}.
\begin{proof}[Proof of Proposition \ref{prop-s=t}]
It suffices to identify the morphisms of (untwisted) constructible complexes \eqref{eq-hecke-operator-as-cohomological-correspondence}, \eqref{eq-verlinde-operator-cohomological-correspondence} over $\Hec_x(\Sht_{G, D}^I)$:
\begin{equation}
\label{eq-s=t-as-isomorphism-of-correspondences}
\widetilde S_{V, x} \simeq \widetilde T(h_{V, x}),
\end{equation}
as it will give rise to \eqref{eq-s=t} by taking compactly supported cohomology.

The isomorphism \eqref{eq-s=t-as-isomorphism-of-correspondences} is our version of \cite[Lemme 6.11]{MR3787407}, and is proved in the same way.\footnote{However, for the reason mentioned in Remark \ref{rem-s=t-gerbe-identifications}, it is not formulated as an identification of cohomological correspondences.} Namely, we may reduce this assertion to the ``vacuum case" $I = \emptyset$ and $W = \mathbf 1$ by rewriting the construction of \eqref{eq-verlinde-operator-cohomological-correspondence} in terms local to $x$ (\emph{cf.}~\cite[\S 6.4]{MR3787407}). For the vacuum case, we argue as in the proof of \cite[Lemme 6.13]{MR3787407}, reducing to the fact that the morphism \eqref{eq-unit-as-cohomological-correspondence-morphism-of-constructible-complexes}, which here is a morphism
$$
\overline{\rationals}_{\ell} \rightarrow (i_{\sharp})^! \mathscr F_{\{2\}, \{1\}, V^*\boxtimes V}
$$
of (untwisted) constructible complexes over $\mathscr Y^{\sharp}$, is uniquely determined by its restriction to the smooth locus of $\mathscr Y^{\sharp}$ (\emph{cf.}~\cite[Lemme 6.15]{MR3787407}).
\end{proof}
\end{void}

\begin{void}
Next, we turn to the Eichler--Shimura relation for the partial Frobenius-equivariance structure on the cohomology of Shtukas.

More precisely, let $I$, $W$ be as in \S\ref{void-verlinde-loop-operator-context}, $x$ be a closed point of $\mathring X$, and $V \in \Rep({}^LH_x)$. Consider the object $\mathscr H_{I\sqcup\{0\}, W \boxtimes V}$ of $\Ind \derived(\mathring X^I \times x)$ equipped with its partial Frobenius-equivariance structure $\varphi_{\{0\}}$ along the $x$-factors (\emph{cf.}~Remark \S\ref{rem-cohomology-of-shtukas-local-form}). Since $F_{\{0\}}^{\deg x}$ is the identity on $\mathring X^I \times x$, we may view $\varphi_{\{0\}}^{\deg x}$ as an endomorphism
\begin{equation}
\label{eq-partial-frobenius-endomorphism-at-extra-point}
\varphi_{\{0\}}^{\deg x} : \mathscr H_{I \sqcup \{0\}, W \boxtimes V} \rightarrow \mathscr H_{I \sqcup \{0\}, W \boxtimes V}.
\end{equation}

The Eicher--Shimura relation is a polynomial equation satisfied by \eqref{eq-partial-frobenius-endomorphism-at-extra-point} with coefficients in the operators $S_{\bigwedge^{\dim V - k}V, x}$ (in its local form \eqref{eq-verlinde-loop-operator-local-form} for the inclusion $\{0\} \subset I \sqcup\{0\}$).
\end{void}

\begin{prop}[The Eichler--Shimura relation]
\label{prop-eichler-shimura-relation}
There is an isomorphism of endomorphisms of $\mathscr H_{I \sqcup \{0\}, W \boxtimes V}$:
\begin{equation}
\label{eq-eichler-shimura-relation}
	\sum_{k = 0}^{\dim V} (-1)^k S_{\bigwedge^{\dim V - k} V, x} \circ (\varphi_{\{0\}}^{\deg x})^k \simeq 0.
\end{equation}
\end{prop}

\begin{proof}
This follows formally from the proof of \cite[Proposition 7.1]{MR3787407}.
\end{proof}

\begin{void}[Hecke-finiteness]
Next, we shall focus on the degree-$0$ part of the cohomology of Shtukas, \emph{i.e.}~the ind-constructible sheaf
$$
\mathscr H_{I, V}^0 := H^0\mathscr H_{I, V} \in \Ind \derived(\mathring X^I)^{\heartsuit}
$$
naturally associated to any finite set $I$ and $V \in \mathscr \Rep({}^L H_{\mathring X})$. Recall that $\mathscr H^0_{\emptyset, \mathbf 1}$ is canonically identified with the vector space of $\zeta$-genuine automorphic forms (\emph{cf.}~\S\ref{void-shtuka-no-leg-automorphic-forms}).

Let $\bar {\eta}$ be a geometric point of $\mathring X^I$ lying over the generic point of each factor. We shall call an element of $\mathscr H_{I, V}^0 |_{\bar {\eta}}$ \emph{Hecke-finite} if it belongs to a finite-dimensional $\overline{\rationals}_{\ell}$-vector space which is stable under actions of the spherical Hecke algebra at any $x \in \mathring X$ (\emph{cf.}~\S\ref{void-hecke-operator-cohomology-of-shtuka-construction}). Write
$$
(\mathscr H^0_{I, V}|_{\bar{\eta}})^{\Hf} \subset \mathscr H^0_{I, V}|_{\bar{\eta}}
$$
for the space of Hecke-finite elements.

As in \cite[Proposition 8.23]{MR3787407}, we identify Hecke-finite elements of $\mathscr H^0_{\emptyset, \mathbf 1}$ with cuspidal $\zeta$-genuine automorphic forms.
\end{void}

\begin{prop}
\label{prop-hecke-finite-cuspidal}
An element of $\mathscr H^0_{\emptyset, \mathbf 1}$ is Hecke-finite if and only if the corresponding $\zeta$-genuine automorphic form is cuspidal.
\end{prop}

\begin{proof}
The fact that cusp forms are Hecke-finite follows from Lemma \ref{lem-cusp-form-finite-dimensional}. For the opposite containment, one applies the proof of \cite[Lemma 8.25]{MR3787407} to each summand in \eqref{eq-weil-uniformization-with-inner-forms-torsors}.
\end{proof}

\begin{void}
\label{void-hecke-finite-specialization}
For any finite set $I$, we write $\eta_I$ for the generic point of $X^I$. Choose a geometric point $\bar{\eta}_I$ lying over $\eta_I$, together with a lift of $\bar{\eta}_I \rightarrow X^I$ to the Henselian local ring of $X^I$ along $\Delta(\bar{\eta})$ (called a \emph{specialization map} $\specialization : \bar{\eta}_I \rightsquigarrow \Delta(\bar{\eta})$).

As in \cite[Proposition 8.27]{MR3787407} (whose proof uses the ``$S = T$" identity and the Eichler--Shimura relation, for which we have the substitutes Proposition \ref{prop-s=t}, Proposition \ref{prop-eichler-shimura-relation}), we have a presentation
$$
(\mathscr H^0_{I, V} |_{\bar{\eta}_I})^{\Hf} \simeq \bigcup_k M_k
$$
where each $M_k$ is a finite-dimensional $\overline{\rationals}_{\ell}$-vector space on which the $\pi_1^{\etale}(\eta_I, \bar{\eta}_I)$-action factors through $\pi_1^{\etale}(\eta, \bar{\eta})^I$ via the partial Frobenius-equivariance structure. Likewise, by \cite[Corollaire 8.34]{MR3787407}, the specialization map induces an isomorphism
\begin{equation}
\label{eq-hecke-finite-cohomology-specialization-isomorphism}
\specialization^* : (\mathscr H^0_{I, V} |_{\Delta(\bar{\eta})})^{\Hf} \simeq (\mathscr H^0_{I, V} |_{\bar{\eta}_I})^{\Hf}.
\end{equation}
\end{void}

\begin{rem}
\label{rem-xue-smoothness}
Xue's work \cite{xue2020smoothness} yields a stronger result: The ind-constructible sheaf $\mathscr H_{I, V}^k := H^k\mathscr H_{I, V}$ is ind-lisse for any $k \in \integers$. Indeed, this is established in the untwisted setting in \cite[Theorem 4.2.3]{xue2020smoothness}, but the argument is general: It uses the family of functors \eqref{eq-cohomology-of-shtukas-local-form}, the ``$S = T$" identity, and the Eichler--Shimura relation.

In particular, we may realize \eqref{eq-hecke-finite-cohomology-specialization-isomorphism} as the Hecke-finite part of the isomorphism
\begin{equation}
\label{eq-cohomology-of-shtukas-specialization-isomorphism}
\specialization^* : \mathscr H^0_{I, V}|_{\Delta(\bar{\eta})} \simeq \mathscr H^0_{I, V}|_{\bar{\eta}_I},
\end{equation}
although \eqref{eq-cohomology-of-shtukas-specialization-isomorphism} (for arbitrary degree $k \in \integers$) is established as a step in the proof that $\mathscr H^k_{I, V}$ is ind-lisse (\emph{cf.}~\cite[Proposition 1.4.3]{xue2020smoothness}).
\end{rem}

\begin{void}
\label{void-spectral-decomposition-proof}
We are now ready to construct the spectral decomposition of $\zeta$-genuine cusp forms \eqref{eq-spectral-decomposition} by repeating \cite[\S12.3.4]{MR3787407}.

\begin{proof}[Proof of Theorem \ref{thm-spectral-decomposition}]
For each finite set $I$ and $V \in \Rep(({}^LH_{\mathring X})^I)$, we denote by $H_{I, V}$ the left-hand-side of \eqref{eq-hecke-finite-cohomology-specialization-isomorphism}, equipped with the continuous $\Gal(\bar F/F)^I$-action induced from the right-hand-side. This defines a family of functors
\begin{align}
\label{eq-lafforgue-functors}
	\Rep(({}^LH_{\mathring X})^I) & \rightarrow \Rep(\Gal(\bar F/F)^I), \\
	\notag
	V & \mapsto H_{I, V},
\end{align}
natural in $I$, together with a canonical identification (\emph{cf.}~Proposition \ref{prop-hecke-finite-cuspidal})
\begin{equation}
\label{eq-lafforgue-functors-vacuum-cusp-forms}
H_{\emptyset, \mathbf 1} \simeq \Fun_{\cusp, \zeta}(\widetilde{\Bun}_{G, D}/\Xi, \overline{\rationals}_{\ell}).
\end{equation}

Note that the naturality of \eqref{eq-lafforgue-functors} in $I$ implies that $H_{I, \mathbf 1}$ is the vector space \eqref{eq-lafforgue-functors-vacuum-cusp-forms} equipped with the trivial $\Gal(\bar F/F)^I$-action.

Given $x \in V$ and $\xi \in V^*$ which are $H_{\bar{\eta}}$-invariant under the diagonal action, and an $I$-tuple $\gamma^I \in \Gal(\bar F/F)^I$, we obtain an endomorphism $S_{I, V, x, \xi, \gamma^I}$ of $H_{\{0\}, \mathbf 1}$ as the composition
\begin{align*}
	H_{\{0\}, \mathbf 1} \xrightarrow{x} H_{\{0\}, V^{H_{\bar{\eta}}}} \rightarrow H_{\{0\}, V} &\simeq H_{I, V} \\
	& \xrightarrow{\gamma^I} H_{I, V} \simeq H_{\{0\}, V} \rightarrow H_{\{0\}, V_{H_{\bar{\eta}}}} \xrightarrow{\xi} H_{\{0\}, \mathbf 1}.
\end{align*}
where the isomorphisms are induced from the naturality of \eqref{eq-lafforgue-functors} with respect to the constant map $I \rightarrow \{0\}$, and $V^{H_{\bar{\eta}}}$ (respectively $V_{H_{\bar{\eta}}}$) denotes the space of $H_{\bar{\eta}}$-invariants (respectively, $H_{\bar{\eta}}$-coinvariants) of $V$ with respect to the diagonal action.

Under the identification of $H_{\{0\}, \mathbf 1}$ with \eqref{eq-lafforgue-functors-vacuum-cusp-forms}, we may view each $S_{I, V, x, \xi, \gamma^I}$ as an endomorphisms of $\Fun_{\cusp, \zeta}(\widetilde{\Bun}_{G, D}/\Xi, \overline{\rationals}_{\ell})$. The proof of \cite[Lemme 10.6]{MR3787407} shows that $S_{I, V, x, \xi, \gamma^I}$ depends only on $I$, $\gamma^I$, and the function
\begin{align}
\label{eq-excursion-operator-biinvariant-function}
H_{\bar{\eta}} \backslash ({}^LH_{\mathring X})^I / H_{\bar{\eta}} & \rightarrow \overline{\rationals}_{\ell} \\
\notag
h^I & \mapsto \langle \xi, h^I\cdot x\rangle,
\end{align}
so one may define $S_{I, f, \gamma^I}$, for any locally constant $H_{\bar{\eta}}$-bi-invariant function $f$ on ${}^LH_{\mathring X}$, to be $S_{I, V, x, \xi, \gamma^I}$ for any $V$, $x$, $\xi$ whose induced function \eqref{eq-excursion-operator-biinvariant-function} equals $f$. We refer to $S_{I, f, \gamma^I}$ as the \emph{excursion operator} associated to $I$, $f$, and $\gamma^I$. They form a mutually commuting family of endomorphisms of $\Fun_{\cusp, \zeta}(\widetilde{\Bun}_{G, D}/\Xi, \overline{\rationals}_{\ell})$, commuting with the Hecke operators at each $x \in \mathring X$ (\emph{cf.}~\cite[Lemme 10.1, Lemme 10.2]{MR3787407}). Furthermore, $S_{I, f, \gamma^I}$ depends only on the image of $\gamma^I$ in $\pi_1^{\etale}(\mathring X, \bar{\eta})^I$ (\emph{cf.}~\cite[Proposition 10.10]{MR3787407}).

Denote by $\mathscr B$ the subalgebra of $\End \Fun_{\cusp, \zeta}(\widetilde{\Bun}_{G, D}/\Xi, \overline{\rationals}_{\ell})$ generated by the excursion operators. It follows from Lemma \ref{lem-cusp-form-finite-dimensional} that $\mathscr B$ is Artinian. Consider the decomposition of \eqref{eq-lafforgue-functors-vacuum-cusp-forms} into generalized $\mathscr B$-eigenspaces:
\begin{equation}
\label{eq-spectral-decomposition-into-excursion-operators}
	\Fun_{\cusp, \zeta}(\widetilde{\Bun}_{G, D}/\Xi, \overline{\rationals}_{\ell}) \simeq \bigoplus_{\chi : \mathscr B \rightarrow \overline{\rationals}_{\ell}} \mathbf H_{D, \chi}
\end{equation}

The desired decomposition \eqref{eq-spectral-decomposition} is defined by coarsening \eqref{eq-spectral-decomposition-into-excursion-operators}, using a map $\chi \mapsto [\sigma]$ from the set of characters of $\mathscr B$ to the set of $L$-parameters obtained from invariant theory. Namely, assigning the function $\gamma^I \mapsto S_{I, f, \gamma^I}$ to $f$ defines a map
\begin{equation}
\label{eq-excursion-operators-pseudo-representation}
C^0(H_{\bar{\eta}} \backslash ({}^LH_{\mathring X})^I / H_{\bar{\eta}}, \overline{\rationals}_{\ell}) \rightarrow C^0(\pi_1^{\etale}(\mathring X, \bar{\eta})^I, \mathscr B)
\end{equation}
where $\mathscr B$ is endowed with the $\ell$-adic topology. Thus, any character $\chi : \mathscr B \rightarrow \overline{\rationals}_{\ell}$ determines a continuous pseudo-representation of $\pi_1^{\etale}(\mathring X, \bar{\eta})$ with values in ${}^LH_{\mathring X}$ by composition with \eqref{eq-excursion-operators-pseudo-representation} (\emph{cf.}~\cite[D\'efinition-Proposition 11.3]{MR3787407}), which corresponds to a semisimple $L$-parameter $[\sigma]$ by \cite[Proposition 11.7]{MR3787407}.

The $L$-parameters arising this way are defined over finite extensions of $\rationals_{\ell}$. The Hecke eigen-property of $\mathbf H_{D, [\sigma]}$ follows from the generalized eigen-property of $\mathbf H_{D, \chi}$ with respect to $\mathscr B$ and the fact that Hecke operators are diagonalizable (\emph{cf.}~\cite[Remarque 11.6]{MR3787407}).
\end{proof}
\end{void}

\begin{rem}
At least for split $G$, one can use \cite{xue2020smoothness} to streamline the proof of Theorem \ref{thm-spectral-decomposition} as follows. We define the \emph{cuspidal subsheaf} $\mathscr H_{I, V, \cusp}^0 \subset\mathscr H^0_{I, V}$ as in \cite[Definition 7.0.3]{xue2020smoothness}. Then by \cite[Proposition 7.0.5]{xue2020smoothness} (\emph{cf.}~Remark \ref{rem-xue-smoothness}), $\mathscr H^0_{I, V, \cusp}$ is a lisse $\overline{\rationals}_{\ell}$-sheaf over $\mathring X^I$. By \cite[Proposition 1.3.4]{xue2020smoothness}, its $\pi_1^{\etale}(\mathring X^I, \bar{\eta}_I)$-action factors through $\pi_1^{\etale}(\mathring X, \bar{\eta})^I$, so by taking fibers at $\bar{\eta}^I$, we obtain a family of functors
\begin{align}
\label{eq-lafforgue-functors-xue}
	\Rep(({}^LH_{\mathring X})^I) & \rightarrow \Rep(\pi_1^{\etale}(\mathring X, \bar{\eta})^I), \\
	\notag
	V & \mapsto H_{I, V} := \mathscr H^0_{I, V, \cusp} |_{\bar{\eta}^I},
\end{align}
natural in $I$. The analogue of \eqref{eq-lafforgue-functors-vacuum-cusp-forms} is now tautological. The rest of the proof of Theorem \ref{thm-spectral-decomposition} proceeds as in \S\ref{void-spectral-decomposition-proof}.
\end{rem}

\newpage

\part{Appendices}

\appendix

\section{Twisting $\infty$-categories}
\label{sec-twisting-infinity-categories}

The goal of this section is to define twisted sections of a sheaf of $\infty$-categories, which is used both for the definition of twisted $\ell$-adic constructible complexes (\emph{cf.}~\S\ref{void-twisted-sheaves}) and for twisting the category of $H$-representations (\emph{cf.}~\S\ref{sec-twisted-representations}).

The material of \S\ref{sec-twisted-sheaves-construction} is similar to \cite[\S1.7]{MR3769731}, but we supply more details than \emph{loc.cit.}~on matters concerning homotopy coherence.

\subsection{Construction}
\label{sec-twisted-sheaves-construction}

\begin{void}
\label{void-twisted-sheaves-context}
Let $X$ be a site. Let $\coeff$ be a ring and $\Perf_{\coeff}$ be the $\infty$-category of perfect complexes of $\coeff$-modules. Let $\mathscr C$ be a sheaf of (small) Karoubi $\Perf_{\coeff}$-module categories over $X$.

Write $\mathscr X$ for the $\infty$-category of $\Spc$-valued sheaves over $X$. Write $\coeff^{\times}$ for the group of units of $\coeff$ and $\deloop^2 \coeff^{\times}$ its double deloop as an object of $\mathscr X$.

In this subsection, we shall explain that given a morphism $\mathscr G : x \rightarrow \deloop^2 \coeff^{\times}$ in $\mathscr X$, we obtain a sheaf $\mathscr C_{\mathscr G}$ of $\Perf_{\coeff}$-module categories over the slice site $X_{/x}$ functorial in $(x, \mathscr G)$. Objects of its global section $\Gamma(x, \mathscr C_{\mathscr G})$ are called \emph{$\mathscr G$-twisted sections of $\mathscr C$ over $x$.}
\end{void}

\begin{rem}
Informally, an object of $\Gamma(x, \mathscr C_{\mathscr G})$ is a compatible family of objects $c_{a, \tau_a} \in \Gamma(a, \mathscr C)$ for every $a \in X_{/x}$ equipped with a neutralization $\tau_a$ of the restriction $\mathscr G_a \in \Gamma(a, \deloop^2\coeff^{\times})$ of $\mathscr G$, such that for any two neutralizations $\tau_a$, $\tau_a'$ related by $\tau_a' \simeq \tau_a \cdot \ell_a$ for some $\ell_a \in */\coeff^{\times}$ (which always exists locally on $a$), there is an isomorphism
$$
c_{a, \tau_a'} \simeq c_{a, \tau_a} \otimes \ell_a,
$$
where $\ell_a$ is viewed as a free $\coeff$-module of rank $1$.
\end{rem}

\begin{void}
\label{void-sheaf-topos-perspective}
In order to construct $\mathscr C_{\mathscr G}$, we need to explain how to glue sheaves in a homotopy-coherent manner. This involves a bit of formalism.

Let $\mathscr O$ be an $\infty$-category admitting limits. Denote by $\Shv(X, \mathscr O)$ the $\infty$-category of $\mathscr O$-valued sheaves over $X$ (\emph{cf.}~\cite[Definition 1.3.1.1]{lurie2018spectral}). Tautologically, we have
$$
\mathscr X \simeq \Shv(X, \Spc).
$$

Given an $\infty$-topos $\mathscr Y$, we write $\Shv(\mathscr Y, \mathscr O)$ for the $\infty$-category of $\mathscr O$-valued sheaves over $\mathscr Y$, \emph{i.e.}~limit-preserving functors $\mathscr Y^{\opposite} \rightarrow \mathscr O$. By \cite[Proposition 1.3.1.7]{lurie2018spectral}, restriction along the sheafification of the Yoneda embedding determines an equivalence of $\infty$-categories:
$$
\Shv(\mathscr X, \mathscr O) \simeq \Shv(X, \mathscr O).
$$

Note that if $\mathscr O$ is presentable, then $\Shv(\mathscr X, \mathscr O)$ is canonically identified with the Lurie tensor product $\mathscr X \otimes \mathscr O$ (\emph{cf.}~\cite[Proposition 4.8.1.17]{lurie2017higher}).
\end{void}

\begin{void}
Denote by $\PrL$ (respectively, $\PrR$) the $\infty$-category of presentable $\infty$-categories with colimit-preserving functors (respectively, limit-preserving functors).

Given a presentable $\infty$-category $\mathscr O$, we shall construct a functor
\begin{equation}
\label{eq-slice-sheaf-functor}
	\Shv(\mathscr X_{/(\cdot)}, \mathscr O) : \mathscr X^{\opposite} \rightarrow \PrR,
\end{equation}
which assigns to $x \in \mathscr X^{\opposite}$ the $\infty$-category $\Shv(\mathscr X_{/x}, \mathscr O)$.

Indeed, the co-Cartesian fibration $\ev_1 : \Fun(\Delta^1, \mathscr X) \rightarrow \mathscr X$ defines a functor
\begin{equation}
\label{eq-slice-category-functor}
\mathscr X \rightarrow \PrL,\quad x\mapsto \mathscr X_{x/}
\end{equation}
Its composition with the endofunctor $(\cdot)\otimes\mathscr O$ of $\PrL$ yields a functor
\begin{equation}
\label{eq-slice-sheaf-functor-opposite}
\mathscr X \rightarrow \PrL,\quad x\mapsto \mathscr X_{/x} \otimes \mathscr O \simeq \Shv(\mathscr X_{/x}, \mathscr O)
\end{equation}
The functor \eqref{eq-slice-sheaf-functor} is defined to be the opposite of \eqref{eq-slice-sheaf-functor-opposite}, using the canonical equivalence $(\PrL)^{\opposite} \simeq \PrR$ (\emph{cf.}~\cite[Corollary 5.5.3.4]{MR2522659}).

The following assertion can be viewed as a precise version of the statement that the $\infty$-category of sheaves itself satisfies gluing. Its proof is contained in \cite[Remark 2.1.0.5]{lurie2018spectral}, but we reproduce the argument for the reader's conveninence.
\end{void}

\begin{lem}
\label{lem-slice-sheaf-functor-descent}
The functor \eqref{eq-slice-sheaf-functor} preserves limits.
\end{lem}

\begin{proof}
It suffices to show that \eqref{eq-slice-sheaf-functor-opposite} preserves colimits. Since the endofunctor $(\cdot) \otimes \mathscr O$ of $\PrL$ preserves colimits, it is enough to prove that \eqref{eq-slice-category-functor} preserves colimits. Passing again to the opposite $\infty$-categories, we obtain the functor
\begin{equation}
\label{eq-slice-category-functor-opposite}
\mathscr X^{\opposite} \rightarrow \PrR,\quad x \mapsto \mathscr X_{x/}
\end{equation}
classifying $\ev_1 : \Fun(\Delta^1, \mathscr X) \rightarrow \mathscr X$ as a \emph{Cartesian} fibration.

Since limits in $\PrL$ and $\PrR$ are both computed by limits in $\Cat$, the fact that \eqref{eq-slice-category-functor-opposite} preserves limits follows from \cite[Theorem 6.1.3.9]{MR2522659}.
\end{proof}

\begin{void}
We now specialize to the case $\mathscr O := \Perf_{\coeff}\Mod$, the $\infty$-category of (small) Karoubi $\Perf_{\coeff}$-module $\infty$-categories.

Given a $\Perf_{\coeff}\Mod$-valued sheaf $\mathscr C$ on $X$, we obtain a $\Perf_{\coeff}\Mod$-valued sheaf over the $\infty$-topos $\mathscr X$ (\emph{cf.}~\S\ref{void-sheaf-topos-perspective}), and consequently a morphism in $\Shv(\mathscr X, \PrR)$:
\begin{equation}
\label{eq-sheaf-of-cat-as-point}
\mathscr C : * \rightarrow \Shv(\mathscr X_{/(\cdot)}, \Perf_{\coeff}\Mod).
\end{equation}
Here, the target is an object of $\Shv(\mathscr X, \PrR)$ thanks to Lemma \ref{lem-slice-sheaf-functor-descent}. Let us view \eqref{eq-sheaf-of-cat-as-point} as a morphism in $\Shv(\mathscr X, \Cat)$ via the forgetful functor $\PrR \rightarrow \Cat$.
\end{void}

\begin{prop}
\label{prop-universal-twisted-category-construction}
The morphism \eqref{eq-sheaf-of-cat-as-point} canonically extends to
\begin{equation}
\label{eq-universal-twisted-category}
\widetilde{\mathscr C} : \deloop^2 \coeff^{\times} \rightarrow \Shv(\mathscr X_{/(\cdot)}, \Perf_{\coeff}\Mod)
\end{equation}
as a morphism in $\Shv(\mathscr X, \Cat)$.
\end{prop}

\begin{void}
In order to construct the extension \eqref{eq-universal-twisted-category}, we note that $\Perf_{\coeff}$ admits a commutative algebra structure. Consequently, the $\infty$-category $\Perf_{\coeff}\Mod$ admits a symmetric monoidal structure, given by Lurie tensor product of the ind-completed $\infty$-categories.\footnote{The functor of ind-completion realizes $\Perf_{\coeff}\Mod$ as the $\infty$-category of $\coeff$-linear presentable $\infty$-categories with functors preserving colimits and compact objects.}

The tensor product functor on $\Perf_{\coeff}\Mod$ gives rise to a functor
$$
\Perf_{\coeff}\Mod \rightarrow \End(\Perf_{\coeff}\Mod),\quad M \mapsto (\cdot)\otimes_{\Perf_{\coeff}} M
$$
by adjunction. It sends $\Perf_{\coeff}$ to $\id_{\Perf_{\coeff}\Mod}$. By taking endomorphisms, we obtain a monoidal morphism
\begin{equation}
\label{eq-commutative-algebra-to-bernstein-center}
	\Perf_{\coeff} \rightarrow \End(\id_{\Perf_{\coeff}}\Mod).
\end{equation}
\end{void}

\begin{rem}
The $\mathbb E_2$-monoidal structure on $\Perf_{\coeff}$ suffices for the construction of \eqref{eq-commutative-algebra-to-bernstein-center}. However, the associative algebra structure does \emph{not} suffice: Informally, \eqref{eq-commutative-algebra-to-bernstein-center} sends $A \in \Perf_{\coeff}$ to the compatible family of endomorphisms of $M$, for each $M \in \Perf_{\coeff}\Mod$, given by multiplication by $A$, which is only $\Perf_{\coeff}$-linear by its $\mathbb E_2$-monoidal structure.
\end{rem}

\begin{void}
\label{void-universal-twist-construction}
We are now ready to construct the morphism \eqref{eq-universal-twisted-category}.

\begin{proof}[Proof of Proposition \ref{prop-universal-twisted-category-construction}]
Since $\deloop^2\coeff^{\times}$ is the constant sheaf associated to the double deloop $K(\coeff^{\times}, 2)$ in the $\infty$-category $\Spc$, it suffices to construct a morphism in $\Cat$
\begin{equation}
\label{eq-universal-twisted-category-global-section}
K(\coeff^{\times}, 2) \rightarrow \Shv(\mathscr X, \Perf_{\coeff} \Mod)
\end{equation}
sending the neutral point to $\mathscr C$.

The pointed morphism \eqref{eq-universal-twisted-category-global-section} is equivalent to a monoidal morphism $*/\coeff^{\times} \rightarrow \End(\mathscr C)$, where the endomorphism monoid is computed in the $\infty$-category $\Fun(\mathscr X^{\opposite}, \Perf_{\coeff}\Mod)$. The desired monoidal morphism is given by the composition
$$
*/\coeff^{\times} \subset \Perf_{\coeff}\rightarrow \End(\id_{\Perf_{\coeff}\Mod}) \rightarrow \End(\mathscr C),
$$
where $*/\coeff^{\times}$ is identified as the full subcategory of $\Perf_{\coeff}$ spanned by $\coeff$, the middle morphism is \eqref{eq-commutative-algebra-to-bernstein-center}, while the last morphism is defined by the functor
$$
(\cdot) \circ \mathscr C : \End(\Perf_{\coeff}\Mod) \rightarrow \Fun(\mathscr X^{\opposite}, \Perf_{\coeff}\Mod)
$$
by evaluating on the endomorphisms of the object $\id_{\Perf_{\coeff}\Mod}$.
\end{proof}

\end{void}

\begin{void}[Construction of $\mathscr C_{\mathscr G}$]
\label{void-twisted-category-construction}
Let us now fix $x \in \mathscr X$ and a morphism $\mathscr G : x \rightarrow \deloop^2\coeff^{\times}$. Having \eqref{eq-universal-twisted-category} at our disposal, we define $\mathscr C_{\mathscr G}$ to be the image of $\mathscr G$ under \eqref{eq-universal-twisted-category}:
$$
\mathscr C_{\mathscr G} := \Gamma(x, \widetilde{\mathscr C})(\mathscr G)
$$

By construction, $\mathscr C_{\mathscr G}$ is an object of $\Shv(\mathscr X_{/x}, \Perf_{\coeff}\Mod)$. We obtain a sheaf on the slice site $X_{/x}$ by restriction along the sheafification of the Yoneda embedding (\emph{cf.}~\S\ref{void-sheaf-topos-perspective}).
\end{void}

\subsection{Monoidality}
\label{sec-monoidal-twist}

\begin{void}
We shall now introduce a variant of the twisting construction for sheaves of \emph{symmetric monoidal} categories, needed for the definition of ``twisted $H$-representations" (\emph{cf.}~\S\ref{sec-twisted-representations}).

Let $\mathscr X$ be an $\infty$-topos and $x\in\mathscr X$. Then $\Shv(\mathscr X_{/x}, \Perf_{\coeff}\Mod)$ inherits a symmetric monoidal structure from that of $\Perf_{\coeff}\Mod$, functorially in $x$.

The construction $\mathscr G, \mathscr C \mapsto \mathscr C_{\mathscr G}$ (\emph{cf.}~\S\ref{void-twisted-category-construction}) may be regarded as a functor
\begin{equation}
\label{eq-universal-twisted-category-as-bifunctor}
\Gamma(x, \deloop^2 \coeff^{\times}) \times \Shv(\mathscr X, \Perf_{\coeff}\Mod) \rightarrow \Shv(\mathscr X_{/x}, \Perf_{\coeff}\Mod).
\end{equation}
\end{void}

\begin{lem}
\label{lem-universal-twist-bifunctor-symmetric-monoidal}
The functor \eqref{eq-universal-twisted-category-as-bifunctor} is naturally symmetric monoidal.
\end{lem}

\begin{proof}
By construction, \eqref{eq-universal-twisted-category-as-bifunctor} is induced via descent from the functor
\begin{equation}
\label{eq-universal-twist-pointwise-bifunctor}
K(\coeff^{\times}, 2) \times \Perf_{\coeff}\Mod \rightarrow \Perf_{\coeff}\Mod
\end{equation}
adjoint to the deloop of $*/\coeff^{\times} \subset \Perf_{\coeff} \rightarrow \End(\id_{\Perf_{\coeff}\Mod})$ (\emph{cf.}~\S\ref{void-universal-twist-construction}). It suffices to endow \eqref{eq-universal-twist-pointwise-bifunctor} with the structure of a symmetric monoidal functor.

However, \eqref{eq-universal-twist-pointwise-bifunctor} is identified with the composition
\begin{align*}
	K(\coeff^{\times}, 2) &\times \Perf_{\coeff}\Mod \\
	& \rightarrow \Perf_{\coeff}\Mod \times \Perf_{\coeff}\Mod \xrightarrow{\otimes} \Perf_{\coeff}\Mod,
\end{align*}
where both functors are symmetric monoidal. For the second functor, we have invoked the fact that the multiplication map of any commutative algebra in a symmetric monoidal $\infty$-category naturally lifts to a map of commutative algebras, which is a consequence of \cite[Proposition 3.2.4.7]{lurie2017higher}.
\end{proof}

\begin{rem}
Informally, the symmetric monoidal structure on \eqref{eq-universal-twisted-category-as-bifunctor} yields a canonical isomorphism of $\Perf_{\coeff}\Mod$-valued sheaves over $\mathscr X_{/x}$:
\begin{equation}
\label{eq-universal-twist-bifunctor-monoidal-explicit}
(\mathscr C_1)_{\mathscr G_1} \otimes (\mathscr C_2)_{\mathscr G_2} \simeq (\mathscr C_1 \otimes \mathscr C_2)_{\mathscr G_1 + \mathscr G_2}
\end{equation}
for any $\Perf_{\coeff}\Mod$-valued sheaves $\mathscr C_1, \mathscr C_2$ over $\mathscr X$ and $A$-gerbes $\mathscr G_1$, $\mathscr G_2$ over $x$ (\emph{i.e.}~morphisms $x\rightarrow \deloop^2\coeff^{\times}$), together with homotopy coherence.
\end{rem}

\begin{void}
\label{void-graded-sheaf-symmetric-monoidal-categories}
Let $\Xi$ be a commutative monoid, whose binary operation is written additively.

By a \emph{$\Xi$-graded sheaf of $\coeff$-linear symmetric monoidal $\infty$-categories over $\mathscr X$}, we shall mean a (right) lax symmetric monoidal functor (\emph{cf.}~\cite[\S2.1.3]{lurie2017higher}):
$$
\mathscr C^{\Xi} : \Xi \rightarrow \Shv(\mathscr X, \Perf_{\coeff}\Mod),\quad \xi \mapsto \mathscr C^{\xi}.
$$

Such functors form an $\infty$-category which we denote by $\Shv^{\Xi, \otimes}(\mathscr X, \Perf_{\coeff}\Mod)$.
\end{void}

\begin{rem}
Informally, the lax symmetric monoidal structure supplies morphisms
\begin{align*}
	\mathbf 1 &: \Perf_{\coeff} \rightarrow \mathscr C^0 \\
	\otimes &: \mathscr C^{\xi_1} \otimes \mathscr C^{\xi_2} \rightarrow \mathscr C^{\xi_1 + \xi_2}
\end{align*}
together with homotopy coherence.

Note that for $\Xi \simeq 0$ the trivial commutative monoid, an object of $\Shv^{0, \otimes}(\mathscr X, \Perf_{\coeff}\Mod)$ is precisely a sheaf of $\coeff$-linear symmetric monoidal categories over $\mathscr X$, with $\mathbf 1$ given by its unit and $\otimes$ given by its monoidal product.
\end{rem}

\begin{void}
\label{void-graded-symmetric-monoidal-category-sum}
Restriction along $\Xi \rightarrow 0$ yields a functor
$$
\Shv^{0, \otimes}(\mathscr X, \Perf_{\coeff}\Mod) \rightarrow \Shv^{\Xi, \otimes}(\mathscr X, \Perf_{\coeff}\Mod),
$$
which admits a left adjoint $\LKE^{\otimes}$ (\emph{cf.}~\cite[Theorem B.1.3]{arinkin2020stack}, which is a special case of \cite[Corollary 3.1.3.5]{lurie2017higher}).

Furthermore, the underlying $\Perf_{\coeff}$-valued sheaf over $\mathscr X$ of $\LKE^{\otimes}(\mathscr C^{\Xi})$ is computed by
$$
\LKE^{\otimes}(\mathscr C^{\Xi}) \simeq \bigoplus_{\xi \in \Xi} \mathscr C^{\xi},
$$
where $\bigoplus$ is the colimit in $\Shv(\mathscr X, \Perf_{\coeff}\Mod)$. Thus, any $\Xi$-graded sheaf of $\coeff$-linear symmetric monoidal categories has an associated sheaf of $\coeff$-linear symmetric monoidal categories, by taking sum over its graded components.

For our purposes, it is convenient to view the datum of $\mathscr C^{\Xi}$ as additional structure on the sheaf of $\coeff$-linear symmetric monoidal categories $\mathscr C := \LKE^{\otimes}(\mathscr C^{\Xi})$, and simply refer to $\mathscr C$ as a $\Xi$-graded sheaf of $\coeff$-linear symmetric monoidal $\infty$-categories.
\end{void}

\begin{void}[Symmetric monoidal twists]
\label{void-symmetric-monoidal-twist}
Let $\mathscr X$ be an $\infty$-topos with $x \in \mathscr X$. We shall consider a variant of \S\ref{void-twisted-category-construction} which takes symmetric monoidal structures into account.

Our input is a triple $(\Xi, \mathscr C^{\Xi}, \nu)$, where
\begin{enumerate}
	\item $\Xi$ is a commutative monoid;
	\item $\mathscr C^{\Xi}$ is a $\Xi$-graded sheaf of $\coeff$-linear symmetric monoidal $\infty$-categories over $\mathscr X$;
	\item $\nu : \Xi \rightarrow \deloop^2 \coeff^{\times}$ is a symmetric monoidal morphism over $\mathscr X_{/x}$.
\end{enumerate}

We define the $\Xi$-graded sheaf of $\coeff$-linear symmetric monoidal $\infty$-categories $\mathscr C^{\Xi}_{\nu}$ over $\mathscr X_{/x}$ as the composition of 
$$
(\nu, \mathscr C^{\Xi}) : \Xi \rightarrow \Gamma(x, \deloop^2 \coeff^{\times}) \times \Shv(\mathscr X, \Perf_{\coeff}\Mod)
$$
with \eqref{eq-universal-twisted-category-as-bifunctor}, which is symmetric monoidal by Lemma \ref{lem-universal-twist-bifunctor-symmetric-monoidal}.
\end{void}

\begin{rem}
By construction, for each $\xi \in \Xi$, the $\xi$-graded component $\mathscr C_{\nu}^{\xi}$ of $\mathscr C^{\Xi}_{\nu}$ is the twist of $\mathscr C^{\xi}$ by $\nu(\xi)$ in the sense of \S\ref{void-twisted-category-construction}.

On the other hand, the lax symmetric monoidal structure on $\mathscr C^{\Xi}_{\nu}$ depends on the symmetric monoidal structure on $\nu$. For example, the monoidal product $\otimes$ may be expressed as the composition
\begin{align*}
 (\mathscr C^{\xi_1})_{\nu(\xi_1)} \otimes (\mathscr C^{\xi_2})_{\nu(\xi_2)} &\simeq (\mathscr C^{\xi_1} \otimes \mathscr C^{\xi_2})_{\nu(\xi_1) + \nu(\xi_2)} \\
 & \xrightarrow{\otimes} (\mathscr C^{\xi_1 + \xi_2})_{\nu(\xi_1) + \nu(\xi_2)} \simeq (\mathscr C^{\xi_1 + \xi_2})_{\nu(\xi_1 + \xi_2)},
\end{align*}
where the first isomorphism is \eqref{eq-universal-twist-bifunctor-monoidal-explicit} and the last isomorphism is supplied by the monoidal structure on $\nu$.
\end{rem}

\medskip

\section{Twisted $\ell$-adic sheaves}
\label{sec-twisted-ell-adic-sheaves}

The goal of this section is to define and study $\ell$-adic constructible complexes twisted by a gerbe, \emph{cf.}~\S\ref{void-twisted-sheaves}. We will extend a number of standard results to the twisted setting. These mainly concern universal local acyclicity (\emph{cf.}~Definition \ref{defn-universal-local-acyclicity}) and hyperbolic localization (\emph{cf.}~Theorem \ref{thm-hyperbolic-localization}).

Although the main body of this article takes place over a field, we work over more general base schemes in this section.

\subsection{Definitions}

\begin{void}
\label{void-sheaf-theory-context}
We work over a base scheme $S$, which is assumed quasi-compact quasi-separated with finitely many irreducible components.\footnote{We use these conditions to ensure the existence of the relative perverse $t$-structure (\emph{cf.}~\S\ref{void-relative-perversity}). For our applications, it is (almost) sufficient to restrict to schemes of finite type over a field.} Denote by $\Sch$ the category of such schemes.

Let $\ell$ be a prime invertible in $S$ and $\coeff$ be a finite extension of $\rationals_{\ell}$. Let $A$ be a finite subgroup of $\coeff^{\times}$ of order invertible in $S$ and write $\zeta : A \rightarrow \coeff^{\times}$ for the inclusion.

Consider the sheaf theory of constructible \'etale sheaves of $\coeff$-vector spaces, which we understand as a functor of $\infty$-categories (\emph{cf.}~\cite[\S3]{MR4609461})
\begin{equation}
\label{eq-sheaf-theory-constructible}
\Sch^{\opposite} \rightarrow \Perf_{\coeff}\Mod,\quad X \mapsto \derived(X) := \derived_{\mathrm{cons}}(X, \coeff).
\end{equation}
To a morphism $f : X_1 \rightarrow X_2$ in $\Sch$, \eqref{eq-sheaf-theory-constructible} assigns the functor $f^* : \derived(X_2) \rightarrow \derived(X_1)$.
\end{void}

\begin{rem}
\label{rem-twisted-six-functor-formalism}
By \cite[Theorem 7.7]{MR4609461}, the $\infty$-category $\derived(X)$ coincides with the classically defined $\infty$-category of $\ell$-adic constructible complexes over $X$.

In particular, we have the functor $f_! : \derived(X_1) \rightarrow \derived(X_2)$, whenever $f$ is separated and of finite presentation: We refer the reader to \cite[\S3.4]{liu2024enhancedoperationsbasechange} for its construction as a functor between $\infty$-categories and to \cite[Expos\'e XVII, Th\'eor\`eme 5.3.6]{SGA4-3} for the preservation of constructibility. The functor $f_!$ admits a right adjoint $f^!$.
\end{rem}

\begin{void}[Twisted $\ell$-adic sheaves]
\label{void-twisted-sheaves}
Denote by $\deloop^2 A$ the double deloop of $A$, viewed as a sheaf on the big \'etale site of $S$.

The assignment \eqref{eq-sheaf-theory-constructible} satisfies \'etale descent (\emph{cf.}~\cite[Corollary 4.7]{MR4609461}), so the twisting construction of \S\ref{sec-twisted-sheaves-construction} yields a functor of $\infty$-categories
\begin{equation}
\label{eq-twisted-sheaves-functor}
(\Sch_{/\deloop^2 A})^{\opposite} \rightarrow \Perf_{\coeff}\Mod, \quad (X, \mathscr G) \mapsto \derived_{\mathscr G, \zeta}(X) := \Gamma(X, \derived(\cdot)_{\mathscr G_{\zeta}}),
\end{equation}
where $\mathscr G_{\zeta}$ denotes the composition of $\mathscr G : X \rightarrow \deloop^2 A$ with the morphism $\zeta : \deloop^2 A \rightarrow \deloop^2\coeff^{\times}$.

We refer to objects of $\derived_{\mathscr G, \zeta}(X)$ as \emph{$(\mathscr G, \zeta)$-twisted constructible complexes over $X$.}
\end{void}

\begin{rem}
\label{rem-etale-local-notion-twisted-complex}
Since $\derived_{\mathscr G, \zeta}(X)$ is identified with $\derived(X)$ once a neutralization of $\mathscr G$ is fixed, and neutralizations of $\mathscr G$ exist \'etale locally, all \'etale-local notions about constructible complexes extend to $(\mathscr G, \zeta)$-twisted constructible complexes.

In particular, given a morphism $(X_1, \mathscr G_1) \rightarrow (X_2, \mathscr G_2)$ in $\Sch_{/\deloop^2 A}$ whose first component $f : X_1 \rightarrow X_2$ is separated and of finite presentation, we have the exceptional adjunction
$$
\begin{tikzcd}
	\derived_{\mathscr G_1, \zeta}(X_1) \ar[r, shift left = 0.5ex, "f_!"] & \derived_{\mathscr G_2, \zeta}(X_2) \ar[l, shift left = 0.5ex, "f^!"]
\end{tikzcd}
$$
constructed by descent along an \'etale cover of $X_2$ trivializing $\mathscr G_2$.
\end{rem}

\begin{void}[The $2$-category $\kernel_S$]
Next, we shall define universal local acyclicity for twisted complexes. Since this is a \'etale-local notion, we can define it using Remark \ref{rem-etale-local-notion-twisted-complex}. However, we will offer a definition intrinsic to twisted complexes, by adapting \cite[\S3]{MR4630128}.

Write $\kernel_S$ for the $2$-category\footnote{It is sufficient for our purposes to construct $\kernel_S$ as an \emph{ordinary} $2$-category.} whose objects are triples $(X, \mathscr G, \mathscr A)$, where $X$ is a separated $S$-scheme of finite presentation, $\mathscr G$ is an $A$-gerbe over $X$, and $\mathscr A \in \Ho \derived_{\mathscr G, \zeta}(X)$.

A morphism
\begin{equation}
\label{eq-category-of-kernels-morphism}
(X_1, \mathscr G_1, \mathscr A_1) \rightarrow (X_2, \mathscr G_2, \mathscr A_2)
\end{equation}
in $\kernel_S$ consists of an object $\mathscr B \in \Ho\derived_{\overset{\leftarrow}\pi{}^*\mathscr G_1 - \overset{\rightarrow}{\pi}{}^*\mathscr G_2, \zeta}(X_1\times_S X_2)$ together with a morphism
$$
\overset{\rightarrow}{\pi}_!(\overset{\leftarrow}{\pi}{}^*\mathscr A_1 \otimes \mathscr B) \rightarrow \mathscr A_2,
$$
where $\overset{\leftarrow}{\pi}$, $\overset{\rightarrow}{\pi}$ denote the projections of $X_1\times_S X_2$ onto $X_1$, respectively $X_2$. Note that morphisms \eqref{eq-category-of-kernels-morphism} form a category in the evident way, and their compositions are defined by convolution product in $\mathscr B$.
\end{void}

\begin{void}
The $2$-category $\kernel_S$ admits a symmetric monoidal structure, with monoidal product
$$
(X_1, \mathscr G_1, \mathscr A_1) \otimes (X_2, \mathscr G_2, \mathscr A_2) := (X_1\times_S X_2, \overset{\leftarrow}{\pi}{}^*\mathscr G_1 + \overset{\rightarrow}{\pi}{}^*\mathscr G_2, \mathscr A_1 \boxtimes \mathscr A_2)
$$
and monoidal unit $(S, \mathbf 0, \coeff)$, where $\mathbf 0$ stands for the trivial $A$-gerbe.

Internal Homs exist in $\kernel_S$ and are given by
\begin{equation}
\label{eq-category-of-kernels-internal-hom}
\SHom_{\kernel_S}((X_1, \mathscr G_1, \mathscr A_1), (X_2, \mathscr G_2, \mathscr A_2)) \simeq (X_1\times_S X_2, -\overset{\leftarrow}{\pi}{}^*\mathscr G_1 + \overset{\rightarrow}{\pi}{}^*\mathscr G_2, \SHom(\overset{\leftarrow}{\pi}{}^*\mathscr A_1, \overset{\rightarrow}{\pi}{}^!\mathscr A_2)).
\end{equation}
\end{void}

\begin{defn}
\label{defn-universal-local-acyclicity}
Given a separated $S$-scheme $X$ of finite presentation and an $A$-gerbe $\mathscr G$ over $X$, an object $\mathscr A \in \derived_{\mathscr G, \zeta}(X)$ is \emph{universally locally acyclic} (ULA) relative to $S$ if $
(X, \mathscr G, \mathscr A)$ is dualizable as an object of $\kernel_S$.
\end{defn}

\begin{rem}
\label{rem-ula-reduction-to-untwisted-notion}
The dualizability of $(X, \mathscr G, \mathscr A)$ can be checked \'etale locally over $X$, where we may assume that $\mathscr G$ is neutral. In that case, Definition \ref{defn-universal-local-acyclicity} coincides with \cite[Definition 3.2]{MR4630128}. Thus, our notion of ULA is equivalent to the classical notion over an \'etale cover of $X$ neutralizing $\mathscr G$.
\end{rem}

\begin{rem}
\label{rem-ula-relative-to-geometric-point}
The analogue of \cite[Theorem 4.4]{MR4630128} holds in our context, by repeating the proof of \emph{loc.cit.}~or by applying Remark \ref{rem-ula-reduction-to-untwisted-notion}.

In particular, the criterion of \cite[Theorem 4.4(iv)]{MR4630128} shows that when $S$ is the spectrum of an algebraically closed field, any $\mathscr A \in \derived_{\mathscr G, \zeta}(X)$ is ULA relative to $S$.
\end{rem}

\begin{void}[Verdier duality]
\label{void-verdier-duality}
Given a separated $S$-scheme $X$ of finite presentation equipped with an $A$-gerbe $\mathscr G$, we define Verdier duality to be the functor
\begin{align*}
\mathbf D_{X/S} : \derived_{\mathscr G, \zeta}(X) &\rightarrow \derived_{-\mathscr G, \zeta}(X) \\
 \mathscr A &\mapsto \SHom(\mathscr A, \pi^! \coeff),
\end{align*}
where $\pi : X \rightarrow S$ denotes the structural morphism.

If $\mathscr A$ is ULA relative to $S$, then the dual of $(X, \mathscr G, \mathscr A)\in\kernel_S$ is computed by \eqref{eq-category-of-kernels-internal-hom} to be $(X, -\mathscr G, \mathbf D_{X/S}(\mathscr A))$. It follows that $\mathbf D_{X/S}(\mathscr A)$ is also ULA relative to $S$ and we have an canonical isomorphism $\mathscr A \simeq \mathbf D_{X/S}(\mathbf D_{X/S} \mathscr A)$.

Since $\kernel_S$ is functorial with respect to ($*$-)pullbacks along $S$, we see that the formation of $\mathbf D_{X/S}(\mathscr A)$ for an ULA object $\mathscr A$ commutes with any base change along $S$.
\end{void}

\begin{rem}
\label{rem-internal-hom-via-duality}
Given two ULA objects $\mathscr A, \mathscr B \in \derived_{\mathscr G, \zeta}(X)$, we may express their internal Hom in terms of duality and tensor product:
\begin{equation}
\label{eq-internal-hom-via-duality}
	\SHom(\mathscr A, \mathscr B) \simeq \mathbf D_{X/S}(\mathscr A \otimes \mathbf D_{X/S} \mathscr B).
\end{equation}

Indeed, the canonical isomorphism $\SHom_{\kernel_S}(a, b) \simeq \SHom_{\kernel_S}(a, \mathbf 1) \otimes b$, for $a := (X, \mathscr G, \mathscr A), b := (X, \mathscr G, \mathscr B) \in \kernel_S$, yields an isomorphism in $\Ho\derived(X)$:
\begin{equation}
\label{eq-internal-hom-via-duality-shriek-tensor}
\SHom(\mathscr A, \mathscr B) \simeq \Delta^!(\mathbf D_{X/S}(\mathscr A) \boxtimes \mathscr B),
\end{equation}
where $\Delta : X \rightarrow X\times_S X$ is the diagonal embedding. Since $\mathbf D_{X/S}(\mathscr A)$ and $\mathscr B$ are both ULA, so is $\mathbf D_{X/S}(\mathscr A) \boxtimes\mathscr B$ (as it corresponds to a monoidal product in $\kernel_S$). Then \eqref{eq-internal-hom-via-duality} follows by rewriting the right-hand-side of \eqref{eq-internal-hom-via-duality-shriek-tensor} using the isomorphism $\Delta^! \mathbf D_{X\times_S X/S} \simeq \mathbf D_{X/S} \Delta^*$.
\end{rem}

\begin{void}[Relative perversity]
\label{void-relative-perversity}
Finally, we recall the relative perverse $t$-structure (\emph{cf.}~\cite[Theorem 6.1]{MR4630128}), which adapts to twisted complexes without change.

Namely, given a finitely presented $S$-scheme $X$ together with an $A$-gerbe $\mathscr G$ over $X$, there is a $t$-structure on $\derived_{\mathscr G, \zeta}(X)$ such that for any $\mathscr A \in \derived_{\mathscr G, \zeta}(X)$,
\begin{align*}
	\mathscr A \in {}^p\derived^{\le 0}_{\mathscr G, \zeta}(X) & \Leftrightarrow (i_{\bar s})^*\mathscr A \in {}^p\derived^{\le 0}_{\mathscr G, \zeta}(X_{\bar s}) \text{ for any geometric point }\bar s \rightarrow S \\
	\mathscr A \in {}^p\derived^{\ge 0}_{\mathscr G, \zeta}(X) & \Leftrightarrow (i_{\bar s})^*\mathscr A \in {}^p\derived^{\ge 0}_{\mathscr G, \zeta}(X_{\bar s}) \text{ for any geometric point }\bar s \rightarrow S
\end{align*}
where $i_{\bar s} : X_{\bar x} := X\times_S \bar x \rightarrow X$ denotes the base change of $\bar s \rightarrow S$.

Here, the fiberwise perverse $t$-structure is defined either by repeating the classical definition, or equivalently by appealing to Remark \ref{rem-ula-reduction-to-untwisted-notion}.
\end{void}

\begin{rem}
\label{rem-ula-perverse-truncation}
The analogue of \cite[Theorem 6.7]{MR4630128} holds in our context. Note that its proof shows that the perverse truncation functors preserve universal local acyclicity.
\end{rem}

\subsection{Vanishing by equivariance}

\begin{void}
We remain in the context of \S\ref{void-sheaf-theory-context}. Let $G$ be a group $S$-scheme of finite type with connected geometric fibers.

Let $\mathscr G$ be an $A$-gerbe over $\deloop G$. Denote by $\derived_{\mathscr G, \zeta}(\deloop G)$ the $\infty$-category of $G$-equivariant objects of $\derived_{\mathscr G, \zeta}(S)$. For any geometric point $\bar s$ of $S$, we write $\mathscr G_{\bar s}$ for the fiber of $\mathscr G$ at $\bar s$, which is an $A$-gerbe over $\deloop G_{\bar s}$.
\end{void}

\begin{lem}
\label{lem-equivariance-vanishing}
Suppose that $\mathscr G_{\bar s}$ is nontrivial for all geometric points $\bar s$ of $S$. Then
$$
\derived_{\mathscr G, \zeta}(\deloop G) \simeq 0.
$$
\end{lem}

\begin{proof}
It suffices to prove that $\mathscr A \in \derived_{\mathscr G, \zeta}(\deloop G)$ has vanishing fiber at any geometric point $\bar s$ of $S$. Thus, we may replace $S$ by $\bar s$ and assume that $S$ is itself a geometric point. Fixing a neutralization of $\mathscr G$ along the base point of $\deloop G$, we may assume that $\mathscr G$ is the deloop of a nontrivial monoidal morphism $\chi : G \rightarrow \deloop A$.

Note that $\chi$ is nontrivial as a plain morphism: Any morphism $G \times G \rightarrow A$ which is trivial along $e\times G$ and $G \times e$ is itself trivial, by connectedness of $G$. Denote by $\mathscr L_{\chi, \zeta}$ the $1$-dimensional $\coeff$-local system over $G$ induced from $\chi$ along $\zeta$. We then have
\begin{equation}
\label{eq-nontrivial-character-vanishing-global-section}
H^0(G, \mathscr L_{\chi, \zeta}) \simeq 0,
\end{equation}
as it is the space of invariants of a nontrivial character of $\pi_1^{\etale}(G, e)$.

Let us describe an object $\mathscr A \in \derived_{\mathscr G, \zeta}(\deloop G)$ by its fiber $\mathscr A_e$ along $e$, which is a constructible complex of $\coeff$-vector spaces, together with (the first piece of) its descent data:
\begin{equation}
\label{eq-classifying-stack-twisted-descent-isomorphism}
\underline{\coeff} \otimes \mathscr A_e \simeq \mathscr L_{\chi, \zeta}\otimes \mathscr A_e \text{ in }\derived(G).
\end{equation}

If $H^i(\mathscr A_e) \neq 0$ for some $i \in \integers$, then by taking $i$th cohomology of \eqref{eq-classifying-stack-twisted-descent-isomorphism} and choosing a surjective character $H^i(\mathscr A_e) \twoheadrightarrow \coeff$, we obtain nonzero global sections of $\mathscr L_{\chi, \zeta}$, contradicting \eqref{eq-nontrivial-character-vanishing-global-section}. Thus, $H^i(\mathscr A_e) \simeq 0$ for all $i\in\integers$ and the complex $\mathscr A$ vanishes.
\end{proof}

\begin{void}
We shall strengthen Lemma \ref{lem-equivariance-vanishing} to a vanishing statement about compactly supported cohomology of twisted constructible complexes.

Let $X$ be a separated $S$-scheme of finite presentation equipped with a $G$-action. Denote by $\derived_{\mathscr G, \zeta}(X/G)$ the $\infty$-category of $G$-equivariant objects of $\derived_{\mathscr G, \zeta}(X)$, formed with respect to the pullback of $\mathscr G$ along the structural morphism $X/G \rightarrow \deloop G$.

The structural morphism $f : X \rightarrow S$ defines a functor
\begin{equation}
\label{eq-compactly-supported-section-twisted-equivariant}
f_! : \derived_{\mathscr G, \zeta}(X) \rightarrow \derived_{\mathscr G, \zeta}(S),
\end{equation}
where the target is formed with respect to the pullback of $\mathscr G$ along $e : S \rightarrow \deloop G$.

The following statement includes Lemma \ref{lem-equivariance-vanishing} as the special case for $X = S$.
\end{void}

\begin{prop}
\label{prop-vanishing-by-equivariance}
Suppose that $\mathscr G_{\bar s}$ is nontrivial for all geometric points $\bar s$ of $S$. Then \eqref{eq-compactly-supported-section-twisted-equivariant} vanishes on the essential image of $\derived_{\mathscr G, \zeta}(X/G)$.
\end{prop}

\begin{proof}
Consider the Cartesian square of stacks
$$
\begin{tikzcd}[column sep = 2em]
	X \ar[r, "p"]\ar[d, "f"] & X/G \ar[d, "g"] \\
	S \ar[r, "e"] & \deloop G
\end{tikzcd}
$$

By descent and base change, there is a functor
$$
g_! : \derived_{\mathscr G, \zeta}(X/G) \rightarrow \derived_{\mathscr G, \zeta}(\deloop G)
$$
related to \eqref{eq-compactly-supported-section-twisted-equivariant} by the isomorphism $e^*g_! \simeq f_! p^*$. Since $\derived_{\mathscr G, \zeta}(\deloop G) \simeq 0$ (\emph{cf.}~Lemma \ref{lem-equivariance-vanishing}), the functor $f_! p^*$ vanishes, as desired.
\end{proof}

\subsection{Hyperbolic localization}
\label{sec-hyperbolic-localization}

\begin{void}
\label{void-hyperbolic-localization-notions}
We remain in the context of \S\ref{void-sheaf-theory-context}. Let $X$ be a quasi-separated $S$-scheme locally of finite presentation equipped with a $\mathbb G_m$-action. Let $\mathscr G$ be an $A$-gerbe over $X / \mathbb G_m$, which we shall often view as a $\mathbb G_m$-equivariant $A$-gerbe over $X$.

The $\mathbb G_m$-action on $X$ defines the \emph{attractor} $X^+$, the \emph{repeller} $X^-$, and the \emph{fix point locus} $X^0$, as $S$-presheaves parametrizing $\mathbb G_m$-equivariant morphisms
\begin{align*}
	X^+ & := \SMaps_{\mathbb G_m}(\mathbb A^{1, +}, X) \\
	X^- & := \SMaps_{\mathbb G_m}(\mathbb A^{1, -}, X) \\
	X^0 & := \SMaps_{\mathbb G_m}(\Spec\base, X)
\end{align*}
were $\mathbb A^{1, +}$ (respectively $\mathbb A^{1, -}$) denotes the affine line over $S$ endowed with the $\mathbb G_m$-action $z\cdot a := za$ (respectively $z\cdot a := z^{-1} a$).

Our assumption on $X$ implies that the $\mathbb G_m$-action on $X$ is \'etale locally linearizable, \emph{i.e.}~it admits a $\mathbb G_m$-equivariant \'etale cover by affine $S$-schemes (\emph{cf.}~\cite[Corollary 10.2]{alper2025etalelocalstructurealgebraic}). It follows that $X^+$, $X^-$, $X^0$ are represented by $S$-schemes locally of finite presentation (\emph{cf.}~\cite[Theorem 1.8]{MR3912059}), with $X^0$ being a closed subscheme of $X$.
\end{void}

\begin{void}
We shall consider the following morphisms of $S$-schemes
\begin{equation}
\label{eq-hyperbolic-localization-diagram}
\begin{tikzcd}[column sep = 1.5em]
	& X^+ \ar[dl, swap, "q^+"]\ar[dr, shift left = 0.5ex, "p^+"] \\
	X & & X^0 \ar[ul, shift left = 0.5ex, "i^+"] \ar[dl, shift right = 0.5ex, swap, "i^-"] \\
	& X^- \ar[ul, "q^-"] \ar[ur, swap, shift right = 0.5ex, "p^-"]
\end{tikzcd}
\end{equation}
where $q^+$, $q^-$ (respectively $p^+$, $p^-$) are defined by evaluation at $1$ (respectively $0$), and $i^+$, $i^-$ are defined by pullback along $\mathbb A^{1, +} \rightarrow S$, $\mathbb A^{1, -} \rightarrow S$.

The morphisms $i^+$, $i^-$ are closed immersions (\emph{cf.}~\cite[Proposition 1.17]{MR3912059}), and $q^+ i^+ = q^- i^-$ is the natural closed immersion of $X^0$ into $X$.
\end{void}

\begin{void}
Writing $\mathscr G^0$ for the restriction of $\mathscr G$ to $X^0$, we have canonical isomorphisms of $A$-gerbes over $X^+$ and $X^-$:
\begin{align}
\label{eq-gerbe-identification-attractor}
	(q^+)^*\mathscr G \simeq (p^+)\mathscr G^0, \\
\label{eq-gerbe-identification-repeller}
	(q^-)^*\mathscr G \simeq (p^-)\mathscr G^0.
\end{align}

Indeed, given an $R$-point of $X^+$ represented by the $\mathbb G_m$-equivariant morphism $\mathbb A^{1, +} \times \Spec R \rightarrow X$, the pullback of $\mathscr G$ to $\mathbb A^{1, +} \times \Spec R$ canonically descends to $\Spec R$, so we obtain an identification of its pullbacks along $\{1\}\times \Spec R$ and $\{0\} \times \Spec R$. This yields \eqref{eq-gerbe-identification-attractor}. The isomorphism \eqref{eq-gerbe-identification-repeller} arises from the same analysis for $X^-$.

Define the full subcategory of \emph{$\mathbb G_m$-monodromic} objects
$$
\derived_{\mathscr G, \zeta}(X)^{\mathbb G_m\mon} \subset \derived_{\mathscr G, \zeta}(X)
$$
to be the essential image of $\derived_{\mathscr G, \zeta}(X/\mathbb G_m)$ under the pullback functor.

The main result of this subsection is the following twisted version of Braden's hyperbolic localization theorem (\emph{cf.}~\cite{MR1996415}). Its proof will be supplied in \S\ref{void-hyperbolic-localization-proof}.
\end{void}

\begin{thm}
\label{thm-hyperbolic-localization}
Suppose that $X \rightarrow S$ is separated and of finite presentation. There is a canonical isomorphism of functors from $\derived_{\mathscr G, \zeta}(X)^{\mathbb G_m\mon}$ to $\derived_{\mathscr G^0, \zeta}(X^0)$:
\begin{equation}
\label{eq-hyperbolic-localization}
(p^+)_!(q^+)^* \simeq (p^-)_*(q^-)^!.
\end{equation}
\end{thm}

\begin{rem}
The formation of $(p^+)_!$ and $(q^-)^!$ require $p^+$ and $q^-$ to be separated and of finite presentation. The morphism $p^+$ is affine (\emph{cf.}~\cite[Corollary 1.12]{MR3912059}), hence separated. Our assumption on $X$ implies that $q^-$ is a monomorphism (\emph{cf.}~\cite[Remark 1.19]{MR3912059}), hence separated as well. Finite presentation follows from \cite[Theorem 1.8(iii)]{MR3912059}.
\end{rem}

\begin{void}
As in the untwisted setting (\emph{cf.}~\cite[\S3.4]{MR3200429}), Theorem \ref{thm-hyperbolic-localization} can be strengthened to a statement concerning quotient stacks.

Namely, we take the quotient of \eqref{eq-hyperbolic-localization-diagram} by $\mathbb G_m$:
\begin{equation}
\label{eq-hyperbolic-localization-diagram-equivariant}
\begin{tikzcd}[column sep = 0em]
	& X^+/\mathbb G_m \ar[dl, swap, "q^+"]\ar[dr, shift left = 0.5ex, "p^+"] \\
	X/\mathbb G_m & & X^0 / \mathbb G_m \ar[ul, shift left = 0.5ex, "i^+"] \ar[dl, shift right = 0.5ex, swap, "i^-"] \\
	& X^- /\mathbb G_m \ar[ul, "q^-"] \ar[ur, swap, shift right = 0.5ex, "p^-"]
\end{tikzcd}
\end{equation}
without changing notation for the morphisms. There is a canonical morphism of functors from $\derived_{\mathscr G, \zeta}(X/\mathbb G_m)$ to $\derived_{\mathscr G^0, \zeta}(X^0/\mathbb G_m)$:
\begin{equation}
\label{eq-hyperbolic-localization-equivariant}
	(p^+)_!(q^+)^* \rightarrow (p^-)_*(q^-)^!,
\end{equation}
where $\derived_{\mathscr G, \zeta}(X/\mathbb G_m)$ denotes the $\infty$-category of $\mathbb G_m$-equivariant objects of $\derived_{\mathscr G, \zeta}(X)$ (and likewise for $\derived_{\mathscr G^0, \zeta}(X^0/\mathbb G_m)$). The construction of \eqref{eq-hyperbolic-localization-equivariant} is a repetition of its untwisted analogue (\emph{cf.}~\cite[Construction 2.2]{MR3912059}).

To prove Theorem \ref{thm-hyperbolic-localization}, it suffices to show that \eqref{eq-hyperbolic-localization-equivariant} is an isomorphism. (More precisely, one constructs the functor \eqref{eq-hyperbolic-localization} in one direction as in \emph{loc.cit.}~and uses its compatibility with \eqref{eq-hyperbolic-localization-equivariant} to achieve this reduction.)
\end{void}

\begin{void}
\label{void-hyperbolic-localization-relevant-locus}
Since $\mathbb G_m$ acts trivially on $X^0$, we have $X^0/\mathbb G_m \simeq X^0\times\deloop\mathbb G_m$, so the quotient map $X^0 \rightarrow X^0/\mathbb G_m$ admits a section $s : X^0/\mathbb G_m \rightarrow X^0$.

The Kummer $\hat{\integers}(1)$-gerbe over $\deloop\mathbb G_m$ (\emph{cf.}~Remark \ref{rem-kummer-class-of-tautological-line-bundle}) induces an isomorphism
\begin{equation}
\label{eq-multiplicative-group-equivariant-gerbe-classification}
\Maps(X^0, A(-1)) \simeq \Maps_{X^0/}(X^0/\mathbb G_m, \deloop^2 A).
\end{equation}
The $A$-gerbe $\mathscr G^0 \otimes s^*(\mathscr G^0|_{X^0})^{-1}$ over $X^0/\mathbb G_m$ defines a morphism $X^0 \rightarrow A(-1)$ under \eqref{eq-multiplicative-group-equivariant-gerbe-classification}, whose vanishing locus $Y^0 \subset X^0$ is an open and closed subscheme.

By definition, the restriction of $\mathscr G^0$ to $Y^0/\mathbb G_m$ is pulled back along $s : Y^0/\mathbb G_m \rightarrow Y^0$. The following observation is an application of \cite[Theorem 1.1]{alper2025etalelocalstructurealgebraic}.
\end{void}

\begin{lem}
\label{lem-gerbe-etale-local-linearization}
Given a point $y \in Y^0$, there exists an affine $S$-scheme $U$ with a $\mathbb G_m$-action and a $\mathbb G_m$-equivariant \'etale morphism $f : U \rightarrow X$ such that
\begin{enumerate}
	\item the image of $f$ contains $y$;
	\item the pullback $f^*\mathscr G$ is $\mathbb G_m$-equivariantly trivial.
\end{enumerate}
\end{lem}

\begin{proof}
The (total space of the) $A$-gerbe $\mathscr G$ over $X/\mathbb G_m$ is an algebraic $S$-stack satisfying the hypotheses of \cite[Theorem 1.1]{alper2025etalelocalstructurealgebraic}. Replacing $y$ by a finite \'etale extension, we may assume that it lifts to a point $\tilde y$ of $\mathscr G$. The residue gerbe $\mathscr G_{\tilde y}$ is non-canonically isomorphic to $\deloop\mathbb G_m \times \deloop A$. Applying \emph{loc.cit.}~to the \'etale morphism
$$
\id \times e : \deloop\mathbb G_m \rightarrow \deloop\mathbb G_m \times \deloop A
$$
we may extend it to a Cartesian square of algebraic stacks
$$
\begin{tikzcd}[column sep = 1.5em]
	\deloop\mathbb G_m \ar[d, "\id\times e"]\ar[r, phantom, "\subset"] & \mathscr U \ar[d, "h"] \\
	\deloop\mathbb G_m \times\deloop A \ar[r, phantom, "\subset"] & \mathscr G
\end{tikzcd}
$$
where $h$ is \'etale and $\mathscr U$ contains $\deloop\mathbb G_m$ as a closed substack.

The composition $\mathscr U \xrightarrow{h} \mathscr G \rightarrow \deloop\mathbb G_m$ restricts to the identity along $\deloop\mathbb G_m \subset \mathscr U$, so by \cite[Proposition 5.3]{alper2025etalelocalstructurealgebraic}, it becomes affine after shrinking $\mathscr U$. The base change of $h$ along the neutral point $S \rightarrow \deloop\mathbb G_m$ then yields a $\mathbb G_m$-equivariant morphism $U \rightarrow \mathscr G|_X$, where $U$ is $S$-affine. The composition $f : U \rightarrow \mathscr G|_X \rightarrow X$ satisfies the desired properties.
\end{proof}

\begin{rem}
The condition $y \in Y^0$ (as opposed to an arbitrary point of $X^0$) in Lemma \ref{lem-gerbe-etale-local-linearization} is crucial: Even for $X$ a geometric point equipped with the trivial $\mathbb G_m$-action, a $\mathbb G_m$-equivariant $A$-gerbe over $X$ may not be $\mathbb G_m$-equivariantly trivial.
\end{rem}

\begin{void}
\label{void-hyperbolic-localization-proof}
We are now ready to prove Theorem \ref{thm-hyperbolic-localization}, or rather, reduce it to hyperbolic localization for \emph{untwisted} constructible complexes (\emph{cf.}~\cite[Theorem 2.6]{MR3912059}).

\begin{proof}[Proof of Theorem \ref{thm-hyperbolic-localization}]
We shall prove that \eqref{eq-hyperbolic-localization-equivariant} is an isomorphism.

Recall the open and closed subscheme $Y^0 \subset X^0$ (\emph{cf.}~\S\ref{void-hyperbolic-localization-relevant-locus}). By Lemma \ref{lem-equivariance-vanishing}, restriction along $Y^0 \subset X^0$ yields an equivalence of $\infty$-categories
\begin{equation}
\label{eq-hyperbolic-localization-relevant-locus-equivalence}
\derived_{\mathscr G^0, \zeta}(X^0/\mathbb G_m) \simeq \derived_{\mathscr G^0, \zeta}(Y^0/\mathbb G_m).
\end{equation}

By Lemma \ref{lem-gerbe-etale-local-linearization}, we find a family of $\mathbb G_m$-equivariant \'etale morphisms $f_i : U_i \rightarrow X$ ($i\in I$) satisfying the following properties:
\begin{enumerate} 
	\item each $U_i$ is $S$-affine;
	\item the pullback $f_i^*\mathscr G$ is $\mathbb G_m$-equivariantly trivial;
	\item the union $Y := \bigcup_{i\in I} f_i(U_i)$ is an open subscheme of $X$ containing $Y^0$.
\end{enumerate}

To prove that \eqref{eq-hyperbolic-localization-equivariant} is an isomorphism, we may replace $X$ by any $\mathbb G_m$-stable open subscheme containing $Y^0$ (for instance, $Y$): In view of \eqref{eq-hyperbolic-localization-relevant-locus-equivalence}, it suffices to show that both sides of \eqref{eq-hyperbolic-localization-equivariant} are unchanged after this replacement. This can be seen from the ``contraction lemma" (\emph{cf.}~\cite[Proposition 5.3.2]{MR3356356}), \emph{i.e.}~isomorphisms of functors
$$
(p^+)_! \simeq (i^+)^!,\quad (p^-)_* \simeq (i^-)^*,
$$
which remain valid for twisted complexes because they are of \'etale local nature over $Y^0$, and $\mathscr G^0$ is $\mathbb G_m$-equivariantly trivial over an \'etale cover of $Y^0$.

After replacing $X$ by $Y$, we may show that \eqref{eq-hyperbolic-localization-equivariant} is an isomorphism over the \'etale cover $f_i : U_i \rightarrow Y$ ($i\in I$), where the $\mathbb G_m$-equivariant triviality of each $f_i^*\mathscr G$ allows us to reduce the statement to the untwisted (affine) case (\emph{cf.}~\cite[\S2.8]{MR3912059}).
\end{proof}

\end{void}

\subsection{Properties of $L_{X/S}$}

\begin{void}
\label{void-hyperbolic-localization-property-context}
We continue in the context of \S\ref{void-sheaf-theory-context} and let $X$ be a separated $S$-scheme of finite presentation equipped with a $\mathbb G_m$-action. Let $\mathscr G$ be an $A$-gerbe over $X/\mathbb G_m$.

Having proved Theorem \ref{thm-hyperbolic-localization}, we denote either side of \eqref{eq-hyperbolic-localization} by
$$
L_{X/S} : \derived_{\mathscr G, \zeta}(X)^{\mathbb G_m\mon} \rightarrow \derived_{\mathscr G^0, \zeta}(X^0)
$$
and call it the \emph{hyperbolic localization} functor. In this subsection, we establish some properties of $L_{X/S}$, which are all standard consequences of Theorem \ref{thm-hyperbolic-localization}.
\end{void}

\begin{lem}
\label{lem-hyperbolic-localization-base-change}
The formation of $L_{X/S}$ commutes with $*$-and $!$-base changes in $S$.
\end{lem}

\begin{proof}
Given a morphism $S' \rightarrow S$ of schemes satisfying the conditions of \S\ref{void-sheaf-theory-context}, with induced morphisms $f : X' := X\times_S S' \rightarrow X$ and $f^0 : (X^0)':= X^0\times_S S' \rightarrow X^0$, we need to construct canonical isomorphisms
\begin{align*}
(f^0)^* L_{X/S} &\simeq L_{X'/S'} f^*, \\
(f^0)^! L_{X/S} &\simeq L_{X'/S'} f^!.
\end{align*}
These are immediate, by presenting $L_{X/S}$ as $(p^+)_!(q^+)^*$, respectively $(p^-)_* (q^-)^!$.
\end{proof}

\begin{prop}
\label{prop-hyperbolic-localization-ula-preservation}
The functor $L_{X/S}$ preserves universal local acyclicity relative to $S$.
\end{prop}

\begin{proof}
Let $\mathscr A$ be an object of $\derived_{\mathscr G, \zeta}(X)^{\mathbb G_m\mon}$ which is ULA relative to $S$. We want to prove that $L_{X/S}(\mathscr A) \in \derived_{\mathscr G^0, \zeta}(X^0)$ is also ULA relative to $S$.

By Lemma \ref{lem-hyperbolic-localization-base-change} and \cite[Theorem 4.4]{MR4630128} (which remains valid in our setting, \emph{cf.}~Remark \ref{rem-ula-relative-to-geometric-point}), we may assume that $S$ is the spectrum of of a rank-$1$ valuation ring with algebraically closed fraction field and prove that the natural map (``cospecialization")
\begin{equation}
\label{eq-hyperbolic-localization-cospecialization}
i^* L_{X/S}(\mathscr A) \rightarrow i^*j_*j^* L_{X/S}(\mathscr A)
\end{equation}
is an isomorphism. Here, $i, j$ are the base changes of the inclusion of the special, respectively generic point of $S$.

Recall that the formation of $L_{X/S}$ commutes with $i^*$ and $j^*$ (\emph{cf.}~Lemma \ref{lem-hyperbolic-localization-base-change}). It also commutes with $j_*$ by presenting $L_{X/S}$ as $(p^-)_*(q^-)^!$. Therefore, \eqref{eq-hyperbolic-localization-cospecialization} is identified with the image of $i^*\mathscr A \rightarrow i^*j_*j^*\mathscr A$ under $L_{X/S}$, which is an isomorphism because $\mathscr A$ is ULA.
\end{proof}

\begin{void}
Next, we study the interaction between $L_{X/S}$ and Verdier duality (\emph{cf.}~\S\ref{void-verdier-duality}).

Denote by $L_{X/S}^-$ the hyperbolic localization functor for $X$ endowed with the \emph{inverted} $\mathbb G_m$-action. Note that $-\mathscr G$ is naturally $\mathbb G_m$-equivariant with respect to the inverted action, and we may present $L_{X/S}^-$ as either of the isomorphic functors
$$
L_{X/S}^- \simeq (p^-)_!(q^-)^* \simeq (p^+)_*(q^+)^!.
$$

The Verdier duality functor $\mathbf D_{X/S} : \derived_{\mathscr G, \zeta}(X) \rightarrow \derived_{-\mathscr G, \zeta}(X)$ restricts to a functor on the full subcategories of $\mathbb G_m$-monodromic objects.
\end{void}

\begin{prop}
\label{prop-hyperbolic-localization-verdier-duality}
The following diagram is canonically commutative:
$$
\begin{tikzcd}[column sep = 1.5em]
	\derived_{\mathscr G, \zeta}(X)^{\mathbb G_m\mon} \ar[r, "L_{X/S}"]\ar[d, "\mathbf D_{X/S}"] & \derived_{\mathscr G^0, \zeta}(X^0) \ar[d, "\mathbf D_{X^0/S}"] \\ 
	\derived_{-\mathscr G, \zeta}(X)^{\mathbb G_m\mon} \ar[r, "L_{X/S}^-"] & \derived_{-\mathscr G^0, \zeta}(X^0)
\end{tikzcd}
$$
\end{prop}

\begin{proof}
For any object $\mathscr A$ of $\derived_{\mathscr G, \zeta}(X)^{\mathbb G_m\mon}$, we have natural isomorphisms
\begin{align*}
	\mathbf D_{X^0/S} L_{X/S}(\mathscr A) & \simeq \SHom((p^+)_!(q^+)^*\mathscr A, \omega_{X^0/S}) \\
	& \simeq (p^+)_*\SHom((q^+)^*\mathscr A, (p^+)^! \omega_{X^0/S}) \\
	& \simeq (p^+)_*\SHom((q^+)^*\mathscr A, (q^+)^! \omega_{X/S}) \\
	& \simeq (p^+)_*(q^+)^!\SHom(\mathscr A, \omega_{X/S}) \simeq L^-_{X/S} \mathbf D_{X/S}(\mathscr A)
\end{align*}
as desired. (Here, $\omega_{Z/S}$ denotes the relative dualizing complex $\pi^!\coeff$ for any separated morphism $\pi : Z \rightarrow S$ of finite presentation.)
\end{proof}

\medskip

\section{Related works}
\label{sec-prior-works-satake}

In this section, we indicate the relation between our version of the twisted geometric Satake equivalence (\emph{cf.}~Theorem \ref{thm-satake-equivalence}) with its predecessors: \cite[Theorem 2.9]{MR2684259}, \cite[Theorem IV.8.3]{MR2956088}, \cite[Theorem 2.1]{lysenko2014twisted}, and \cite[\S9.2]{MR3769731}.

Since the works \cite{MR2684259, lysenko2014twisted} use the K-theoretic parametrization of covers of \cite{MR1896177}, we shall first recall this formalism in \S\ref{sec-integral-vs-etale-levels} in order to set the stage. Much of this material already appears in \cite[\S3.4]{MR3769731}. In \S\ref{sec-finkelberg-lysenko}-\ref{sec-reich}, we discuss the aforementioned works on the twisted geometric Satake equivalence in relation to our version. We also point out some shortcomings in them.\footnote{Some of these shortcomings are known to experts. However, in putting them in writing, the author assumes sole responsibility for any potential misrepresentation of these works.} Finally, in \S\ref{sec-gaitsgory-semisimplicity}, we sketch another proof of the semisimplicity of the Satake category (\emph{cf.}~Corollary \ref{cor-semisimplicity}), communicated to the author by Gaitsgory.\footnote{Any deficiency in the presentation of this proof is due to the author.}

\subsection{Brylinski--Deligne covers}
\label{sec-integral-vs-etale-levels}

\begin{void}
Let $S$ be a smooth scheme over a field. Denote by $\Ktheory_2$ the Zariski sheafification of the $S$-presheaf sending $\Spec R \rightarrow S$ to the second algebraic K-group $K_2(R)$.

Given a smooth affine group $S$-scheme $G$, a \emph{Brylinski--Deligne cover} (or an \emph{integral level}) of $G$ is a central extension
\begin{equation}
\label{eq-brylinski-deligne-cover}
1 \rightarrow \Ktheory_2 \rightarrow E \rightarrow G \rightarrow 1
\end{equation}
of sheaves of groups on the (big) Zariski site of $S$.

Equivalently, a Brylinski--Deligne cover of $G$ is a morphism of pointed Zariski stacks $\deloop_{\Zar} G \rightarrow \deloop_{\Zar}^2 \Ktheory_2$, where $\deloop_{\Zar}$ is the deloop functor for Zariski stacks (\emph{cf.}~\S\ref{void-bar-construction}).
\end{void}

\begin{void}
Let $G$ be a reductive group $S$-scheme and $N$ be an integer invertible over $S$.

Gaitsgory constructed an \emph{\'etale realization} functor (\emph{cf.}~\cite[\S6]{MR4117995}, \cite[\S2.3]{zhao2022metaplectic})
\begin{equation}
\label{eq-brylinski-deligne-covers-etale-realization}
R_{\etale} : \Maps_*(\deloop_{\Zar} G, \deloop^2_{\Zar}\Ktheory_2) \rightarrow \Maps_*(\deloop G, \deloop^4\mu_N^{\otimes 2}),
\end{equation}
where $\deloop$ denotes the deloop functor for \'etale stacks, as in the main body of the text. In particular, when $\mu_N$ is constant over $S$, every Brylinski--Deligne cover defines an \'etale level valued in the abelian group $A := \mu_N(S)$.

Notably, when $S$ is the spectrum of a local or global field $F$ containing $N$ distinct $N$th roots of unity, $\mu_N$ is constant over $S$. This brings us to the setting of the Langlands--Weissman program (\emph{cf.}~\cite{MR3802417}). In this setting, the construction of topological covers factors through the \'etale realization functor \eqref{eq-brylinski-deligne-covers-etale-realization} (\emph{cf.}~\cite[Proposition 2.3.13]{zhao2022metaplectic}).
\end{void}

\begin{rem}
In the setting of the Langlands--Weissman program, it is natural to compare Weissman's construction of the $L$-group (\emph{cf.}~\cite{MR3802418}, assuming that $G$ is quasi-split) with our construction for its \'etale realization (formulated for function fields in \S\ref{void-ell-group}).

As in \cite{MR3802418}, our $L$-group is obtained from a Baer sum of extensions (\emph{cf.}~Remark \ref{rem-ell-group-baer-sum}), where the ``first twist" is identified with Weissman's meta-Galois twist (\emph{cf.}~Theorem \ref{thm-meta-galois-twist}) and the ``second twist" is obtained from the $\integers$-linear morphism \eqref{void-metaplectic-dual-morphism-linear-component}, whose construction is valid over an arbitrary base scheme $S$ over which $\mu_N$ is constant.

In particular, our construction carries over to number fields, by replacing the ``first twist" by Weissman's meta-Galois twist. What remains to be done is thus a comparison of \eqref{void-metaplectic-dual-morphism-linear-component} with the construction of \cite[\S3]{MR3802418} over a general base. We leave this for future work.
\end{rem}

\begin{void}
\label{void-factorization-line-bundle}
Let us now work over a smooth curve $X$ over a field $\base$ (\emph{cf.}~\S\ref{void-local-hecke-stack-context}). Let $N$ be an integer invertible in $\base$. Denote by $\Psi$ the reduction mod $N$ of the Kummer map (\emph{cf.}~Remark \ref{rem-kummer-class-of-tautological-line-bundle}), assigning to a line bundle the $\mu_N$-gerbe parametrizing its $N$th roots:
$$
\Psi : \deloop\mathbb G_m \rightarrow \deloop^2 \mu_N.
$$

We fix a Brylinski--Deligne cover \eqref{eq-brylinski-deligne-cover} of a reductive group $X$-scheme $G$. According to \cite[Corollary 3.4.7]{zhao2023halfintegral}, there is an associated factorization central extension
\begin{equation}
\label{eq-factorization-central-extension-loop-group}
1 \rightarrow \mathbb G_{m, \Ran} \rightarrow \widetilde{LG} \rightarrow LG \rightarrow 1
\end{equation}
equipped with a splitting over $L^+G$. Here, $L^+G$, $LG$ are the arc, respectively loop groups over the Ran space (\emph{cf.}~\S\ref{void-local-hecke-stack-definition}).

Taking quotients of $\widetilde{LG}$ by the left and right $L^+G$-actions defined by the above section, we obtain a factorization line bundle
\begin{equation}
\label{eq-hecke-stack-factorization-line-bundle}
\mathscr L_{\Hec_G} := L^+G \backslash \widetilde{LG} / L^+G
\end{equation}
over the local Hecke stack $\Hec_G$, equipped with a multiplicative structure (in the sense of \S\ref{void-local-hecke-stack-gerbe-multiplicative}, after replacing $\deloop^2 A$ by $\deloop\mathbb G_m$) induced from $\widetilde{LG}$.\footnote{We remark that the construction of $\widetilde{LG}$, and consequently that of $\mathscr L_{\Hec_G}$, is significantly less direct than that of $\mathscr G_{\Hec_G}$ (\emph{cf.}~\S\ref{void-local-hecke-stack-gerbe}). This makes it much harder to prove properties about $\mathscr L_{\Hec_G}$.}
\end{void}

\begin{void}
\label{void-satake-category-from-line-bundle}
Let us now assume further that $\mu_N$ is constant over $\base$ and fix an injective character $\zeta : \mu_N(\base) \hookrightarrow \overline{\rationals}{}^{\times}_{\ell}$, where $\ell$ is a prime invertible in $\base$.

Applying the construction of the Satake category (\emph{cf.}~\S\ref{void-satake-subcategory}) with $\mathscr G_{\Hec_G}$ replaced by $\Psi(\mathscr L_{\Hec_G})$, we obtain a monoidal category
\begin{equation}
\label{eq-satake-category-from-line-bundle}
\Sat_{\Psi(\mathscr L), \zeta}(\Hec_{G, \underline x})
\end{equation}
associated to any $S$-point $\underline x$ of $\Ran$.

Concretely, objects of \eqref{eq-satake-category-from-line-bundle} are perverse ULA sheaves over the total space of $\mathscr L_{\Hec_{G, \underline x}}$ (relative to $S$) which are $\mathbb G_m$-equivariant against the character local system induced from $\zeta$. The monoidal structure on \eqref{eq-satake-category-from-line-bundle} is given by convolution with respect to the multiplicative structure on $\mathscr L_{\Hec_G}$.

Using the factorization structure on $\mathscr L_{\Hec_G}$, we may lift the monoidal category \eqref{eq-satake-category-from-line-bundle}, for $\underline x$ the canonical map $X^I \rightarrow \Ran$ corresponding to a finite set $I$, to a symmetric monoidal category following \S\ref{void-fusion-product}--\S\ref{void-normalized-satake-category}. This results in a symmetric monoidal category
\begin{equation}
\label{eq-satake-category-from-line-bundle-normalized}
{}^+\Sat_{\Psi(\mathscr L), \zeta}(\Hec_{G, I})
\end{equation}
naturally associated to any finite set $I$.
\end{void}

\begin{void}
In order to compare the Satake category \eqref{eq-satake-category-from-line-bundle-normalized} with the one introduced in \S\ref{sec-satake-category} (for the \'etale realization of the given Brylinski--Deligne cover), we need to compare $\Psi(\mathscr L_{\Hec_G})$ with the $\mu_N$-gerbe $\mathscr G_{\Hec_G}$ (\emph{cf.}~\S\ref{void-local-hecke-stack-gerbe}).

The following version of \cite[Conjecture 3.4.5]{MR3769731} is expected to be valid without the assumption that $\mu_N$ be constant, as long as $N$ is invertible in $\base$.\footnote{Although we constructed $\mathscr G_{\Hec_G}$ for coefficients in the constant sheaf $A$, the construction generalizes \emph{verbatim} for the locally constant sheaf of coefficients $A := \mu_N$.}
\end{void}

\begin{conj}
\label{conj-factorization-line-bundle-kummer-gerbe}
There is a canonical isomorphism of $\mu_N$-gerbes over $\Hec_G$
$$
\Psi(\mathscr L_{\Hec_G}) \simeq \mathscr G_{\Hec_G}
$$
compatibly with the factorization and multiplicative structures (\emph{cf.}~Proposition \ref{prop-local-hecke-stack-gerbe-factorization}).
\end{conj}

\begin{void}
More functorially, we expect the following diagram to commute:
\begin{equation}
\label{eq-etale-realization-kummer-compatibility}
\begin{tikzcd}[column sep = 1.5em]
	\Maps_*(\deloop_{\Zar} G, \deloop^2_{\Zar} \Ktheory_2) \ar[r, "R_{\etale}"]\ar[d] & \Maps_*(\deloop G, \deloop^4\mu_N^{\otimes 2}) \ar[d] \\
	\Maps_{\fact}(\Hec_G^{[\cdot]}, \deloop_{\Ran}\mathbb G_m^{[\cdot]})\ar[r, "\Psi"] & \Maps_{\fact}(\Hec_G^{[\cdot]}, \deloop_{\Ran}^2 \mu_N^{[\cdot]}) 
\end{tikzcd}
\end{equation}
where the lower row consists of the spaces of factorization multiplicative line bundles, respectively $\mu_N$-gerbes over $\Hec_G$, and the vertical arrows are given by \cite[Corollary 3.4.7]{zhao2023halfintegral}, respectively Proposition \ref{prop-local-hecke-stack-gerbe-factorization}.

Let us indicate how one may obtain the commutativity of \eqref{eq-etale-realization-kummer-compatibility} by somewhat indirect means. Note that the vertical arrows in \eqref{eq-etale-realization-kummer-compatibility} are equivalences: For the left vertical arrow, this follows from \cite[Corollary 3.4.7]{zhao2023halfintegral}, while for the right vertical arrow, this follows from \cite[Proposition 3.1.9]{MR3769731} and \cite[Proposition 7.2.5]{MR3769731}.

Therefore, in order to prove that \eqref{eq-etale-realization-kummer-compatibility} commutes, we may replace the vertical arrows by their inverses. Then we may replace the top row by their classification data using \cite[Theorem 7.2]{MR1896177}, respectively \cite[Theorem 5.1.13]{zhao2022metaplectic}. It then remains to show that the resulting square commutes, which one reduces to the case of tori.
\end{void}

\begin{rem}
The reason for the indirectness of this argument is that the left vertical arrow in \eqref{eq-etale-realization-kummer-compatibility} lacks a ``natural" definition (\emph{cf.}~the discussion in \cite[\S4.5]{zhao2023halfintegral}). Namely, if it can be constructed in the same spirit as our construction of $\mathscr G_{\Hec_G}$, then the commutativity of \eqref{eq-etale-realization-kummer-compatibility} would be tautological.
\end{rem}

\begin{void}
Let us return to the context of \S\ref{void-satake-category-from-line-bundle} and assume the validity of Conjecture \ref{conj-factorization-line-bundle-kummer-gerbe}. We then have an equivalence of symmetric monoidal categories:
\begin{equation}
\label{eq-satake-category-line-bundle-gerbe-comparison}
	{}^+\Sat_{\Psi(\mathscr L), \zeta}(\Hec_{G, I}) \simeq {}^+ \Sat_{\mathscr G, \zeta}(\Hec_{G, I})
\end{equation}
naturally in the finite set $I$.

In particular, one can reformulate our geometric Satake equivalence (\emph{cf.}~Theorem \ref{thm-satake-equivalence}) in terms of the factorization line bundle $\mathscr L_{\Hec_G}$:
\begin{equation}
\label{eq-satake-equivalence-line-bundle}
	{}^+\Sat_{\Psi(\mathscr L), \zeta}(\Hec_{G, I}) \simeq \Rep_{H^{\boxtimes I}, (\nu + \vartheta)^{\boxplus I}}
\end{equation}
by composing \eqref{eq-satake-equivalence} with \eqref{eq-satake-category-line-bundle-gerbe-comparison}.
\end{void}

\subsection{\cite{MR2684259} and \cite{lysenko2014twisted}}
\label{sec-finkelberg-lysenko}

\begin{void}
The general context for \cite{MR2684259, lysenko2014twisted} is that of \S\ref{void-factorization-line-bundle}-\S\ref{void-satake-category-from-line-bundle}. However, these texts present more assumptions. Let us first briefly indicate their roles.

In \cite{MR2684259}, two further assumptions are imposed:
\begin{enumerate}
	\item $G$ is almost simple and simply connected;
	\item the characteristic $p$ of $\base$ does not divide $2\check h/d$, where $\check h$ is the dual Coxeter number of $G$ and $d$ is defined as in \cite[\S2.1]{MR2684259}.
\end{enumerate}

The reason for (1) is that \cite{MR2684259} constructs the central extension \eqref{eq-factorization-central-extension-loop-group} only in this context (\emph{cf.}~\cite[Proposition 2.2]{MR2684259}). Note that \emph{op.cit.}~does not construct the factorization structure on \eqref{eq-hecke-stack-factorization-line-bundle}, but rather uses the factorization structure of its $(2\check h/d)$-power to construct the commutativity constraint on the Satake category. This is responsible for assumption (2) which, however, can be removed following \cite[Remark 2.10]{MR2684259}.

In \cite{lysenko2014twisted}, one encounters a series of input data concerning central extensions of $LG$ by $\mathbb G_{m, \Ran}$ (\emph{cf.}~\cite[\S2.3]{lysenko2014twisted}). The reason for this additional complication again has to do with \eqref{eq-factorization-central-extension-loop-group}, which was not constructed in general at the time \cite{lysenko2014twisted} was written.
\end{void}

\begin{void}[Fiber functor]
\label{void-finkelberg-lysenko-fiber-functor}
Next, we mention an issue in the proof of \cite[Theorem 2.9]{MR2684259}, which is pointed out by Reich (\emph{cf.}~\cite[\S V.1]{MR2956088}). Namely, the twisted Satake equivalence for tori stated in \cite[\S4.2]{MR2684259} is \emph{not} a symmetric monoidal equivalence.

Let us explain the issue more precisely. Assume that $G$ is split simple and simply connected, endowed with a pinning $T \subset B \subset G$, $\mathfrak g_{\check{\alpha}} \simeq \mathbb G_a$ ($\check{\alpha}\in\check{\Delta}$) and a Brylinski--Deligne cover \eqref{eq-brylinski-deligne-cover}. Let us also fix a square root $\Omega^{1/2}$ of the canonical line bundle of $X$.

In this context, it is effectively stated in \cite[\S4.2]{MR2684259} that we have an equivalence
\begin{equation}
\label{eq-finkelberg-lysenko-equivalence}
	\Sat_{\Psi(\mathscr L), \zeta}(\Hec_{T, \{1\}}) \simeq \Rep_{T_H},
\end{equation}
where $T_H$ is the canonical maximal torus of $H$ (\emph{i.e.}~the Langlands dual of $T^{\sharp}$), given by restrictions along $\varpi^{\lambda} : X \rightarrow \Hec_{T, \{1\}}$ for $\lambda \in \Lambda^{\sharp}$ (\emph{cf.}~\S\ref{void-schubert-cell-notation}).\footnote{We use the notation $T^{\sharp}$, $T_H$, and $H$ as in \S\ref{sec-metaplectic-dual-data}, but here these objects are obtained directly from the combinatorial classification of Brylinski--Deligne covers. Namely, to \eqref{eq-brylinski-deligne-cover}, one may attach an integral Weyl-invariant quadratic form $Q$ on $\Lambda$ (\emph{cf.}~\cite[Theorem 6.2]{MR1896177}). Applying the constructions of \S\ref{sec-metaplectic-dual-data} to $Q$ mod $N$ then yields $T^{\sharp}$, $T_H$, and $H$.

Note that $Q$ mod $N$ coincides with the strict Weyl-invariant quadratic form associated to the \'etale realization of \eqref{eq-brylinski-deligne-cover} (\emph{cf.}~\S\ref{void-etale-level-fiber-sequence}), in view of \cite[\S6.1]{MR1896177} and \cite[\S5.2]{zhao2022metaplectic}.}

The same statement appears in \cite[\S3.2]{lysenko2014twisted} with weaker assumptions on $G$.
\end{void}

\begin{void}
\label{void-finkelberg-lysenko-counterexample}
Let us specialize to the example $G = \SL_2$ equipped with the unique Brylinski--Deligne cover whose quadratic form $Q$ takes value $1$ on the simple coroot $\alpha \in \Lambda$.

The factorization line bundle $\mathscr L_{\Hec_G}$ is then identified with the \emph{inverse} of the determinant line bundle. Namely, the pullback of $\mathscr L_{\Hec}$ along an $R$-point $P^0 \overset{\underline x}{\sim} P^1$ of $\Hec_G$ is the line bundle $\det(P^1 \mid P^0)^{\otimes -1}$, where $P^0$, $P^1$ are viewed as rank-$2$ vector bundles over $D_{\underline x}$ identified over $\mathring D_{\underline x}$ (\emph{cf.}~\cite[\S2]{MR1961134}). Here, the notation $\det(P^1 \mid P^0)$ means the relative determinant
$$
\det(P^1 \mid P^0) := \det(L^1 / L) \otimes \det(L^0 / L)^{-1},
$$
where $L^0$, $L^1$ are the lattices in the Tate $R$-module $\Gamma(\mathring D_{\underline x}, P^0) \simeq \Gamma(\mathring D_{\underline x}, P^1)$ defined by $P^0$, respectively $P^1$, and $L$ is a lattice contained in $L^0 \cap L^1$.

The restriction $\mathscr L_{\Hec_T} := \mathscr L_{\Hec_G} |_{\Hec_T}$ (where $T := \mathbb G_m$ is the diagonal maximal torus) admits a canonical square-root as a factorization \emph{super} line bundle:
\begin{equation}
\label{eq-determinant-restriction-square-root}
\mathscr L_{\Hec_T} \simeq (\mathscr L_{\Tate})^{\otimes 2},
\end{equation}
where $\mathscr L_{\Tate}$ is the super line bundle over $\Hec_{\mathbb G_m}$ assigning $\det(\mathscr L^1 \otimes \Omega^{1/2} \mid \mathscr L^0 \otimes \Omega^{1/2})^{\otimes -1}$ with its natural $\integers/2$-grading to an $R$-point $\mathscr L^0 \overset{\underline x}{\sim} \mathscr L^1$ of $\Hec_{\mathbb G_m}$ (\emph{cf.}~\cite[\S5.2.4]{MR3769731}). The latter has grading $\lambda$ mod $2$ over the component of $\Hec_{\mathbb G_m, \{1\}}$ corresponding to $\lambda \in \Lambda \simeq \integers$.

Let us take $N = 2$, so $\Lambda^{\sharp} = \Lambda$, which we will identify with $\integers$. It follows from \eqref{eq-determinant-restriction-square-root} that $\Psi(\mathscr L_{\Hec_T})$ is canonically isomorphic to the pullback of the sign gerbe along the reduction mod $2$ map $\epsilon : \Lambda \rightarrow \integers/2$ (\emph{cf.}~\cite[\S4.9.1]{MR3769731}). Consequently, we obtain an equivalence of symmetric monoidal categories
\begin{equation}
\label{eq-finkelberg-lysenko-counterexample}
\Sat_{\Psi(\mathscr L), \zeta}(\Hec_{T, \{1\}}) \simeq {}^{\epsilon}\Rep_{T_H},
\end{equation}
where ${}^{\epsilon}\Rep_{T_H}$ is equivalent to $\Rep_{T_H}$ as a $\Lambda$-graded monoidal category, but its commutativity constraint is given by the isomorphism
\begin{align*}
	V^{\lambda_1} \otimes V^{\lambda_2} & \simeq V^{\lambda_2} \otimes V^{\lambda_1}\\
	a \otimes b & \mapsto (-1)^{\epsilon(\lambda_1)\epsilon(\lambda_2)} b\otimes a
\end{align*}
for any two objects $V^{\lambda_1}, V^{\lambda_2}$ with $\Lambda$-gradings $\lambda_1$, respectively $\lambda_2$.

There are no symmetric monoidal equivalences between ${}^{\epsilon} \Rep_{T_H}$ and $\Rep_{T_H}$ compatible with the $\Lambda$-gradings, so the same applies to the two sides of \eqref{eq-finkelberg-lysenko-equivalence}.
\end{void}

\begin{rem}
Assuming the validity of Conjecture \ref{conj-factorization-line-bundle-kummer-gerbe}, we may also obtain an equivalence \eqref{eq-finkelberg-lysenko-counterexample} from \eqref{eq-satake-equivalence-line-bundle} (for $G = T$ and $I = \{1\}$).

Indeed, let $\mu$ be the \'etale realization of the Brylinski--Deligne cover of \S\ref{void-finkelberg-lysenko-counterexample}. The restriction of $\mu$ to $T \simeq \mathbb G_m$ is canonically identified with the \'etale level of \S\ref{void-hilbert-central-extension} (\emph{cf.}~\cite[Lemma 5.1.19]{zhao2022metaplectic}). It follows that the $\mathbb E_{\infty}$-monoidal morphism $\nu : \Lambda \rightarrow \deloop^2 \{\pm 1\}$ admits an $\mathbb E_1$-monoidal trivialization, in terms of which the commutativity constraint of its fiber is given by the pairing $\lambda_1 \otimes \lambda_2 \mapsto (-1)^{\epsilon(\lambda_1)\epsilon(\lambda_2)}$ on $\Lambda \otimes \Lambda$ (\emph{cf.}~\cite[\S2.2.9]{shi2025extendedpureinnerforms}).

Thus, we obtain a symmetric monoidal equivalence
$$
\Rep_{T_H, \nu} \simeq {}^{\epsilon}\Rep_{T_H},
$$
and consequently \eqref{eq-finkelberg-lysenko-counterexample}, by composing with \eqref{eq-satake-equivalence-line-bundle} and identifying $\nu + \vartheta$ with $\nu$ using the choice of $\Omega^{1/2}$ (\emph{cf.}~Remark \ref{rem-theta-characteristic-trivializes-first-twist}). We have not checked that this construction of \eqref{eq-finkelberg-lysenko-counterexample} coincides with the one of \S\ref{void-finkelberg-lysenko-counterexample}.
\end{rem}

\begin{void}
While the issue of \S\ref{void-finkelberg-lysenko-fiber-functor} may seem innocuous, it causes a gap in the proof of \cite[Theorem 2.9]{MR2684259}. Namely, because \eqref{eq-finkelberg-lysenko-equivalence} is not symmetric monoidal, it is not possible to apply the Tannakian formalism to ${}^+\Sat_{\Psi(\mathscr L), \zeta}(\Hec_{G, \{1\}})$.

We do not see a way to fill this gap without invoking the results of this text.
\end{void}

\begin{void}[Semisimplicity]
\label{void-lysenko-semisimplicity}
The work \cite{MR2684259} refers to the proof of \cite[Proposition 11]{MR2265675} for the semisimplicity of the twisted Satake category.

However, the proof of \cite[Proposition 11]{MR2265675} essentially states the the $A$-gerbe $\Psi(\mathscr L_{\Hec_G})$ (at least in the special case $G = \Sp_{2n}$ and $N = 2$) is trivial over the Demazure resolutions of closures of Bruhat cells in the affine flag variety. This is \emph{not} true.

For a concrete example, we may again consider the one of \S\ref{void-finkelberg-lysenko-counterexample}. Denote by $\Fl_{G, x}^w$ the Bruhat cell in the affine flag variety $\Fl_{G, x}$ at a $\base$-point $x$ of $X$, associated to an element $w$ of the (extended) affine Weyl group $W^{\aff}$ (\emph{cf.}~\S\ref{sec-parity-vanishing}). Write $\overline{\Fl}{}_{G, x}^w$ for its closure.

If $w$ is a simple reflection, then $\overline{\Fl}{}_{G, x}^w$ coincides with its own Demazure resolution, and is isomorphic to $\mathbb P^1$ (\emph{cf.}~\cite[\S3]{MR1961134}). If we take $w$ to be the simple reflection $s$ associated to the \emph{affine} simple root, then the pullback of the determinant line bundle $\mathscr L_{\Hec_G}$ to $\overline{\Fl}{}_{G, x}^s$ is isomorphic to $\mathscr O_{\mathbb P^1}(1)$ (\emph{cf.}~\cite[Theorem 7]{MR1961134}). This line bundle does not have a square root, so $\Psi(\mathscr L_{\Hec_G})$ is nontrivial over $\overline{\Fl}{}_{G, x}^s$. 
\end{void}

\begin{rem}
The above computation is consistent with Lemma \ref{lem-gerbe-dimension-one-bruhat-variety}, as the $\{\pm 1\}$-gerbe $\mathscr G^{\aff}$ there coincides with the pullback of $\Psi(\mathscr L_{\Hec_G})$. (This is a special case of Conjecture \ref{conj-factorization-line-bundle-kummer-gerbe}, but for the example at hand, it can be established directly.)
\end{rem}

\subsection{\cite{MR2956088} and \cite{MR3769731}}
\label{sec-reich}

\begin{void}
We now comment on Reich's treatment of the twisted geometric Satake equivalence \cite{MR2956088}, along with further developments due to Gaitsgory and Lysenko \cite{MR3769731}. In what follows, we shall take $G$ to be a split reductive group.

We first remark that Reich's text \cite{MR2956088} is written in the complex analytic context, so ``gerbes" in \emph{op.cit.}~refers to $\complexes^{\times}$-gerbes in the analytic topology.

The text \cite{MR3769731} is written in the same context as ours (\emph{i.e.}~\'etale cohomology with torsion coefficients), but its treatment of the Satake equivalence requires the characteristic $p$ of the ground field $\base$ to be coprime to the order of $\pi_1 (G_{\der})$ (\emph{cf.}~Remark \ref{rem-gaitsgory-lysenko-characteristic-assumption} below).
\end{void}

\begin{void}
The text \cite{MR2956088} contains three main results:
\begin{enumerate}
	\item the classification of factorization gerbes over $\Gr_G$ (\emph{cf.}~\cite[Theorem II.7.3]{MR2956088});
	\item the proof that any factorization gerbe over $\Gr_G$ admits a canonical $L^+G$-equivariance structure, satisfying certain properties (\emph{cf.}~\cite[Theorem III.2.10]{MR2956088});
	\item the construction of the geometric Satake equivalence \cite[Theorem IV.8.3]{MR2956088} associated to each factorization gerbe over $\Gr_G$. 
\end{enumerate}

Note that in both \cite{MR2956088} and \cite{MR3769731}, results (1) \& (2) are needed for the formulation of the geometric Satake equivalence, because their authors use factorization gerbes over $\Gr_G$ to parametrize covering groups (as opposed to \'etale levels, \emph{cf.}~\S\ref{void-local-hecke-stack-context}).\footnote{As mentioned in the introduction, the idea of using \'etale levels to parametrize covering groups dates much earlier: It first appeared in Deligne's work \cite{MR1441006}.}

In what follows, we will discuss the status of these three results in order.
\end{void}

\begin{void}[Classification of factorization gerbes]
\label{void-classification-factorization-gerbes}
The result \cite[Theorem II.7.3]{MR2956088} is incorrect as stated, because the fiber sequence of \emph{loc.cit.}~does not split. Since its proof relies on an incorrectly defined splitting, we do not see a way to fix it.

The correct classification of factorization gerbes over $\Gr_G$ appears as \cite[Proposition 3.1.9]{MR3769731}: This is the assertion that factorization gerbes over $\Gr_G$ are \emph{equivalent} to \'etale levels. In particular, it implies that the input data of our version of the Satake equivalence (\emph{cf.}~Theorem \ref{thm-satake-equivalence}) are equivalent to those of \cite[\S9]{MR3769731}, except that we do \emph{not} need the condition that $p$ is coprime to $\pi_1(G_{\der})$.
\end{void}

\begin{rem}
\label{rem-gaitsgory-lysenko-characteristic-assumption}
In \cite{MR3769731}, the condition on $p$ is needed for the construction of the metaplectic dual data (\emph{cf.}~\cite[\S6.2, Appendix A]{MR3769731}).

It is unnecessary for us since our construction of the metaplectic dual data is different from that of \cite{MR3769731} and does not make use of the affine Grassmannian (\emph{cf.}~\S\ref{sec-metaplectic-dual-data}).
\end{rem}

\begin{void}[$L^+G$-equivariance structure]
\label{void-arc-equivariance}
The result \cite[Theorem III.2.10]{MR2956088} does not appear to be adequately justified.

More precisely, in the proofs of \cite[Proposition III.2.8]{MR2956088} and \cite[Theorem III.2.10]{MR2956088}, the author applies the notion of ``orders" of gerbes (\emph{cf.}~\cite[\S I.5]{MR2956088}) along certain Cartier divisors in indschemes to justify their extension properties. This relies on cohomological purity of the Cartier divisors in question, which is not justified.

However, \cite[\S7.3]{MR3769731} offers two constructions of the $L^+G$-equivariance structure: One of them is essentially restated in \S\ref{sec-local-hecke-stack-gerbe} of this text, but the other one uses different ideas.
\end{void}

\begin{void}[Geometric Satake equivalence]
\label{void-reich-satake-equivalence}
Taking into account the corrections \S\ref{void-classification-factorization-gerbes}-\ref{void-arc-equivariance} offered in \cite{MR3769731}, one may state Reich's version of the geometric Satake equivalence as in \cite[\S9.2]{MR3769731}. However, \emph{op.cit.}~refers to \cite[Theorem IV.8.3]{MR2956088} for its proof and there are at least two issues with Reich's proof:
\begin{enumerate}
	\item The proof of the semisimplicity of the pointwise Satake category (\emph{cf.}~\cite[Proposition IV.6.13]{MR2956088}) contains an error.
	\item It does not offer a correct construction of the fiber functor.
\end{enumerate}

The issue with semisimplicity occurs in \cite[Lemma IV.6.9]{MR3769731}, where the criterion for equality there does \emph{not} hold in general.

The issue with the fiber functor occurs in \cite[Lemma IV.7.8]{MR3769731}, whose proof indicates that the $\hat Z_H$-grading (denoted by $X^*(Z(\check G_Q))$ in \emph{loc.cit.})~on the twisted Satake category is induced from its $\pi_1 G$-grading. We cannot make sense of this statement since the natural map $\hat Z_H \rightarrow \pi_1 G$ is \emph{not} injective in general.

On the other hand, assuming the semisimplicity of the Satake category, one can obtain a $\hat Z_H$-grading by coarsening its $\Lambda^{\sharp}$-grading. However, the $\hat Z_H$-grading obtained this way is not obviously compatible with the monoidal structure. Thus, we do not know how to fix Reich's construction of the fiber functor without invoking the results of this text.
\end{void}

\begin{rem}
As mentioned in \S\ref{void-finkelberg-lysenko-fiber-functor} and \S\ref{void-lysenko-semisimplicity}, these two issues of the proof are already present in \cite{MR2684259}, and at least one of them is explicitly pointed out by Reich.

Let us also acknowledge that Reich first observed the need to construct a $\hat Z_H$-grading on the Satake category compatible with the monoidal structure, even though this may not have been achieved in \cite{MR2956088}.
\end{rem}

\subsection{Gaitsgory's proof of semisimplicity}
\label{sec-gaitsgory-semisimplicity}

\begin{void}
Finally, let us sketch an unpublished proof due to Gaitsgory of the semisimplicity of the twisted Satake category at a point (\emph{cf.}~Corollary \ref{cor-semisimplicity}). It avoids parity vanishing (\emph{cf.}~Proposition \ref{prop-parity-vanishing}) and instead uses the $\hat Z_H$-grading of \S\ref{sec-virtual-connected-components}, combined with the method of the proof of \cite[Theorem 3]{MR1836286}.

Let us return to the context of \S\ref{sec-pointwise-studies}. Furthermore, we shall fix a $\base$-point $x$ of $X$ and consider the local Hecke stack $\Hec_G$ at $x$, suppressing $x$ from the notation. Our goal is to prove that $\Sat_{\mathscr G, \zeta}(\Hec_G)$ is semisimple.
\end{void}

\begin{void}
\label{void-ext-vanishing-same-coweight}
First, we establish the following vanishing result: Given $\lambda \in \Lambda^{\sharp, +}$ and $(\mathscr G_{\varpi^{\lambda}}, \zeta)$-twisted lisse sheaves $\mathscr E_1$, $\mathscr E_2$ over $x$, we have
\begin{equation}
\label{eq-ext-vanishing-same-coweight}
\Ext^1(\IC^{\lambda}(\mathscr E_1), \IC^{\lambda}(\mathscr E_2)) \simeq 0,
\end{equation}
where $\IC^{\lambda}$ is the functor of \S\ref{void-affine-grassmannian-ic-sheaf}.

To prove \eqref{eq-ext-vanishing-same-coweight}, we may choose a trivialization of $\mathscr G_{\varpi^{\lambda}}$ at $x$ and reduce to the case where $\mathscr E_1$, $\mathscr E_2$ are both the constant sheaf $\coeff$. Then the argument of \cite[Proposition 1]{MR1826370} applies. (Details may also be found in \cite[Lemma IV.6.6]{MR2956088}).

Next, we need the following consequence of the Decomposition Theorem.
\end{void}

\begin{prop}[Reich]
\label{prop-convolution-preserves-semisimplicity}
Given $\lambda_1, \lambda_2 \in \Lambda^{\sharp, +}$ and $\mathscr G_{\varpi^{\lambda_1}}$, $\mathscr G_{\varpi^{\lambda_2}}$-twisted lisse sheaves $\mathscr E_1$, $\mathscr E_2$ over $x$, the convolution product $\IC^{\lambda_1}(\mathscr E_1) \star \IC^{\lambda_2}(\mathscr E_2)$ is a direct sum of objects of the form $\IC^{\lambda}(\mathscr E)$, for $\lambda \in \Lambda^{\sharp, +}$ and $\mathscr G_{\varpi^{\lambda}}$-twisted lisse sheaves $\mathscr E$ over $x$.
\end{prop}

\begin{void}
\label{void-convolution-semisimple-summand}
Proposition \ref{prop-convolution-preserves-semisimplicity} is a reformulation of \cite[Lemma IV.6.11]{MR2956088}. Note, however, that the proof of \emph{loc.cit.}~uses a corollary of the problematic \cite[Lemma IV.6.9]{MR2956088} (\emph{cf.}~\S\ref{void-reich-satake-equivalence}), but the problematic part (the criterion for equality) is not used.

Let us write the conclusion of Proposition \ref{prop-convolution-preserves-semisimplicity} as an isomorphism
\begin{equation}
\label{eq-convolution-semisimple-decomposition}
	\IC^{\lambda_1}(\mathscr E_1) \star \IC^{\lambda_2}(\mathscr E_2) \simeq \bigoplus_{\lambda} \IC^{\lambda}(\mathscr E)
\end{equation}
in $\Sat_{\mathscr G, \zeta}(\Hec_G)$. We now need a consequence of the $\hat Z_H$-grading on the Satake category: By Corollary \ref{cor-virtual-connected-components-convolution-compatibility} (and the proof of Proposition \ref{prop-virtual-connected-components}), all elements $\lambda \in \Lambda^{\sharp, +}$ appearing in the sum \eqref{eq-convolution-semisimple-decomposition} share the same equivalence class in $\hat Z_H$ as $\lambda_1 + \lambda_2$.
\end{void}

\begin{void}
As in \S\ref{sec-fiber-functor}, we construct the untwisted Satake category ${}^+\Sat_{\mathscr G, \zeta}(\Hec_G)_{-(\nu + \vartheta)}$ together with the fiber functor $\omega$ valued in finite-dimensional $\coeff$-vector spaces.

Applying Tannakian formalism, we obtain a symmetric monoidal equivalence
\begin{equation}
\label{eq-tannaka-dual-without-semisimplicity}
	{}^+\Sat_{\mathscr G, \zeta}(\Hec_G)_{-(\nu + \vartheta)} \simeq \Rep_{\check G},
\end{equation}
for some affine monoid $\coeff$-scheme $\check G$. The rigidity of the Satake category (\emph{cf.}~Proposition \ref{prop-convolution-monoidal-dual}) implies that $\check G$ is an affine group $\coeff$-scheme.

The construction of \S\ref{void-tannaka-dual-maximal-torus} yields a group subscheme $T_H \subset \check G$, where $T_H$ is the canonical maximal torus of $H$. Write $\hat B$ for the group subscheme of $\check G$ consisting of automorphisms of $\omega$ preserving the filtration defined by $2\check{\rho}$, \emph{i.e.}~$F^{\ge n}\omega \subset \omega$ ($n \in \integers$) is the sum of components of $T_H$-weights $\lambda$ satisfying $\langle 2\check{\rho}, \lambda\rangle \ge n$.
\end{void}

\begin{void}
Write $\mathscr C$ for the full subcategory of ${}^+\Sat_{\mathscr G, \zeta}(\Hec_G)_{-(\nu + \vartheta)}$ generated by the images of (the untwisted forms of) $\IC^{\lambda}$ ($\lambda \in \Lambda^{\sharp, +}$) under \emph{direct sums}.

By \eqref{eq-ext-vanishing-same-coweight}, $\mathscr C$ is semisimple. Since ${}^+\Sat_{\mathscr G, \zeta}(\Hec_G)_{-(\nu + \vartheta)}$ is generated by the images of $\IC^{\lambda}$ ($\lambda \in \Lambda^{\sharp, +}$) under extensions, $\mathscr C$ coincides with its \emph{maximal} semisimple subcategory. By Proposition \ref{prop-convolution-preserves-semisimplicity}, $\mathscr C$ inherits a symmetric monoidal structure.

By Tannakian formalism, we have a symmetric monoidal equivalence $\mathscr C \simeq \Rep_{\check G_{\red}}$ for some affine group scheme $\check G_{\red}$, equipped with a natural homomorphism
\begin{equation}
\label{eq-tannaka-group-reductive-quotient}
\check G \rightarrow \check G_{\red}.
\end{equation}

The morphism \eqref{eq-tannaka-group-reductive-quotient} is faithfully flat by \cite[Proposition 2.21(a)]{Deligne1982}. Furthermore, the proof of Lemma \ref{lem-tannaka-dual-reductive} shows that $\check G_{\red}$ is reductive and the proof of Lemma \ref{lem-tannaka-group-maximal-torus} shows that the composition $T_H \hookrightarrow \check G \rightarrow \check G_{\red}$ realizes $T_H$ as a maximal torus of $\check G_{\red}$. Denote by $\check B_{\red}$ the image of $\check B$ in $\check G_{\red}$, which is a Borel subgroup.

The following implication is the main insight of the proof.
\end{void}

\begin{lem}
\label{lem-root-data-implies-semisimplicity}
Assuming that $T_H \subset \check B_{\red} \subset \check G_{\red}$ and $T_H \subset B_H \subset H$ define the same root data, then \eqref{eq-tannaka-group-reductive-quotient} is an isomorphism.
\end{lem}

\begin{proof}
It suffices to proving the vanishing statement
\begin{equation}
\label{eq-ext-vanishing-different-coweight}
\Ext^1(\IC^{\lambda_1}(\coeff), \IC^{\lambda_2}(\coeff)) \simeq 0,
\end{equation}
for any $\lambda_1, \lambda_2 \in \Lambda^{\sharp, +}$, where the Ext group is computed in ${}^+\Sat_{\mathscr G, \zeta}(\Hec_G)_{-(\nu + \vartheta)}$. Note that the special case for $\lambda_1 = \lambda_2$ is known by \S\ref{void-ext-vanishing-same-coweight}. Moreover, by Proposition \ref{prop-virtual-connected-components}, \eqref{eq-ext-vanishing-different-coweight} also holds when $\lambda_1$, $\lambda_2$ have distinct images in $\hat Z_H$, so let us assume that $\lambda_1$, $\lambda_2$ have the same image in $\hat Z_H$ in what follows.

Note that we have an isomorphism
\begin{equation}
\label{eq-ext-vanishing-different-coweight-adjoint}
\Ext^1(\IC^{\lambda_1}(\coeff), \IC^{\lambda_2}(\coeff)) \simeq \Ext^1(\IC^0(\coeff), \IC^{\lambda_1}(\coeff)^{\vee} \otimes \IC^{\lambda_2}(\coeff)),
\end{equation}
where the monoidal dual $\IC^{\lambda_1}(\coeff)^{\vee}$ is a simple object of ${}^+\Sat_{\mathscr G, \zeta}(\Hec_G)_{-(\nu + \vartheta)}$ with $\hat Z_H$-grading the image of $-\lambda_1$ (\emph{cf.}~Corollary \ref{cor-virtual-connected-components-convolution-compatibility}). This implies that $\IC^{\lambda_1}(\coeff)^{\vee} \otimes \IC^{\lambda_2}(\coeff)$ is a direct sum of objects of the form $\IC^{\lambda}(\coeff)$, where $\lambda \in \Lambda^{\sharp, +}$ has vanishing image in $\hat Z_H$ (\emph{cf.}~\S\ref{void-convolution-semisimple-summand}).

It thus remains to prove \eqref{eq-ext-vanishing-different-coweight} when $\lambda_1 = 0$ and $\lambda_2$ belonging to the root lattice of $H$. Since $T_H \subset \check B_{\red} \subset \check G_{\red}$ and $T_H \subset B_H \subset H$ are assumed to have the same root data, $\lambda_2$ belongs to the root lattice of $\hat G_{\red}$. Applying \cite[Lemma 6.1.3]{MR1836286} to $\check G_{\red}$, we may realize $\IC^{\lambda_2}(\coeff)$ as a summand of $\IC^{\lambda}(\coeff)^{\vee} \star \IC^{\lambda}(\coeff)$, for some $\lambda \in \Lambda^{\sharp, +}$. Applying \eqref{eq-ext-vanishing-different-coweight-adjoint} again, we reduce to the assertion \eqref{eq-ext-vanishing-different-coweight} for $\lambda_1 = \lambda_2$, which is established in \S\ref{void-ext-vanishing-same-coweight}.
\end{proof}

\begin{void}
We shall now show that the assumption of Lemma \ref{lem-root-data-implies-semisimplicity} holds unconditionally.

This requires a minor modification of the proof of Proposition \ref{prop-tannaka-group-root-data}, because the construction of $\check G_{\red}$ is \emph{a priori} incompatible with constant term functors.
\end{void}

\begin{lem}
\label{lem-tannaka-dual-reductive-quotient-root-data}
The triples $T_H \subset \check B_{\red} \subset \check G_{\red}$, $T_H \subset B_H \subset H$ define the same root data.
\end{lem}

\begin{proof}
The case where $G$ has a unique simple root is handled as in the proof of Proposition \ref{prop-tannaka-group-root-data}. For the general case, we fix a simple root $\check{\alpha}$ of $G$ and write $P$ for the corresponding subminimal parabolic with Levi quotient $M$. The constant term functor induces a morphism of affine group $\coeff$-schemes
\begin{equation}
\label{eq-tannaka-group-levi-subgroup-without-semisimplicity}
\check M \rightarrow \check G.
\end{equation}

Applying Lemma \ref{lem-root-data-implies-semisimplicity} to $M$, we find that $\check M$ is a reductive group $\coeff$-scheme with unique simple root $\alpha^{\sharp}$ and simple coroot $\check{\alpha}{}^{\sharp}$. Then, the argument of Proposition \ref{prop-tannaka-group-root-data} shows that the composition $\check M \rightarrow \check G_{\red}$ of \eqref{eq-tannaka-group-levi-subgroup-without-semisimplicity} and \eqref{eq-tannaka-group-reductive-quotient} is a closed immersion, and consequently yields the computation of the root data of $T_H \subset \check B_{\red} \subset \check G_{\red}$.
\end{proof}

\begin{prop}
\label{prop-tannaka-group-equals-reductive-quotient}
The morphism \eqref{eq-tannaka-group-reductive-quotient} is an isomorphism.
\end{prop}

\begin{proof}
This follows by combining Lemma \ref{lem-root-data-implies-semisimplicity} and Lemma \ref{lem-tannaka-dual-reductive-quotient-root-data}.
\end{proof}

\begin{void}
Proposition \ref{prop-tannaka-group-equals-reductive-quotient} implies that $\check G$ is reductive, thereby giving a new proof of the semisimplicity of the Satake category (\emph{cf.}~Corollary \ref{cor-semisimplicity}).
\end{void}

\medskip

\bibliographystyle{amsalpha}
\bibliography{bibliography.bib}

\end{document}